\documentclass[a4paper, 10pt]{article}
\usepackage{srcltx,graphicx}
\usepackage{amsmath, amssymb, amsthm}

\usepackage{color, xcolor}
\usepackage{array}

\usepackage{lscape}
\usepackage{algorithm}      
\usepackage[noend]{algpseudocode}
\usepackage{multirow}

\usepackage{psfrag}
\usepackage{pifont}
\usepackage{bbding}

\usepackage{hyperref}
\usepackage{makeidx}

\usepackage{caption}
\usepackage{subcaption} 

\usepackage[section]{placeins}

\newtheorem{conclusion}{Conclusion}

{\theoremstyle{remark} }
\numberwithin{figure}{subsection}
\numberwithin{equation}{section}
 
\hypersetup{linktoc = all}                
\hypersetup{hidelinks}
\hypersetup{bookmarksnumbered}
\newcommand\sign     {\text{sgn}}
\newcommand\showfigures[1]{} 
\newcommand \be { \begin{equation}}
\newcommand \ee { \end{equation}}
\newcommand \bel {\begin{equation}\label}
\newcommand \del \partial 

\newcommand \red {}

\begin{document}

\title{\red A class of well-balanced algorithms for relativistic fluids 
\\
on a Schwarzschild background\footnotetext{
$^{1}$Laboratoire Jacques-Louis Lions \& Centre National de la Recherche Scientifique,
Sorbonne Universit\'e, 4 Place Jussieu, 75252 Paris, France.
Email: {\sl contact@philippelefloch.org}
\newline
$^2$Departamento de An\'alisis Matem\'atico, Estad\'istica e Investigaci\'on Operativa, y Matem\'atica aplicada, Universidad de M\'alaga, Bulevar Louis Pasteur, 31, 29010, M\'alaga, Spain.
Email: {\sl pares@uma.es, erpigar@uma.es}
\newline 
\red
{\it Keywords and Phrases}. Compressible fluid; relativistic flow; Schwarzschild black hole; well-balanced algorithm; asymptotic behavior.}
}

\author{Philippe G. LeFloch$^1$, Carlos Par\'es$^2$, and Ernesto Pimentel-Garc\'ia$^2$}

\date{}

\maketitle

\begin{abstract} 
\red
For the evolution of a compressible fluid in spherical symmetry on a Schwarzschild curved background, we design a class of well-balanced numerical algorithms with first-order or second-order of accuracy. We treat both the relativistic Burgers-Schwarzschild model and the relativistic Euler-Schwarzschild model and take advantage of the explicit or implicit forms available for the stationary solutions of these models. Our schemes follow the finite volume methodology and preserve the stationary solutions. Importantly, they allow us to investigate the global asymptotic behavior of such flows and determine the asymptotic behavior of the mass density and velocity field of the fluid. 
\end{abstract} 
 

{\small 

\setcounter{secnumdepth}{3}
\setcounter{tocdepth}{1} 
\tableofcontents
}

\section{Introduction} 
\label{section--1}
 
{\red

We are interested in the numerical approximation and the long time behavior of relativistic compressible fluid flows on a Schwarzschild black hole background. 
The flow is assumed to enjoy spherical symmetry and therefore we deal with nonlinear hyperbolic systems of partial differential equations (PDEs) in one space variable. 
This paper is two-fold: on the one hand, designing and testing  numerically finite volume algorithms that are well-balanced; on the other to perform 
a thorough investigation of the behavior of the solutions and numerically infer definite conclusions about the long-time behavior of such flows.  
Our study should provide first and useful insights for, on the one hand, further development concerning the mathematical analysis of the models and, on the other hand, further 
investigations to the same problem in higher dimensions without symmetry restriction. 

Two models are of interest in the present paper. 
We treat first the {\sl relativistic Burgers-Schwarzschild model} (as it is called in \cite{PLF-SX1,PLF-SX2}): 
\begin{subequations}
\label{eq:burguers}
\be
v_t + F(v,r)_r = S(v,r), \qquad t \geq 0, \quad r>2M, \qquad
\ee
where $v=v(t,r) \in[-1,1]$ is the unknown function and the flux and source terms read 
\be
 F(v,r) =  \Big( 1-\frac{2M}{r} \Big) \frac{v^2-1}{2}, \qquad S(v,r)=  \frac{2M}{r^2} \, (v^{2}-1),
\ee
\end{subequations} 
while the constant $M>0$ represents the mass of the black hole. Obviously, the speed of propagation for this scalar balance law reads 
\be
\del_v F(v,r) =  \Big( 1-\frac{2M}{r} \Big) v, 
\ee
which vanishes at the boundary $r=2M$, so that no boundary condition is required in order to pose the Cauchy problem. 

Next, we consider the {\sl relativistic Euler-Schwarzschild model} (as it is called in \cite{PLF-SX1,PLF-SX2}):  
\begin{subequations}
\label{eq:Euler}
\be
V_{t}+ F(V,r)_{r} = S(V,r), \qquad t \geq 0, \quad r>2M, \qquad
\ee
whose unknowns are the fluid density $\rho=\rho(t,r) \geq 0$ and the normalized velocity $v=v(t,r) \in (-1,1)$. 
These functions are defined for all $r>2M$ and the limiting values $v=\pm 1$ being reached at the boundary $r=2M$ only, 
and 
\be
V=\left(\begin{array}{c}
V^{0}\\
V^{1}
\end{array}\right)=
\left(\begin{array}{c}
\displaystyle \frac{1+k^{2}v^{2}}{1-v^{2}}\rho\\
\displaystyle \frac{1+k^{2}}{1-v^{2}}\rho v\\
\end{array}\right), 
\qquad 
F(V,r)= \left(\begin{array}{c}
\displaystyle \Big( 1-\frac{2M}{r} \Big)\frac{1+k^{2}}{1-v^{2}}\rho v \\
\displaystyle \Big( 1-\frac{2M}{r} \Big)\frac{v^{2}+k^{2}}{1-v^{2}}\rho
\end{array}\right),
\ee
\be
S(V,r)=\left(\begin{array}{c}
\displaystyle -\frac{2}{r}\Big( 1-\frac{2M}{r} \Big)\frac{1+k^{2}}{1-v^{2}}\rho v \\
\displaystyle \frac{-2r+5M}{r^{2}}\frac{v^{2}+k^{2}}{1-v^{2}}\rho-\frac{M}{r^{2}}\frac{1+k^{2}v^{2}}{1-v^{2}}\rho+2\frac{r-2M}{r^{2}}k^{2}\rho
\end{array}\right), 
\ee
with
\bel{eq:primitive_variables}
v=\frac{1+k^{2}-\sqrt{(1+k^{2})^{2}-4k^{2}\left(\frac{V^{1}}{V^{0}}\right)^{2}}}{2k^{2}\frac{V^{1}}{V^{0}}}, \qquad
 \rho=\frac{V^{1}(1-v^{2})}{v(1+k^{2})}.
\ee
\end{subequations} 
Here, $k\in(-1,1)$ denotes the (constant) speed of sound, and this second model can be checked to be 
a strictly hyperbolic system. The eigenvalues of the corresponding Jacobian of the flux function read
\bel{eigen}
\mu_{\pm} = \Big( 1-\frac{2M}{r} \Big)\frac{v \pm k}{1 \pm k^2 v}. 
\ee
 As usual, a state $(\rho, v)$, by definition, is said to be  \textit{sonic} if one of the eigenvalues vanishes, i.e.~if $|v| = k$, 
\textit{supersonic}  if both eigenvalues have the same sign, i.e.~if $ |v| > k$, 
or \textit{subsonic}  if the eigenvalues have different signs, i.e.~if $ |v| < k$.
We will need to distinguish between these regimes in order to design a robust scheme for this model. Both eigenvalues $\mu_{\pm}$ 
vanish at the boundary $r=2M$, so that no boundary condition is required in order to pose the Cauchy problem. 
 

 In order to be able of running reliable and accurate numerical simulations for these two models, in this paper we design 
 shock-capturing, high-order, and well-balanced finite volume methods of first- and second-order of accuracy (and even third-order accurate for \eqref{eq:burguers}). 
 Specifically, we extend to the present problem the well-balanced methodology proposed recently by Castro and Par\'es \cite{Castro-Pares} for nonlinear hyperbolic systems of balance laws. For earlier work on well-balanced schemes we also refer to \cite{CGLP,Russo1,Russo2} and, concerning the design of geometry-preserving schemes, we refer for instance to \cite{BLM,bouchut2004nonlinear,AB-PLF,castro2017well,DongPLF,DKM,GiesselmanPLF,Rossmanith-2004,Vazquez} and the references therein. 
 
The properties of the stationary solutions play a fundamental role in the design of well balanced schemes, as well as in 
the study of the long time behavior of solutions. We thus also built here upon earlier investigations by LeFloch and collaborators~\cite{PLF-HM,PLF-SX1,PLF-SX2} on the theory and approximation of the relativistic Burgers- and Euler-Schwarzschild model \eqref{eq:burguers} and \eqref{eq:Euler}. Remarkably, the stationary solutions to both models are available in explicit or implicit form.

An outline of the content of this paper is as follows. 

In Section \ref{section--2} we describe the methodology for this paper and indicate the challenges met with the two models above.
The actual design of the schemes is the content of Section \ref{section--3} (Burgers equation) and Section \ref{section--5} (Euler equations). 
The proposed well-balanced algorithms, by construction, preserve all of the steady state solutions, which is 
an essential property since numerical methods without this property may lead to wrong conclusions. 
Furthermore, for each model we investigate the efficiency, accuracy, and robustness of the proposed algorithms 
first in Section \ref{section--4} (Burgers equation) and Section \ref{section--6} (Euler equations). 
Our numerical experiments below make comparisons between well-balanced and standard methods, and 
 we demonstrate that the proposed schemes are numerically well-balanced and we emphasize the importance of this property in order to reach reliable results. 
Furthermore, we study the late time behavior of solutions to both models and discuss the role of the choice of the value of the initial data at the boundary. 
Finally, we also describe how steady shocks behave under small (or large) perturbations. 
In short, we demonstrate that the global dynamics can be accurately determined by the proposed algorithms and we reach some definite conclusions in Sections \ref{subsec_conclusions_Burgers} and \ref{subsec_conclusions_Euler}, respectively.  
Further details concerning our methodology and conclusions are found in the corresponding sections for each model. 

}


\section{Well-balanced methodology and strategy proposed in this paper}\label{section--2}

\paragraph{General methodology}

Both problems of interest are of the form
\bel{sbla}
V_t + F(V,r)_r = S(V,r), \qquad r > 2M, 
\ee
with unknown $V =V(t,r)\in \mathbb{R}^N$ and $N = 1$ or $2$. Systems of this form have non-trivial stationary solutions, which satisfy the ODE 
\bel{sbl}
F(V,r)_r = S(V,r).
\ee
Our goal is to introduce a family of numerical methods that are well-balanced, i.e.~that preserve the stationary solutions  in a sense to be specified.  We follow the strategy in \cite{Castro-Pares} to which we refer for further details and arguments of proof.

We consider semi-discrete finite volume numerical methods of the form 
\bel{eq:nummet}
\frac{dV_{i}}{dt}= -\frac{1}{\Delta r}
\Big(
F_{i+\frac{1}{2}}-F_{i-\frac{1}{2}} - \int_{r_{i-\frac{1}{2}}}^{r_{i+\frac{1}{2}}} S(\mathbb{P}^t_{i}(r), r)\,dr
\Big),
\ee
and the following notation is used. 
\begin{itemize}
\item  $I_i= [r_{i-\frac{1}{2}},r_{i+\frac{1}{2}}]$ denote the computational cells, whose length $\Delta r$ is taken to be a constant for the sake simplicity in the presentation.

\item $V_i(t)$ denotes the approximate average of the exact solution in the $i$th cell at the time $t$, that is, 
\be
V_i(t) \cong \frac{1}{\Delta r} \int_{r_{i-\frac{1}{2}}}^{r_{i+\frac{1}{2}}}  V(r,t) \, dr. 
\ee

\item $\mathbb{P}^t_{i}(r)$ denotes the approximate solution in the $i$th cell, as given by a high-order reconstruction operator based on the cell averages $\{ V_i(t) \}$, that is,  
$
\mathbb{P}^t_{i}(r) = \mathbb{P}^t_i(r; \{V_j(t)\}_{j \in \mathcal{S}_i})$. 
Here, $\mathcal{S}_i$ denotes the set of cell indices associated with the stencil of the $i$th cell.

\item The flux terms are denoted by $F_{i+\frac{1}{2}} = \mathbb{F}(V_{i+\frac{1}{2}}^{t, -}, V_{i+\frac{1}{2}}^{t, +}, r_{i+\frac{1}{2}})$, where $V_{i+\frac{1}{2}}^{t, \pm}$ are the reconstructed states at the interfaces, i.e.
\be
V_{i+\frac{1}{2}}^{t, -}=\mathbb{P}_i^t(r_{i+\frac{1}{2}}), \quad V_{i+\frac{1}{2}}^{t, +}=\mathbb{P}_{i+1}^t(r_{i+\frac{1}{2}}). 
\ee
Here, $\mathbb{F}$ is a consistent numerical flux{{{{}}}{, i.e.~a continuous function $\mathbb{F}: \mathbb{R}^N \times \mathbb{R}^N \times (2M, + \infty) \to \mathbb{R}^N$ satisfying 
$\mathbb{F}(V,V,r) = F(V,r)$ for all $V, r$. 
}}

\end{itemize}

Furthermore, given a stationary solution $V^*$ of (\ref{sbl}), we use the following terminology. 

\begin{itemize}

\item The numerical method (\ref{eq:nummet}) is said to be well-balanced for $V^*$ if the vector of cell averages of $V^*$ is an equilibrium of the ODE system  (\ref{eq:nummet}).

\item The reconstruction operator is said to be well-balanced for $V^*$ if we have $
\mathbb{P}_i(r)=V^*(r)$ for all $r \in [r_{i-\frac{1}{2}}, r_{i+\frac{1}{2}}]$, 
where $\mathbb{P}_i$ is the approximation of $V^*$ obtained by applying the reconstruction operator to the vector of cell averages of $V^*$.
\end{itemize}
It is easily checked that, if the reconstruction operator is well-balanced for a continuous stationary solution $V^*$ of \eqref{sbl} then the numerical method is also well-balanced for $V^*$.  
The following strategy to design a well-balanced reconstruction operator $\mathbb{P}_i$ on the basis of a standard operator $\mathbb{Q}_i$ was introduced in \cite{CGLP}: 

{{{{}}}{Given}} a family of cell values $\{V_i\}$, in every cell $I_i=[r_{i-\frac{1}{2}}, r_{i+\frac{1}{2}}]$ we proceed as follows. 

\begin{enumerate}
\item Seek, {{{}}}{(whenever possible),} a stationary solution $V_i^*(x)$ {{{{}}}{defined in the stencil of cell $I_i$ ($\cup_{j\in\mathcal{S}
			_{i}}I_j$)}} such that 
\begin{equation} \label{step1}
\frac{1}{\Delta r} \int_{r_{i-\frac{1}{2}}}^{r_{i+\frac{1}{2}}} V_i^* (r) \, dr = V_i.
\end{equation}
{{{}}}{If such a solution does not exist,  take $V^*_i \equiv 0$.}

\item Apply the reconstruction operator to the cell values $\{W_j\}_{j \in \mathcal{S}_i}$ given by
\be 
W_j= V_j - \frac{1}{\Delta r} \int_{r_{j-\frac{1}{2}}}^{r_{j+\frac{1}{2}}} V_i^* (r) \, dr,
\qquad j \in \mathcal{S}_i, 
\ee
in order to obtain
$
\mathbb{Q}_i(r)=\mathbb{Q}_i(r;\{W_j\}_{j \in \mathcal{S}_i}).
$

\item Define finally 
\begin{equation} \label{step3}
\mathbb{P}_i(r)=V_i^*(r)+\mathbb{Q}_i(r).
\end{equation}
\end{enumerate}

It can be then easily shown that the reconstruction operator $\mathbb{P}_i$ in (\ref{step3}) is well-balanced for every stationary solution provided that the reconstruction operator
 $\mathbb{Q}_i$ is exact for the zero function. Moreover, if $\mathbb{Q}_i$ is conservative then $\mathbb{P}_i$ is conservative, in the sense that 
\be 
\frac{1}{\Delta r} \int_{r_{i-\frac{1}{2}}}^{r_{i+\frac{1}{2}}}\mathbb{P}_i (r) \, dr = V_i,
\ee
and $\mathbb{P}_i$ has the same accuracy as $\mathbb{Q}_i$ if  the stationary solutions are sufficiently regular.

{{{}}}{Observe that, if there does not exists a stationary solution defined in the stencil and satisfying \eqref{step1} then the standard reconstruction is used. 
Observe that this choice does not spoil the well-balanced character of the numerical method since, in this case, the cell values in the stencil cannot be the averages of a stationary solution (otherwise, there would be at least one solution $V^*_i$) and therefore there exists no local equilibrium which would be required to preserve. 
On the other hand, if there exist more than one stationary solution defined on the stencil and satisfying \eqref{step1}, 
a criterion is needed in order to select one of them and indeed 
this sometimes happen for the Euler-Schwarzschild model of interest in the present paper. 
The relevant criterion in this regime depends on the particular problem and, in the case of the Euler-Schwarzschild model, 
we propose to distinguish between the (subsonic or supersonic) regimes of the flow, as we will explain later in Section \ref{section--5}.
}

If a quadrature formula (whose order of accuracy must be greater or
equal to the one of the reconstruction operator)
$$
\int_{r_{i-\frac{1}{2}}}^{r_{i+\frac{1}{2}}} f(x) \,dx \approx \Delta r\sum_{l=0}^q \alpha_l f(r_{i,l})
$$
where $\alpha_0, \dots, \alpha_q$, $r_{i,0}, \dots, r_{i,q}$ represent the weights and the nodes of the formula, is used to compute the averages of the initial condition, namely 
$
V_{i,0} = \sum_{l=0}^q \alpha_l  V_0 (r_{i,l})$, 
the reconstruction procedure has to be modified to preserve the well-balanced property:.
Then Steps 1 and 2 have to be replaced by the following ones. 

\begin{enumerate}
\item Seek, {{{{}}} if possible,} the stationary solution $V_i^*(x)$ {{{{}}}{defined in the stencil of cell $I_i$ ($\cup_{j\in\mathcal{S}
			_{i}}I_j$)}} such that 
\begin{equation} \label{step1qf}
\sum_{l=0}^q \alpha_l V_i^* (r_{i,l}) = V_i.
\end{equation}
{{{}}}{If this solution does not exist, take $V^*_i \equiv 0$.}

\item Apply the reconstruction operator to the cell values $\{W_j\}_{j \in \mathcal{S}_i}$ given by
\begin{equation*}
W_j= V_j - \sum_{l=0}^q \alpha_l V_i^* (r_{j,l}),  \qquad j \in \mathcal{S}_i.
\end{equation*}
\end{enumerate}
For first- or second-order methods, if  the midpoint rule is selected to compute the initial averages, i.e.
$
V_{i,0}  = V_0(r_i)$, 
then at the first step of the reconstruction procedure, the problem \eqref{step1qf} reduces to finding the stationary solution satisfying
\bel{step1mp}
V^*_i(r_i) = V_i.
\ee

The well-balanced property of the method can be lost if the quadrature formula is used to compute the integral appearing at the  right-hand side of \eqref{eq:nummet}. In order to circumvent this difficulty,  in \cite{Castro-Pares} it is proposed to rewrite the methods as follows: 
\begin{equation}\label{eq:nummet2}
\begin{split}
\frac{dV_i}{dt} &=-\frac{1}{\Delta r} \left( F_{i+\frac{1}{2}}-F\left( V_i^{t,*}(r_{i+\frac{1}{2}}), r_{i+\frac{1}{2}}\right) - F_{i-\frac{1}{2}} +F\left( V_i^{t,*}(r_{i-\frac{1}{2}}), r_{i-\frac{1}{2}}\right) \right) \\
& \quad 
+ \frac{1}{\Delta r} \int_{r_{i-\frac{1}{2}}}^{r_{i+\frac{1}{2}}} \left( S(\mathbb{P}_i^t (r), r) -S(V_i^{t,*}(r), r) \right) \, dr,
\end{split}
\end{equation}
where $V_i^{t,*}$ is the {{{}}}{function selected in Step 1 for the $i$th cell at time $t$}. 
In this equivalent  form,  a quadrature formula can be applied  to the integral without losing the well-balanced property, and this leads to a numerical method of the form:
\begin{equation}\label{eq:nummet3}
\begin{split}
\frac{dV_i}{dt} &=-\frac{1}{\Delta r} \left( F_{i+\frac{1}{2}}-F\left( V_i^{t,*}(r_{i+\frac{1}{2}}), r_{i+\frac{1}{2}}\right) - F_{i-\frac{1}{2}} +F\left( V_i^{t,*}(r_{i-\frac{1}{2}}), r_{i-\frac{1}{2}}\right) \right) \\
& \quad + \sum_{l =0}^q \alpha_l \left( S(\mathbb{P}_i^t (r_{i,l} ), r_{i,l} ) -S(V_i^{t,*}(r_{i,l}), r_{i,l}) \right). 
\end{split}
\end{equation}
First-order well-balanced methods are obtained by selecting the trivial constant piecewise reconstruction operator as the standard one, i.e.
\be
\mathbb{Q}_i(r, V_i)=V_i, \qquad 
 r \in [r_{i-\frac{1}{2}}, r_{i+\frac{1}{2}}].
\ee
It can be easily checked that the numerical method then reduces to 
\begin{equation}\label{eq:nummetfo}
\frac{dV_i}{dt} =-\frac{1}{\Delta r} \left( F_{i+\frac{1}{2}}-F\left( V_i^{t,*}(r_{i+\frac{1}{2}}), r_{i+\frac{1}{2}}\right) - F_{i-\frac{1}{2}} +F\left( V_i^{t,*}(r_{i-\frac{1}{2}}), r_{i-\frac{1}{2}}\right)\right),
\end{equation}
where
$
 F_{i+\frac{1}{2}} = \mathbb{F}\left(V^*_i( r_{i+\frac{1}{2}}), V^*_{i+1}(r_{i+\frac{1}{2}}), r_{i+\frac{1}{2}} \right).
$

\paragraph{Further strategy for this paper} 

{{{{}}}{When applying the well-balanced reconstruction methodology, a main difficulty 
 comes from the first step, in which we need to find a stationary solution (i.e.~a solution of the ODE system \eqref{sbl}) satisfying \eqref{step1}, \eqref{step1qf}, or
\eqref{step1mp}, 
and this depends upon the way that the relevant integrals are computed. 
This amounts, in general, to a non-trivial, nonlinear problem whose analysis and solution depend on the system of balance laws under consideration.

For the Burgers-Schwarzschild model, the explicit form of the general solution of the ODE \eqref{sbl} is available (as we recall in Section \ref{section--3}) and the nonlinear problem to be solved in Step 1 has
 always a unique solution $v^*_i$. 
 Nevertheless, it may not be possible to extend this solution to the whole stencil and,  in such a case, the standard reconstruction operator must then be used.

For the Euler-Schwarzschild model, only first and second-order methods will be considered in the present paper, 
and the mid-point rule will be used in order to numerically compute the relevant integrals. 
As a consequence, only solutions of the ODE system \eqref{sbl} satisfying \eqref{step1mp} (i.e.~solutions of a standard Cauchy problem) should be found in Step 1. 
Now, the general solution of the ODE system is available in implicit form (as we will recall in Section \ref{section--5})
and, for this system,
 the Cauchy problem \eqref{sbl}, \eqref{step1mp} may have a single solution or no solutions as well as two or more solutions. If there exists no solution, the standard reconstruction is used. 
 If there are two or more solutions, a criterion based on the (subsonic, supersonic) nature of the fluid flow will be used 
 in order to select one of them, as mentioned above.

This methodology can be extended to systems of balance laws for which the solutions to the ODE system 
\eqref{sbl} are not available neither in explicit or implicit form.
For such system, the nonlinear problems arising in Step 1 must be solved numerically. 
For instance, in \cite{gomez2021high} a control-based strategy combined with a standard ODE solver is used 
and solutions
of \eqref{sbl} are computed that satisfy averaging conditions like \eqref{step1} or \eqref{step1qf}. 

Concerning the extension to multidimensional problems, the main challenge is again the problem in Step 1, 
which now consists in finding a solution to the {\sl nonlinear system of PDEs}
 satisfied by the stationary solutions with prescribed average in the cell 
or from values at a point.  
Such a problem, clearly, is  much more difficult to solve (either exactly or numerically) than an ODE system. Moreover, there may exist infinitely many stationary solutions satisfying the average or point value conditions. 
Nevertheless, the methodology can be still applied and allows one to design numerical methods that preserve a certain family of stationary solutions. 
In Step 1, the stationary solution belonging to the prescribed family that satisfies the imposed condition would be selected or the one that is closer in some sense to the cell values in the stencil. For instance, it is possible to apply this methodology 
and design numerical methods for the Euler-Schwarzschild model in three spatial dimension and 
impose spherical-symmetric stationary solutions are preserved. }}

\section{Burgers-Schwarzschild model: designing the numerical algorithm}\label{section--3}

\subsection{Preliminaries}

\begin{figure}[H]
		\centering
		\includegraphics[width=.7\linewidth]{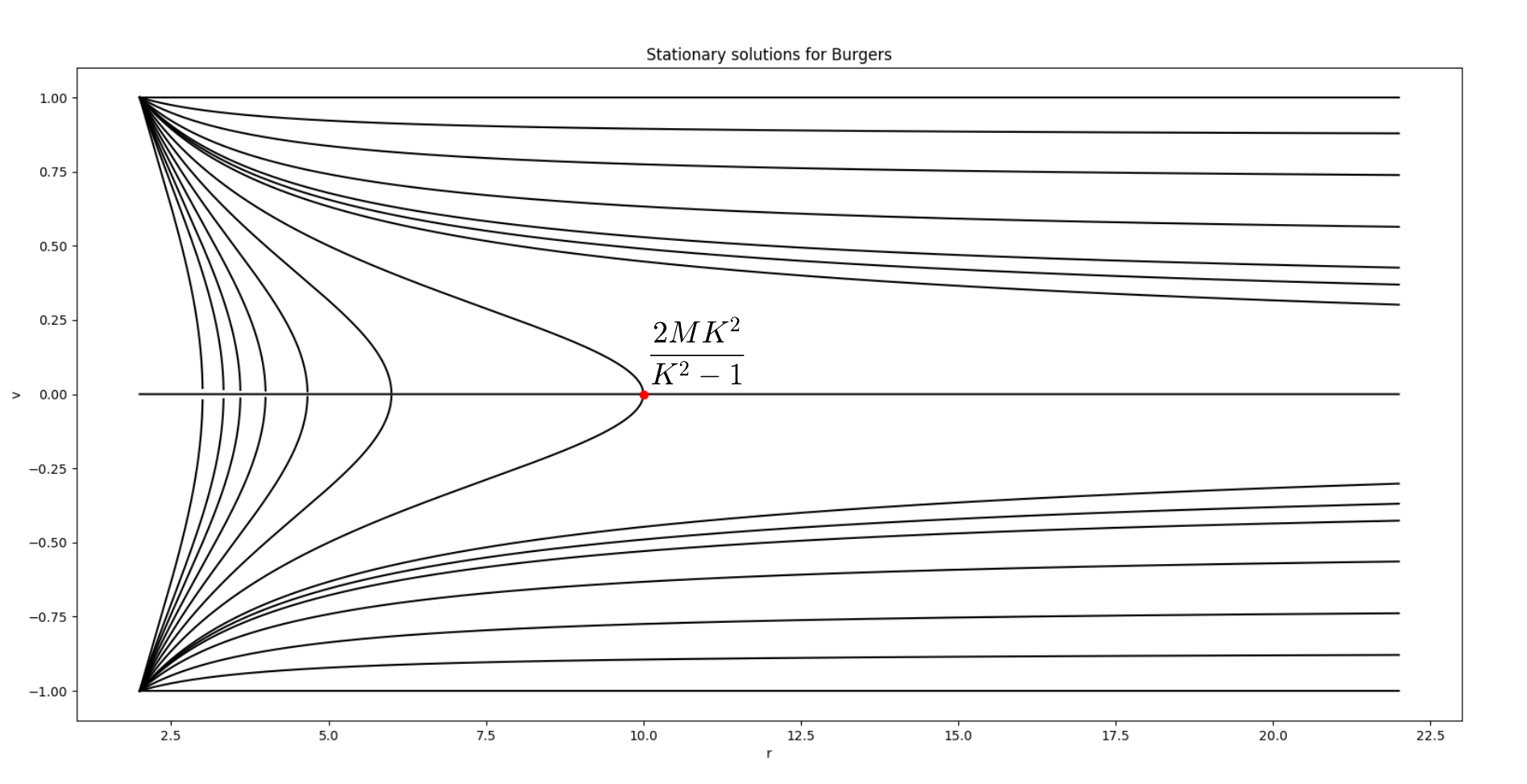}
		\caption{Steady solutions to the Burgers model.}
		\label{fig:steadysolutionsforburgerslabelled}
\end{figure}

For the Burgers-Schwarzschild equation \eqref{eq:burguers}, the steady state solutions are of the form
\be
	v^{*}(r) = \pm \sqrt{1-K^{2}\Big( 1-\frac{2M}{r} \Big)}, \ \ K>0. 
\ee
In Figure \ref{fig:steadysolutionsforburgerslabelled} we plot the steady solutions for several values of $K^{2}$.
The domain of definition of these stationary solutions is
\be
D_{K} = \begin{cases}
[2M, \infty),
 & \text{ $K^{2}\leq 1$,} 
 \\
\displaystyle \left[2M, \frac{2MK^{2}}{K^{2}-1}\right],
 & \text { $K^{2}>1$}.
\end{cases}
\ee	
It can be easily checked that, given a pair $(K, r^* )$ such that $r^* \in D_K$,
 the discontinuous function defined in $D_K$   by
\begin{equation}\label{steadydisc}
w^*(r) = \begin{cases} 
\sqrt{1-K^{2}\Big( 1-\frac{2M}{r} \Big)},
 & \text{$r \leq r^*$,} 
 \\
- \sqrt{1-K^{2}\Big( 1-\frac{2M}{r} \Big)}, & \text{ otherwise,}
\end{cases}
\end{equation}
is an entropy weak stationary solution of  the Burgers-Schwarzschild model. (For further properties, see the study in \cite{PLF-SX1,PLF-SX2}.) 

\subsection{First-order method}

If the midpoint rule is used to compute the initial averages, at the first step of the reconstruction procedure one has to search for $K_{i}^{2}$ such that
	\bel{eq:step1Burger}
	\sqrt{1-K_{i}^{2}\left(1-\frac{2M}{r_{i}}\right)} = |v_{i}|.
	\ee
	There is always a unique solution $\widetilde K_{i}$ given by
	\bel{tildeK}
	\widetilde K_{i}^{2} = \frac{1-v_{i}^{2}}{1-\displaystyle \frac{2M}{r_{i}}},
	\ee
	so that the stationary solution is
	\bel{eq:stadystate1}
	v_{i}^{*}(r) =\sign(v_i)\sqrt{1-\widetilde K_{i}^{2}\Big( 1-\frac{2M}{r} \Big)}.
	\ee
In order to apply the numerical method \eqref{eq:nummet2}, this stationary solution has to be computed at $r_{i\pm \frac{1}{2}}$ and this is only possible if 
$r_{i+\frac{1}{2}} \in D_{\widetilde K_i}$, that is, 
 provided 
\bel{eq:condexfo}
\widetilde K_{i}^{2} \leq 1 \text{ or } \left(\widetilde K_{i}^{2}>1 \text{ and } r_{i+\frac{1}{2}} \leq  \frac{2M\widetilde K_{i}^{2}}{\widetilde K_{i}^{2}-1}\right).
\ee
If this condition is satisfied, then the numerical method \eqref{eq:nummetfo} can be used.

If this condition is not satisfied, then the standard trivial reconstruction is considered, i.e.~
 $
 \mathbb{Q}_{i}(r) = v_{i}. 
 $
 The numerical method writes then as follows:
		\be
		\frac{dv_{i}}{dt}=  -\frac{1}{\Delta r}\left(F_{i+\frac{1}{2}}-F_{i-\frac{1}{2}} - S(v_{i}, r_i)\right),
		\ee
where
$
F_{i+\frac{1}{2}} = \mathbb{F}(v_{i}, v_{i+1}, r_{i+\frac{1}{2}}).
$

Summing up, the expression of the semi-discrete first-order method reads 
\bel{eq:semi_disc-merge}
\frac{dv_{i}}{dt}= -\frac{1}{\Delta r}\left(F_{i+\frac{1}{2}}-F_{i-\frac{1}{2}} - S_i \right),
\ee
where 
\bel{eq:source-merge}
S_ i 
= \begin{cases}
F(v_{i}^{*}(r_{i+\frac{1}{2}}), r_{i+\frac{1}{2}}) - F(v_{i}^{*}(r_{i-\frac{1}{2}}), r_{i-\frac{1}{2}})), & \text{ if \eqref{eq:condexfo} holds,}
\\
S(v_i, r_i), & \text{ otherwise.} 
\end{cases}
\ee
The forward Euler method is used for the time discretization.
Furthermore, we emphasize that, if \eqref{eq:condexfo} is not satisfied, then $v_i$ cannot be the point value of a stationary solution defined in the computational domain, so that the use of the standard reconstruction does not destroy the well-balanced property of the method, since in this case there is no stationary solution to preserve.


\subsection{Second-order method}

Let us suppose again that the midpoint rule is used to compute the cell averages and that the minmod reconstruction operator is considered: see \cite{VanLeer}. The stationary solution $v^*_i$  selected at the first stage of the reconstruction procedure is again \eqref{eq:stadystate1} with $\tilde K_i$ given by \eqref{tildeK}. In order to compute the reconstructions, this
stationary solution has  to be computed at the points $r_{i-1}$, $r_{i -\frac{1}{2}}$, $r_{i+\frac{1}{2}}$, $r_{i+1}$ so that the following condition has to be satisfied $r_{i+1} \in
D_{\widetilde{K}_i}$, i.e.
\bel{eq:condexso}
\widetilde K_{i}^{2} \leq 1 \text{ or } \left(\widetilde K_{i}^{2}>1 \text{ and } r_{i+1} \leq  \frac{2M\widetilde K_{i}^{2}}{\widetilde K_{i}^{2}-1}\right).
\ee
If this condition is satisfied, the following step of the reconstruction procedure consists in computing the fluctuations:
\be
\label{rec1}
\aligned 
& w_{i-1} = v_{i-1} - \sign(v_{i-1})\sqrt{1-\widetilde{K}_{i}^{2}\left(1-\frac{2M}{r_{i-1}}\right)},
\\
&w_{i} = v_{i} - \sign(v_{i})\sqrt{1-\widetilde{K}_{i}^{2}\left(1-\frac{2M}{r_{i}}\right)} = 0, 
\qquad\qquad
w_{i+1} = v_{i+1} - \sign(v_{i+1})\sqrt{1-\widetilde{K}_{i}^{2}\left(1-\frac{2M}{r_{i+1}}  \right)}.
\endaligned
\ee
Then the reconstruction is defined as 
\be
\mathbb{P}_{i}(r) = v_{i}^{*}(r) + \text{minmod}\left(\frac{w_{i+1}-w_{i}}{\Delta r}, \frac{w_{i+1}-w_{i-1}}{2\Delta r},\frac{w_{i}-w_{i-1}}{\Delta r}\right)(r-r_{i}),
\ee
where
\be
\text{minmod}(a,b,c) = \begin{cases}
\min\{a,b,c\}, & \text{  } \ \ a,b,c>0, \\
\max\{a,b,c\}, & \text{ } \ \ a,b,c <0, \\
0, & \text{ otherwise.}
\end{cases}
\ee
Once the well-balanced reconstruction operator has been computed, the form \eqref{eq:nummet2} of the numerical method is considered and the midpoint rule is used again to approximate the integral. Observe however that, in this case:
$$
\int_{r_{i-\frac{1}{2}}}^{r_{i+\frac{1}{2}}} \left( S(\mathbb{P}_i^t (r), r) -S(V_i^{t,*}(r), r) \right) \, dr \cong \Delta r  \left( S(\mathbb{P}_i^t (r_i), r_i) -S(V_i^{t,*}(r_i), r_i)  \right)= 0.
$$
Therefore, the expression \eqref{eq:nummet3} reduces again to \eqref{eq:nummetfo} with
$
F_{i+\frac{1}{2}} = \mathbb{F}(v_{i+\frac{1}{2}}^{t, -}, v_{i+\frac{1}{2}}^{t, +}, r_{i+\frac{1}{2}}).
$

If \eqref{eq:condexso} is not satisfied, then the standard MUSCL reconstruction is applied, namely 
\be
\mathbb{Q}_{i}(r) = v_i + \text{minmod}\left(\frac{v_{i+1}-v_{i}}{\Delta r}, \frac{v_{i+1}-v_{i-1}}{2\Delta r},\frac{v_{i}-v_{i-1}}{\Delta r}\right)(r-r_{i}).
\ee
The expression of the numerical method is given again by \eqref{eq:semi_disc-merge}-\eqref{eq:source-merge} with the difference that the second-order reconstructions are used now to compute
the numerical fluxes. The TVDRK2 method \cite{Gottlieb-S} is used in order to discretize the equations in time. 

Observe that, according to the  well-balanced reconstruction procedure described in the previous section, the fluctuations to be reconstructed should be in this case
$$
w_j = v_j -  v^*_i(r_j) = v_j -  \sign(v_{i})\sqrt{1-\widetilde{K}_{i}^{2}\left(1-\frac{2M}{r_{j}}\right)},
\qquad j=i-1,i,i+1,
$$
but  in \eqref{rec1} the sign of $v_i$ has been replaced by that of $v_j$. This modification allows one to preserve not only the continuous stationary solutions
solution but also the discontinuous stationary solutions of the family \eqref{steadydisc}.


\subsection{Third-order method}

In order to design a third-order method the \textit{CWENO} reconstruction of order 3 (see \cite{Cravero-S}, \cite{LPR}) will be considered and the two-point Gauss quadrature will be used to
compute the initial averages and the integrals coming  from the source term, namely 
$
v_{i}^{0} = \frac{1}{2} (v_{0}(r_{i,0})+v_{0}(r_{i,1})),
$
where
$$
r_{i,0} = r_{i-\frac{1}{2}} + \frac{\Delta r}{2}\left(-\sqrt{\frac{1}{3}}+1\right),\ \
r_{i,1} = r_{i-\frac{1}{2}} + \frac{\Delta r}{2}\left(\sqrt{\frac{1}{3}}+1\right).
$$
Therefore, at the first step of the reconstruction procedure one has to look for $K_{i}^{2}$ such that:
	\bel{eq:constant}
\frac{1}{2}\left(\sqrt{1-K_{i}^{2}\left(1-\frac{2M}{r_{i,0}}\right)}+ \sqrt{1-K_{i}^{2}\left(1-\frac{2M}{r_{i,1}}\right)}\right)=|v_{i}|.
	\ee
If we define the function
$$
g(x) = \frac{1}{2}\left(\sqrt{1-x\left(1-\frac{2M}{r_{i,0}}\right)}+ \sqrt{1-x\left(1-\frac{2M}{r_{i,1}}\right)}\right), \qquad x \geq 0,
$$
it is easily checked that $g$ is a positive decreasing function defined in the interval $[0, K_{i,c}^2]$ with
	$
	K_{i,c}^{2} = \big( 1-\displaystyle \frac{2M}{r_{i,1}}\big)^{-1}.
	$
 Therefore there are two possibilities:
	\begin{itemize}

		\item If $|v_{i}| \in [g(K_{i,c}^{2}), 1]$, there is a unique  $\widetilde K^{2}_{i}$ satisfying (\ref{eq:constant}).
		
                     \item In other case, \eqref{eq:constant} has no solution. 

             \end{itemize}

If \eqref{eq:constant} is satisfied, the corresponding stationary solution
		$
		v_{i}^{*}(r) = \sign(v_{i})\sqrt{1-\widetilde K_{i}^{2}\Big( 1-\frac{2M}{r} \Big)}
		$
has to be computed in the points $\{r_{i-1,0}, r_{i-1,1}, r_{i+1,0}, r_{i+1,1}\}$ in the reconstruction procedure. Therefore, these points have to be in the interval of definition of $v^*_i$,
and this happens if $r_{i+1,1} \in D_{\widetilde K_i}$. Therefore, the well-balanced reconstruction can be computed only if the following condition is satisfied:
\bel{eq:condexto}
|v_{i}| \in [g(K_{i,c}^{2}), 1] \text{ and } \left( \widetilde K_{i}^{2} \leq 1 \text{ or } \left(\widetilde K_{i}^{2}>1 \text{ and } r_{i+1,1} \leq  \frac{2M\widetilde K_{i}^{2}}{\widetilde K_{i}^{2}-1}\right)\right).
\ee

If this condition is satisfied, the fluctuations can be then computed:
$$
w_j = v_j - \sign(v_{j})\frac{1}{2}\left(\sqrt{1-\widetilde K_{i}^{2}\left(1-\frac{2M}{r_{j,0}}\right)}+ 
\sqrt{1-\widetilde K_{i}^{2}\left(1-\frac{2M} {r_{j,1}}\right)}\right),\qquad j=i-1,i,i+1,
$$
and the well-balanced reconstruction is given by
$
\mathbb{P}_{i}(r) = v_{i}^{*}(r) + \mathbb{Q}_{i}(r; w_{i-1}, w_i, w_{i+1}),
$
where $\mathbb{Q}$ represents the CWENO approximation.
On the other hand, if \eqref{eq:condexto} is not satisfied, the standard CWENO reconstruction is applied, namely 
$
\mathbb{Q}_{i}(r) = \mathbb{Q}_{i}(r; v_{i-1}, v_i, v_{i+1}).
$
The expression of the semi-discrete method is finally \eqref{eq:semi_disc-merge} where the numerical fluxes are computed at the reconstructed states and
\be
S_i = \begin{cases} 
F(v_{i}^{*}(r_{i+\frac{1}{2}}), r_{i+\frac{1}{2}}) - F(v_{i}^{*}(r_{i-\frac{1}{2}}), r_{i-\frac{1}{2}})  
+ 
  \displaystyle \frac{\Delta r}{2}\sum_{j=0,1} \left( S(\mathbb{P}_{i}(r_{i,j}), r_{i,j})-S(v_{i}^{*}(r_{i,j}), r_{i,j})  \right) & \text{ if \eqref{eq:condexto} holds,}
\\
\displaystyle \frac{\Delta r}{2}\sum_{j=0,1} S(\mathbb{Q}_{i}(r_{i,j}), r_{i,j}),
&\text{ otherwise.}
\end{cases}
\ee
The TVDRK3 method of order 3 \cite{Gottlieb-S} will be used for the time discretization.


\subsection{Preserving the exact averages of the stationary solutions}

The three methods presented above can be modified to preserve the exact averages of the stationary solutions instead of its approximation computed with a quadrature formula. To do this, the problem to be solved at the first stage of the well-balanced reconstruction procedure is the following one: find $K_i^{2}$ such that:
	\bel{eq:step1Burger_explicit}
	\frac{1}{\Delta r}\int_{r_{i-\frac{1}{2}}}^{r_{i+\frac{1}{2}}} \sqrt{1-K_i^{2}\Big( 1-\frac{2M}{r} \Big)}dr = |v_{i}|.
	\ee
	If we define the function
	$$
	g(x) = \frac{1}{\Delta r}\int_{r_{i-\frac{1}{2}}}^{r_{i+\frac{1}{2}}} \sqrt{1-x\Big( 1-\frac{2M}{r} \Big)}dr = |v_{i}|,
	$$
	it can be easily checked that it is a decreasing function defined in  $[0, K_{e,i}^{2}]$ where
	$
	K_{e,i}^{2} = \left(1-2M / r_{i+\frac{1}{2}} \right)^{-1} 
	$
	and $g(0) = 1.$ Therefore, \eqref{eq:step1Burger_explicit} has a unique solution $\widetilde K_i^2$ if
	\bel{eq:condexex}
	|v_i| \leq g(K_{e,i}).
	\ee
	The explicit expression of $g$ can be obtained:
$
	g(x) = \frac{1}{\Delta r}\left(f(x, r_{i+\frac{1}{2}})- f(x, r_{i-\frac{1}{2}})\right),
$
	where 
	$$ 
f(x, r) = \begin{cases}
	r\sqrt{1-x\Big( 1-\frac{2M}{r} \Big)}
	+\displaystyle \frac{xM}{\sqrt{1-x}}\log \left(x(M-r) + r + \sqrt{1-x}r\sqrt{1-x\Big( 1-\frac{2M}{r} \Big)}\right),
	 & \text{$0 \leq x < 1,$}\\
	2r\sqrt{\displaystyle \frac{2M}{r}},
	 & \text{ $x = 1$, }\\
	r\sqrt{1-x}\Big( 1-\frac{2M}{r} \Big)
	-\displaystyle\frac{2xM}{\sqrt{x-1}}\tan^{-1}\left(\frac{\sqrt{1-x\Big( 1-\frac{2M}{r} \Big)}}{\sqrt{x-1}}\right),
	 & \text{ $x>1$.}
	\end{cases}
	$$
is a primitive function of $\sqrt{1-x\Big( 1-\frac{2M}{r} \Big)}$.	Therefore  $g(K_{e,i})$ can be explicitly computed. 

{{{{}}}{The well-balanced reconstruction can thus be computed if \eqref{eq:condexex} is satisfied and  the cells of the stencil $\mathcal{S}_i$ are contained in the domain of definition $D_{\widetilde K_i}$ of the corresponding stationary solution. Otherwise, the standard reconstruction is applied. } }The expression of the numerical methods is the same as the ones above.


\section{Burgers-Schwarzschild model: a numerical study}
\label{section--4}

\subsection{Preliminaries}

In this section several numerical tests are applied to check the performance of the well-balanced numerical methods introduced in the previous section. We consider the spatial  interval $[2M, L]$ with $M=1$ and $L=4$, a 256-point uniform mesh and the CFL number is set to 0.5. At $r=2M$ we impose $F_{-\frac{1}{2}}=0$ as boundary condition since $\Big( 1-\frac{2M}{r} \Big)= 0$. At $r=L$ we will use a transmissive boundary condition using ghost-cells if the initial condition is not a stationary solution; otherwise,  the stationary solution is imposed in the ghost-cells. The following numerical flux will be used:
$$
F_{i+\frac{1}{2}} = \mathbb{F}( v_{i}, v_{i+1},r_{i+\frac{1}{2}}) = \left(1-\frac{2M}{r_{i+\frac{1}{2}}}\right)\frac{q^{2}(0;v_{i}, v_{i+1})-1}{2},
$$
where $q(\cdot; v_L, v_R)$ is the self-similar solution of the Riemann problem for the standard Burgers equation with the initial condition 
$$
v_{0}(r)=\left\{\begin{array}{l}
v_{L},
 \quad \ r<{0},\\
v_{R}, \quad \ \ r>{0}.
\end{array}\right.
$$ 
In order to check the relevance of the well-balanced property, these methods will be compared with those  based on the same numerical fluxes and the standard first-, second-, or third-order reconstructions.


\subsection{Stationary solutions}

\paragraph{{{\red Positive stationary solution}}}
 
\begin{figure}[h]
	\begin{subfigure}[h]{0.32\textwidth}
		\centering
		\includegraphics[width=1\linewidth]{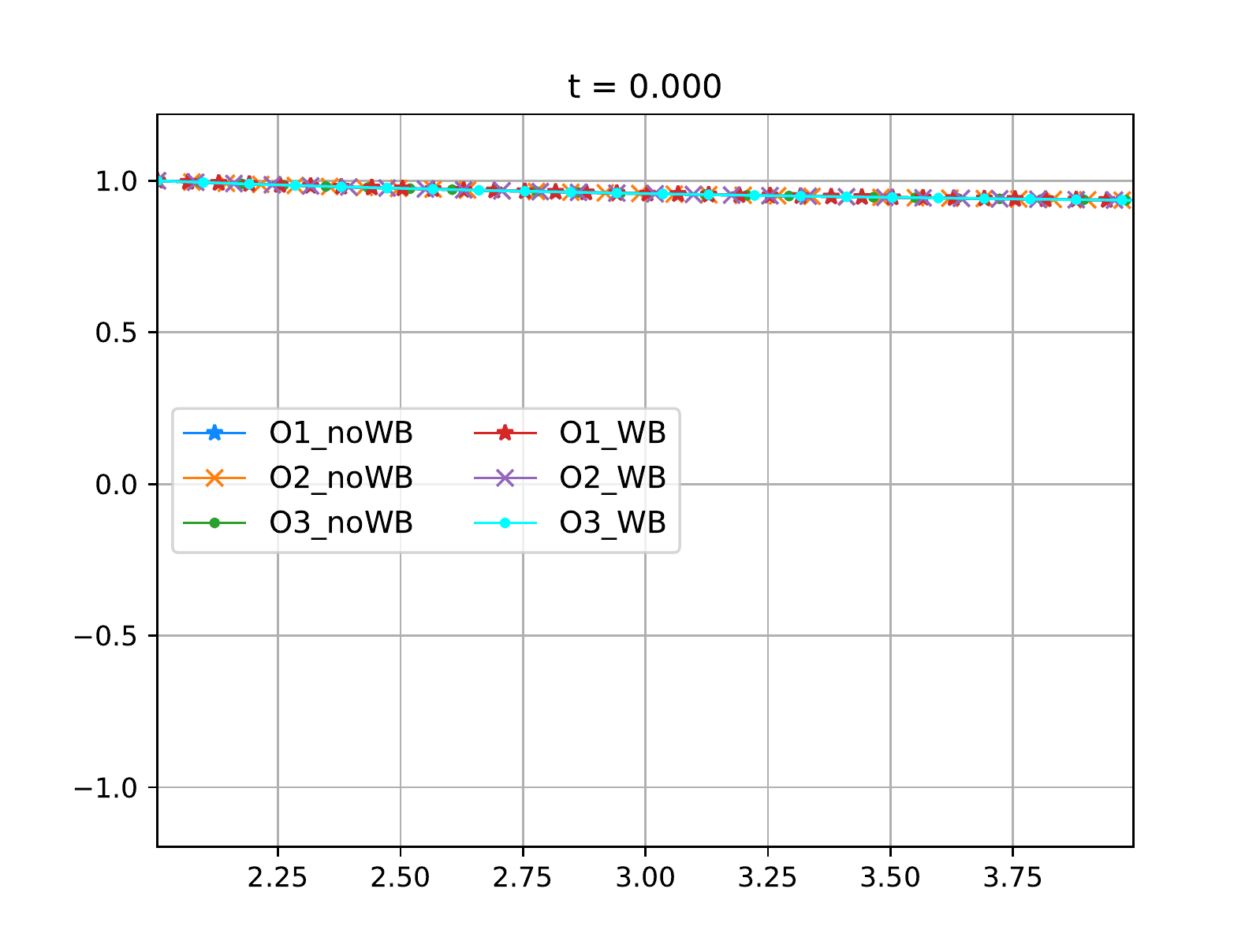}
		\label{fig:ko1_ko2_ko3_WB_vs_noWB_testWB1_t_0}
	\end{subfigure}
	\begin{subfigure}[h]{0.32\textwidth}
		\centering
		\includegraphics[width=1\linewidth]{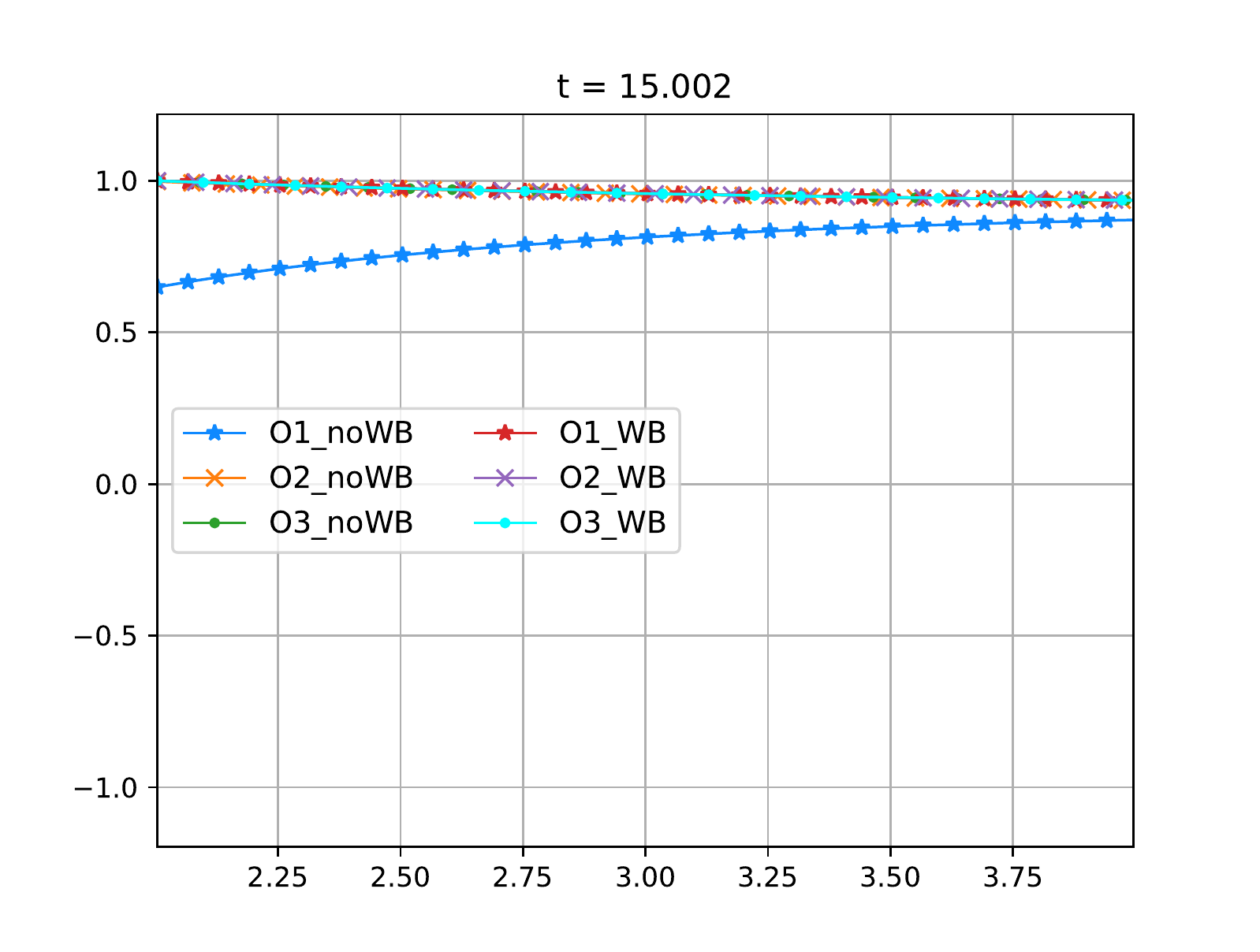}
		\label{fig:ko1_ko2_ko3_WB_vs_noWB_testWB1_t_15}
	\end{subfigure}
	\begin{subfigure}[h]{0.32\textwidth}
		\centering
		\includegraphics[width=1\linewidth]{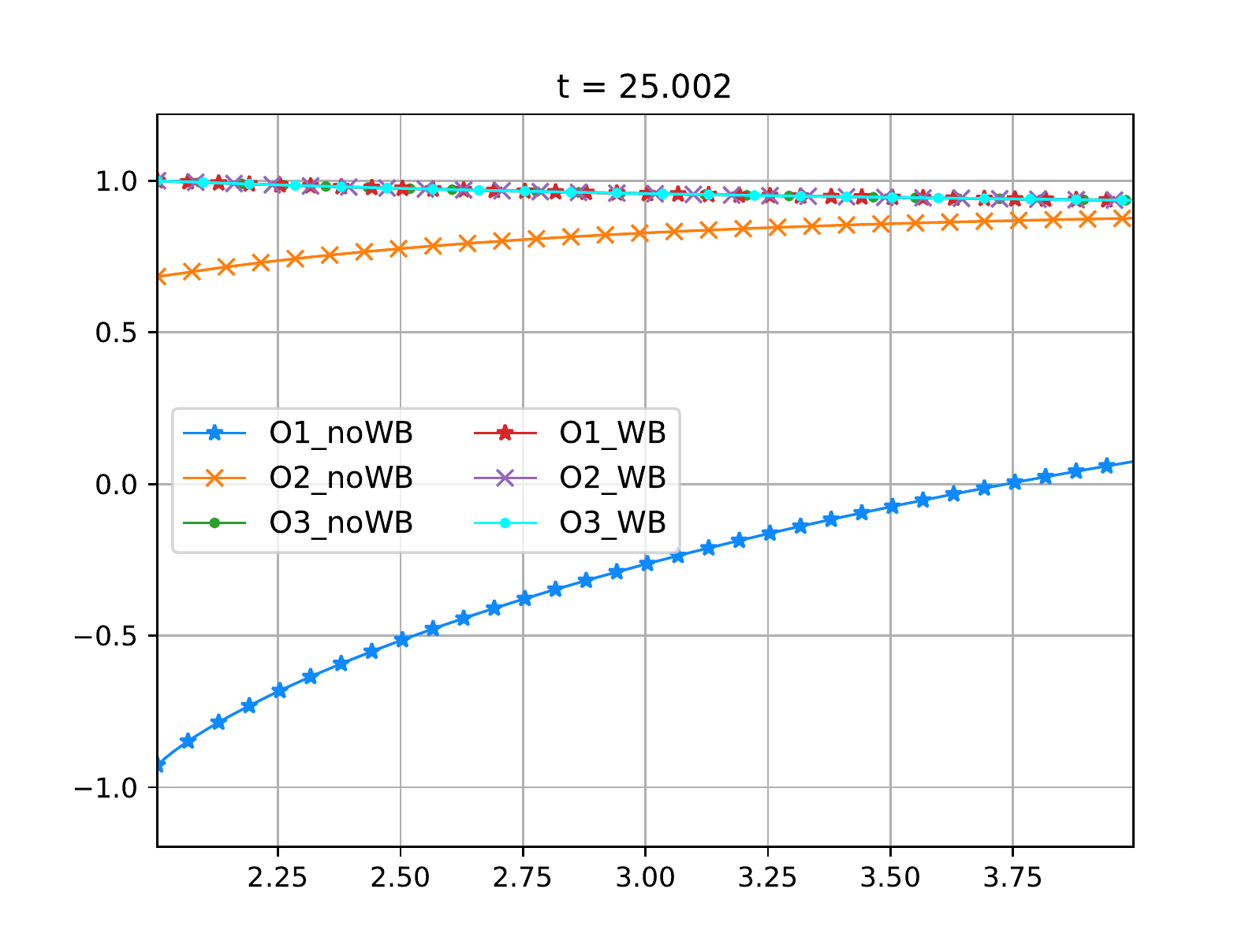}
		\label{fig:ko1_ko2_ko3_WB_vs_noWB_testWB1_t_25}
	\end{subfigure}
	\begin{subfigure}[h]{0.32\textwidth}
		\centering
		\includegraphics[width=1\linewidth]{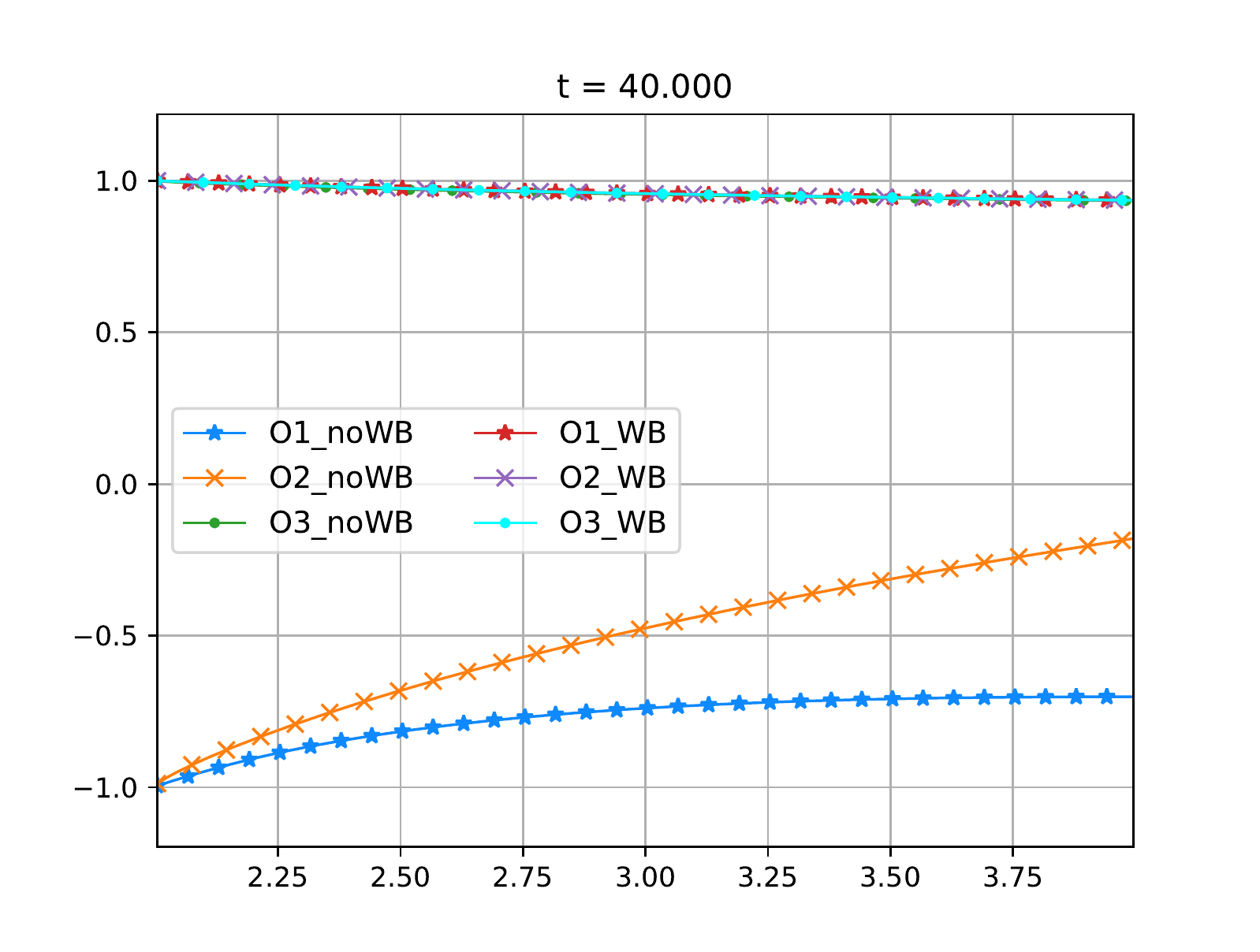}
		\label{fig:ko1_ko2_ko3_WB_vs_noWB_testWB1_t_40}
	\end{subfigure}
\begin{subfigure}[h]{0.32\textwidth}
	\centering
	\includegraphics[width=1\linewidth]{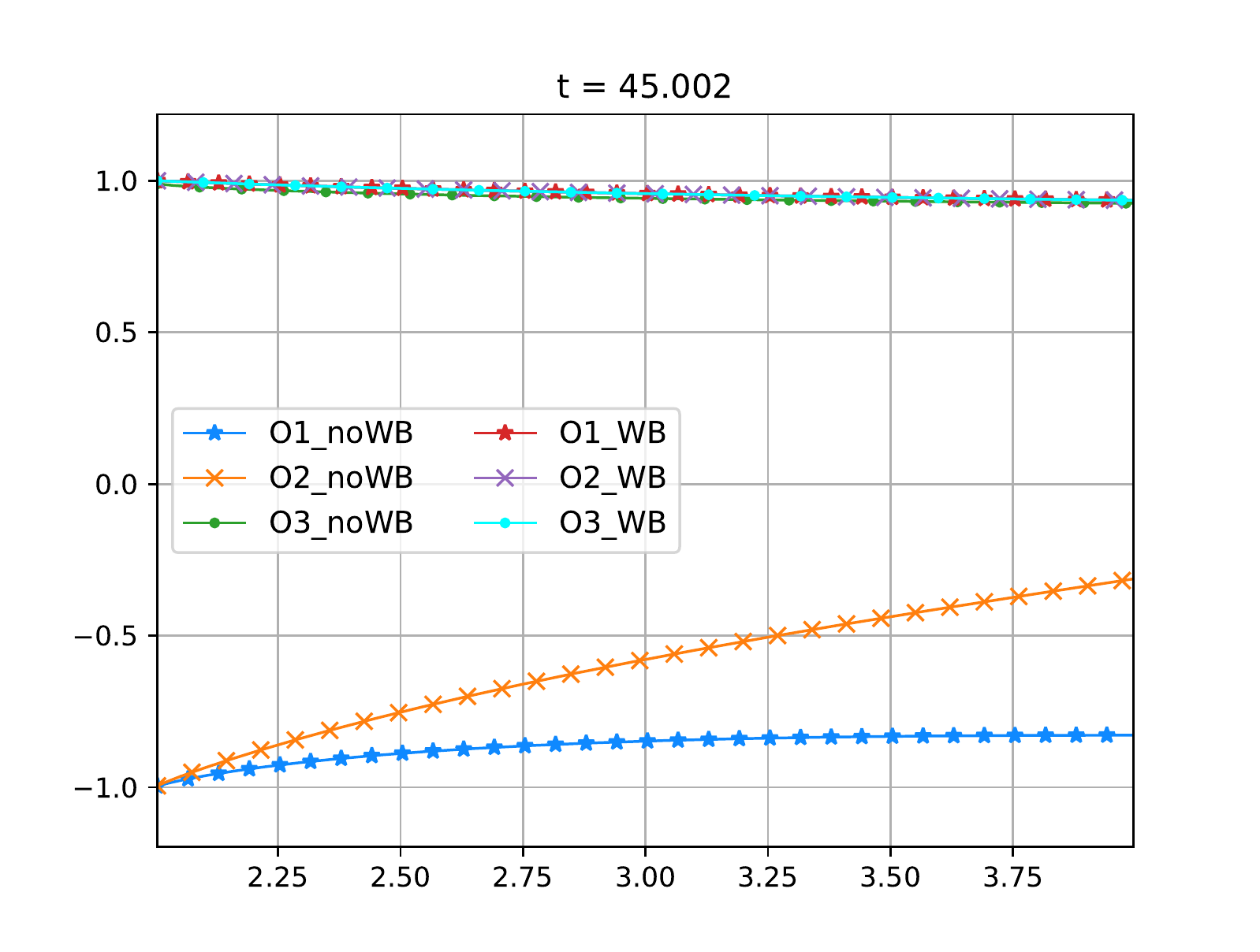}
	\label{fig:ko1_ko2_ko3_WB_vs_noWB_testWB1_t_45}
\end{subfigure}
\begin{subfigure}[h]{0.32\textwidth}
	\centering
	\includegraphics[width=1\linewidth]{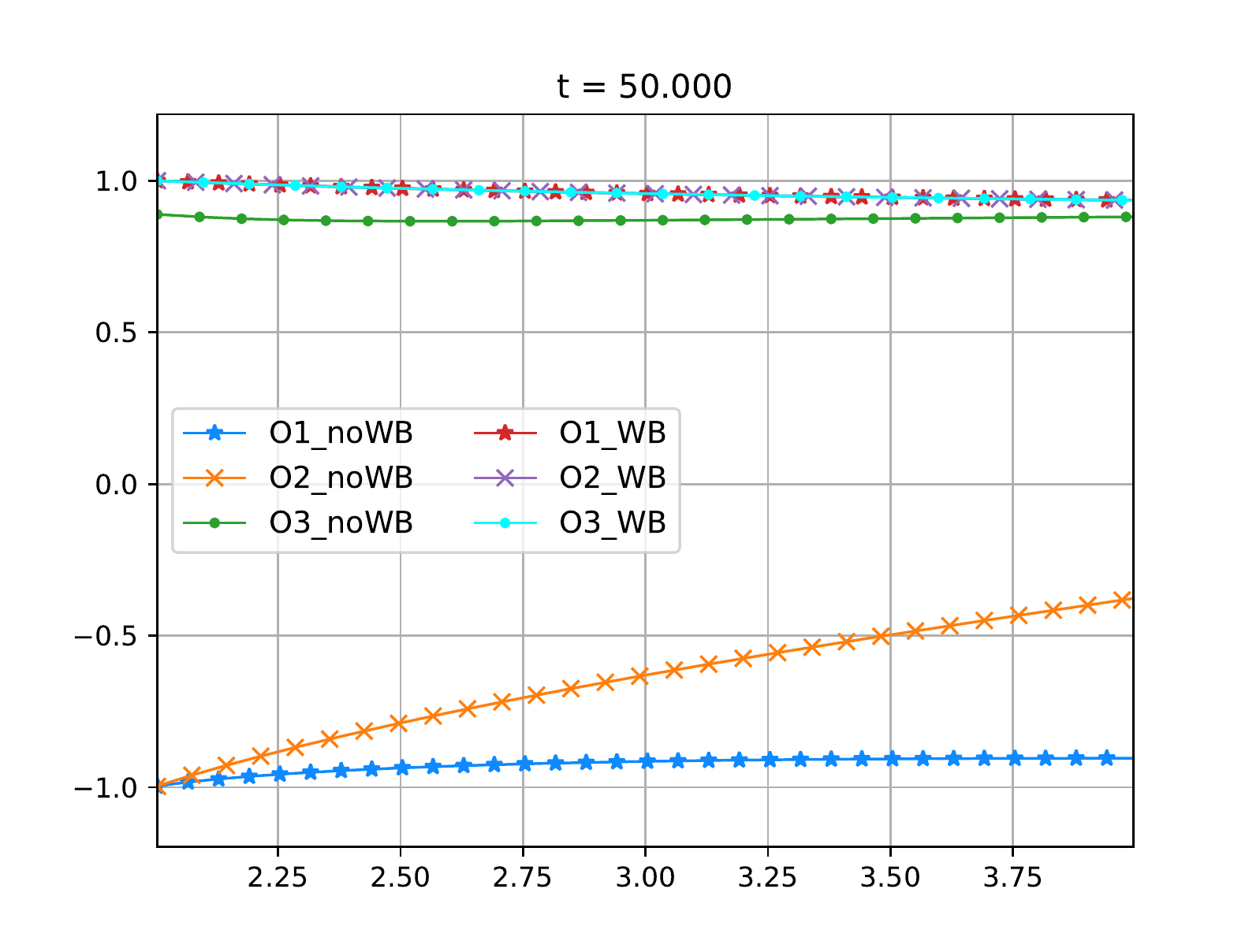}
	\label{fig:ko1_ko2_ko3_WB_vs_noWB_testWB1_t_50}
\end{subfigure}
	\caption{Burgers-Schwarzschild model with the initial condition \eqref{testB1}: first-, second-, and third-order well-balanced and not-well-balanced methods at various times.}
	\label{fig:comparison_ko1_ko2_ko3_WB_vs_noWB_testWB1}
\end{figure}

We consider the initial condition
\bel{testB1}
v_{0}(r) = \sqrt{\displaystyle \frac{3}{4} + \frac{1}{2r}}
\ee
corresponding to the positive stationary solution with $K = \frac{1}{2}$.   Table \ref{tab:Error_TestB1} shows the error in $L^1$ norm between the initial condition and the numerical solution
at time $ t = 50$.
Figure \ref{fig:comparison_ko1_ko2_ko3_WB_vs_noWB_testWB1} compares the numerical solutions obtained with the well-balanced and the non-well-balanced methods: it can be seen how the latter are unable to capture the stationary solution. 
After a time that decreases with the order, the numerical solutions depart from the steady state.

\begin{table}[ht]
	\centering
	\begin{tabular}{|c|c|c|c|}
		\hline 
		Scheme (256 cells) & Error (1st) & Error (2nd) & Error (3rd) \\ 
		\hline 
		Well-balanced & 1.13E-14 & 8.72Ee-17 & 7.22E-14 \\ 
		\hline 
		Non well-balanced & 1.89 & 1.61 & 8.78E-02 \\ 
		\hline 
	\end{tabular} 
	\caption{Well-balanced versus non-well-balanced schemes: $L^{1}$ errors at  $t=50$ for the Burgers model with the initial condition (\ref{testB1}).}
\label{tab:Error_TestB1}
\end{table}

\paragraph{{{\red Negative stationary solution}}}

\begin{figure}[h]
	\begin{subfigure}[h]{0.32\textwidth}
		\centering
		\includegraphics[width=1\linewidth]{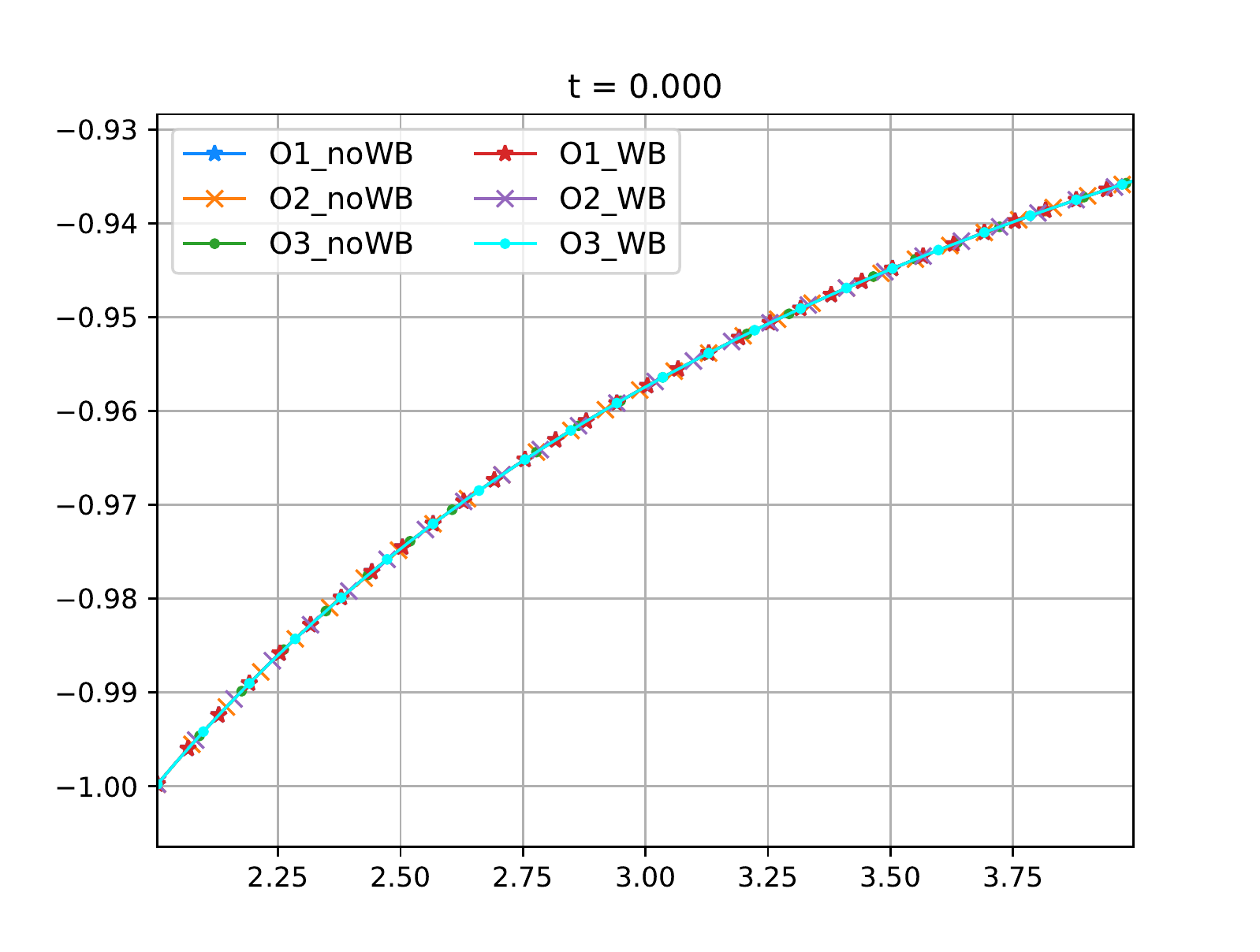}
		\label{fig:ko1_ko2_ko3_WB_vs_noWB_testWB2_t_0}
	\end{subfigure}
	\begin{subfigure}[h]{0.32\textwidth}
		\centering
		\includegraphics[width=1\linewidth]{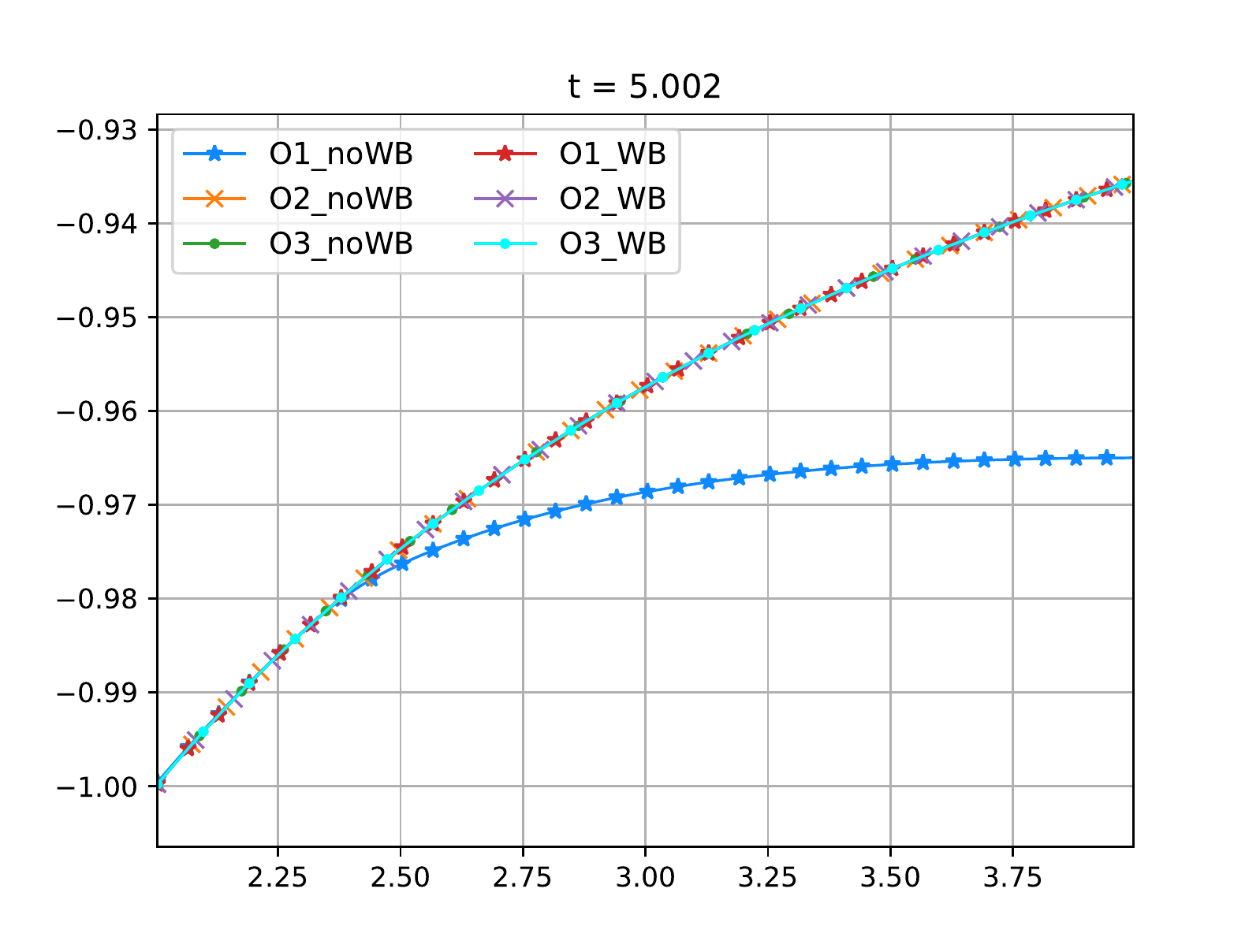}
		\label{fig:ko1_ko2_ko3_WB_vs_noWB_testWB2_t_5}
	\end{subfigure}
	\begin{subfigure}[h]{0.32\textwidth}
		\centering
		\includegraphics[width=1\linewidth]{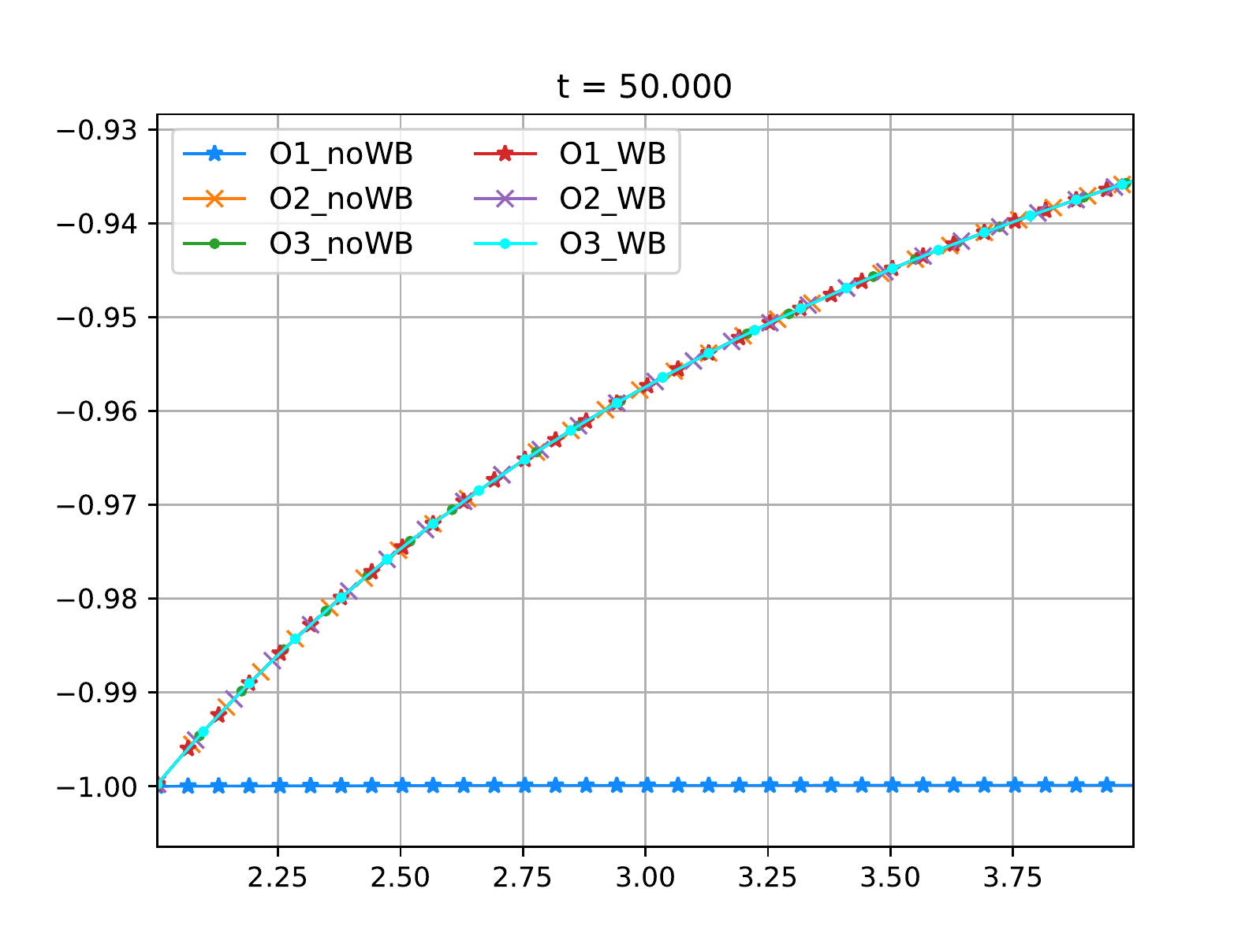}
		\label{fig:ko1_ko2_ko3_WB_vs_noWB_testWB2_t_50}
	\end{subfigure}
	\caption{Burgers-Schwarzschild model with the initial condition \eqref{testB2}: first-, second-, and third-order well-balanced and non-well-balanced methods at selected times.}
	\label{fig:comparison_ko1_ko2_ko3_WB_vs_noWB_testWB2}
\end{figure}

Let us consider now as initial condition the negative stationary condition corresponding to $K=\frac{1}{2}$:
\bel{testB2}
v_{0}(r) = -\sqrt{\displaystyle \frac{3}{4} + \frac{1}{2r}}.
\ee
{{{}}}{Figure \ref{fig:comparison_ko1_ko2_ko3_WB_vs_noWB_testWB2} shows the numerical solutions obtained with the different numerical methods. 
Observe that the scale of the vertical axis is not the same as the one in Figure \ref{fig:comparison_ko1_ko2_ko3_WB_vs_noWB_testWB1}: it has been changed  so that the difference between the numerical solutions can be better seen.} Table \ref{tab:Error_TestB2} shows the error in $L^1$ norm between the initial condition and the numerical solution
at time $ t = 50$.
According to Figure \ref{fig:comparison_ko1_ko2_ko3_WB_vs_noWB_testWB2} and Table \ref{tab:Error_TestB2} we need more time to see the differences between the well-balanced and non-well-balanced schemes of order 1, 2 and 3 but the errors are again much smaller with the well-balanced schemes for this test. {{{{}}}{In this case we need more time to see these differences since this negative stationary solution is close to the constant state $v(r) = -1$ where it seems that the non-well-balanced schemes converge.}}

\begin{table}[ht]
	\centering
	\begin{tabular}{|c|c|c|c|}
		\hline 
		Scheme (256 cells) & Error (1st) & Error (2nd) & Error (3rd) \\ 
		\hline 
		Well-balanced & 6.98E-16 & 1.24E-16 & 4.03E-16 \\ 
		\hline 
		Non well-balanced & 3.92E-02 & 3.20E-07 & 1.63E-10 \\ 
		\hline 
	\end{tabular} 
	\caption{Well-balanced versus non-well-balanced schemes:  $L^{1}$ errors at $t=50$ for the Burgers model with the initial condition (\ref{testB2}).}
\label{tab:Error_TestB2}
\end{table}


\paragraph{{{\red Discontinuous stationary entropy weak solution}}}

\begin{figure}[h]
	\begin{subfigure}[h]{0.32\textwidth}
		\centering
		\includegraphics[width=1\linewidth]{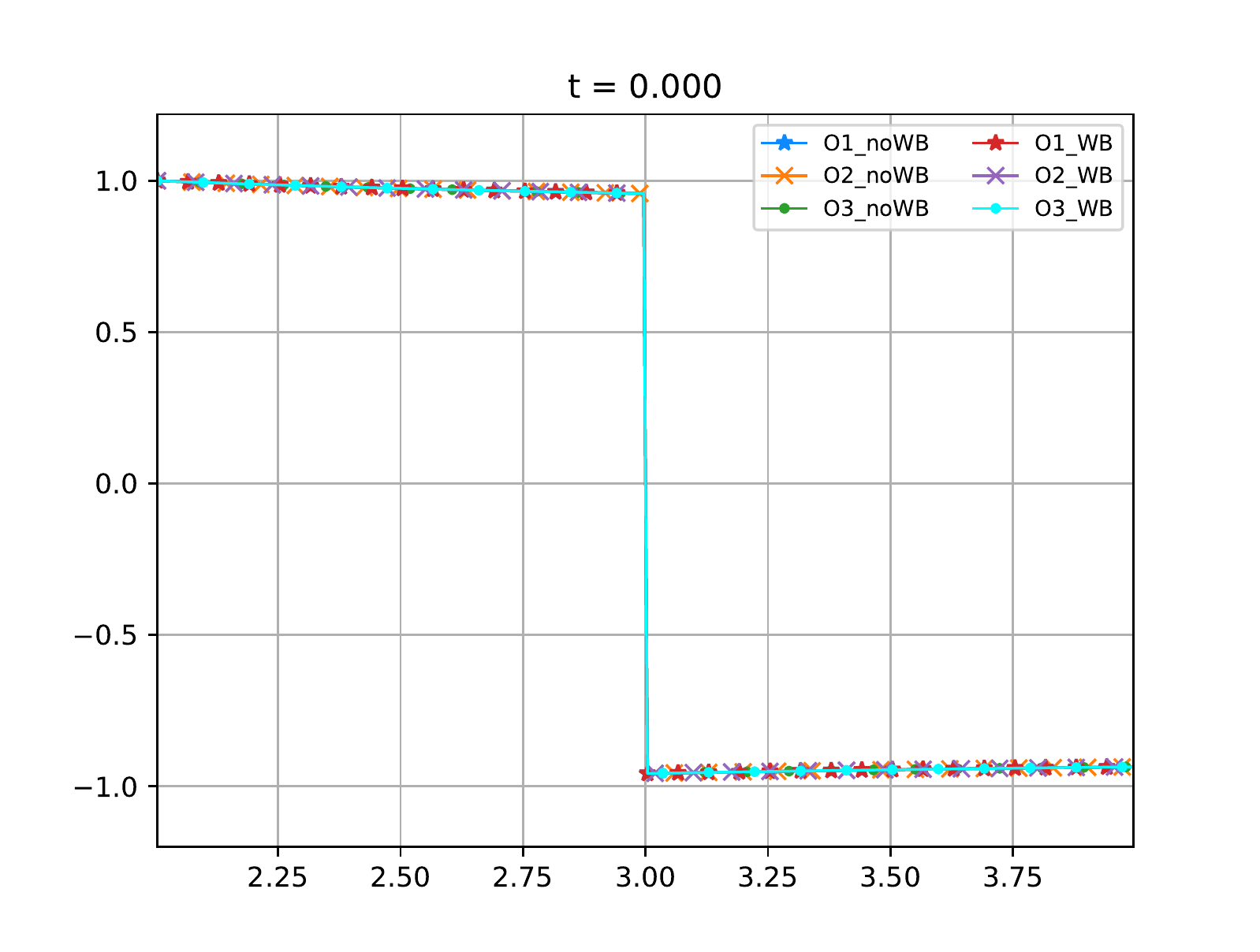}
		\label{fig:ko1_ko2_ko3_WB_vs_noWB_testWB3_t_0}
	\end{subfigure}
	\begin{subfigure}[h]{0.32\textwidth}
		\centering
		\includegraphics[width=1\linewidth]{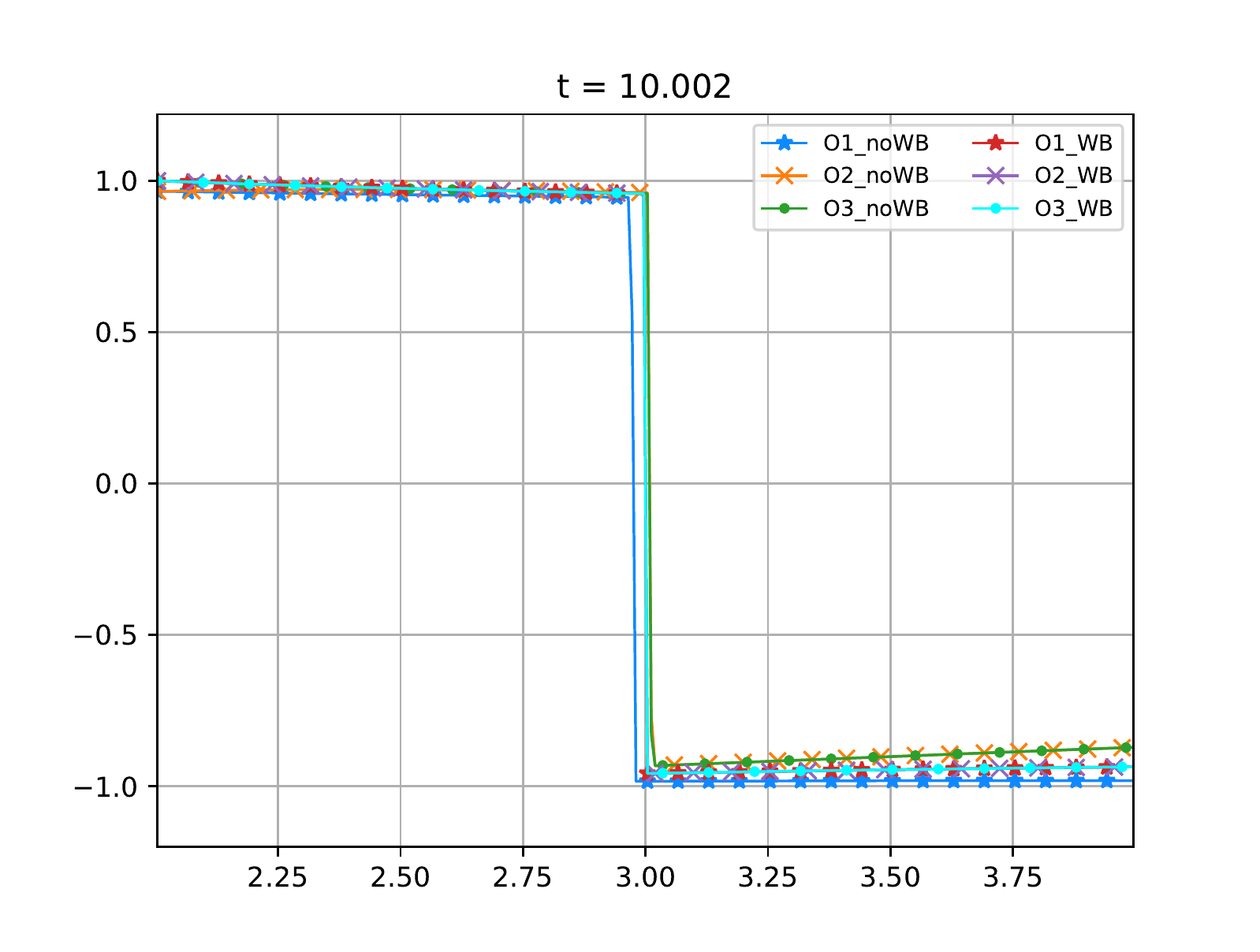}
		\label{fig:ko1_ko2_ko3_WB_vs_noWB_testWB3_t_10}
	\end{subfigure}
	\begin{subfigure}[h]{0.32\textwidth}
		\centering
		\includegraphics[width=1\linewidth]{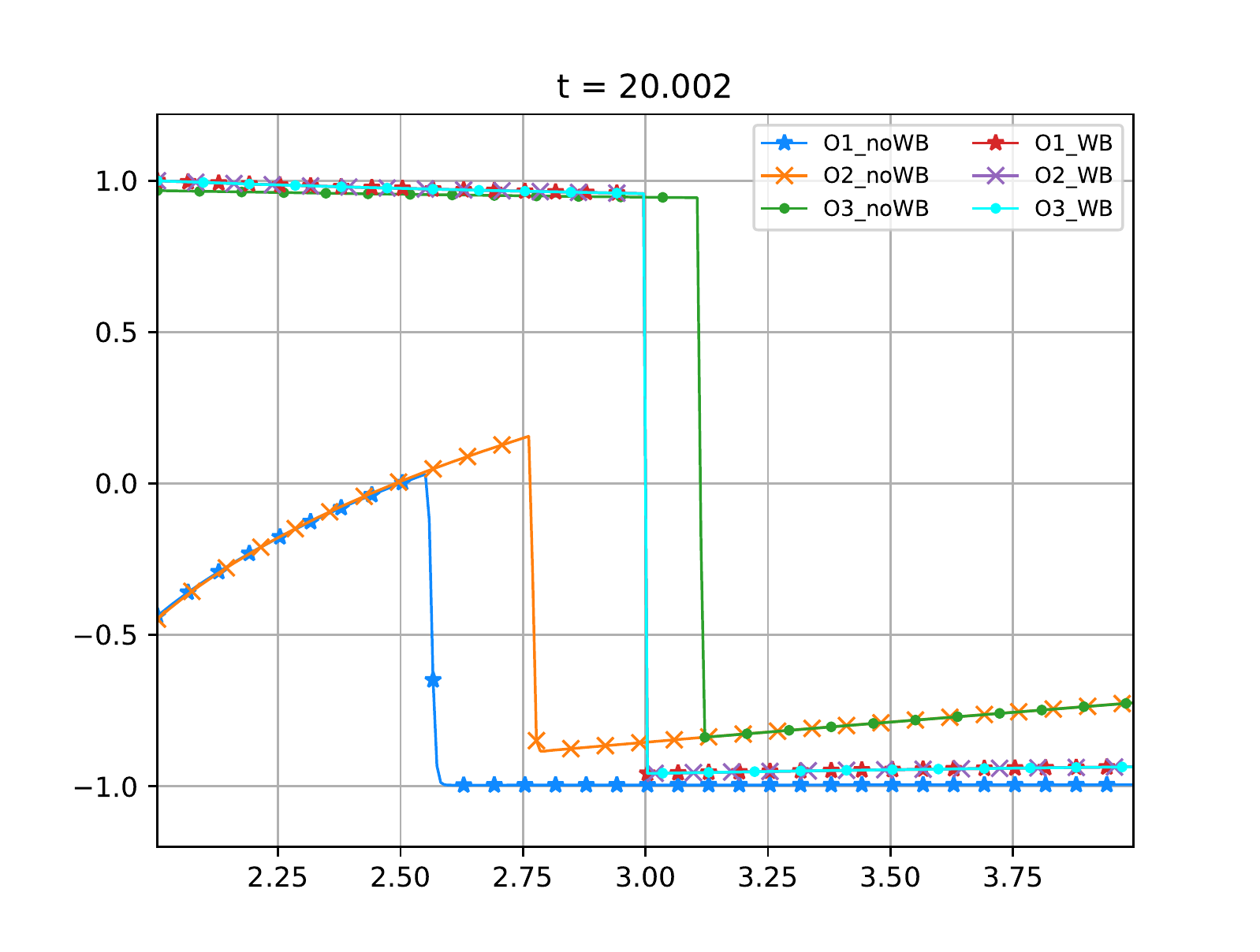}
		\label{fig:ko1_ko2_ko3_WB_vs_noWB_testWB3_t_20}
	\end{subfigure}
	\begin{subfigure}[h]{0.32\textwidth}
		\centering
		\includegraphics[width=1\linewidth]{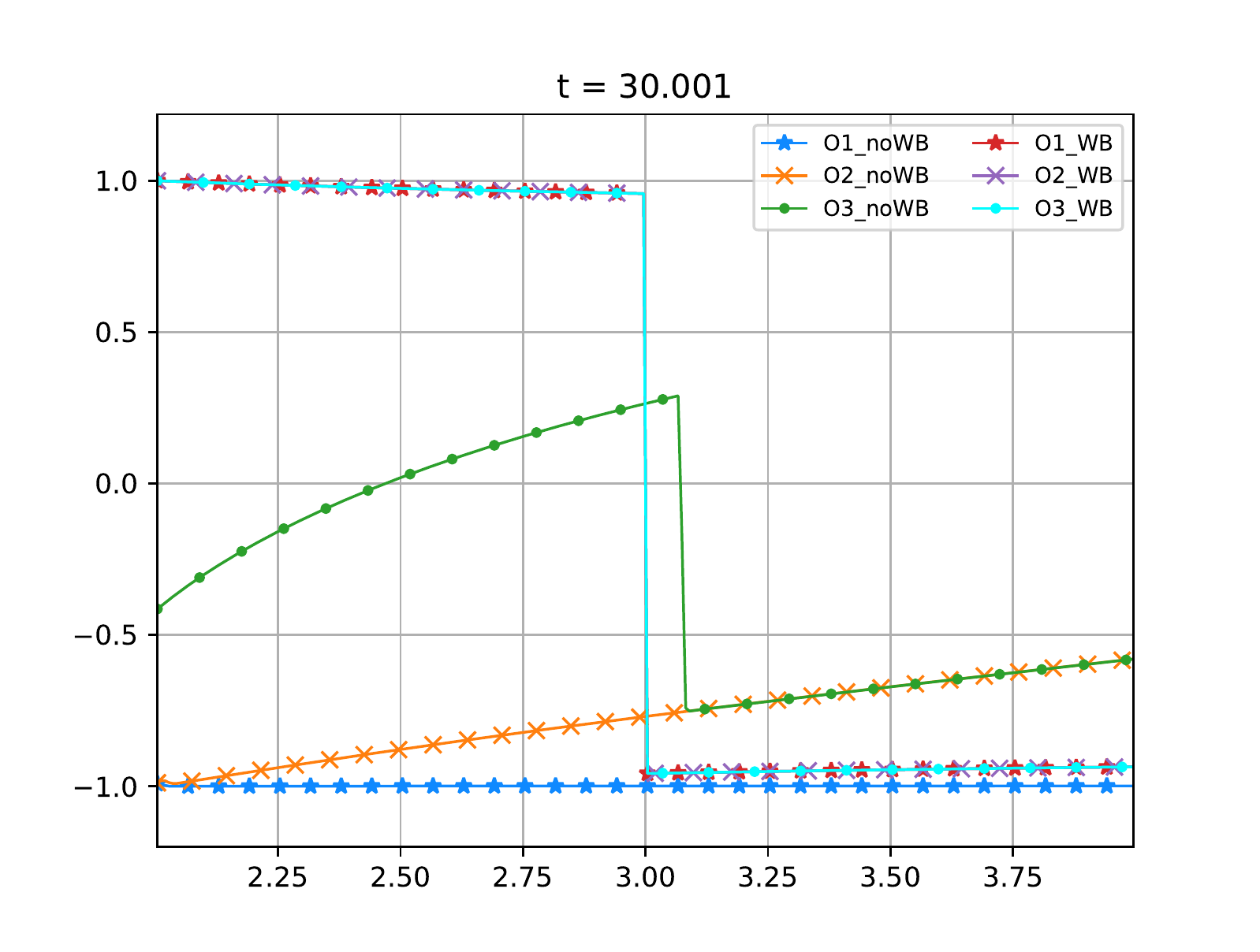}
		\label{fig:ko1_ko2_ko3_WB_vs_noWB_testWB3_t_30}
	\end{subfigure}
	\begin{subfigure}[h]{0.32\textwidth}
		\centering
		\includegraphics[width=1\linewidth]{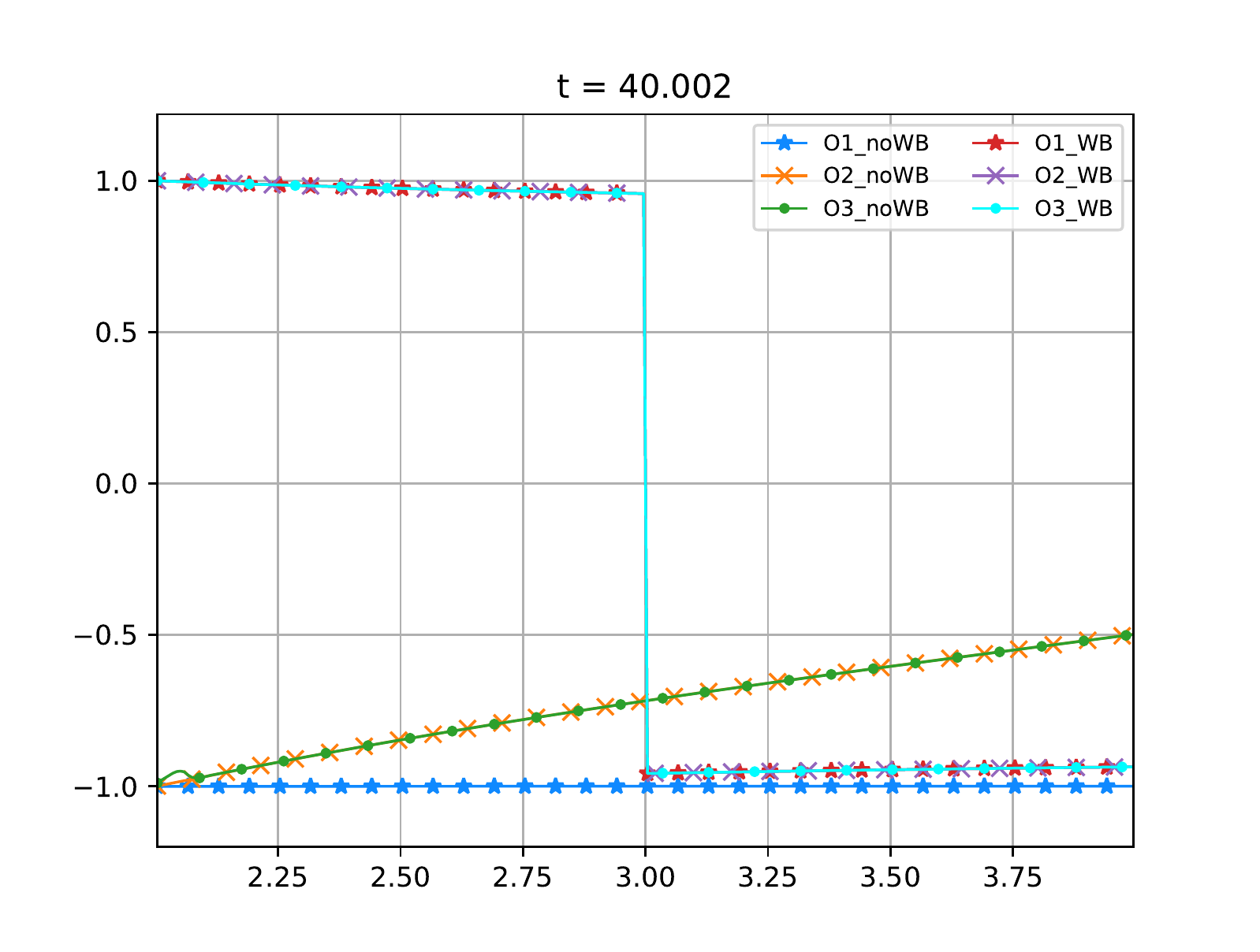}
		\label{fig:ko1_ko2_ko3_WB_vs_noWB_testWB3_t_40}
	\end{subfigure}
	\begin{subfigure}[h]{0.32\textwidth}
		\centering
		\includegraphics[width=1\linewidth]{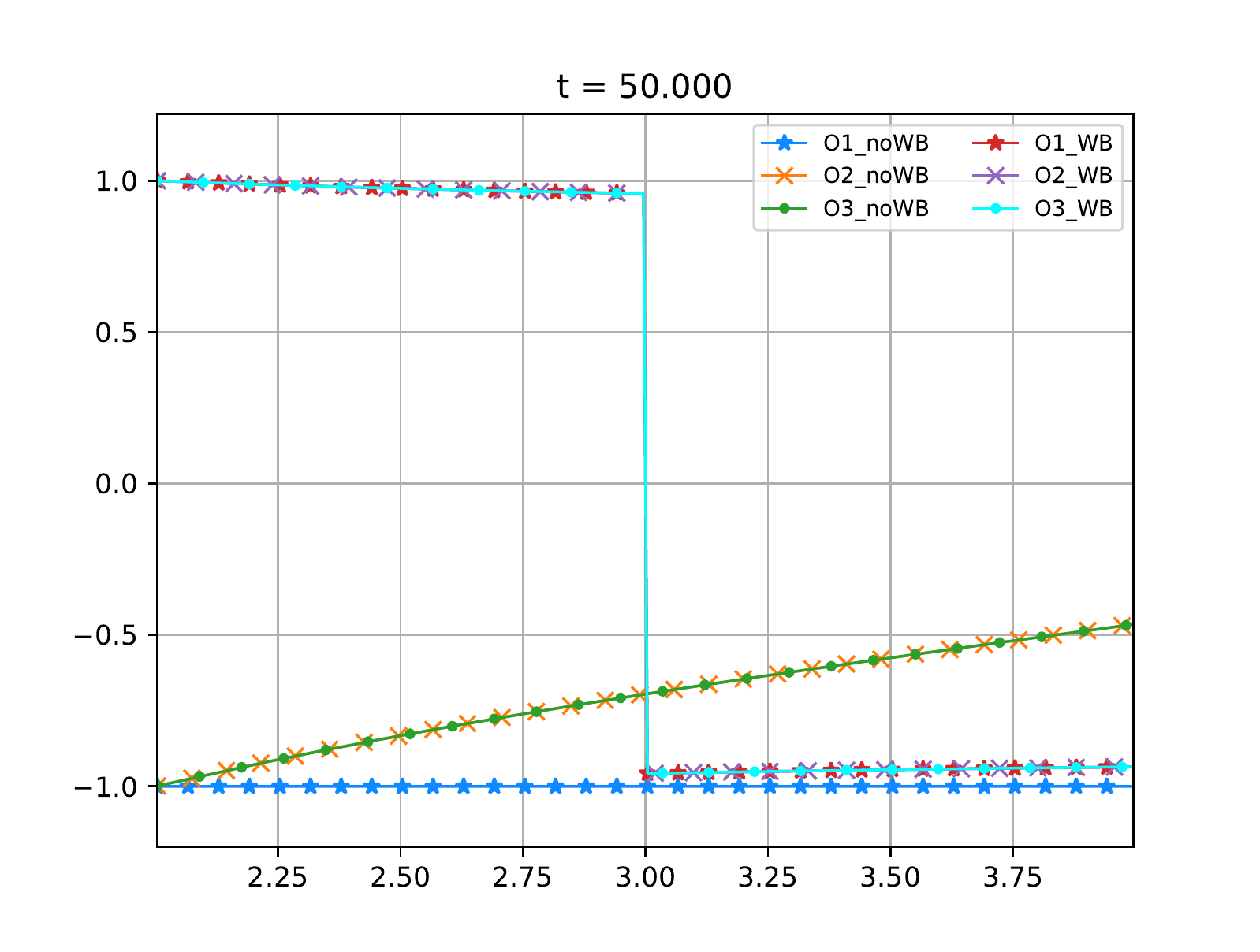}
		\label{fig:ko1_ko2_ko3_WB_vs_noWB_testWB3_t_50}
	\end{subfigure}
	\caption{Burgers-Schwarzschild model with the initial condition \eqref{testB3}: first-, second-, and third-order well-balanced and non-well-balanced methods at selected times.}
	\label{fig:comparison_ko1_ko2_ko3_WB_vs_noWB_testWB3}
\end{figure}

Let us consider finally the discontinuous initial condition
\be
v_{0}(r) = \begin{cases}\label{testB3}
	\sqrt{\displaystyle \frac{3}{4} + \frac{1}{2r}}, & \text{ $2<r<3$},\\
	-\sqrt{\displaystyle \frac{3}{4} + \frac{1}{2r}}, & \text{ otherwise},
\end{cases}
\ee
that is a stationary entropy weak solution of the family \eqref{steadydisc}.
Table \ref{tab:Error_TestB3} shows the error in $L^1$ norm between the initial condition and the numerical solution
at time $ t = 50$.
Figure \ref{fig:comparison_ko1_ko2_ko3_WB_vs_noWB_testWB3} shows the differences between the numerical solutions obtained with well-balanced and non-well-balanced methods: again the latter depart from the stationary solution at time that decrease with the order.  The numerical results obtained for the equation with initial conditions \eqref{testB1}, \eqref{testB2}, and \eqref{testB3} clearly  show the need of using well-balanced methods for this equation.

\begin{table}[ht]
	\centering
	\begin{tabular}{|c|c|c|c|}
		\hline 
		Scheme (256 cells) & Error (1st) & Error (2nd) & Error (3rd) \\ 
		\hline 
		Well-balanced & 8.68E-15 & 8.54E-17 & 7.90E-14 \\ 
		\hline 
		Non well-balanced & 1.02 & 1.09 & 1.09 \\ 
		\hline 
	\end{tabular} 
\caption{Well-balanced versus non-well-balanced schemes:  $L^{1}$ errors at $t=50$ for the Burgers model with the initial condition (\ref{testB3}).}
\label{tab:Error_TestB3}
\end{table}

\subsection{Moving shocks connecting two steady profiles}

\paragraph{{{\red Right-moving shock}}}

\begin{figure}[h]
	\begin{subfigure}[h]{0.5\textwidth}
		\centering
		\includegraphics[width=1\linewidth]{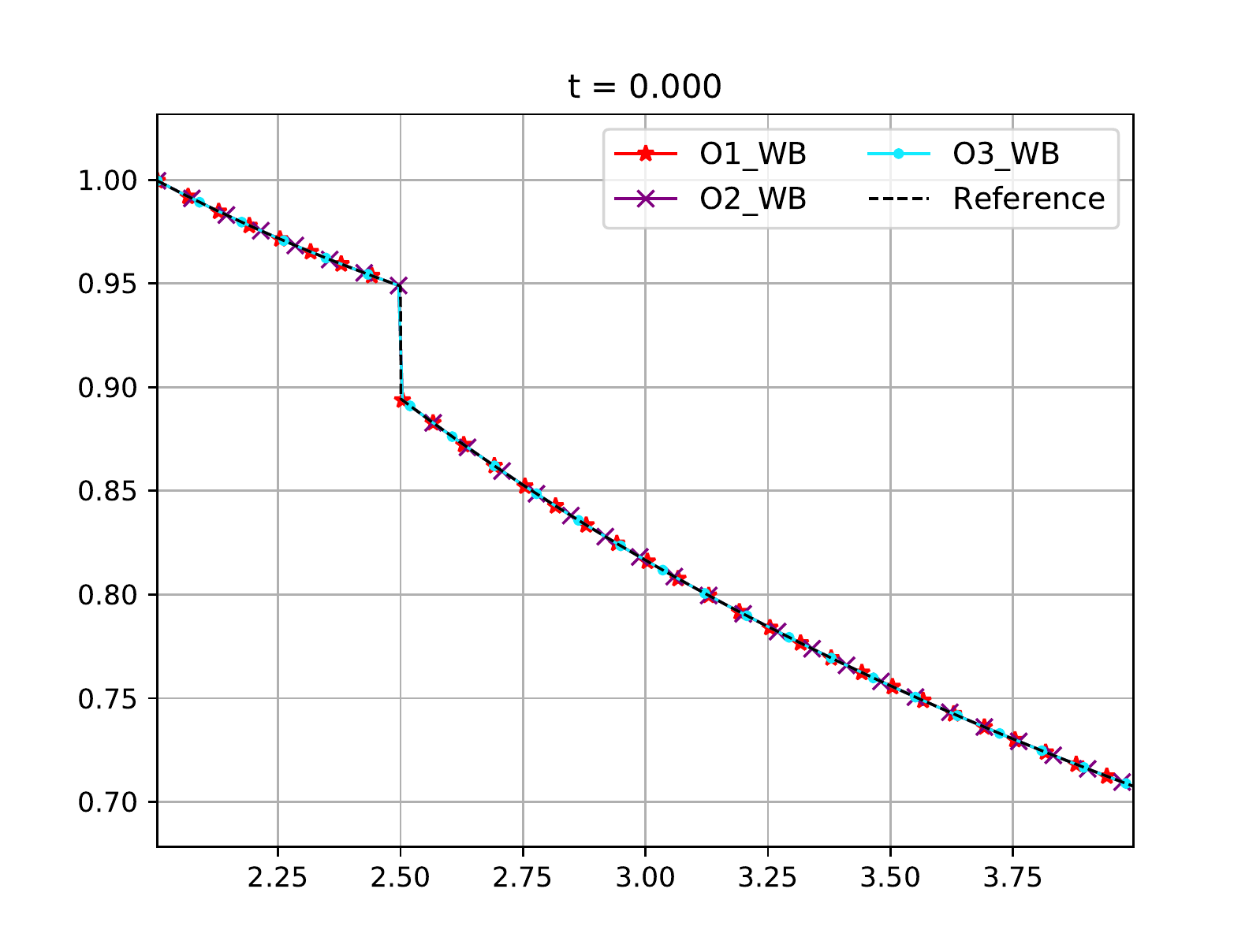}
		\label{fig:ko1_ko2_ko3_wb_testPositiveShock_with_reference_t_0}
	\end{subfigure}
	\begin{subfigure}[h]{0.5\textwidth}
		\centering
		\includegraphics[width=1\linewidth]{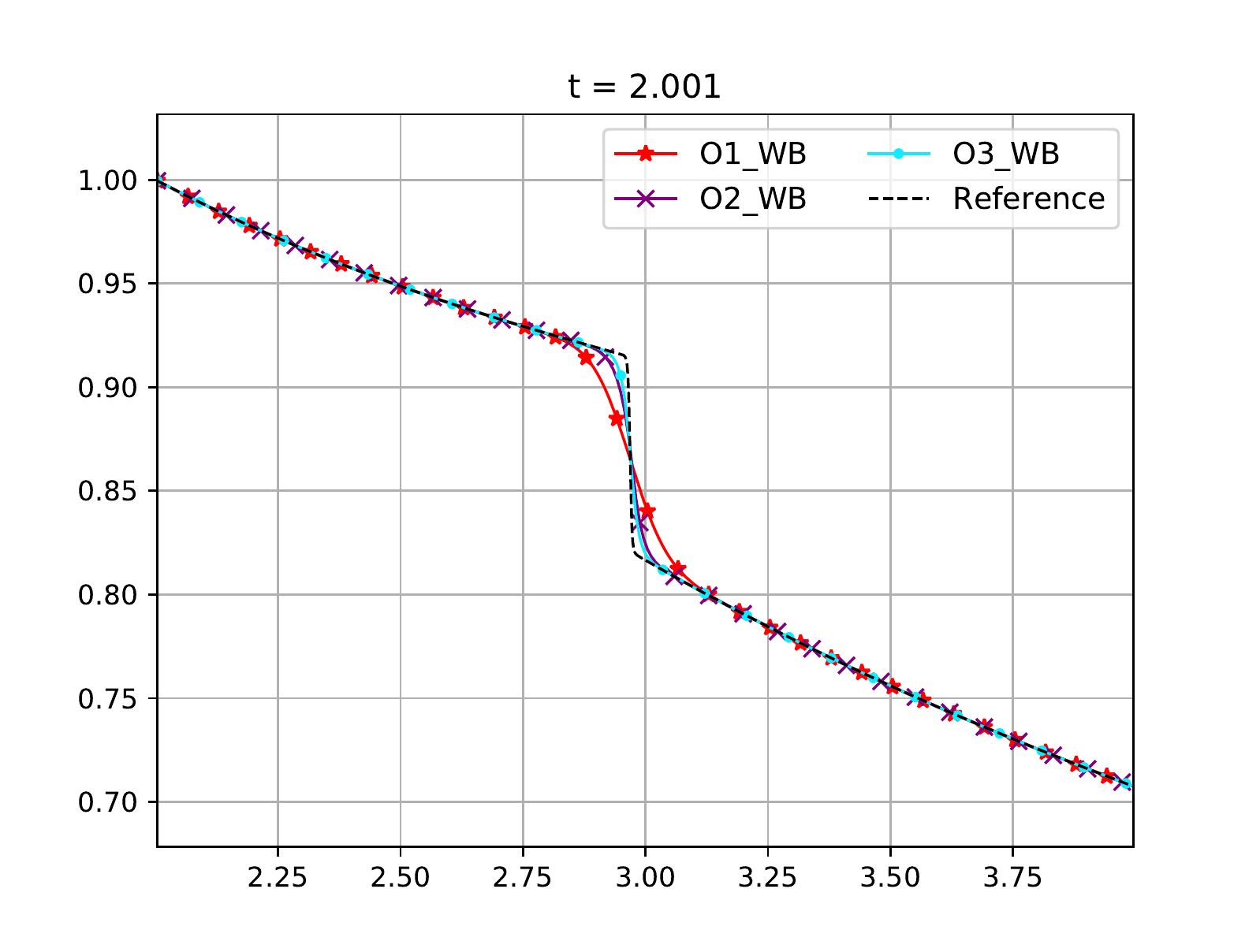}
		\label{fig:ko1_ko2_ko3_wb_testPositiveShock_with_reference_t_2}
	\end{subfigure}
	\begin{subfigure}[h]{1\textwidth}
		\centering
		\includegraphics[width=0.5\linewidth]{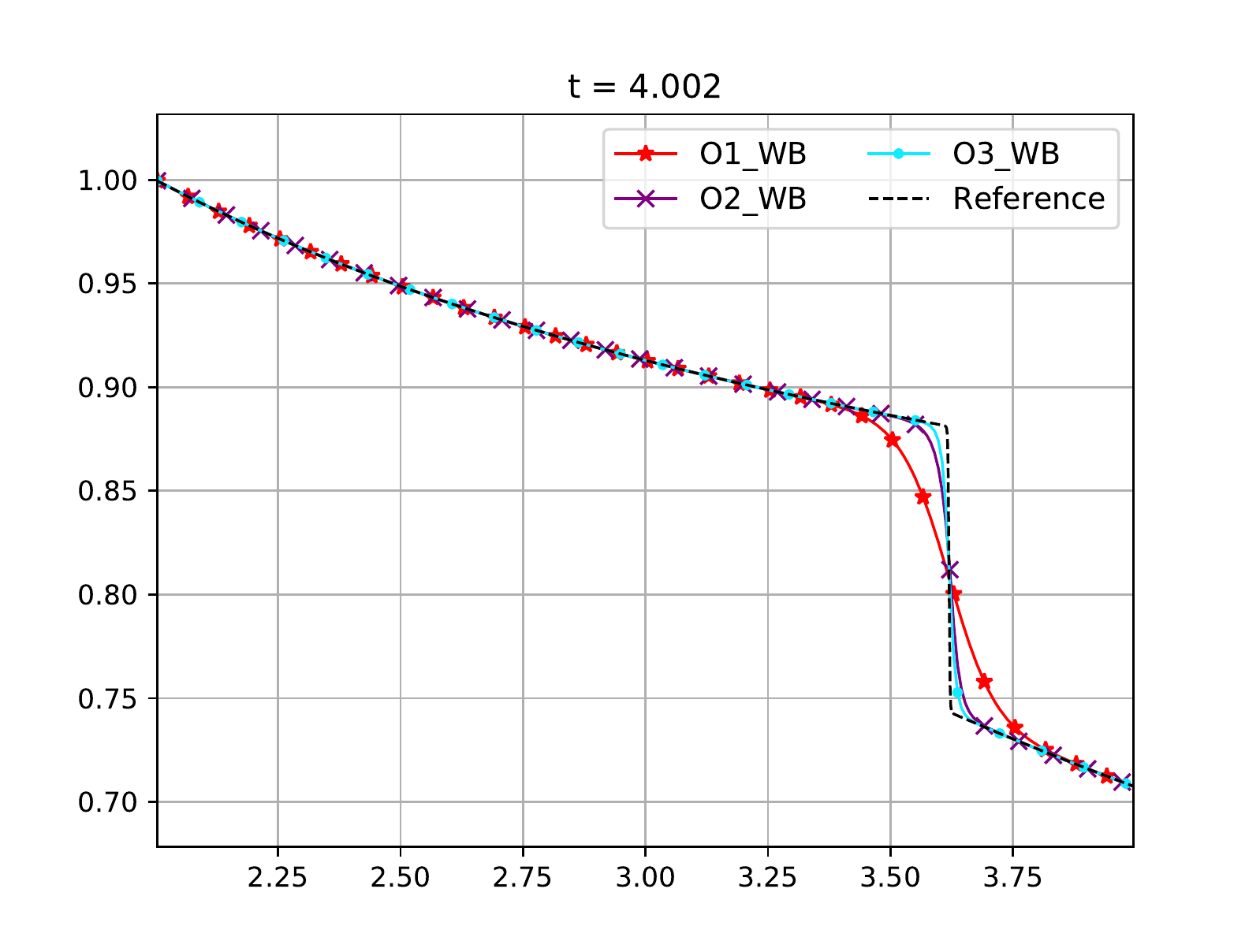}
		\label{fig:ko1_ko2_ko3_wb_testPositiveShock_with_reference_t_4}
	\end{subfigure}
	\caption{Burgers-Schwarzschild model with the initial condition \eqref{testB4}: first-, second-, and third-order well-balanced  methods at selected times.}
	\label{fig:comparison_ko1_ko2_k03_wb_testPositiveshock}
\end{figure}

We consider now the  initial condition
\bel{testB4}
v_{0}(r) = \begin{cases}
	\sqrt{\displaystyle \frac{1}{2} + \frac{1}{r}}, & \text{ $2<r<2.5$},\\
	\sqrt{\displaystyle \frac{2}{r}}, & \text{ otherwise}.
\end{cases}
\ee
The corresponding solution consists of a right-moving shock connecting two branches of stationary solutions. Figure \ref{fig:comparison_ko1_ko2_k03_wb_testPositiveshock} shows the numerical solutions obtained with the first-, second-, and third-order well-balanced methods and a reference solution computed with the first-order standard method using a mesh of 10000 cells. As  it can be seen, the well-balanced methods capture correctly the shock with a resolution that increases with the order as expected.


\paragraph{{{\red Left-moving shock}}}

\begin{figure}[h]
	\begin{subfigure}[h]{0.5\textwidth}
		\centering
		\includegraphics[width=1\linewidth]{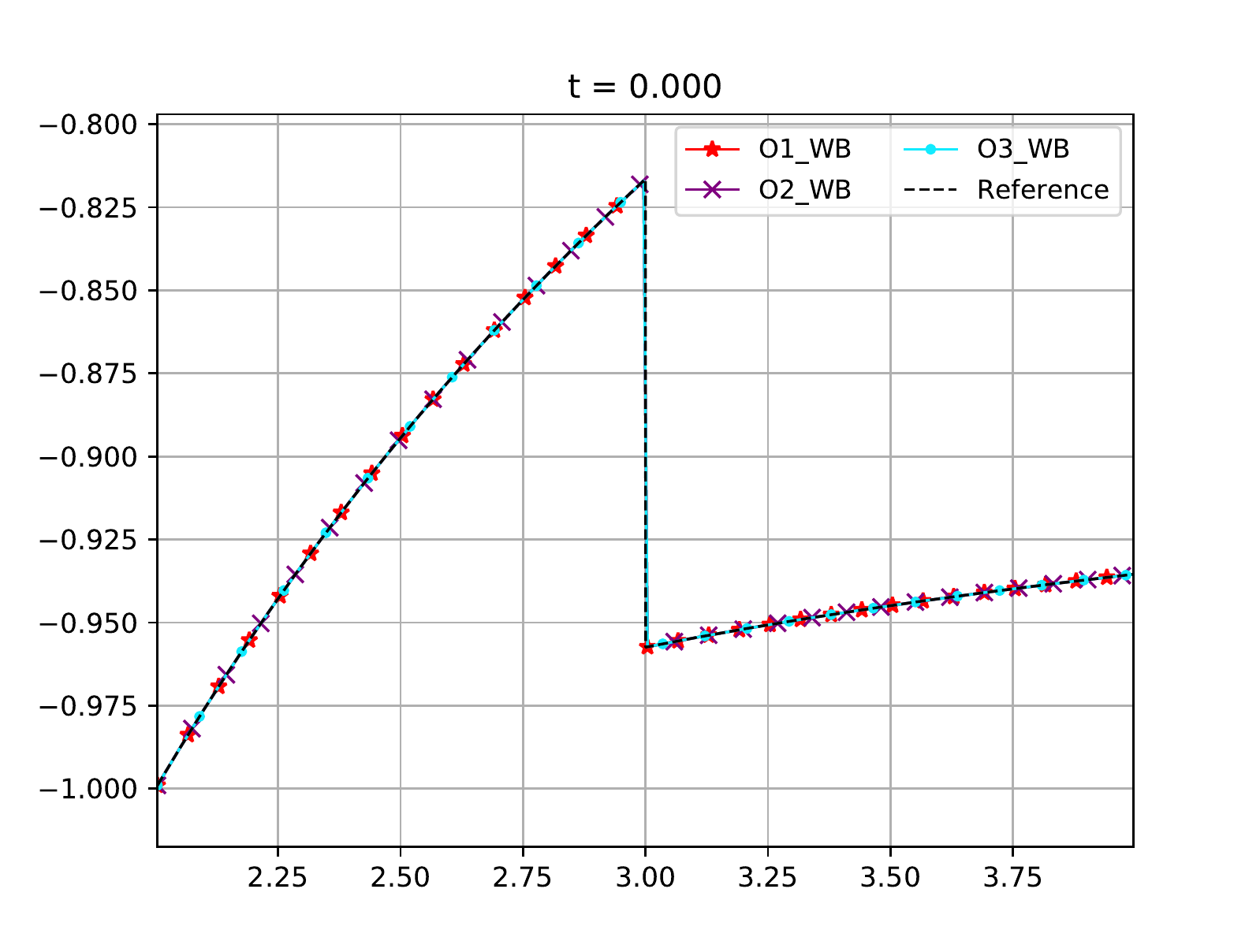}
		\label{fig:ko1_ko2_ko3_wb_testNegativeShock_t_0}
	\end{subfigure}
	\begin{subfigure}[h]{0.5\textwidth}
		\centering
		\includegraphics[width=1\linewidth]{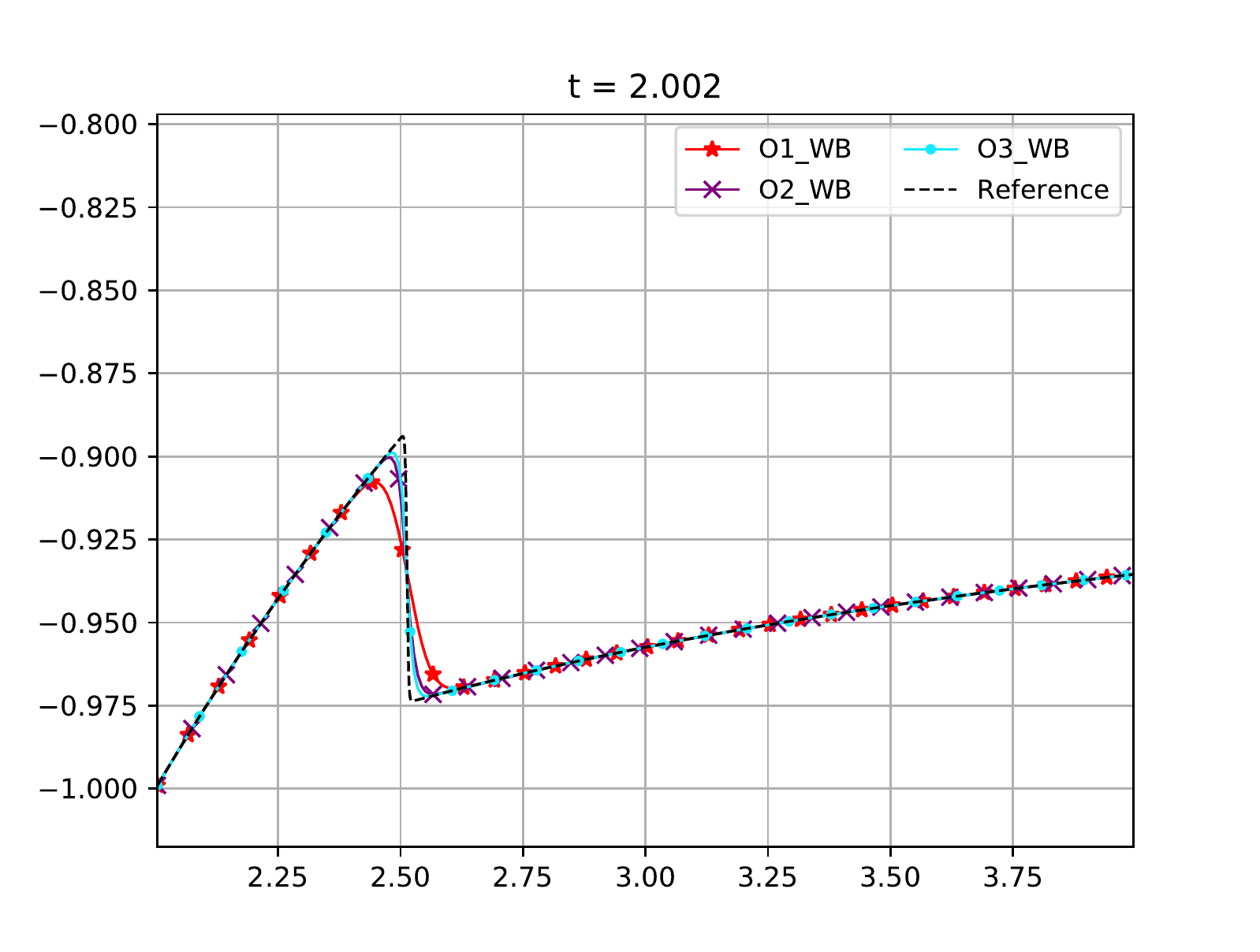}
		\label{fig:ko1_ko2_ko3_wb_testNegativeShock_t_2}
	\end{subfigure}
	\begin{subfigure}[h]{1\textwidth}
		\centering
		\includegraphics[width=0.5\linewidth]{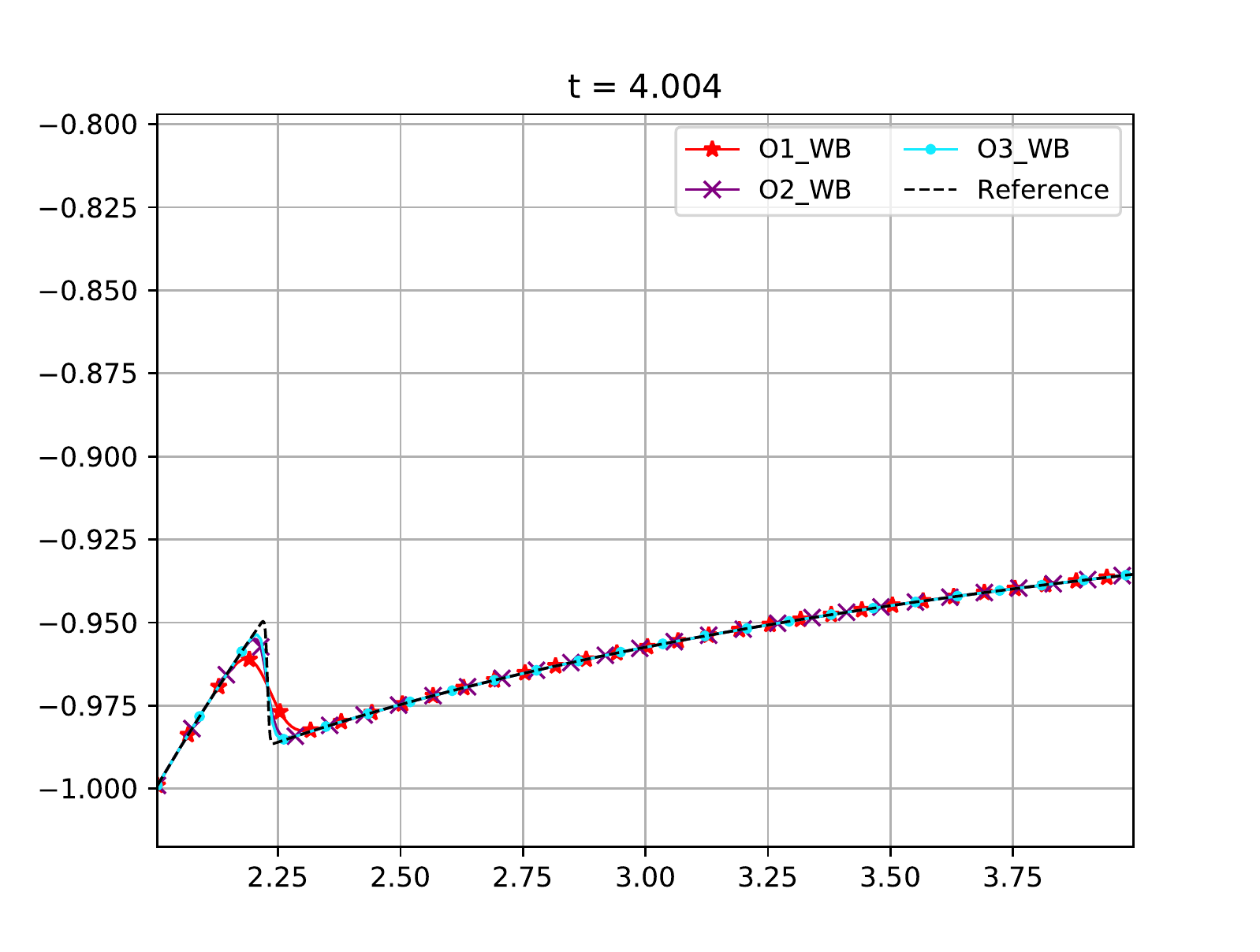}
		\label{fig:ko1_ko2_ko3_wb_testNegativeShock_t_4}
	\end{subfigure}
	\caption{Burgers-Schwarzschild model with the initial condition \eqref{testB5}: first-, second-, and third-order well-balanced  methods at selected times.}
	\label{fig:comparison_o1_o2_03_wb_testNegativeshock}
\end{figure}

Similar conclusions can be drawn for the left-moving shock linking two branches of stationary solutions that generates from the 
initial condition:
\bel{testB5}
v_{0}(r) = \begin{cases}
	-\sqrt{\displaystyle \frac{2}{r}}, & \text{ $2<r<2.5$,}\\
	-\sqrt{\displaystyle \frac{3}{4} + \frac{1}{2r}}, & \text{ otherwise,}
\end{cases}
\ee
see Figure \ref{fig:comparison_o1_o2_03_wb_testNegativeshock}. A reference solution computed with the first-order standard method has been computed again using a mesh of 10000 cells.

\subsection{Perturbation of a steady shock solution}

\paragraph{{{\red Left-hand perturbation}}}

\begin{figure}[h]
	\begin{subfigure}[h]{0.5\textwidth}
		\centering
		\includegraphics[width=1\linewidth]{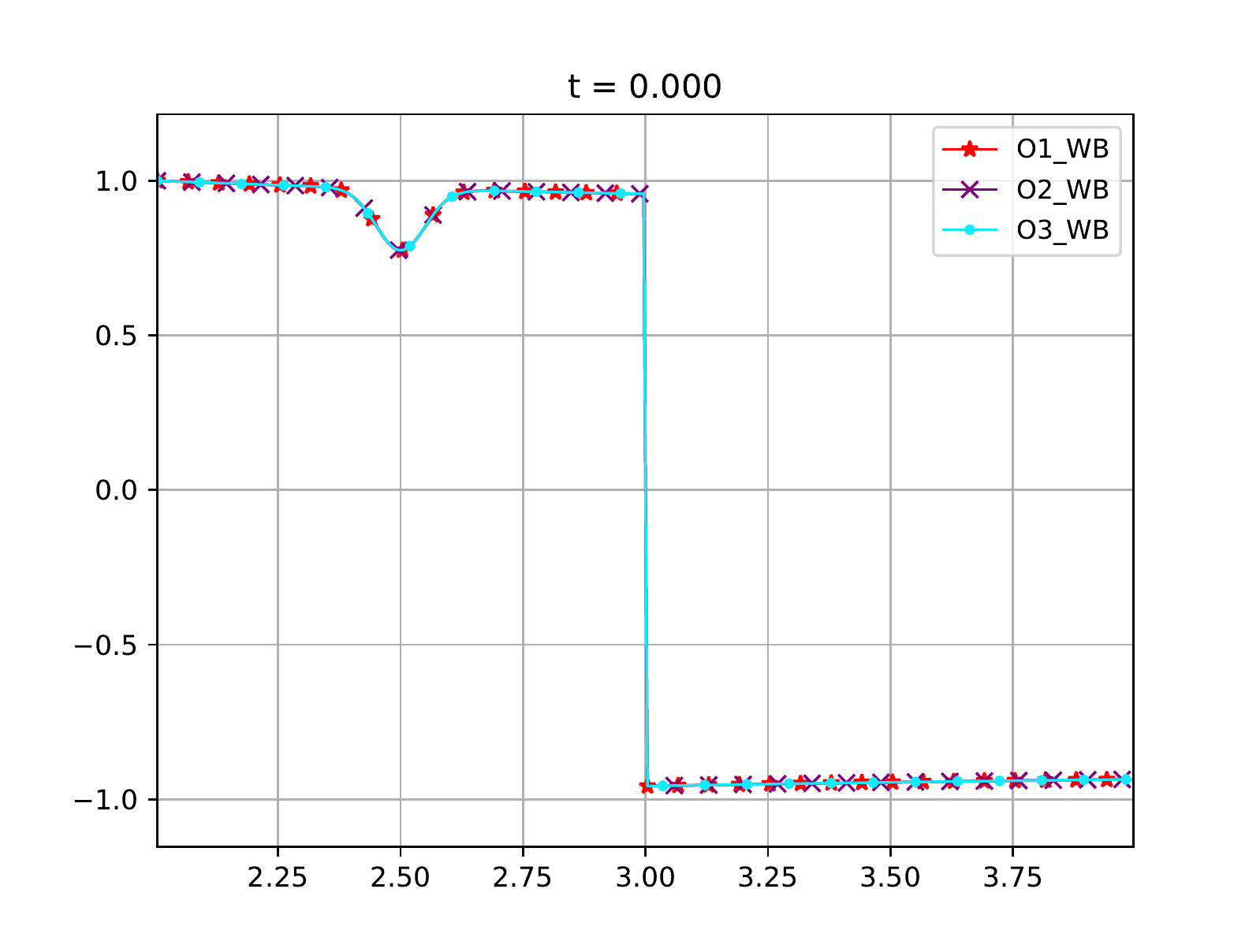}
		\label{fig:ko1_ko2_ko3_wb_testSteadyShockLeftPerturbated_t_0}
	\end{subfigure}
	\begin{subfigure}[h]{0.5\textwidth}
		\centering
		\includegraphics[width=1\linewidth]{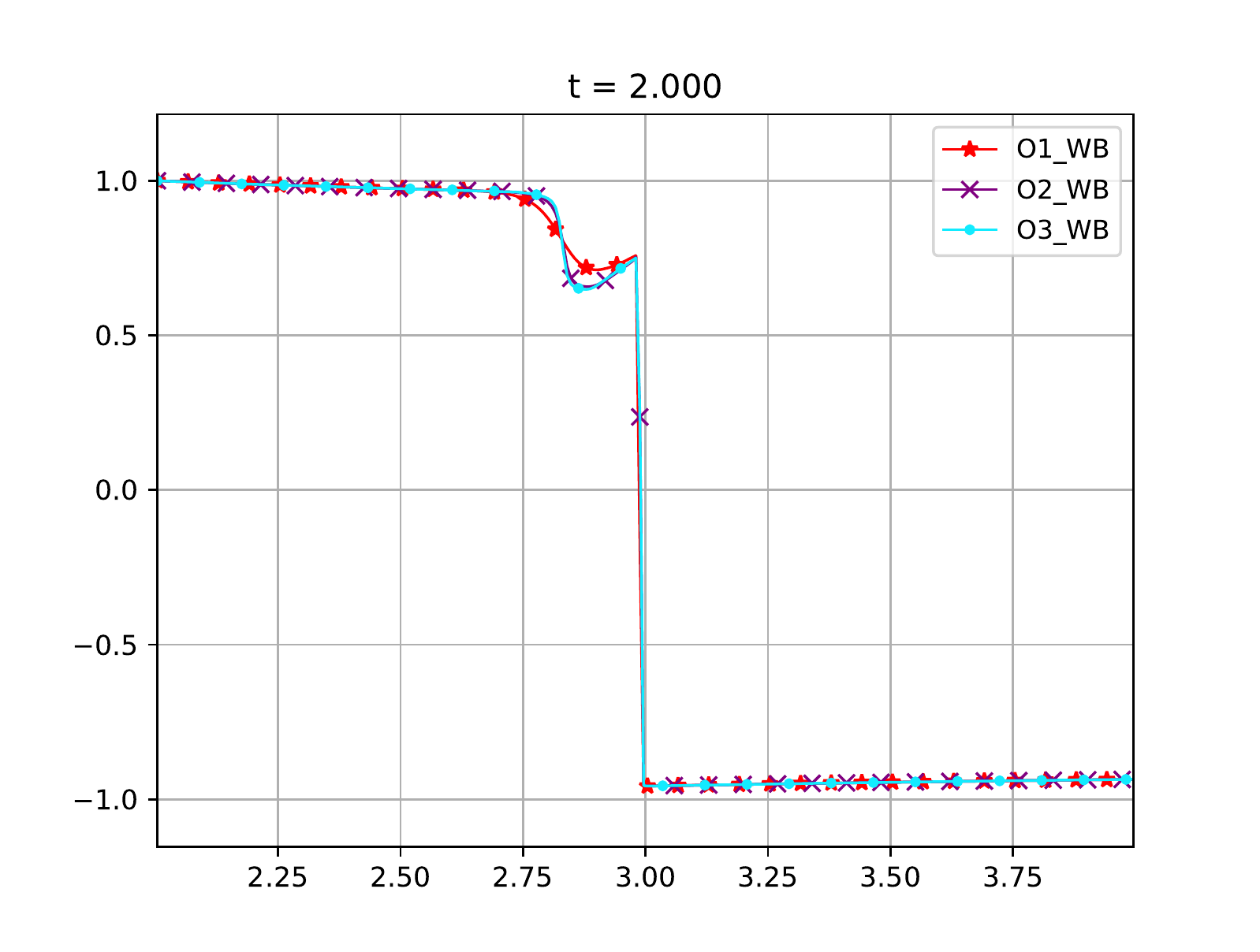}
		\label{fig:ko1_ko2_ko3_wb_testSteadyShockLeftPerturbated_t2}
	\end{subfigure}
	\begin{subfigure}[h]{1\textwidth}
		\centering
		\includegraphics[width=0.5\linewidth]{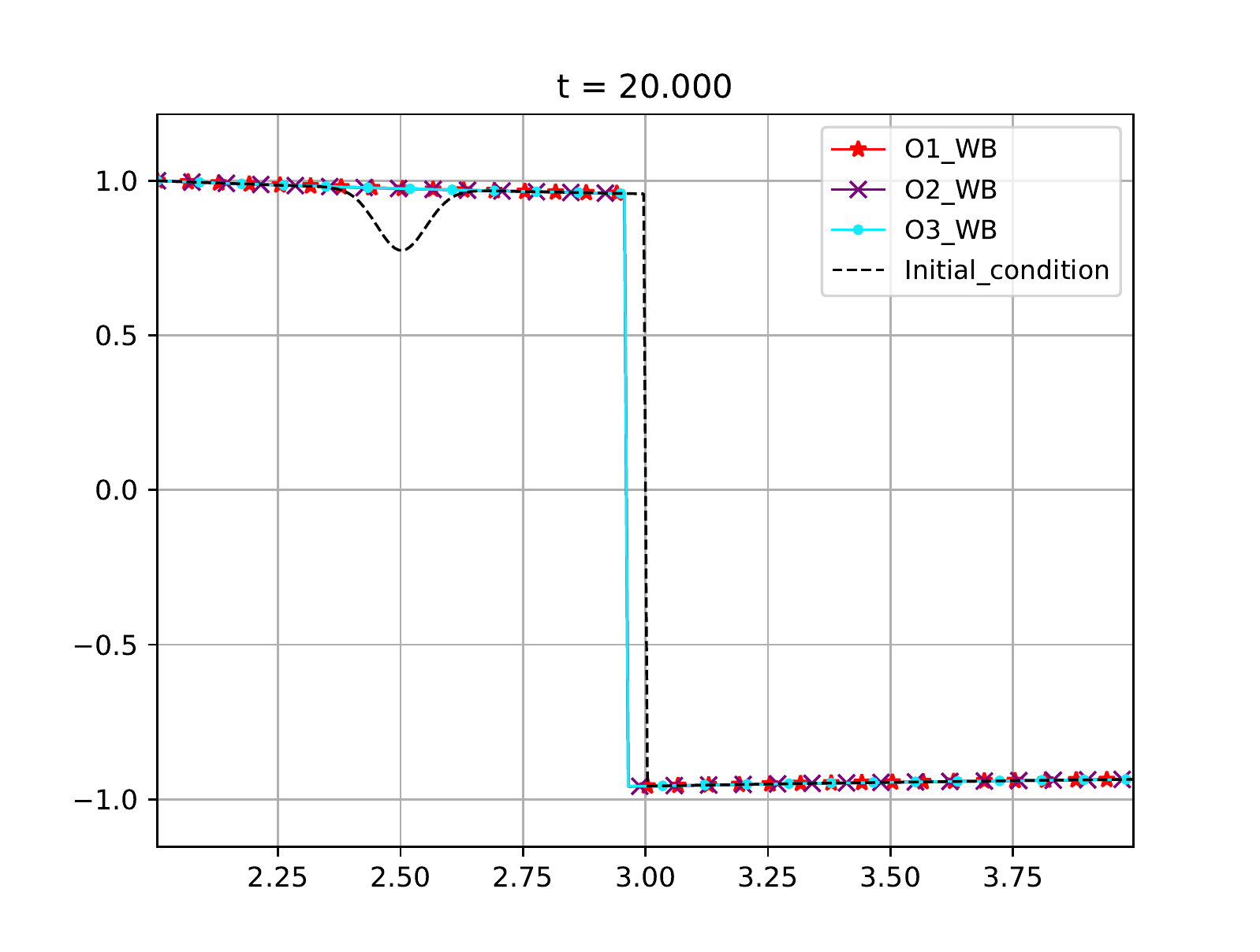}
		\label{fig:ko1_ko2_ko3_wb_testSteadyShockLeftPerturbated_t_20}
	\end{subfigure}
	\caption{Burgers-Schwarzschild model with the initial condition \eqref{testB6}-\eqref{testB3}-\eqref{pl}: 
	 first-, second-, and third-order well-balanced methods at selected times.}
	\label{fig:ko1_ko2_ko3_wb_testSteadyShockLeftPerturbated}
\end{figure}

In this test case  we consider the initial condition:
\bel{testB6}
\tilde v_{0}(r) = v_0(r) + p_L(r),
\ee
where $v_0$ is the steady shock solution given by \eqref{testB3} and
\bel{pl}
p_L(r) =
\begin{cases}
 \displaystyle -\frac{1}{5}e^{-200(r-2.5)^{2}},  & \text{ $2.2<r<2.8$,}\\
0, & \text{ otherwise.} 
\end{cases}
\ee
The first-, second-, and third-order well-balanced methods have been applied to this problem.
In  Figure \ref{fig:ko1_ko2_ko3_wb_testSteadyShockLeftPerturbated} 
it can be observed that, after the wave generated by the initial perturbation leaves the computational domain, the stationary solution 
\eqref{testB3} is not recovered: a different stationary solution of the family \eqref{steadydisc} is obtained whose shock is placed at a different location. Observe that all the three methods capture the same stationary solution.


\paragraph{{{\red Right side perturbation}}}

\begin{figure}[h]
	\begin{subfigure}[h]{0.5\textwidth}
		\centering
		\includegraphics[width=1\linewidth]{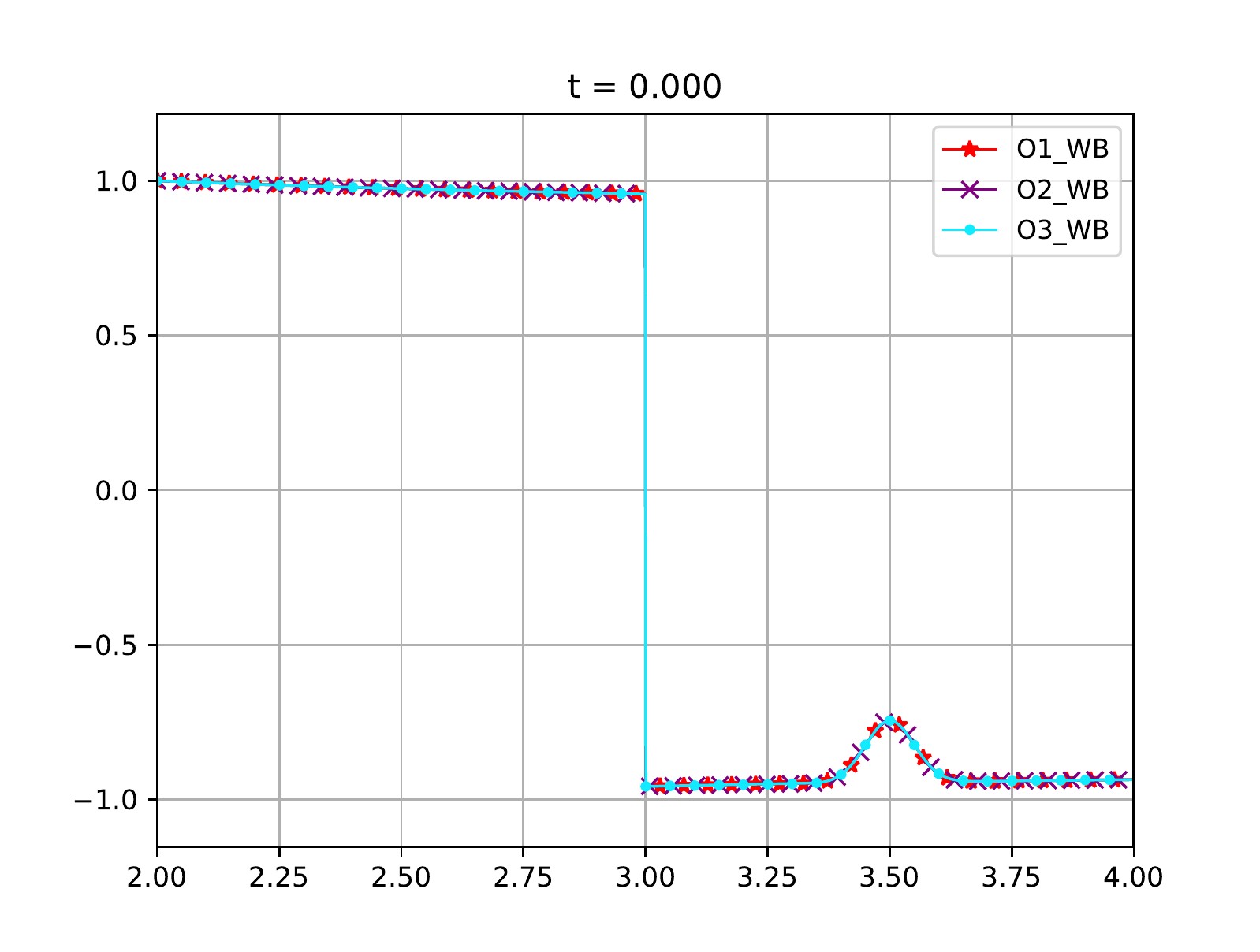}
		\label{fig:ko1_ko2_ko3_wb_testSteadyShockRightPerturbated_t_0_2000}
	\end{subfigure}
	\begin{subfigure}[h]{0.5\textwidth}
		\centering
		\includegraphics[width=1\linewidth]{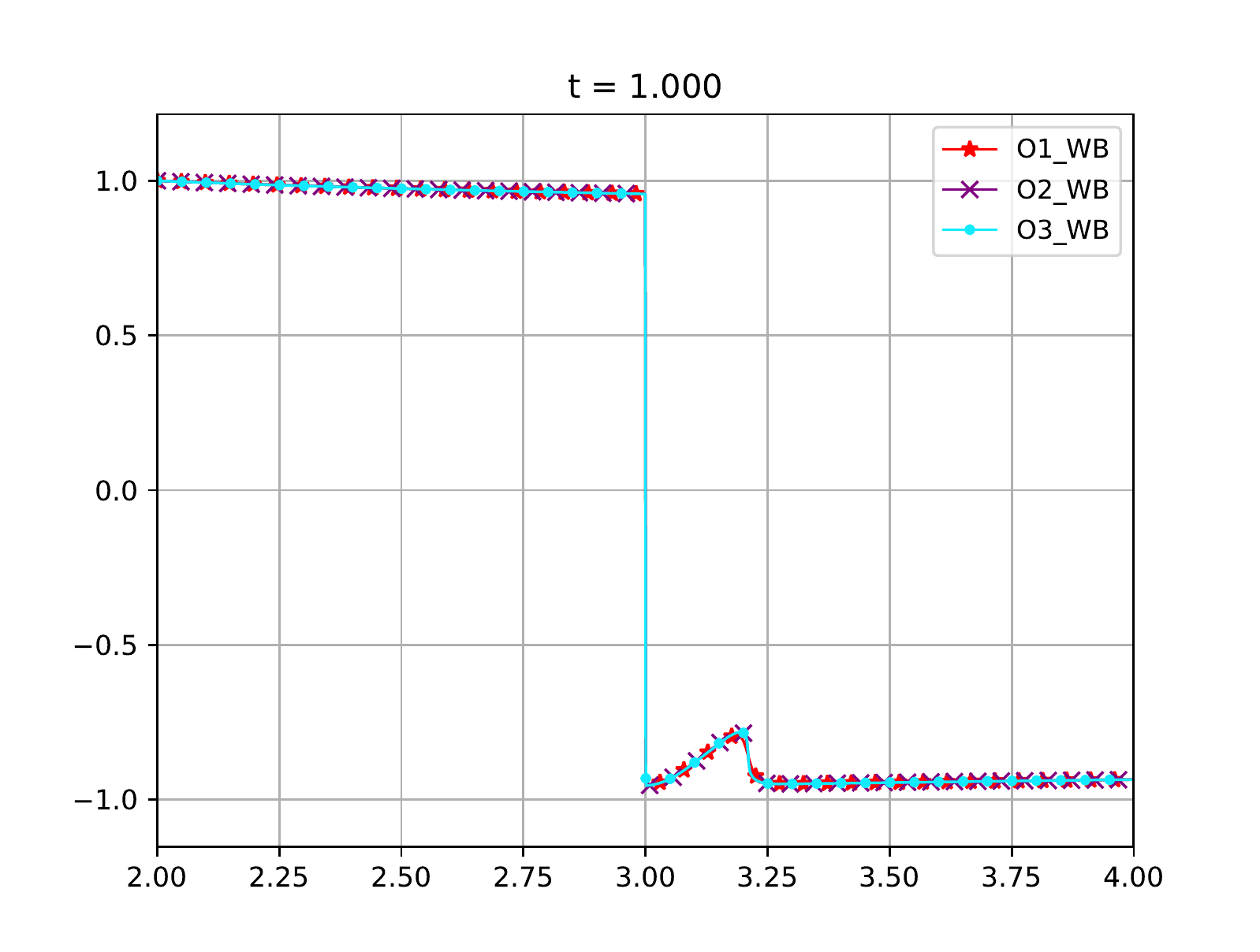}
		\label{fig:ko1_ko2_ko3_wb_testSteadyShockRightPerturbated_t_1_2000}
	\end{subfigure}
	\begin{subfigure}[h]{0.5\textwidth}
		\centering
		\includegraphics[width=1\linewidth]{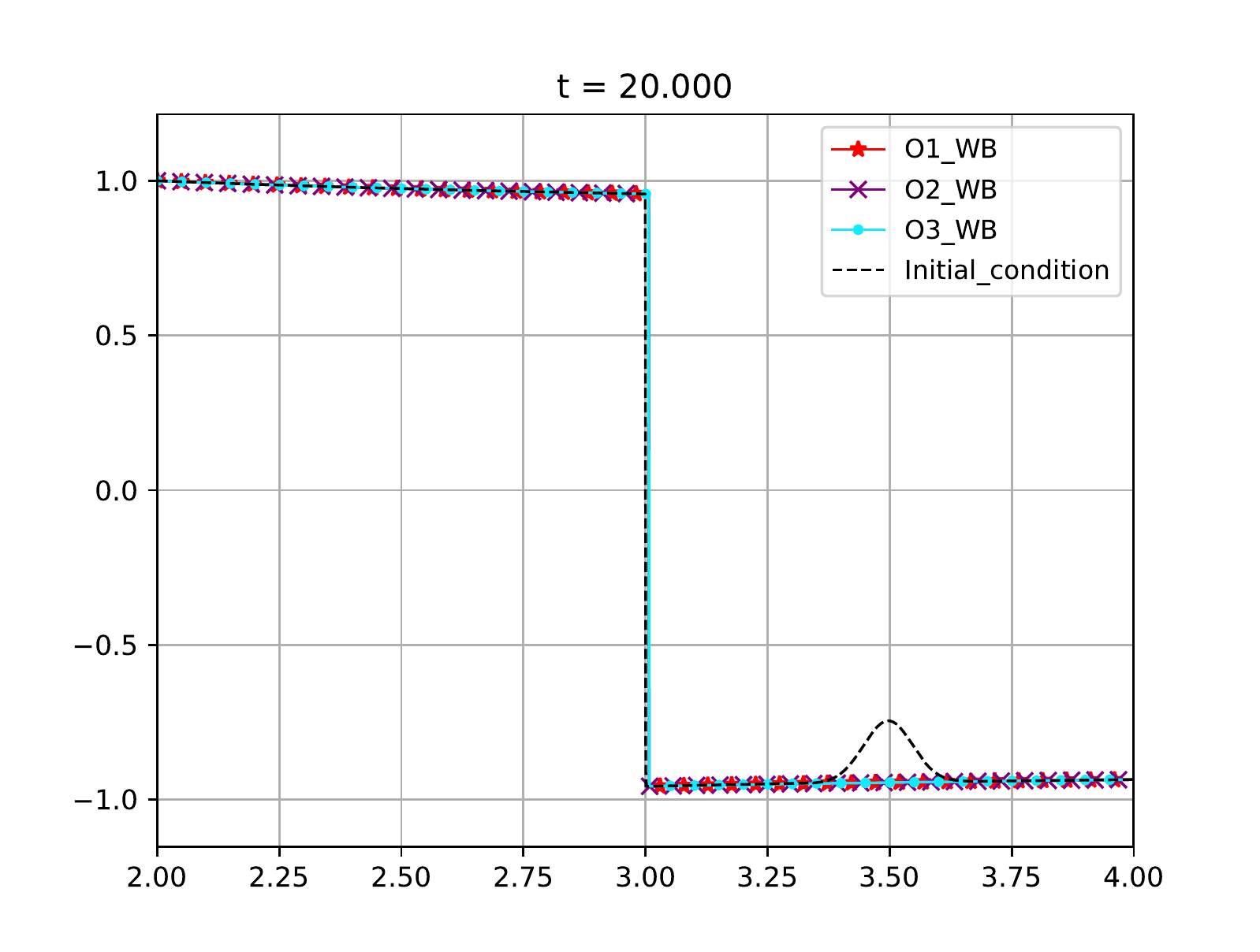}
		\label{fig:ko1_ko2_ko3_wb_testSteadyShockRightPerturbated_t_20_2000}
	\end{subfigure}
\begin{subfigure}[h]{0.5\textwidth}
	\centering
	\includegraphics[width=1\linewidth]{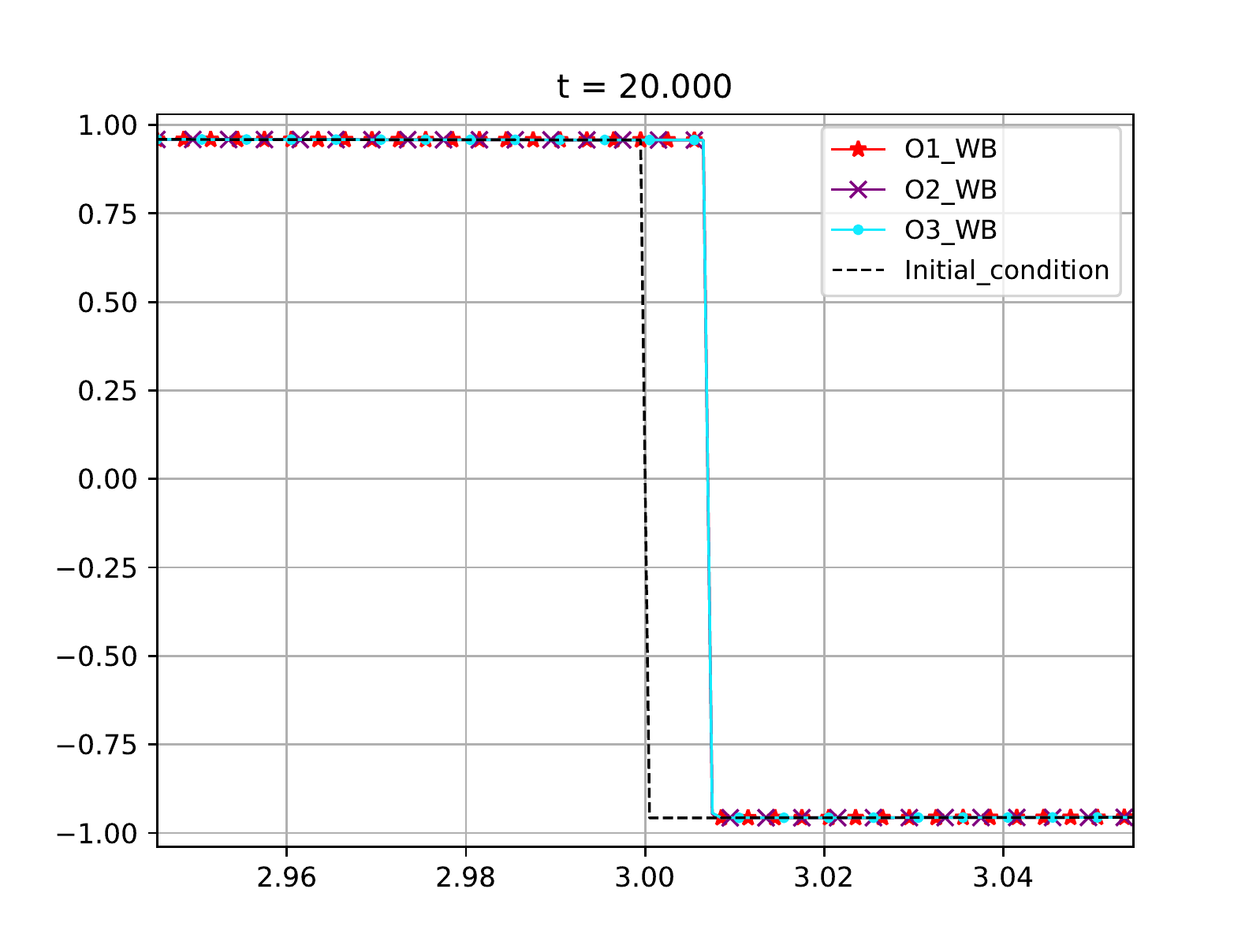}
	\caption*{Zoom of the solution.}
	\label{fig:ko1_ko2_ko3_wb_testSteadyShockRightPerturbated_t_20_2000_zoom}
\end{subfigure}
	\caption{Burgers-Schwarzschild model with the initial condition \eqref{testB7}-\eqref{testB3}-\eqref{pr}: 
	first-, second-, and third-order well-balanced methods at selected times.}
	\label{fig:ko1_ko2_ko3_wb_testSteadyShockRightPerturbated}
\end{figure}

Similar conclusions can be drawn if a perturbation at the right side of the shock is superposed to the stationary solution $v_0$ given by
\eqref{testB3}:
\bel{testB7}
\tilde v_0(r) = v_0(r) + p_R(r),
\ee
with
\bel{pr}
p_R(r) =  \begin{cases}
 \displaystyle \frac{1}{5}e^{-200(r-3.5)^{2}},  & \text{   $3.2<r<3.8$},\\
0, & \text{ otherwise,}
\end{cases}
\ee
as displayed in Figure \ref{fig:ko1_ko2_ko3_wb_testSteadyShockRightPerturbated}. In this case we have used a 2000-point uniform mesh since the displacement of the shock is
 smaller in this case and more points in the mesh are needed in order to see that the steady shock is not recovered.


\paragraph{{{\red Left-hand perturbation with zero average}}}

\begin{figure}[h]
	\begin{subfigure}[h]{0.5\textwidth}
		\centering
		\includegraphics[width=1\linewidth]{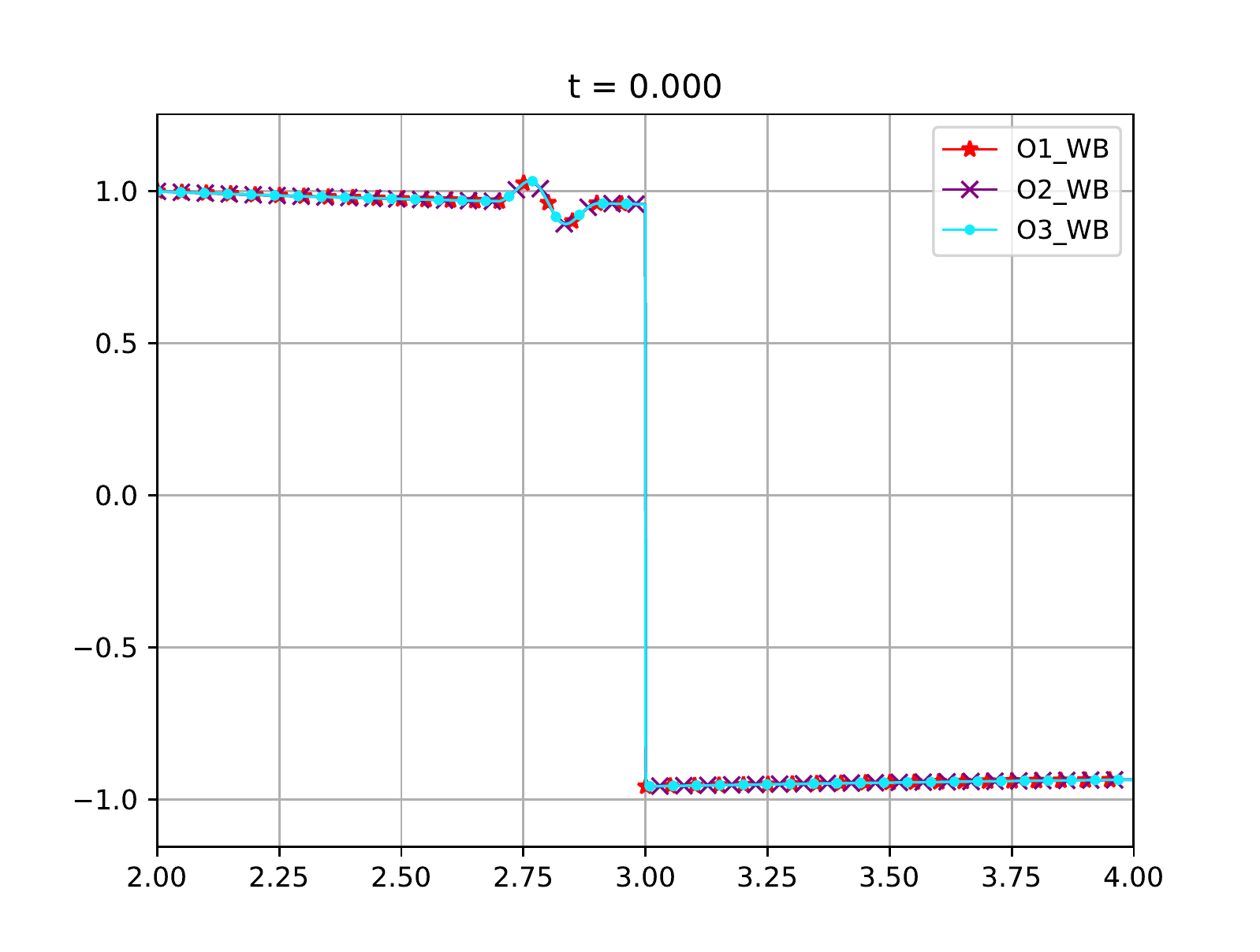}
		\label{fig:ko1_ko2_ko3_wb_testSteadyShockLeftPerturbatedSumZero_t_0_2000}
	\end{subfigure}
	\begin{subfigure}[h]{0.5\textwidth}
		\centering
		\includegraphics[width=1\linewidth]{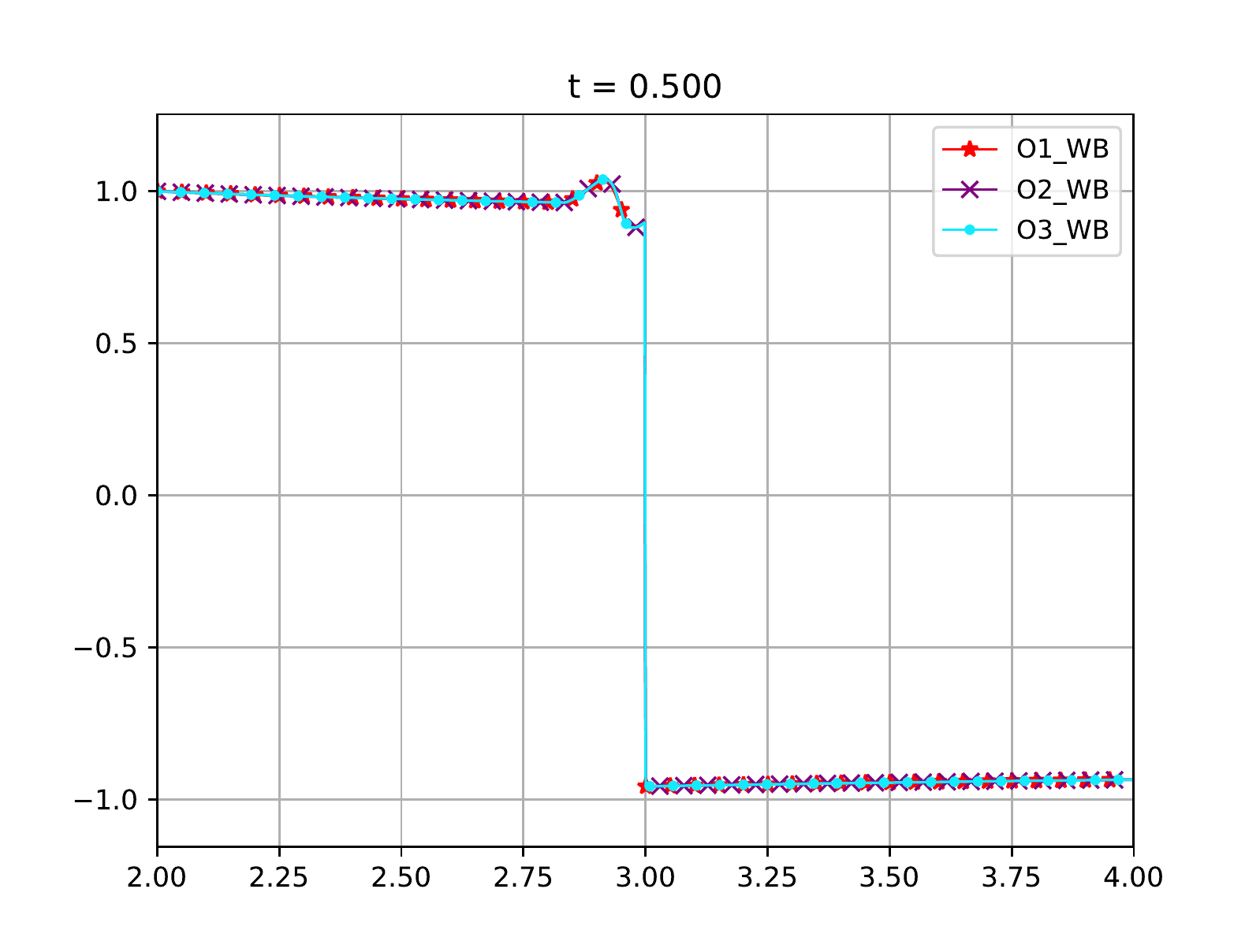}
		\label{fig:ko1_ko2_ko3_wb_testSteadyShockLeftPerturbatedSumZero_t_0_5_2000}
	\end{subfigure}
	\begin{subfigure}[h]{0.5\textwidth}
		\centering
		\includegraphics[width=1\linewidth]{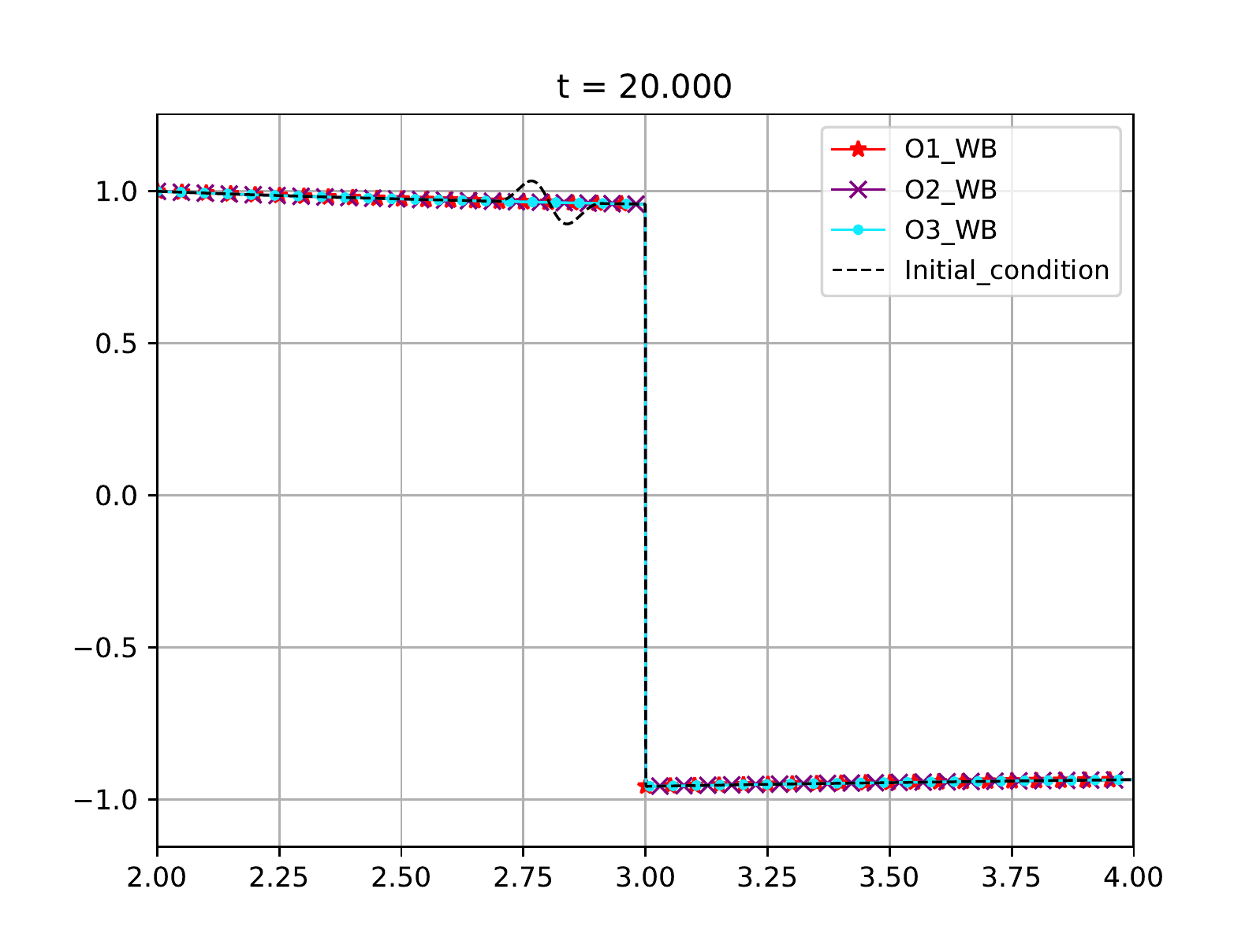}
		\label{fig:ko1_ko2_ko3_wb_testSteadyShockLeftPerturbatedSumZero_t_20_2000}
	\end{subfigure}
	\begin{subfigure}[h]{0.5\textwidth}
		\centering
		\includegraphics[width=1\linewidth]{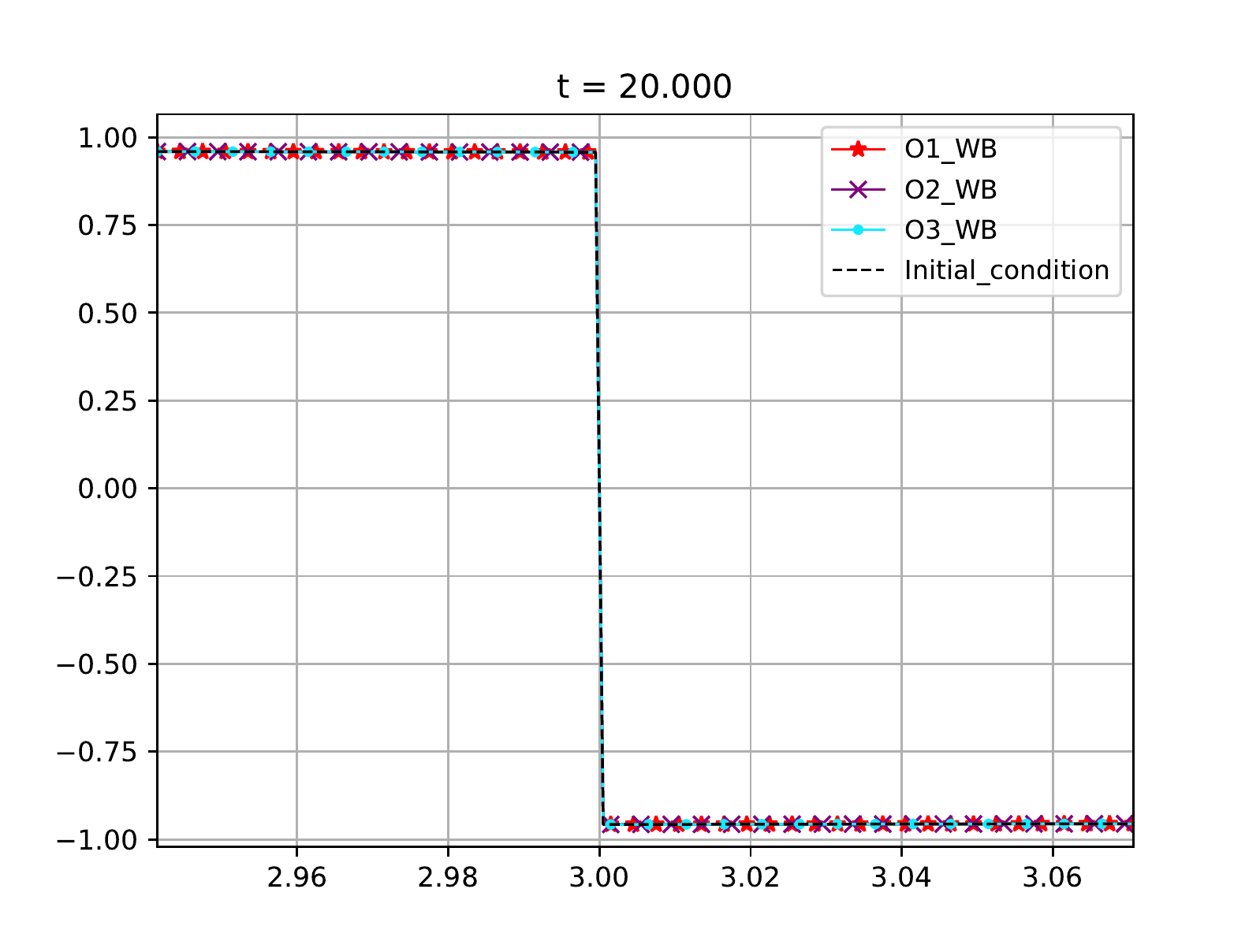}
		\label{fig:ko1_ko2_ko3_wb_testSteadyShockLeftPerturbatedSumZero_t_20_2000_zoom}
	\end{subfigure}
	\caption{Burgers-Schwarzschild model with the initial condition \eqref{testB6}-\eqref{testB3}-\eqref{pl_SumZero}: e first-, second-, and third-order well-balanced methods at selected times, and zoom of the initial and final  stationary shocks (right-down).}
	\label{fig:ko1_ko2_ko3_wb_testSteadyShockLeftPerturbatedSumZero}
\end{figure}

Now we consider an  initial condition of the form (\ref{testB6}) with a  perturbation  $p_L$ such that $\int p_L(r)dr = 0$, in particular:
\bel{pl_SumZero}
p_L(r) =
\begin{cases}
	\displaystyle 0.1\text{cos}(-25.5 \pi + 10\pi x)e^{-200(x-2.8)(x-2.8)},  & \text{ $2.7<r<2.9$,}\\
	0, & \text{ otherwise.} 
\end{cases}
\ee
In Figure \ref{fig:ko1_ko2_ko3_wb_testSteadyShockLeftPerturbatedSumZero} it can be observed that now, after the wave generated by the initial perturbation leaves the computational domain, the stationary solution 
\eqref{testB3} is recovered. Here we have used again a 2000-point uniform mesh to verify that the steady state does not move.

\paragraph{{{\red Right side perturbation with zero average}}}

\begin{figure}[h]
	\begin{subfigure}[h]{0.5\textwidth}
		\centering
		\includegraphics[width=1\linewidth]{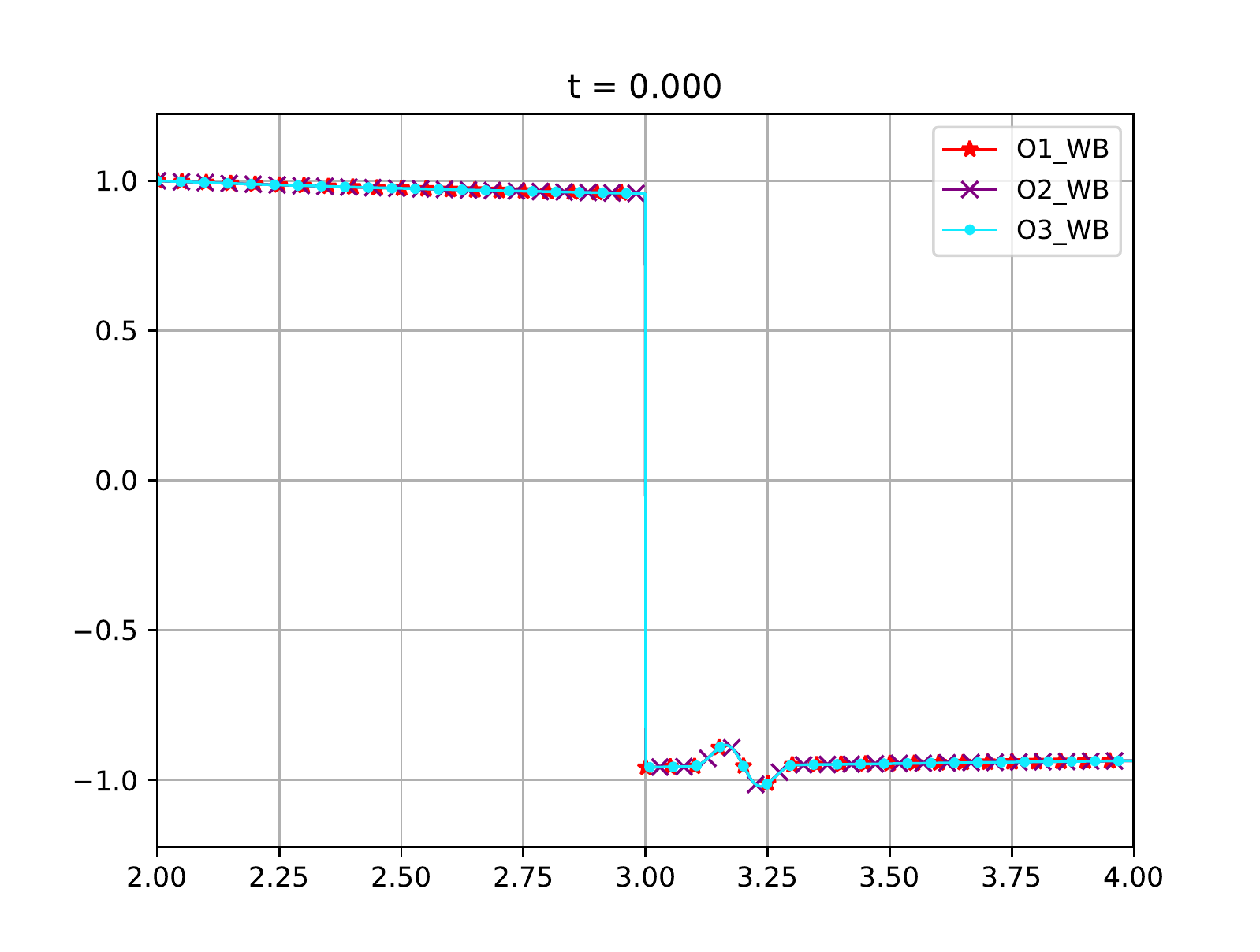}
		\label{fig:ko1_ko2_ko3_wb_testSteadyShockRightPerturbatedSumZero_t_0_2000}
	\end{subfigure}
	\begin{subfigure}[h]{0.5\textwidth}
		\centering
		\includegraphics[width=1\linewidth]{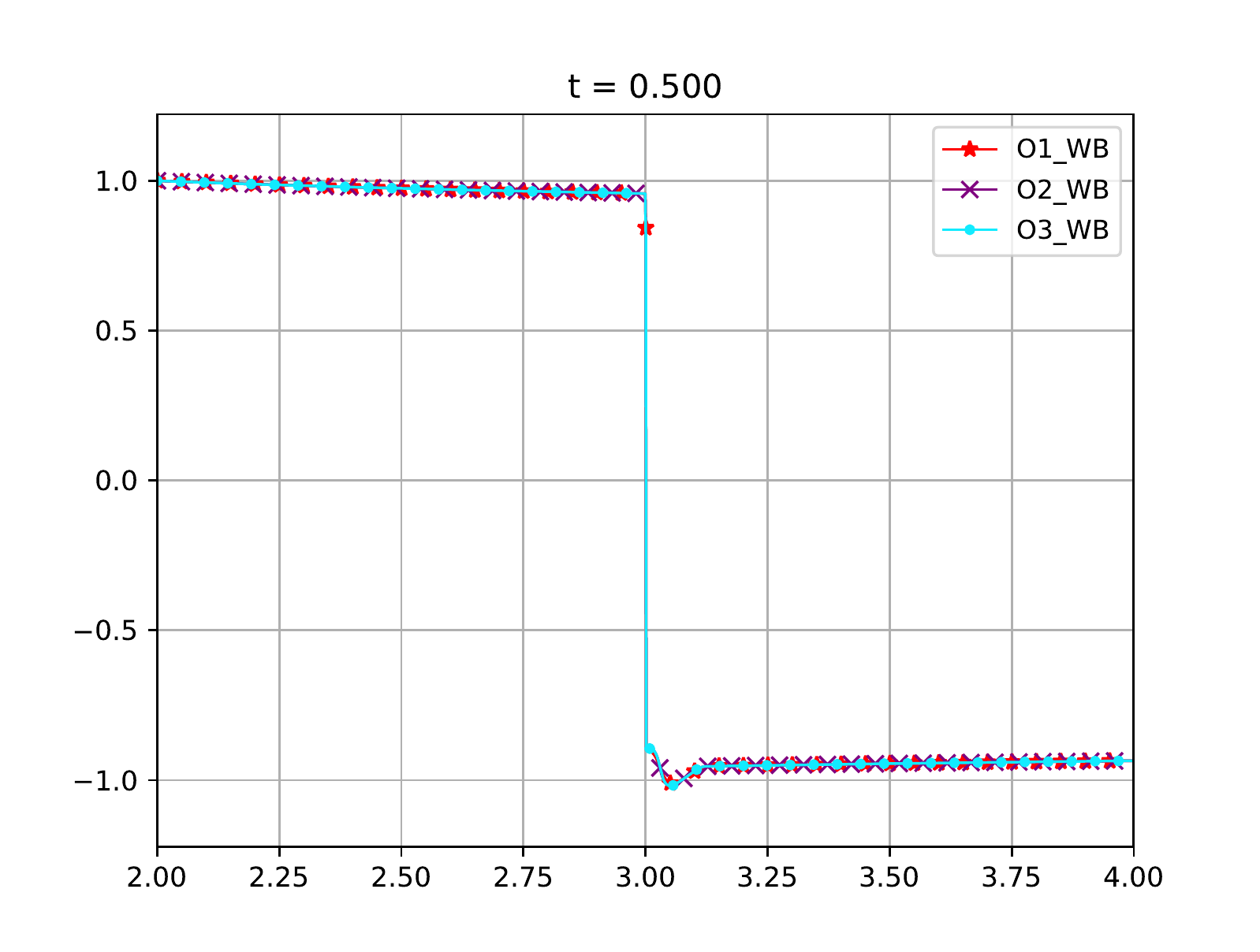}
		\label{fig:ko1_ko2_ko3_wb_testSteadyShockRightPerturbatedSumZero_t_0_5_2000}
	\end{subfigure}
	\begin{subfigure}[h]{0.5\textwidth}
		\centering
		\includegraphics[width=1\linewidth]{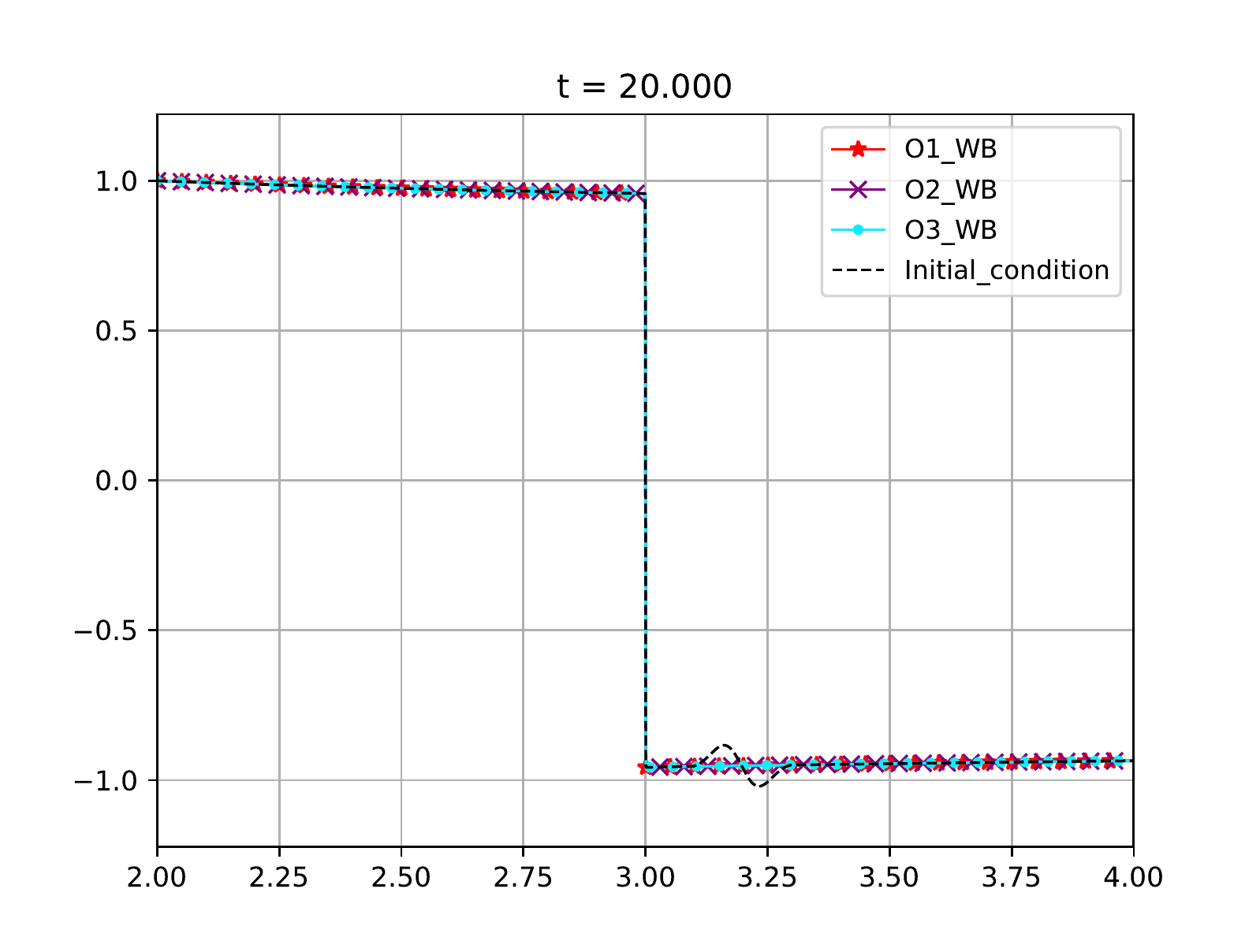}
		\label{fig:ko1_ko2_ko3_wb_testSteadyShockRightPerturbatedSumZero_t_20_2000}
	\end{subfigure}
	\begin{subfigure}[h]{0.5\textwidth}
		\centering
		\includegraphics[width=1\linewidth]{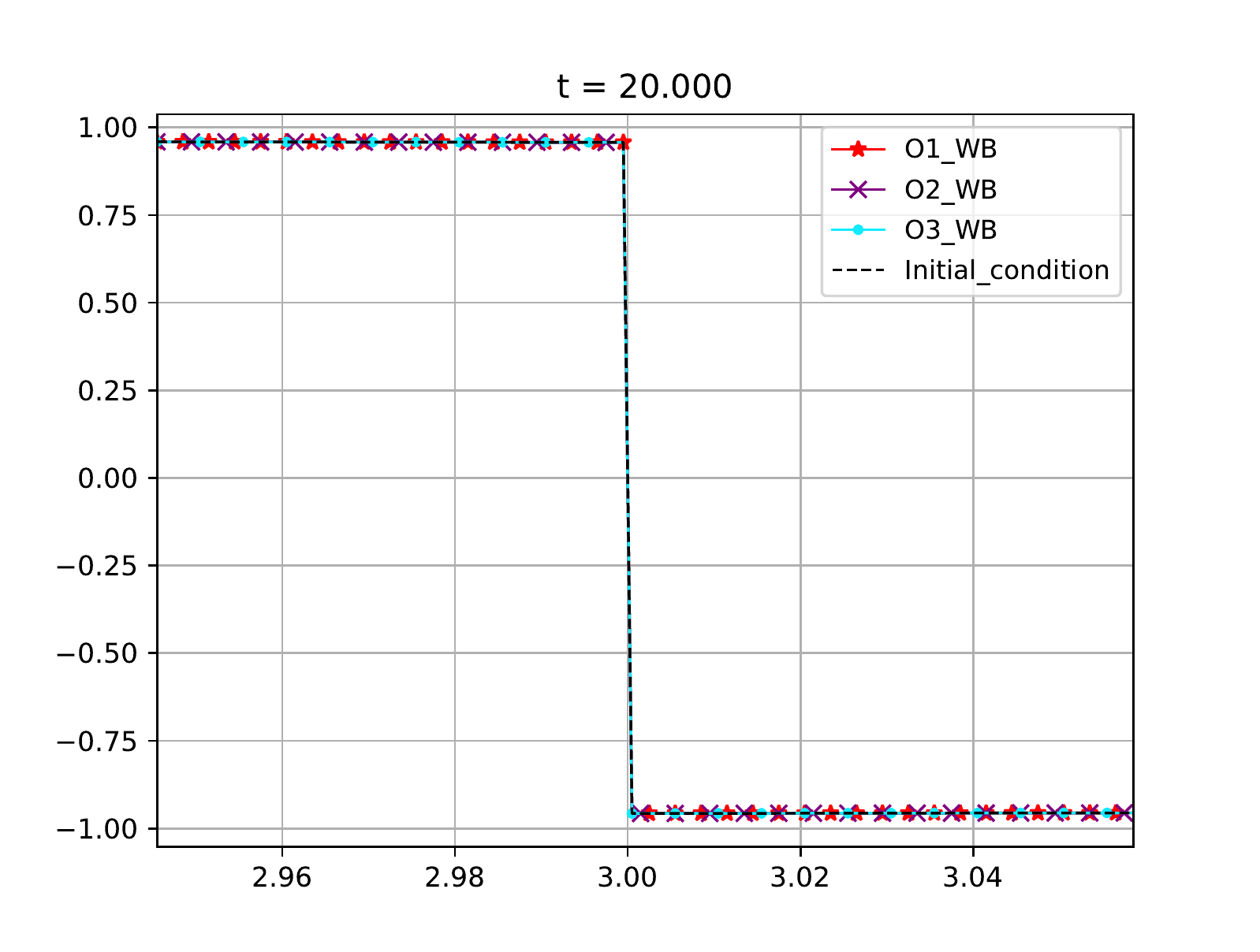}
		\caption*{Zoom of the solution.}
		\label{fig:ko1_ko2_ko3_wb_testSteadyShockRightPerturbatedSumZero_t_20_2000_zoom}
	\end{subfigure}
	\caption{Burgers-Schwarzschild model with the initial condition \eqref{testB7}-\eqref{testB3}-\eqref{pr_SumZero}:  first-, second-, and third-order well-balanced  methods at selected times, and zoom of the initial and final  stationary shocks (right-down).}
	\label{fig:ko1_ko2_ko3_wb_testSteadyShockRightPerturbatedSumZero}
\end{figure}

Similar conclusions can be drawn if we consider an initial condition of the form (\ref{testB7}) with $\int p_R(r)dr = 0$. In particular we take
\bel{pr_SumZero}
p_R(r) =
\begin{cases}
	\displaystyle 0.1\text{cos}(-29.5 \pi + 10\pi x)e^{-200(x-3.2)(x-3.2)},  & \text{ $3.1<r<3.3$,}\\
	0, & \text{ otherwise,} 
\end{cases}
\ee
see Figure \ref{fig:ko1_ko2_ko3_wb_testSteadyShockRightPerturbatedSumZero}. Here we have used again a 2000-point uniform mesh to verify that the steady state does not move.

\paragraph{{{\red Left-hand and right-hand perturbations with zero average}}}

\begin{figure}[h]
	\begin{subfigure}[h]{0.5\textwidth}
		\centering
		\includegraphics[width=1\linewidth]{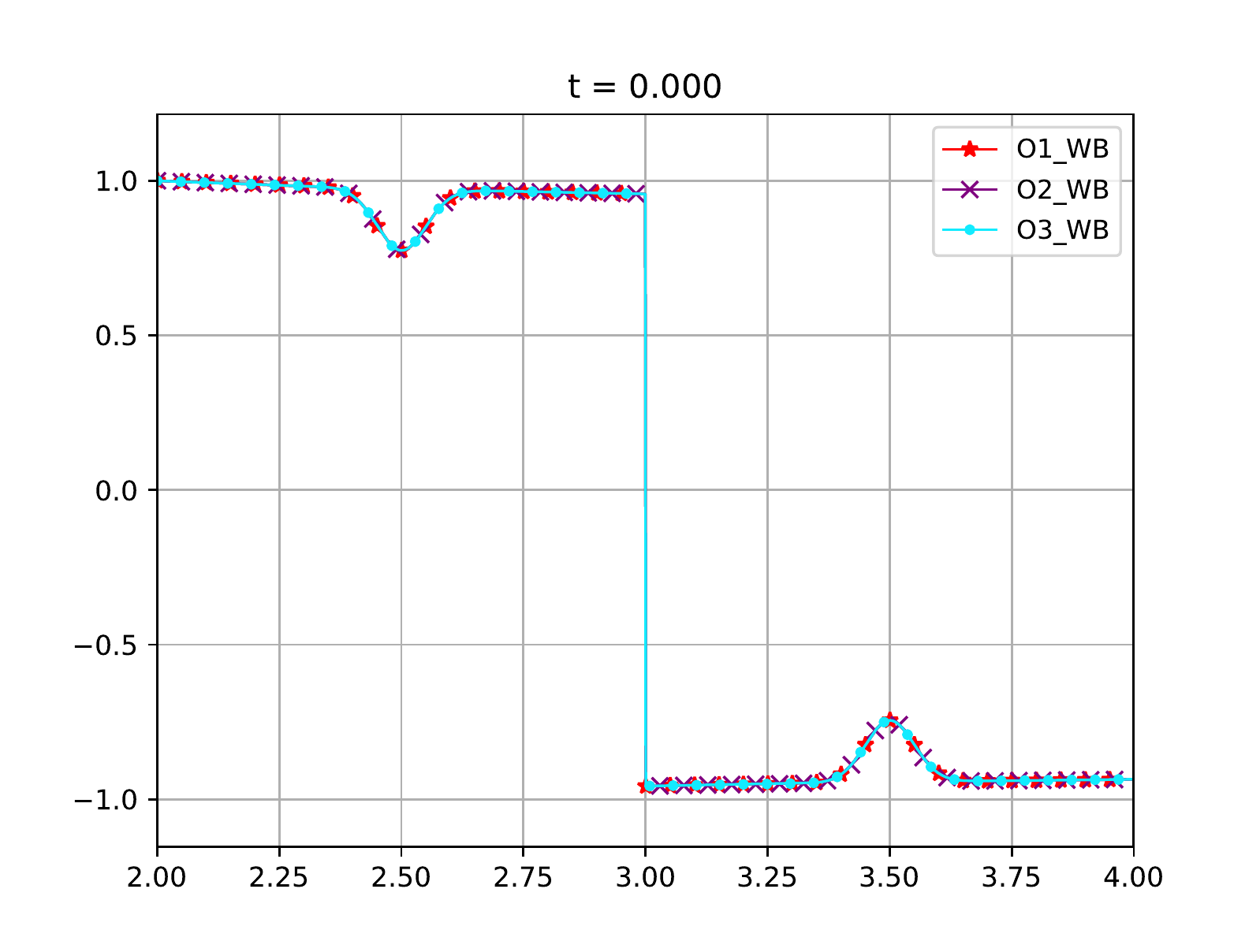}
		\label{fig:ko1_ko2_ko3_wb_testSteadyShockBothPerturbatedSumZero_t_0_2000}
	\end{subfigure}
	\begin{subfigure}[h]{0.5\textwidth}
		\centering
		\includegraphics[width=1\linewidth]{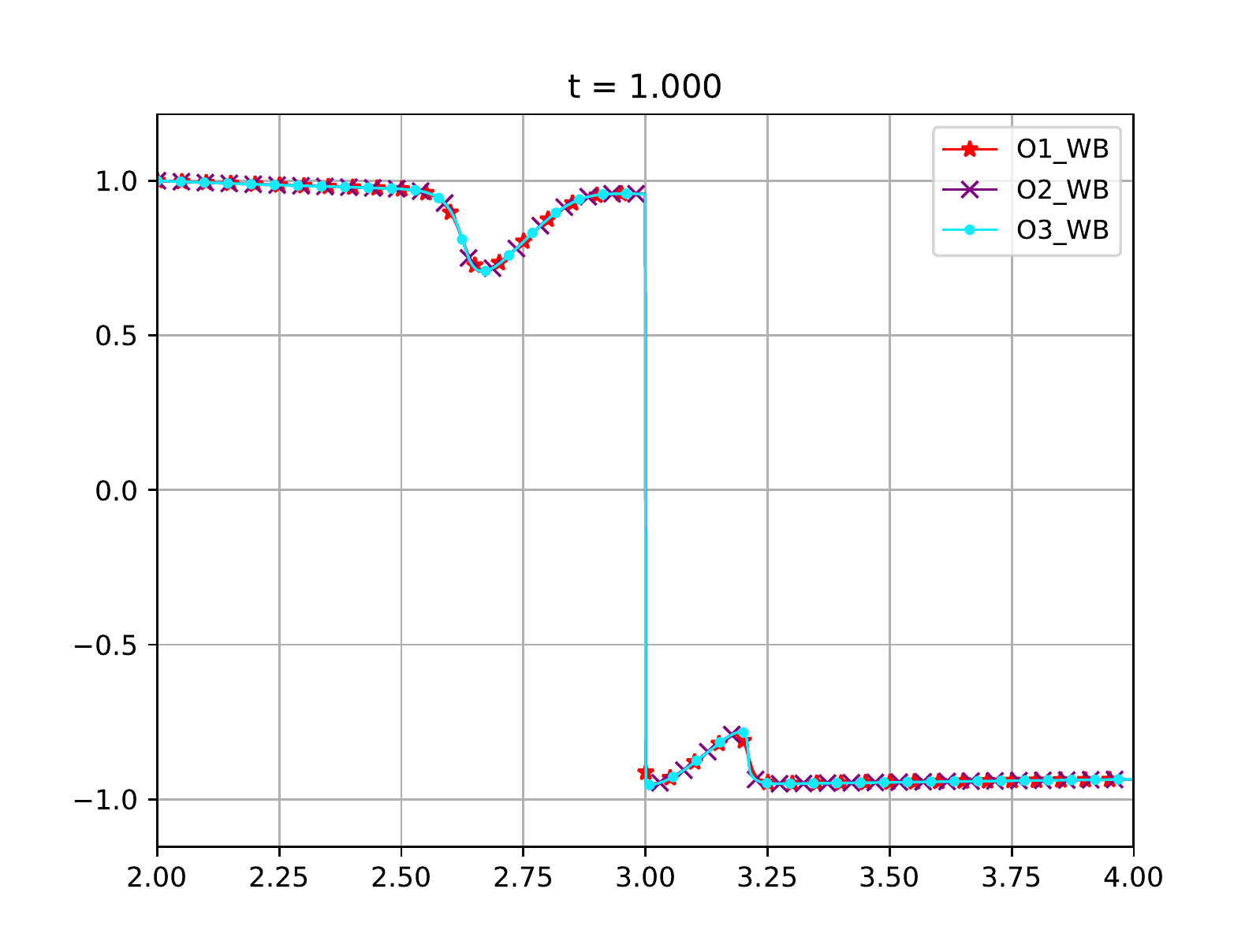}
		\label{fig:ko1_ko2_ko3_wb_testSteadyShockBothPerturbatedSumZero_t_1_2000}
	\end{subfigure}
	\begin{subfigure}[h]{0.5\textwidth}
		\centering
		\includegraphics[width=1\linewidth]{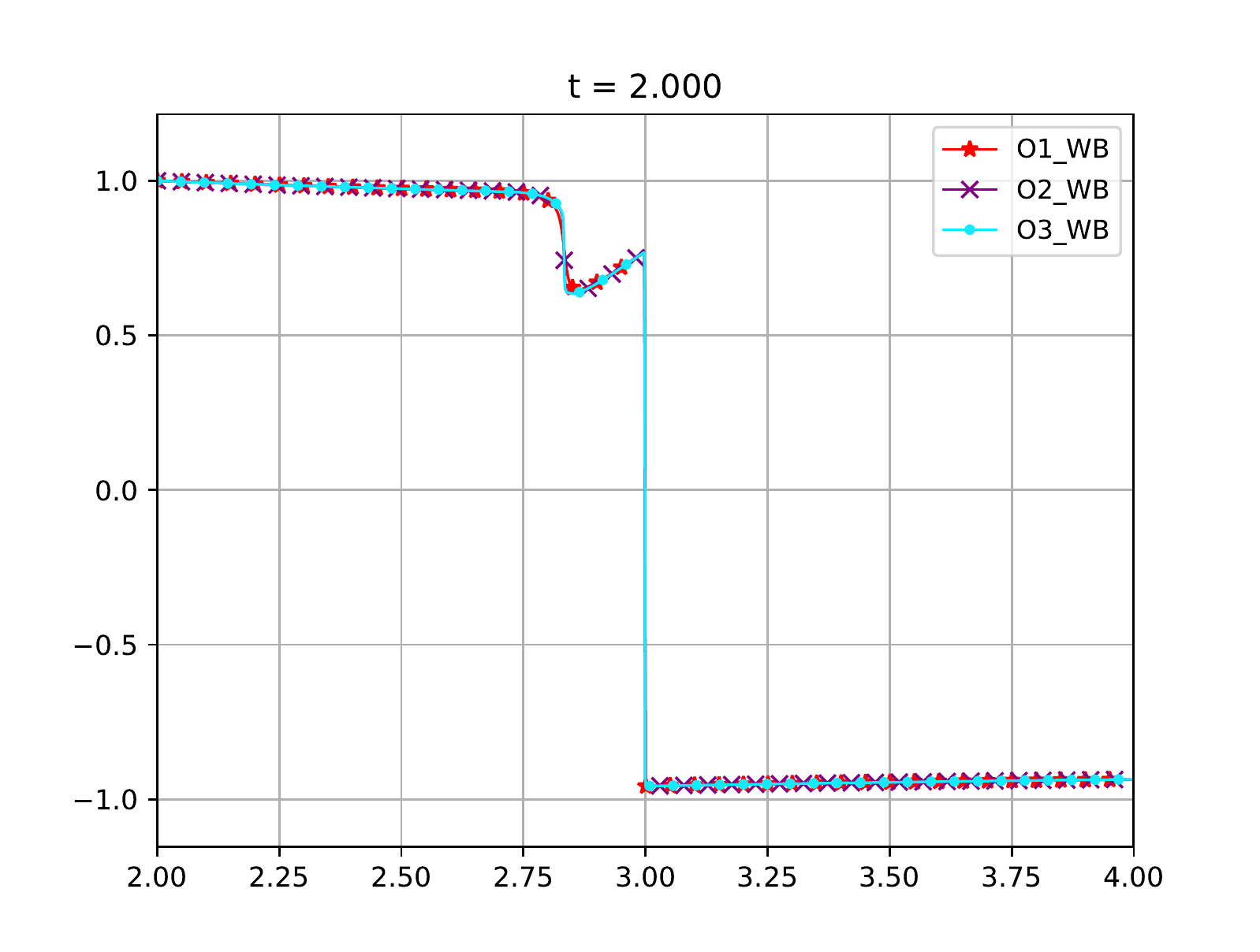}
		\label{fig:ko1_ko2_ko3_wb_testSteadyShockBothPerturbatedSumZero_t_2_2000}
	\end{subfigure}
	\begin{subfigure}[h]{0.5\textwidth}
		\centering
		\includegraphics[width=1\linewidth]{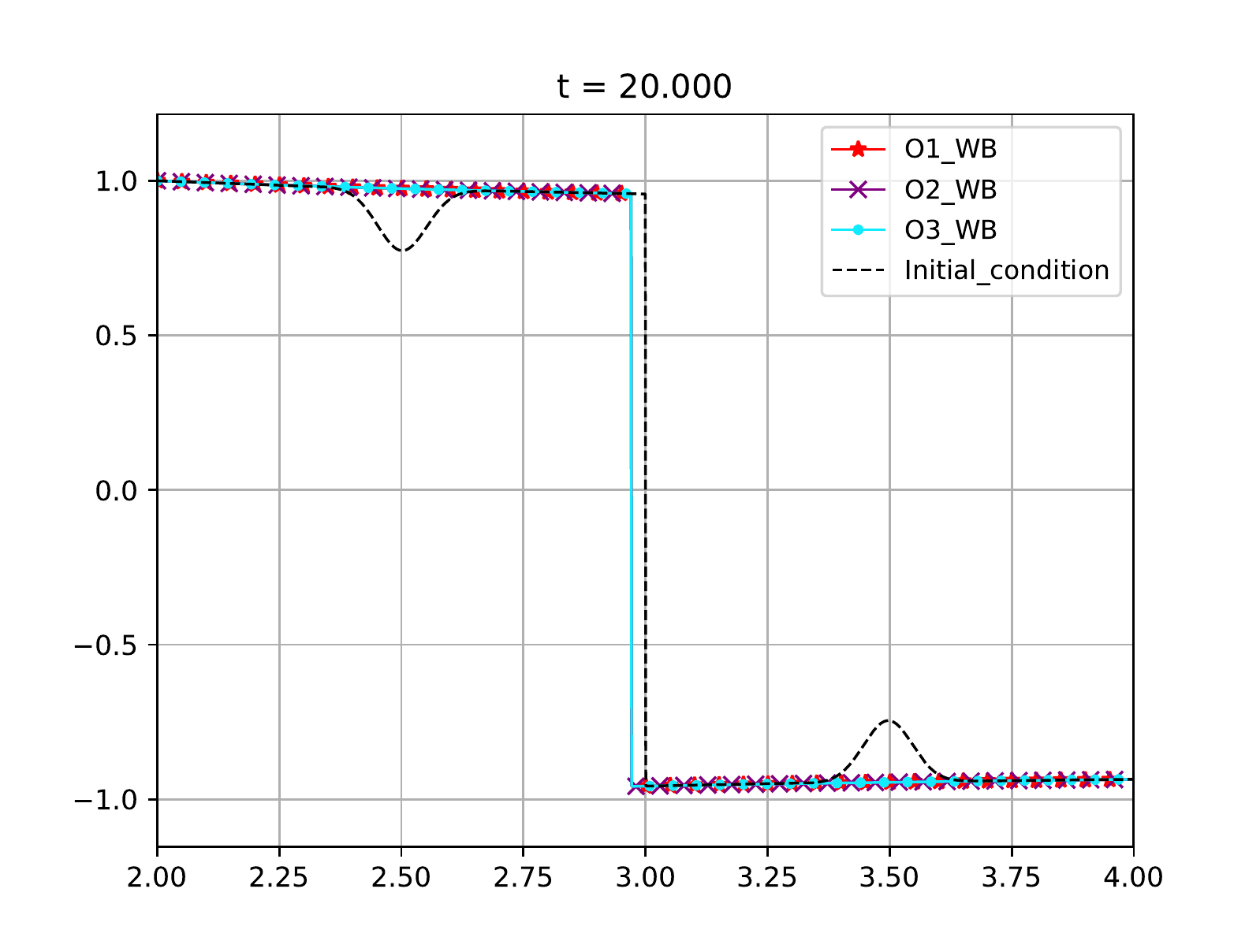}
		\label{fig:ko1_ko2_ko3_wb_testSteadyShockBothPerturbatedSumZero_t_20_2000}
	\end{subfigure}
	\caption{Burgers-Schwarzschild model with the initial condition \eqref{testB7}-\eqref{testB3}-\eqref{pr_SumZero}: comparison between the first-, second-, and third-order well-balanced  methods at selected times.}
	\label{fig:ko1_ko2_ko3_wb_testSteadyShockBothPerturbatedSumZero}
\end{figure}

In order to study the relation between the amplitude of the perturbation and the distance between the initial and the final stationary shocks, we consider the initial condition:
\bel{testB6+7}
\tilde v_{0}(r) = v_0(r) + p_L(r) + p_R(r),
\ee
where $v_0$ is the steady shock solution given by \eqref{testB3} and $\int (p_L(r)+p_R(r))dr =0$. In particular we take $p_L(r)$ as in ($\ref{pl}$) and $p_R(r)$ as in ($\ref{pr}$). In Figure \ref{fig:ko1_ko2_ko3_wb_testSteadyShockBothPerturbatedSumZero} it can be observed that, after the wave generated by the initial perturbation leaves the computational domain, the stationary solution 
\eqref{testB3} is not recovered: a different stationary solution of the family \eqref{steadydisc} with the shock is placed at a different location. This is a natural result since, as we saw before, the right perturbation creates a lower displacement than the left-hand perturbation. Here we have used again a 2000-point uniform mesh.

\paragraph{{{\red Relation between the perturbation and the displacement of the shock}}}

\begin{figure}[h]
	\begin{subfigure}[h]{0.5\textwidth}
		\centering
		\includegraphics[width=0.8\linewidth]{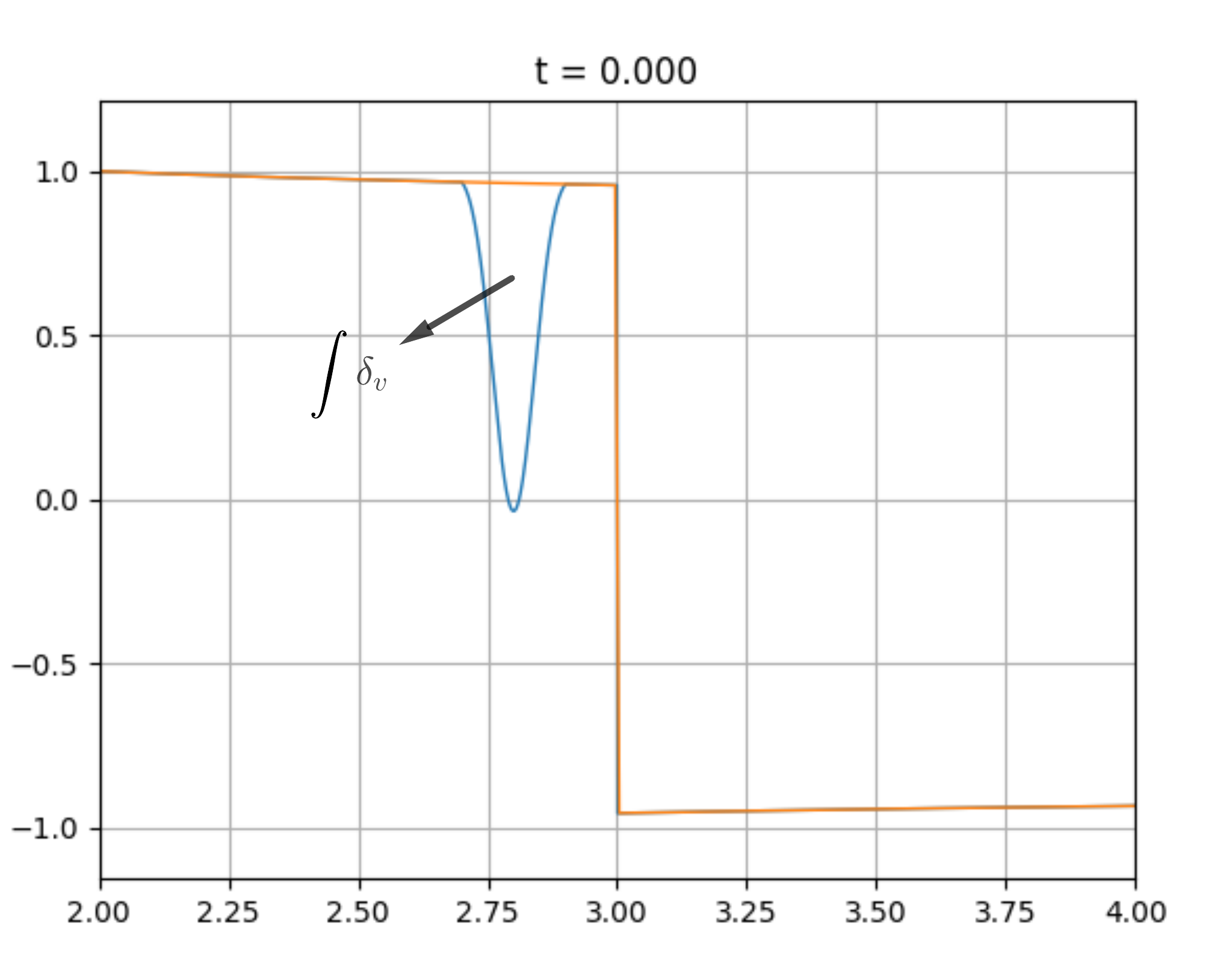}
		\caption{Area of the perturbation for $\alpha=1$.}
		\label{fig:ko1wbteststeadyshockperturbatedalphaintegralperturbationlabelled}
	\end{subfigure}
	\begin{subfigure}[h]{0.5\textwidth}
		\centering
		\includegraphics[width=0.8\linewidth]{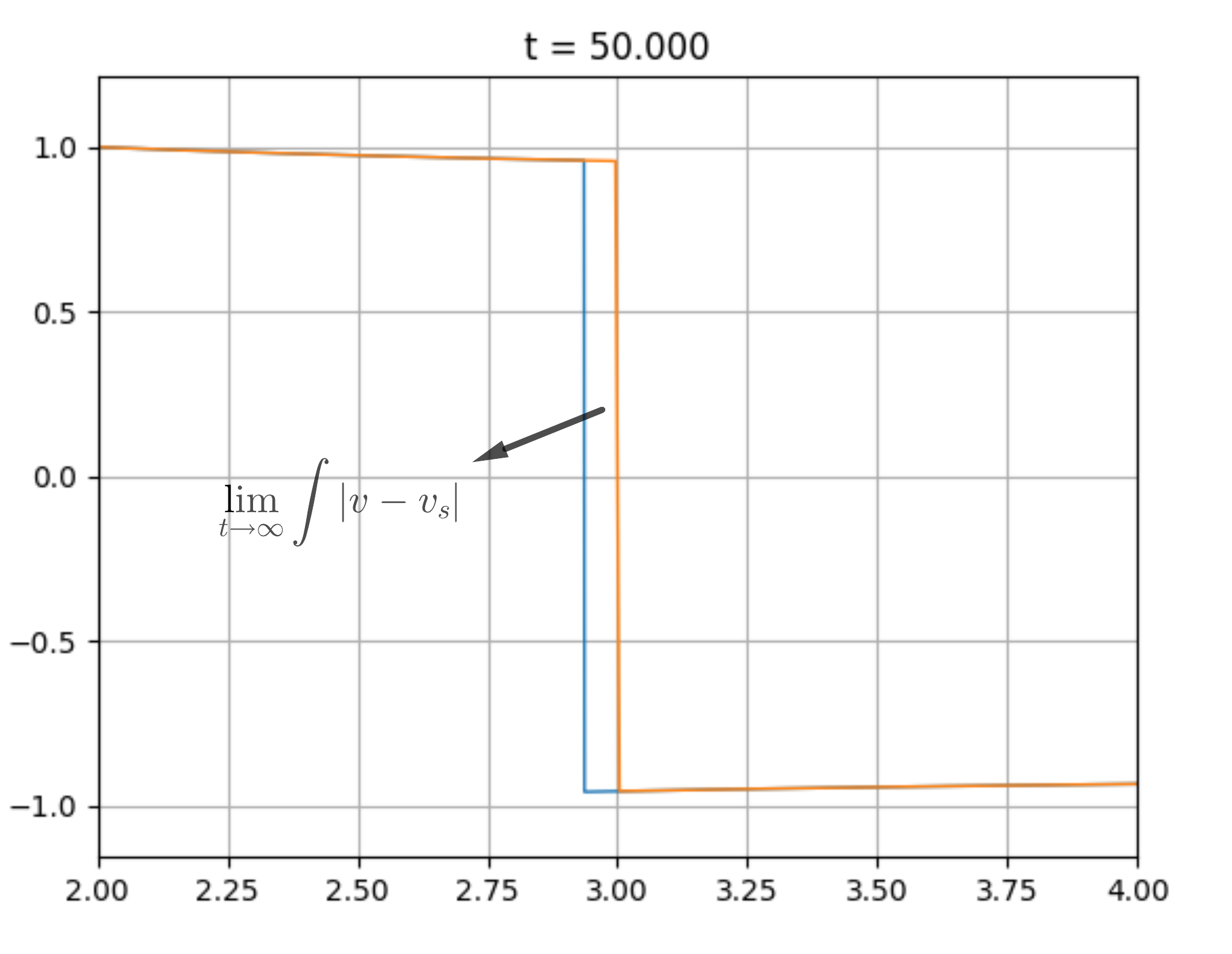}
		\caption{Area between the initial  and the final steady shock and the final steady shocks for $\alpha=1$.}
		\label{fig:ko1wbteststeadyshockperturbatedalphaintegralsteadyshockperturbated}
	\end{subfigure}
	\caption{Burgers-Schwarzschild model with the initial condition \eqref{testB8}-\eqref{testB3}-\eqref{delta}: measures of the perturbation and the shock displacement for 
$\alpha = 1$.}
	\label{fig:Areas_with_alpha_1}
\end{figure}

\begin{figure}[h]
	\centering
	\includegraphics[width=0.5\linewidth]{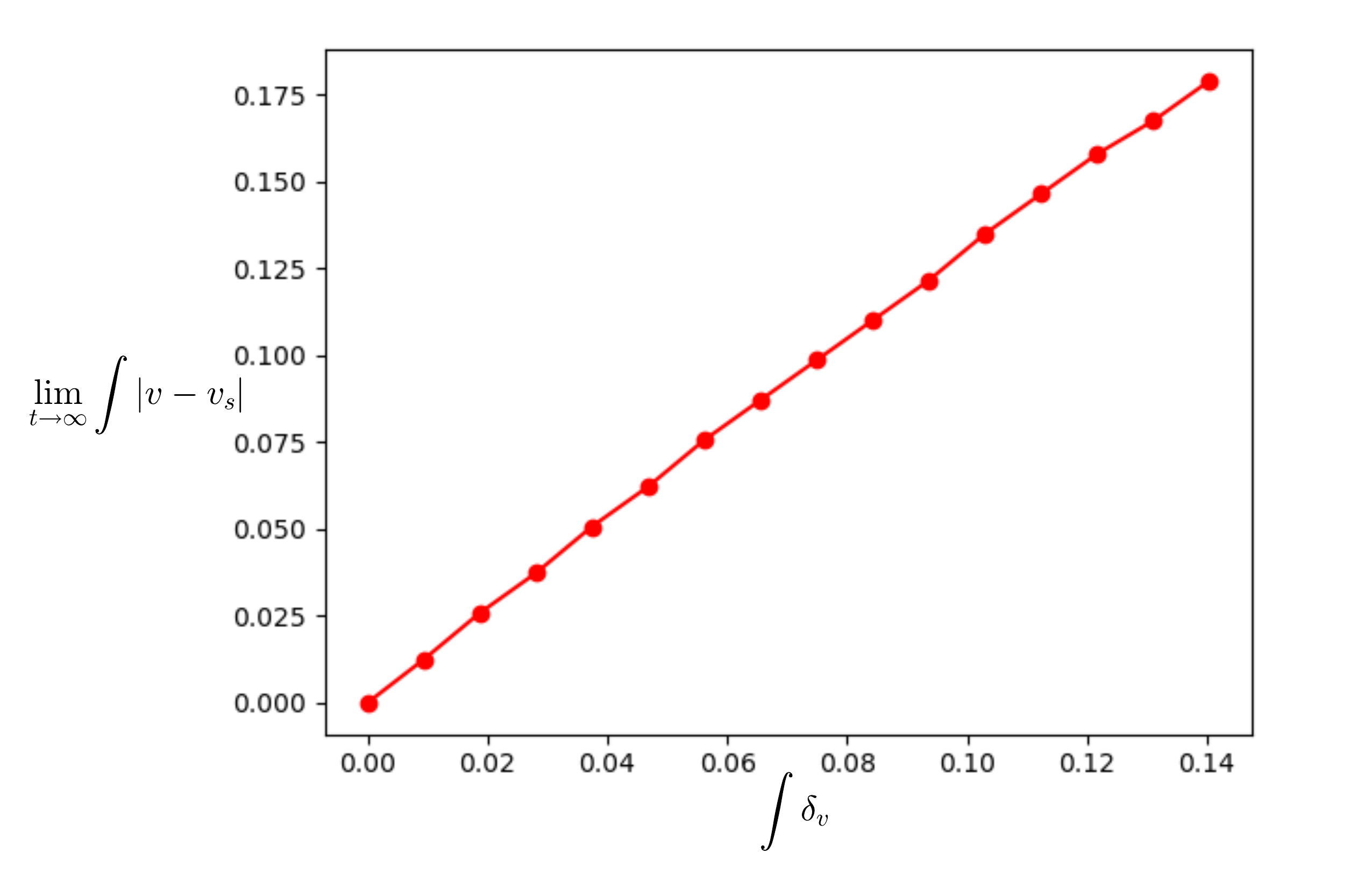}
	\caption{Burgers-Schwarzschild model with the initial condition \eqref{testB8}-\eqref{testB3}-\eqref{delta}: values of $\lim_{t\to\infty} \int |v -v_{s}|$ as a function of
 $\int \delta_{v}$.}
	\label{fig:integralperturbationvsintegralsteadyshockperturbedlabelled}
\end{figure}

In order to study the relationship between the amplitude of the perturbation and the distance between the initial and the final shock locations,
we consider the family of initial conditions:
\bel{testB8}
\tilde v_0 (r) = v_0(r) + \delta_v(\alpha, r),
\ee
where $v_0$ is given again by \eqref{testB3} and
\bel{delta}
\delta_v(\alpha, r) = \begin{cases}
\alpha \cos(5\pi r-12\pi)e^{-200(r-2.8)^{2}},  & \text{ $2.7<r<2.9$,}\\
0, & \text{ otherwise},
\end{cases}
\ee
with $\alpha>0$. 
The amplitude of the perturbation is measured by 
$
\int \delta_v(\alpha, r) \,dr
$
and the distance between the shocks are measured by
$
\lim_{t \to \infty} \int |v(r,t) - v_0(r)| \, dr.
$
See Figure \ref{fig:Areas_with_alpha_1}. Table \ref{tab:Areas_for_different_values_of_alpha}  and Figure \ref{fig:integralperturbationvsintegralsteadyshockperturbedlabelled} show the relationship between those magnitudes that is clearly linear.

\begin{table}[h]
	\centering
	\begin{tabular}{||c|c|c||} 
		\hline
		$\alpha$ & $\int \delta_{v}$ & $\lim_{t\to\infty} \int |v -v_{s}|$\\ 
		\hline
		0.0 & 0.00000 & 0.00000 \\ 
		\hline
		0.1 & 0.00936 & 0.01245 \\ 
		\hline
		0.2 & 0.01873 & 0.02586 \\
		\hline
		0.3 & 0.02809 & 0.03735 \\
		\hline
		0.4 & 0.03745 & 0.05076  \\
		\hline
		0.5 & 0.04682 & 0.06225 \\ 
		\hline
		0.6 & 0.05618 & 0.07566 \\ 
		\hline
		0.7 & 0.06554 & 0.08715 \\
		\hline
		0.8 & 0.07491 & 0.09864 \\
		\hline
		0.9 & 0.08427 & 0.11013  \\
		\hline
		1.0 & 0.09364 & 0.12163 \\ 
		\hline
		1.1 & 0.10300 & 0.13503 \\ 
		\hline
		1.2 & 0.11236 & 0.14653 \\
		\hline
		1.3 & 0.12172 & 0.15802 \\
		\hline
		1.4 & 0.13109 & 0.16760  \\
		\hline
		1.5 & 0.14045 & 0.17909 \\
		\hline
	\end{tabular}
	\caption{Burgers-Schwarzschild model with the initial condition \eqref{testB8}-\eqref{testB3}-\eqref{delta}: measures of the perturbation and the shock displacement for 
different values of $\alpha$.}
	\label{tab:Areas_for_different_values_of_alpha}
\end{table}

%


\subsection{Long-time behavior of the solutions}

In  this section we consider different  initial conditions and investigate the long-time behavior of the corresponding solutions using the first-order well-balanced scheme. A large number of tests have been performed with the first-order methods (that is the less costly one) considering different initial conditions,  different meshes,  and different lengths of the computational domain:  the observed behavior of the numerical solutions have been always one of the four ones shown here depending on the value at $2M$ (1 or lower) and at the right boundary (positive or negative).

\begin{enumerate}

\item{Initial condition satisfying $v_{0}(2M)=1$ and $v_{0}(L)\geq 0$:}
let us consider the  initial condition
\bel{testB9}
v_{0}(r) = \begin{cases}
1,  & \text{   $2<r<2.1$},\\
\cos(30r)e^{\frac{-1}{(x-2.5)^{2}}}, & \text{ otherwise},
\end{cases}
\ee
that takes value 1 in a neighborhood of $2M =2$ and a positive value at the right boundary of the computational domain $x = 4$. As it can be observed in Figure \ref{fig:ko1_wb_testOnePositiveNonSteady} after a transient regime,  the numerical solution takes the form of a right-moving shock  linking the stationary solution
$v \equiv 1$ with the negative stationary solution that takes value $-1$ at $x= 2M$ and value 0 at $x = 4$. Once this shock leaves the domain, the stationary solution
$v\equiv 1$ is reached in the whole computational domain.

\begin{figure}[h]
	\begin{subfigure}[h]{0.5\textwidth}
		\centering
		\includegraphics[width=1\linewidth]{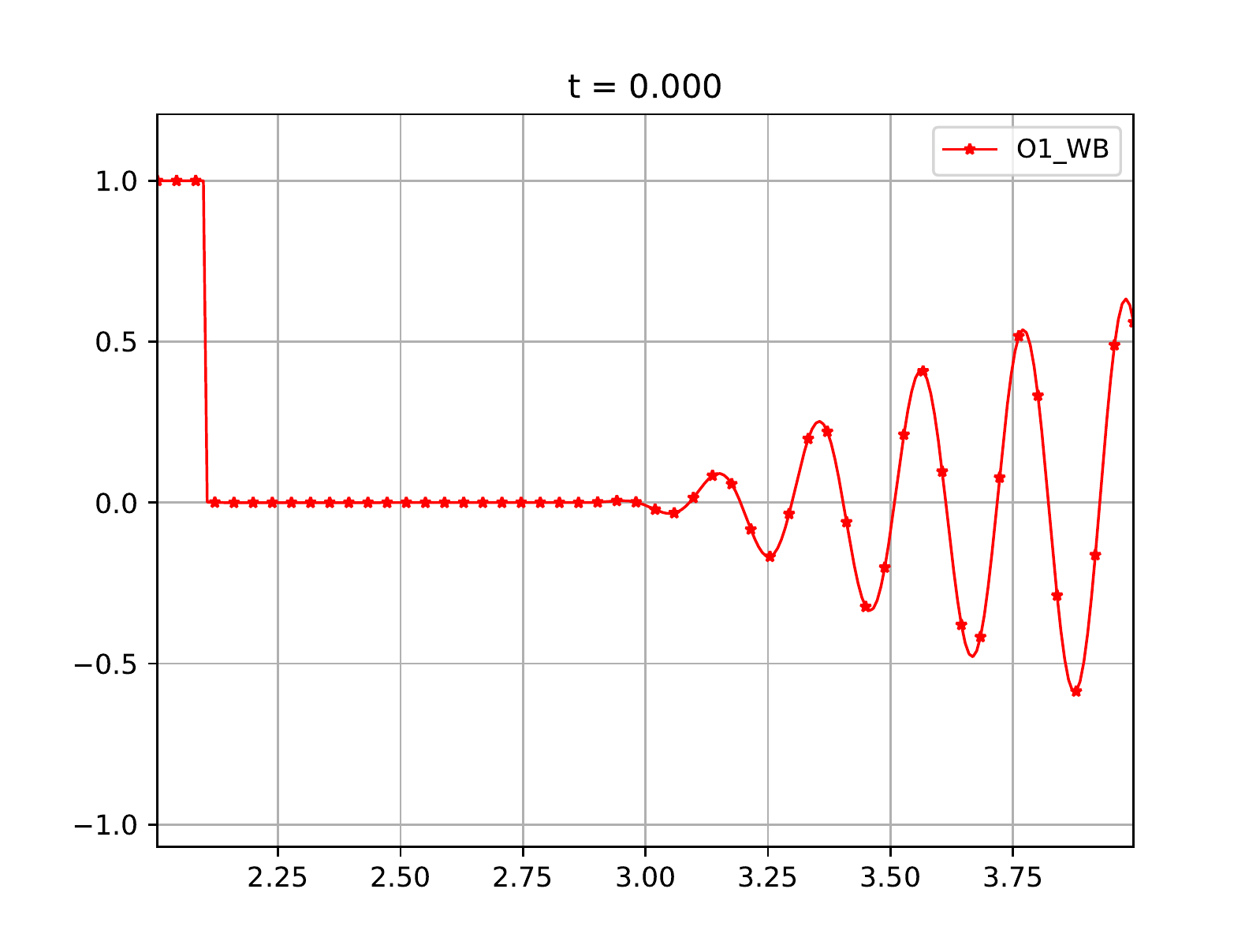}
		\label{fig:ko1_wb_testOnePositiveNonSteady_t_0}
	\end{subfigure}
	\begin{subfigure}[h]{0.5\textwidth}
		\centering
		\includegraphics[width=1\linewidth]{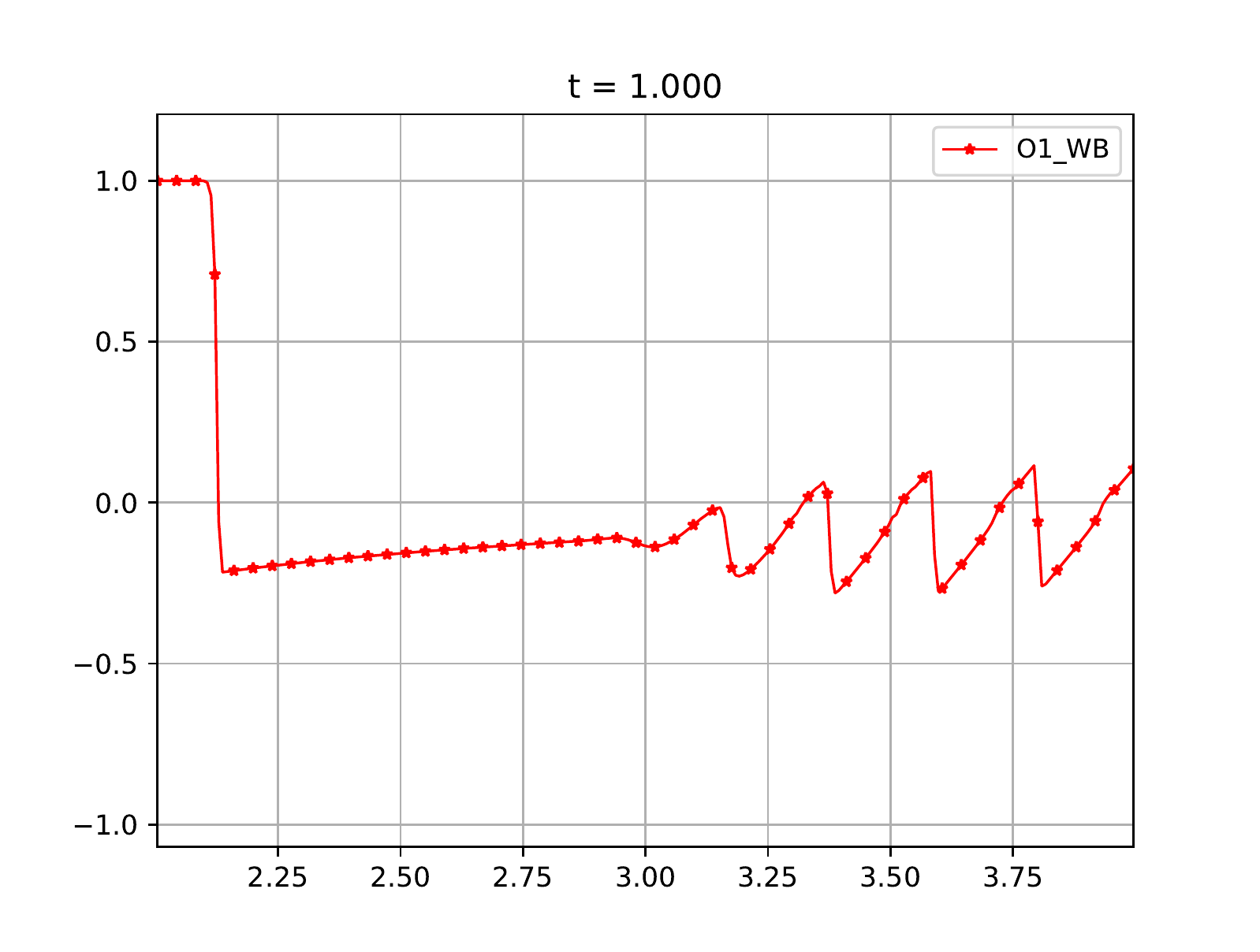}
		\label{fig:ko1_wb_testOnePositiveNonSteady_t_1}
	\end{subfigure}
	\begin{subfigure}[h]{0.5\textwidth}
		\centering
		\includegraphics[width=1\linewidth]{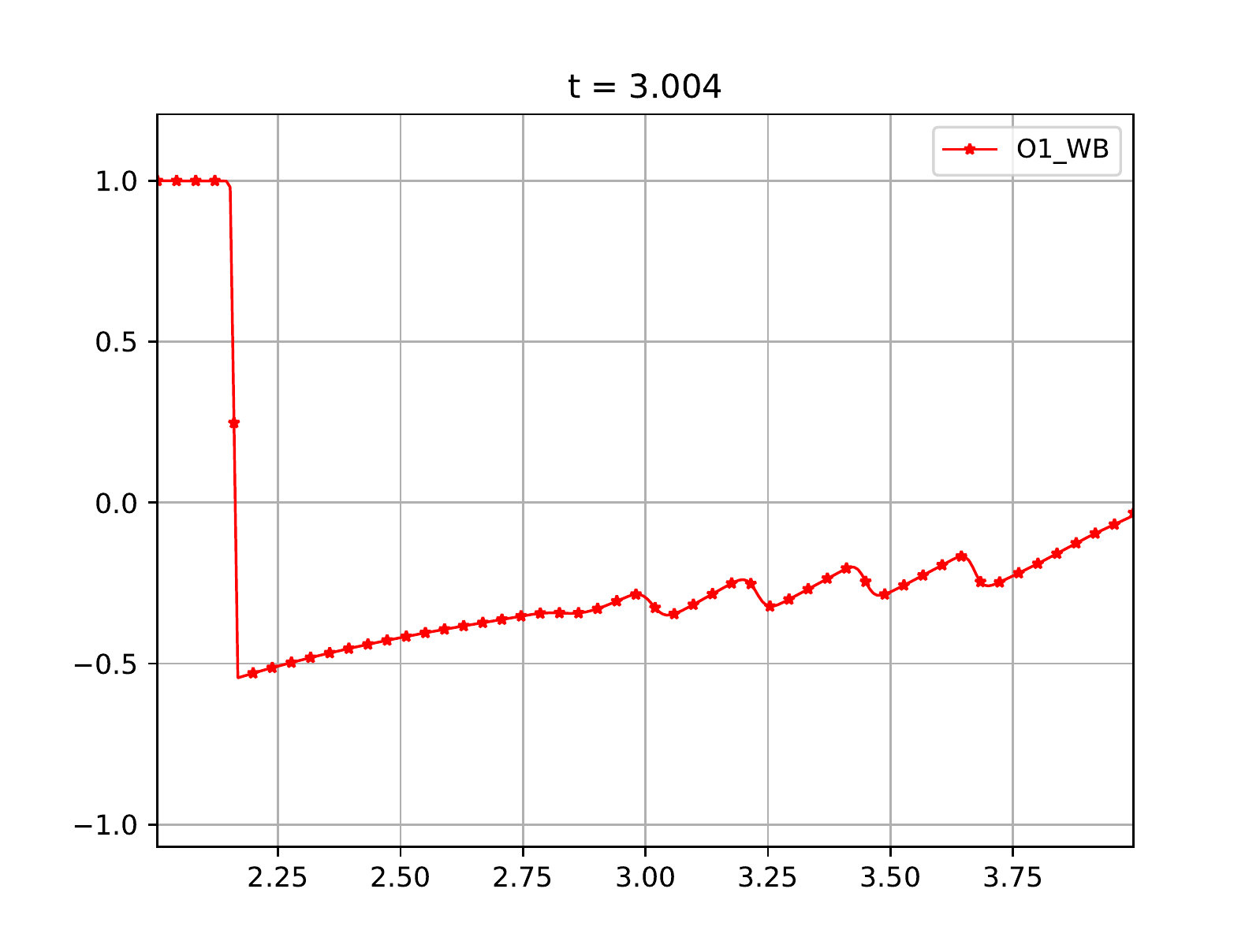}
		\label{fig:ko1_wb_testOnePositiveNonSteady_t_3}
	\end{subfigure}
	\begin{subfigure}[h]{0.5\textwidth}
		\centering
		\includegraphics[width=1\linewidth]{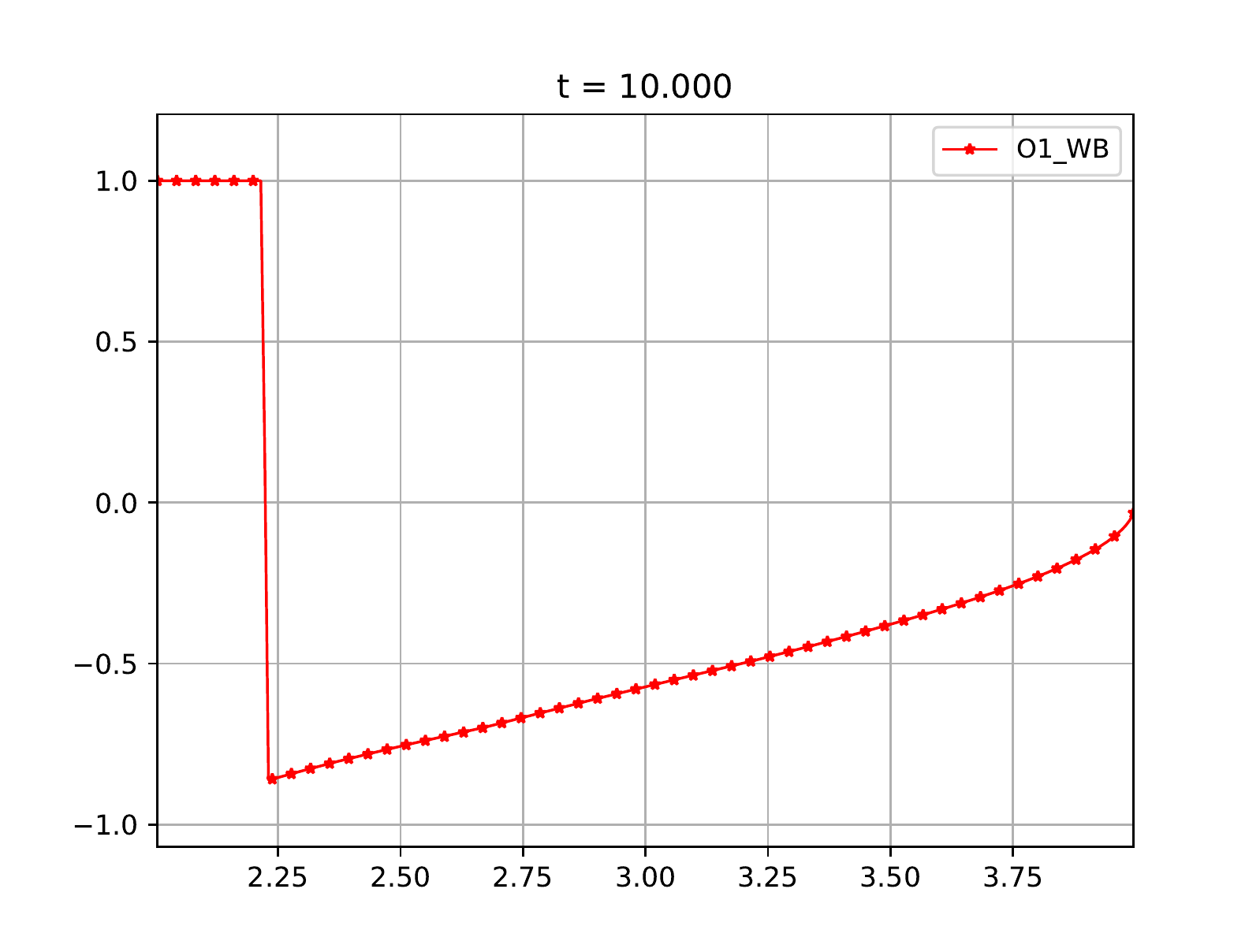}
		\label{fig:ko1_wb_testOnePositiveNonSteady_t_10}
	\end{subfigure}
	\begin{subfigure}[h]{0.5\textwidth}
		\centering
		\includegraphics[width=1\linewidth]{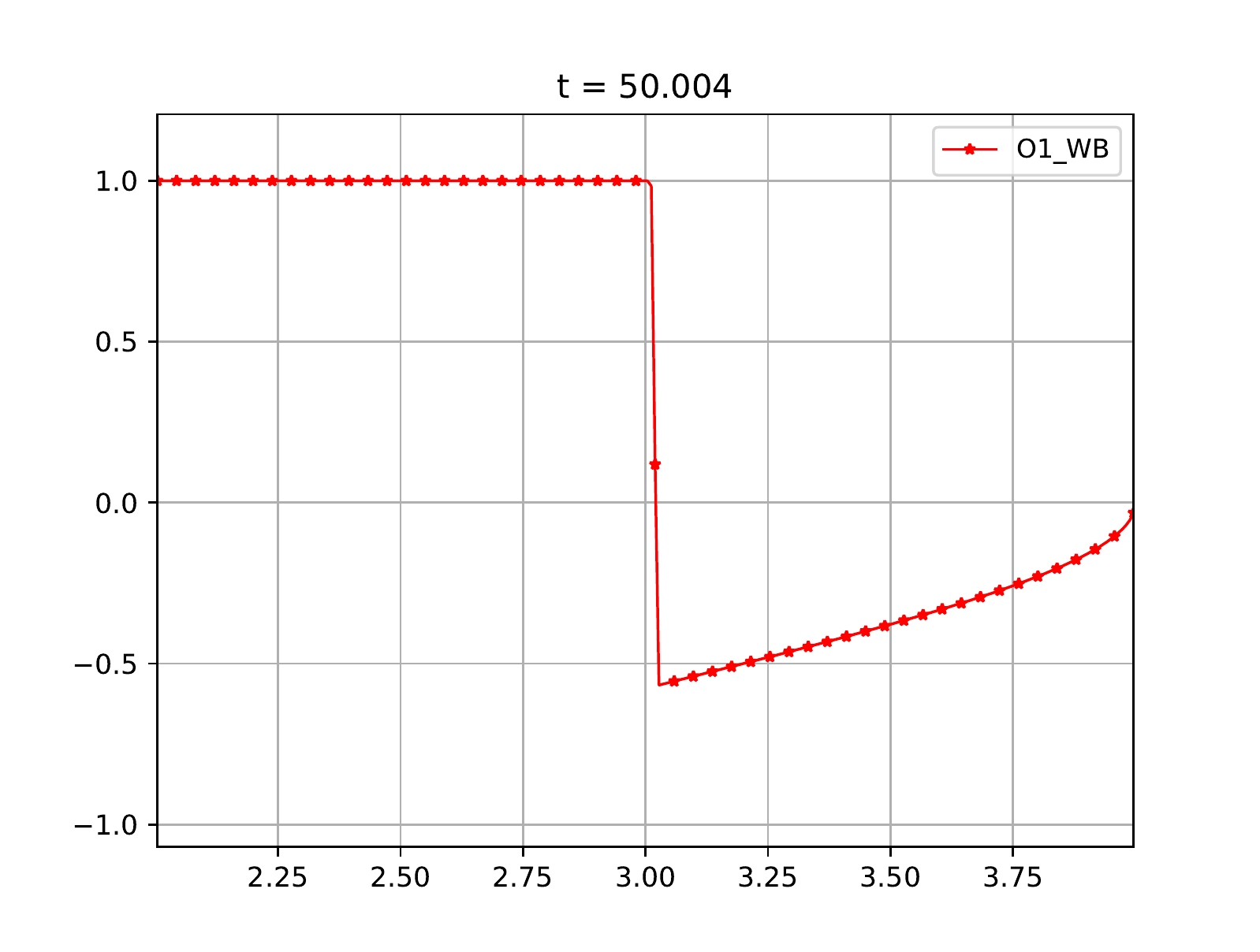}
		\label{fig:ko1_wb_testOnePositiveNonSteady_t_50}
	\end{subfigure}
	\begin{subfigure}[h]{0.5\textwidth}
		\centering
		\includegraphics[width=1\linewidth]{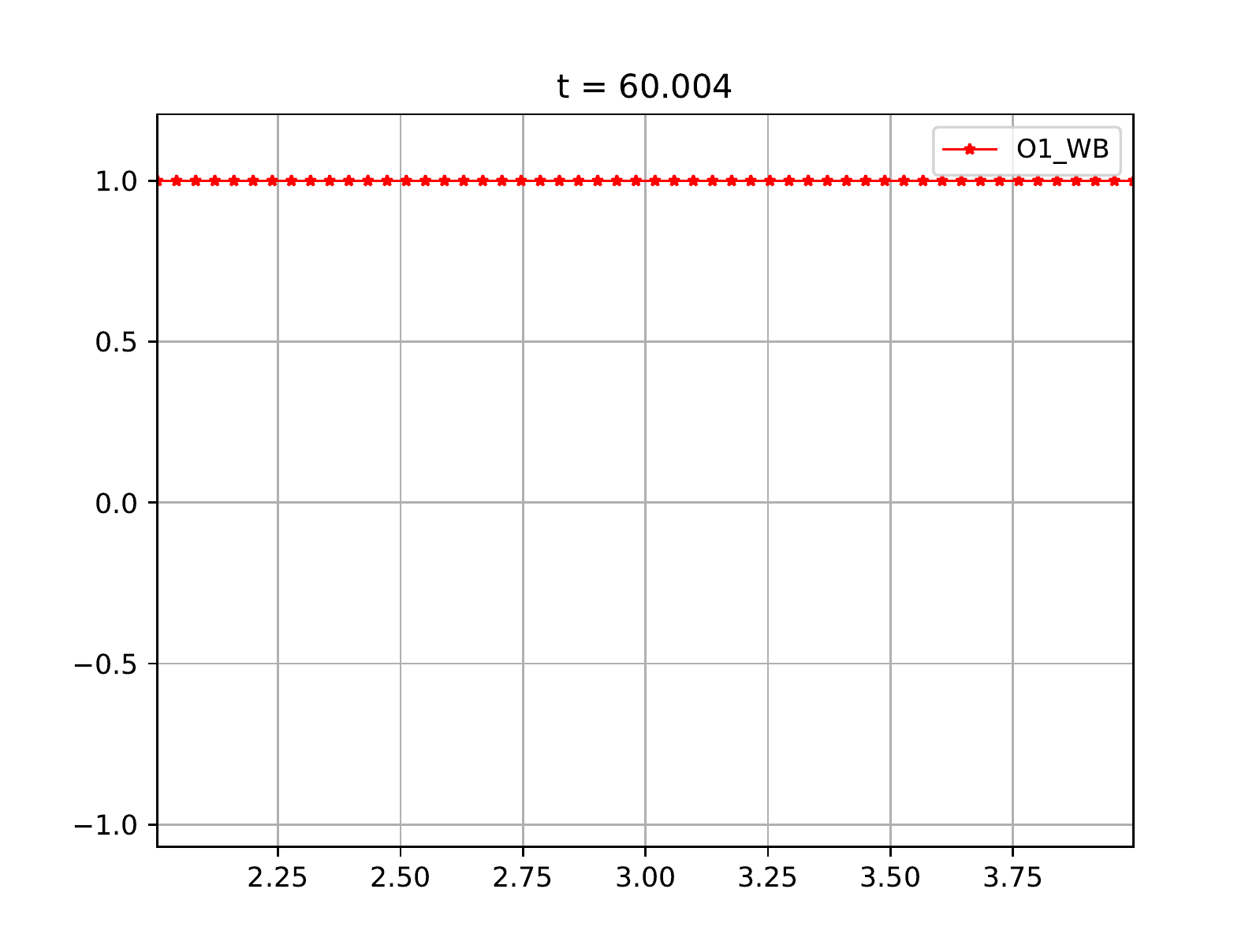}
		\label{fig:ko1_wb_testOnePositiveNonSteady_t_60}
	\end{subfigure}
	\caption{Burgers-Schwarzschild model with the initial condition \eqref{testB9}: first-order well-balanced scheme at selected times.}
	\label{fig:ko1_wb_testOnePositiveNonSteady}
\end{figure}


\item  Initial condition satisfying $v_{0}(2M)=1$ and $v_{0}(L)< 0$: we consider now the
\bel{testB10}
v_{0}(r) = \begin{cases}
1,  & \text{ $2<r<2.1$,}\\
\cos(20r)e^{\frac{-1}{(x-2.5)^{2}}}, & \text{ otherwise,}
\end{cases}
\ee
that takes value 1 in a neigborhood of $2M =2$ and negative value at the right boundary of the computational domain $x = 4$. As it can be observed in Figure  \ref{fig:ko1_wb_testOneNegativeNonSteady}  after a transient period,  the numerical solution takes the form of a right-moving shock  linking the stationary solution
$v \equiv 1$ with the negative stationary solution that takes value $-1$ at $x= 2M$ and value $v_0(4)$ at $x = 4$. Once this shock leaves the domain, the stationary solution
$v\equiv 1$ is reached in the whole computational domain. 
\begin{figure}[h]
	\begin{subfigure}[h]{0.5\textwidth}
		\centering
		\includegraphics[width=1\linewidth]{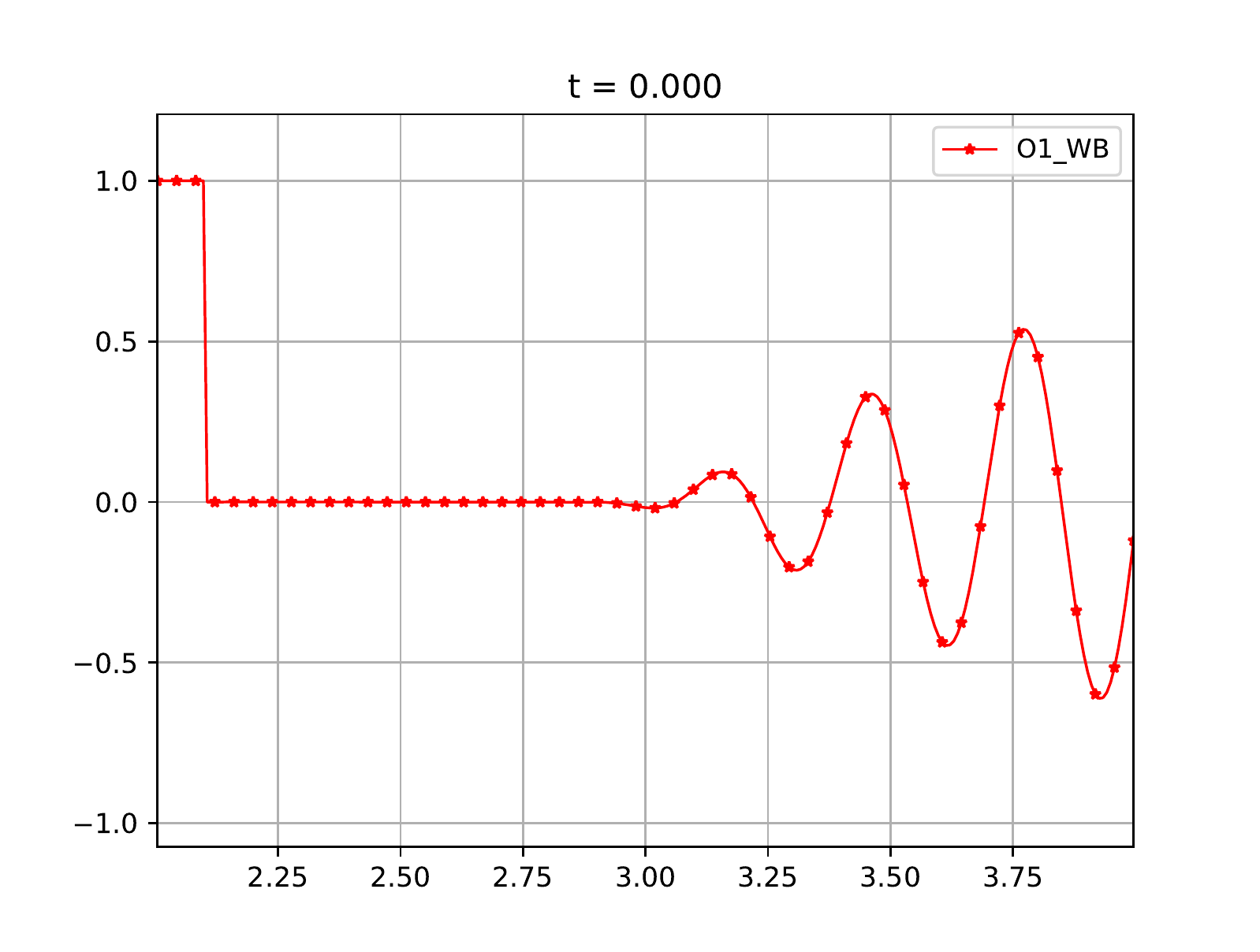}
		\label{fig:ko1_wb_testOneNegativeNonSteady_t_0}
	\end{subfigure}
	\begin{subfigure}[h]{0.5\textwidth}
		\centering
		\includegraphics[width=1\linewidth]{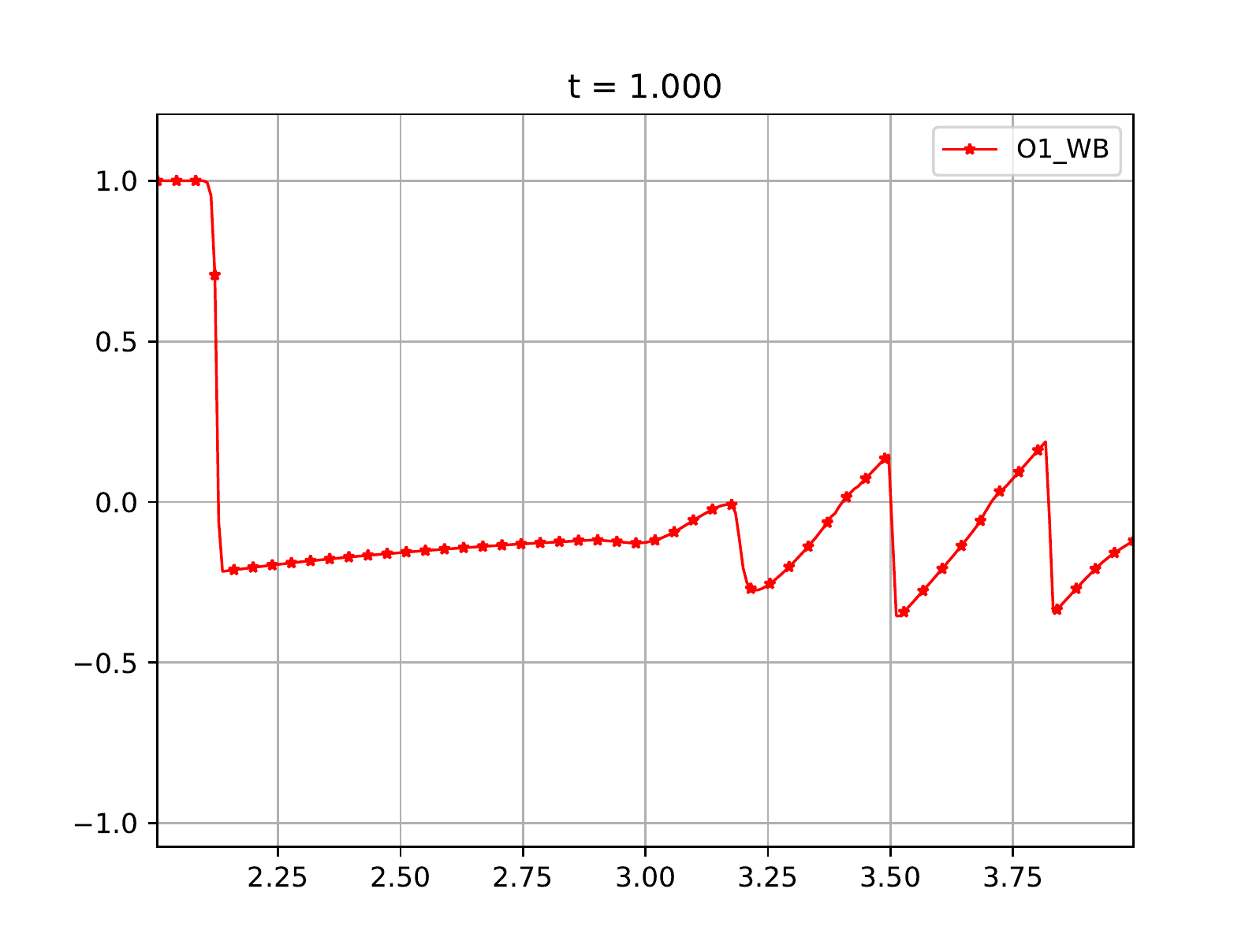}
		\label{fig:ko1_wb_testOneNegativeNonSteady_t_1}
	\end{subfigure}
	\begin{subfigure}[h]{0.5\textwidth}
		\centering
		\includegraphics[width=1\linewidth]{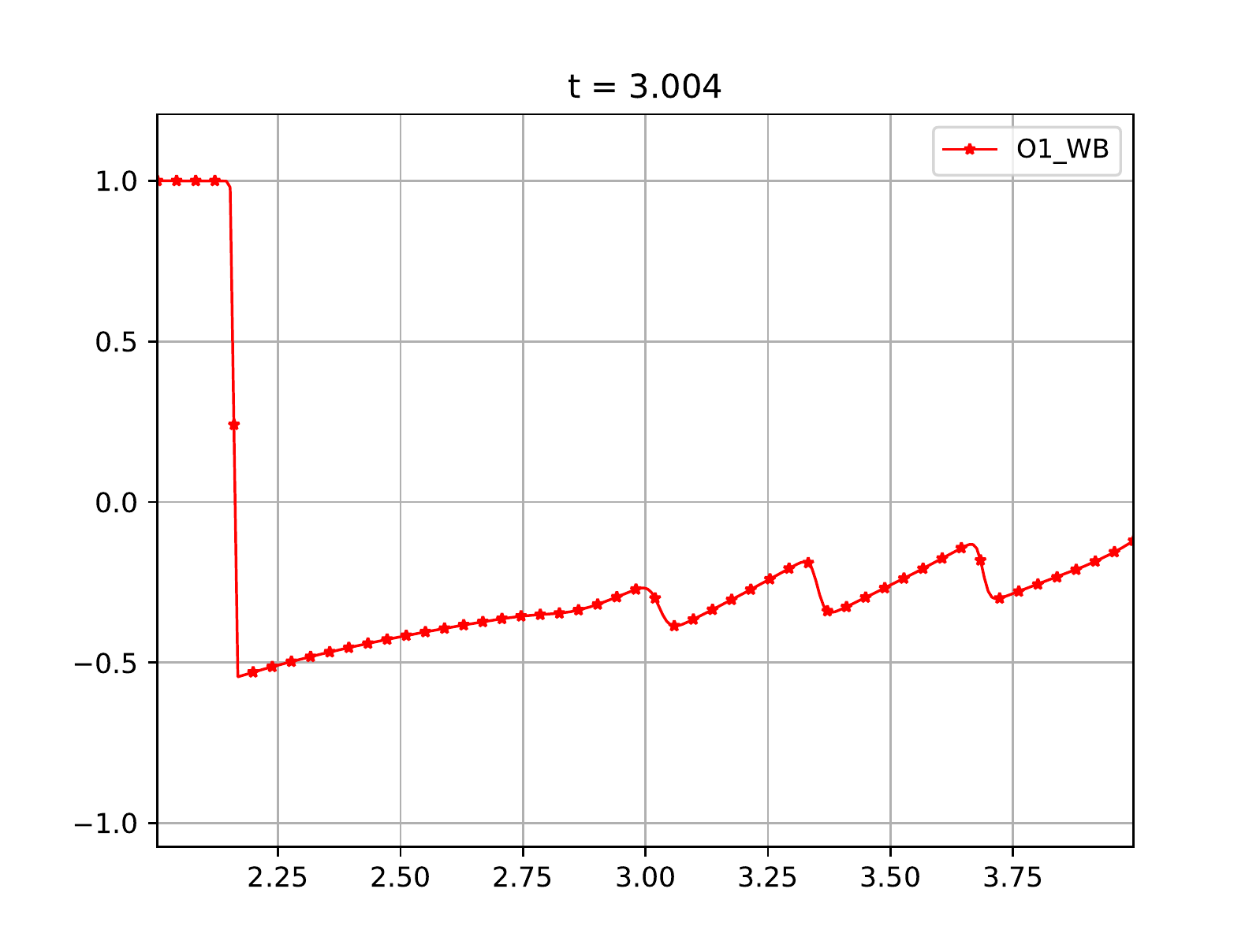}
		\label{fig:ko1_wb_testOneNegativeNonSteady_t_3}
	\end{subfigure}
	\begin{subfigure}[h]{0.5\textwidth}
		\centering
		\includegraphics[width=1\linewidth]{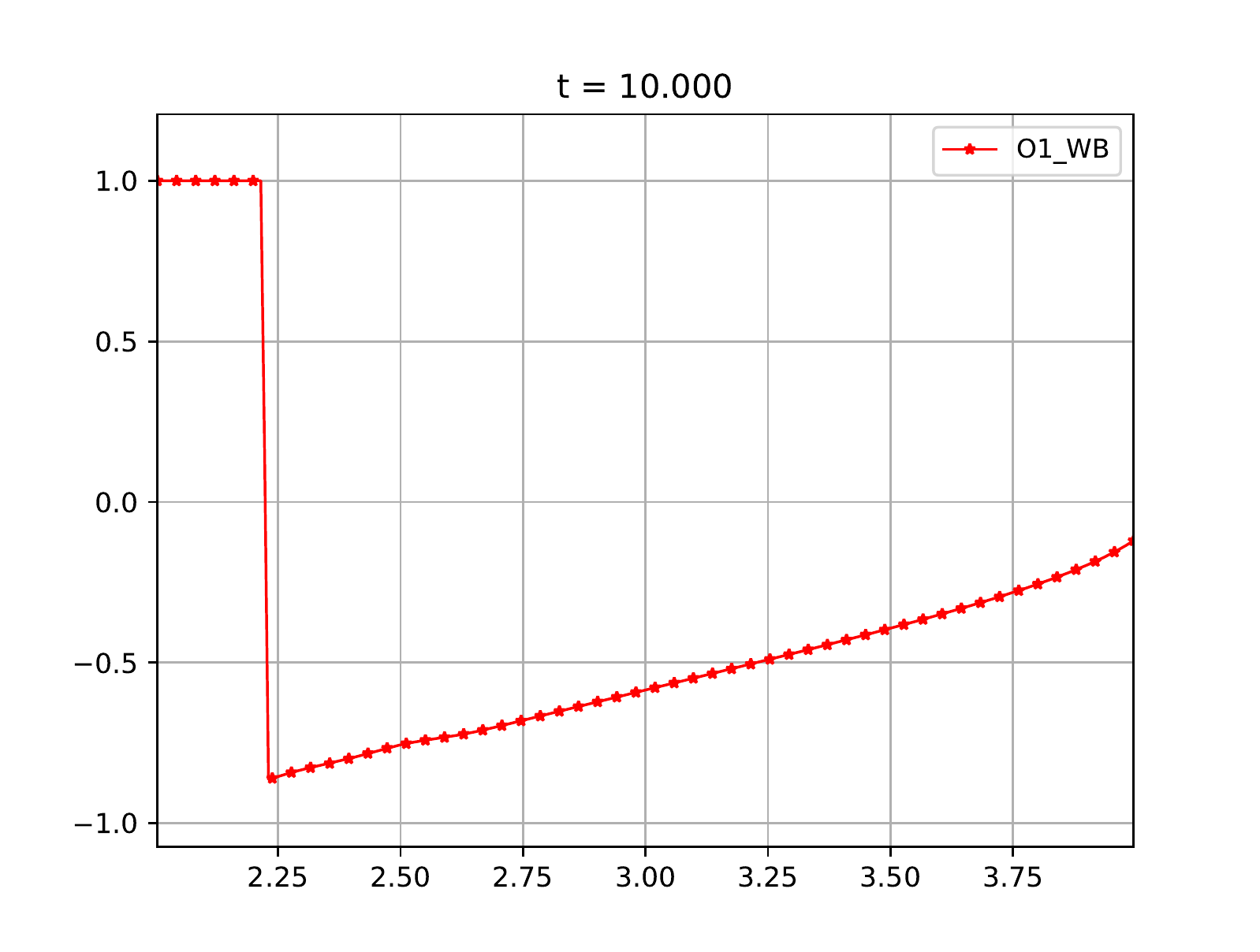}
		\label{fig:ko1_wb_testOneNegativeNonSteady_t_10}
	\end{subfigure}
	\begin{subfigure}[h]{0.5\textwidth}
		\centering
		\includegraphics[width=1\linewidth]{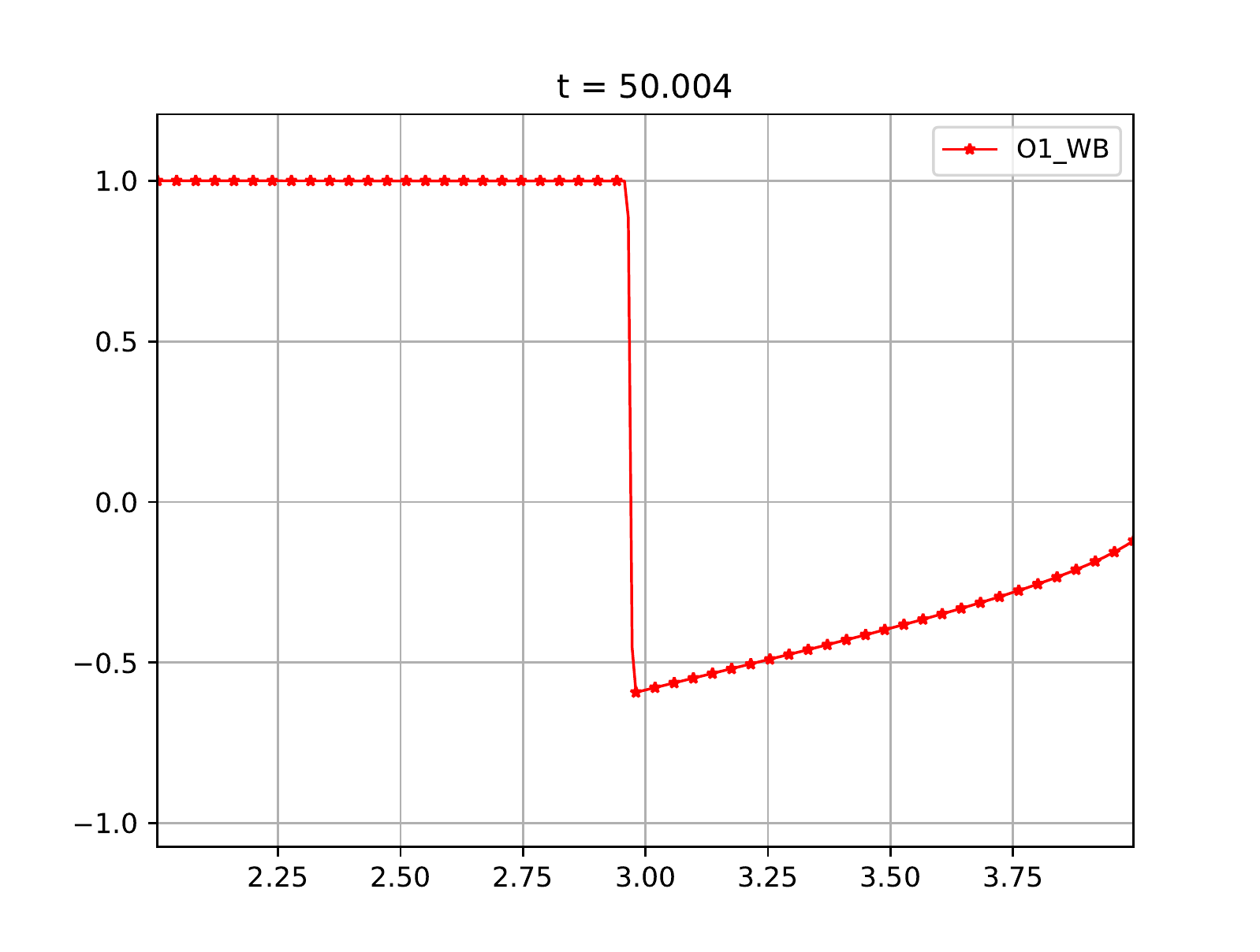}
		\label{fig:ko1_wb_testOneNegativeNonSteady_t_50}
	\end{subfigure}
	\begin{subfigure}[h]{0.5\textwidth}
		\centering
		\includegraphics[width=1\linewidth]{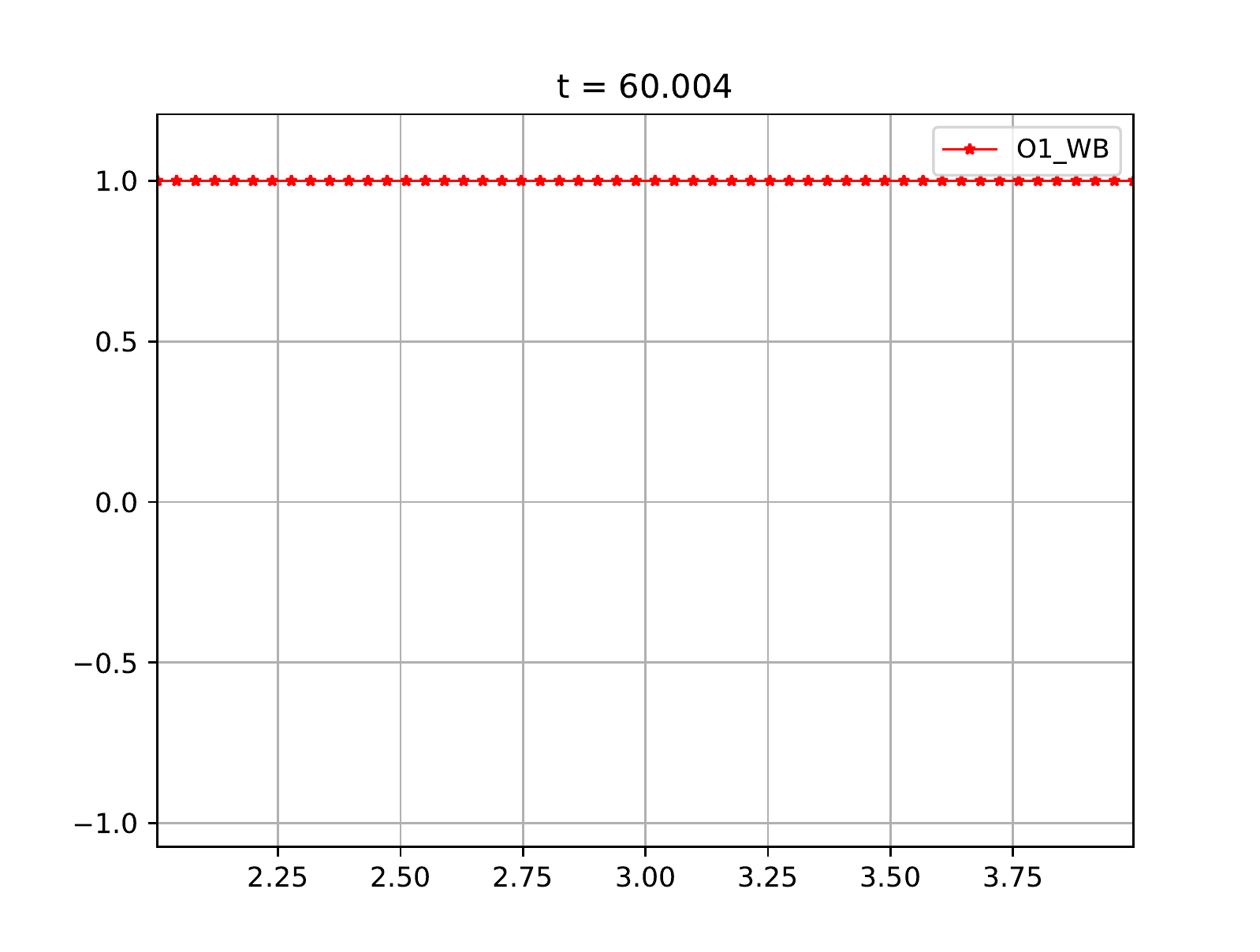}
		\label{fig:ko1_wb_testOneNegativeNonSteady_t_60}
	\end{subfigure}
	\caption{Burgers-Schwarzschild model with the initial condition \eqref{testB10}: first-order well-balanced scheme at selected times.}
	\label{fig:ko1_wb_testOneNegativeNonSteady}
\end{figure}


\item Initial condition satisfying $v_{0}(2M)<1$ and $v_{0}(L)\geq 0$:  we consider now the  initial condition
\bel{testB11}
v_{0}(r) = \begin{cases}
0.8, & \text{  $ 2<r<2.1$,}\\
\cos(30r)e^{\frac{-1}{(x-2.5)^{2}}}, & \text{ otherwise}.
\end{cases}
\ee
In this case the numerical solution reaches in finite time  the negative stationary solution $v^*$ such that $v^*(2) = -1$ and $v^*(4) = 0$: see  Figure \ref{fig:ko1_wb_testNonOnePositiveNonSteady}. 
\begin{figure}[h]
	\begin{subfigure}[h]{0.45\textwidth}
		\centering
		\includegraphics[width=1\linewidth]{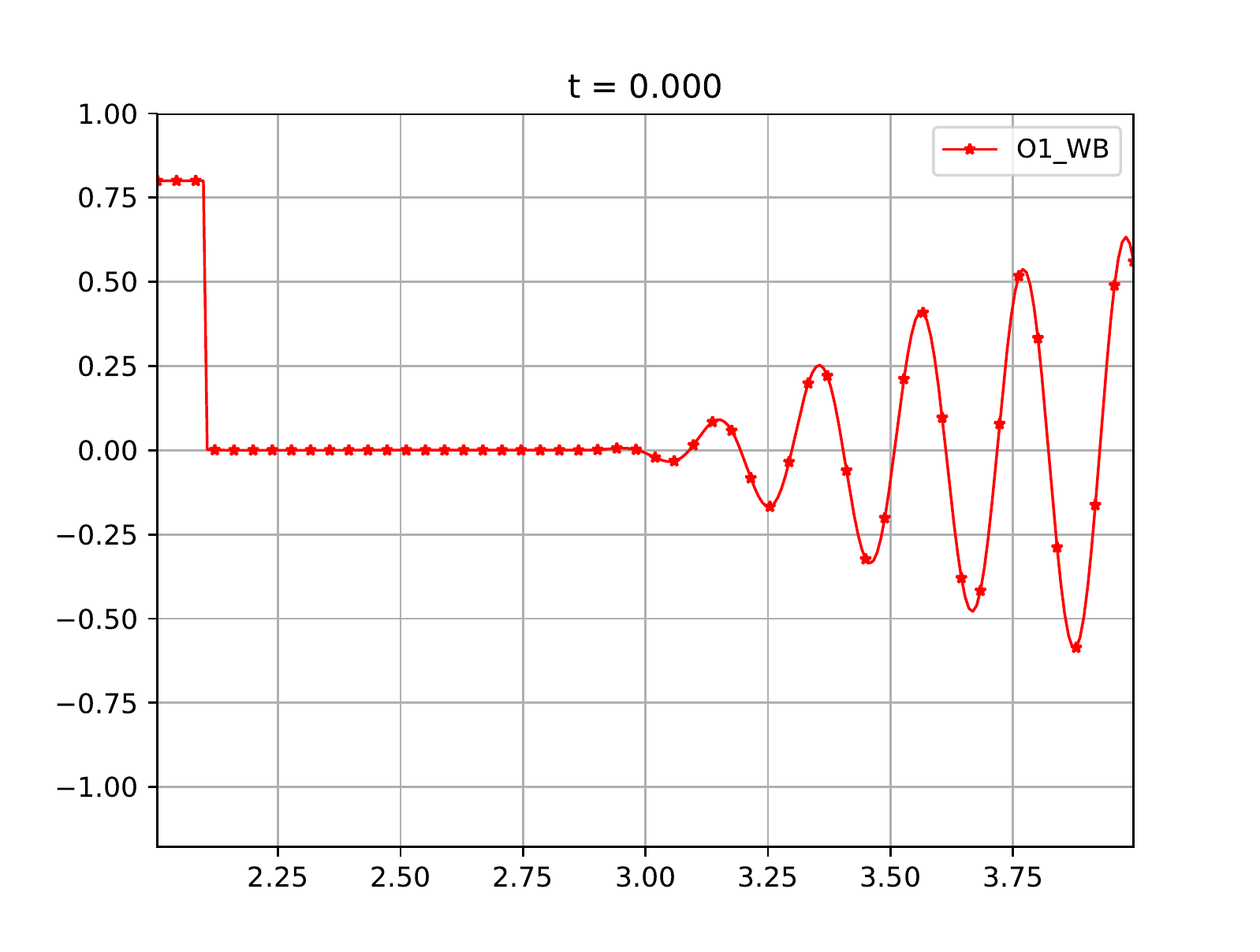}
		\label{fig:ko1_wb_testNonOnePositiveNonSteady_t_0}
	\end{subfigure}
	\begin{subfigure}[h]{0.45\textwidth}
		\centering
		\includegraphics[width=1\linewidth]{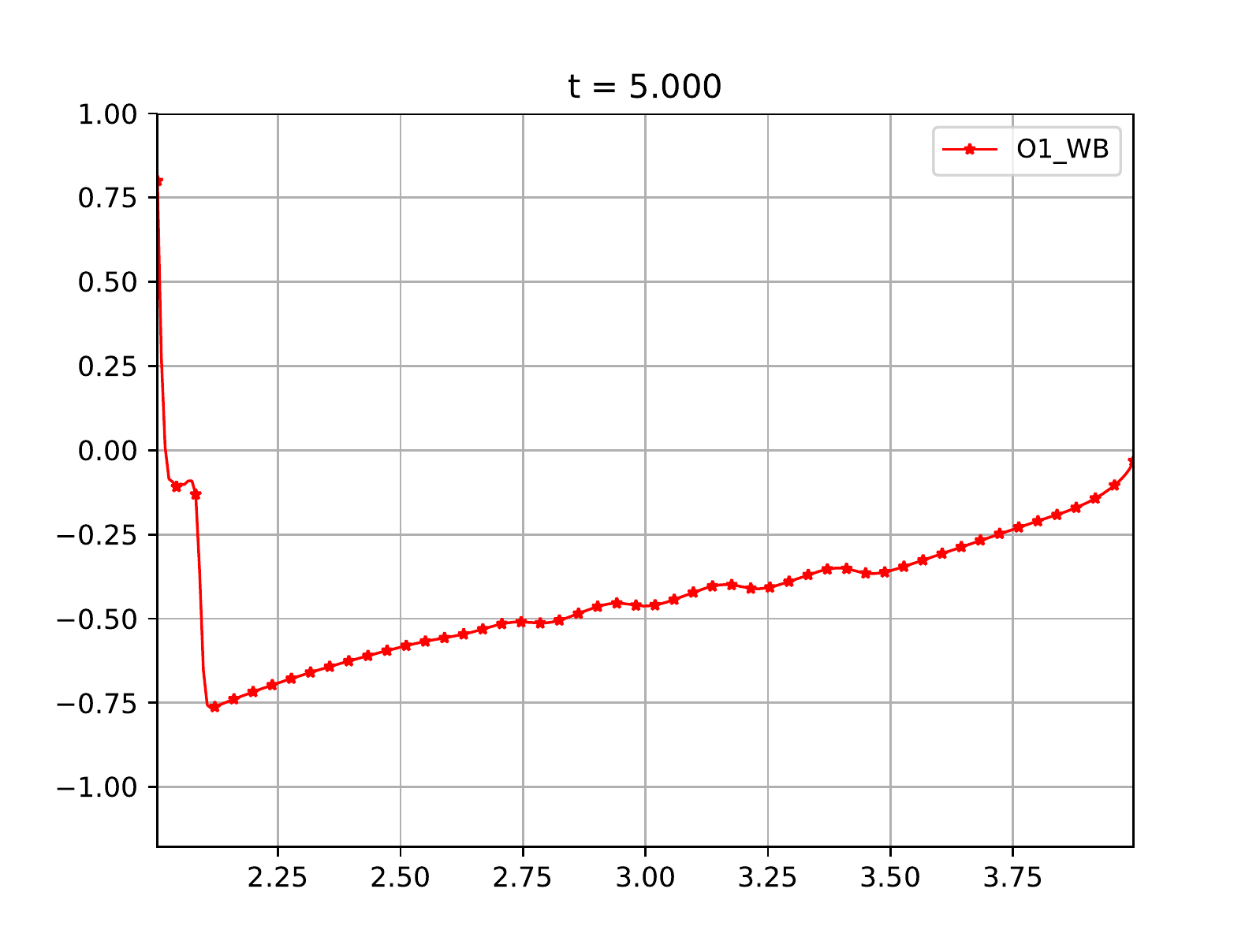}
		\label{fig:ko1_wb_testNonOnePositiveNonSteady_t_5}
	\end{subfigure}
	\begin{subfigure}[h]{0.45\textwidth}
		\centering
		\includegraphics[width=1\linewidth]{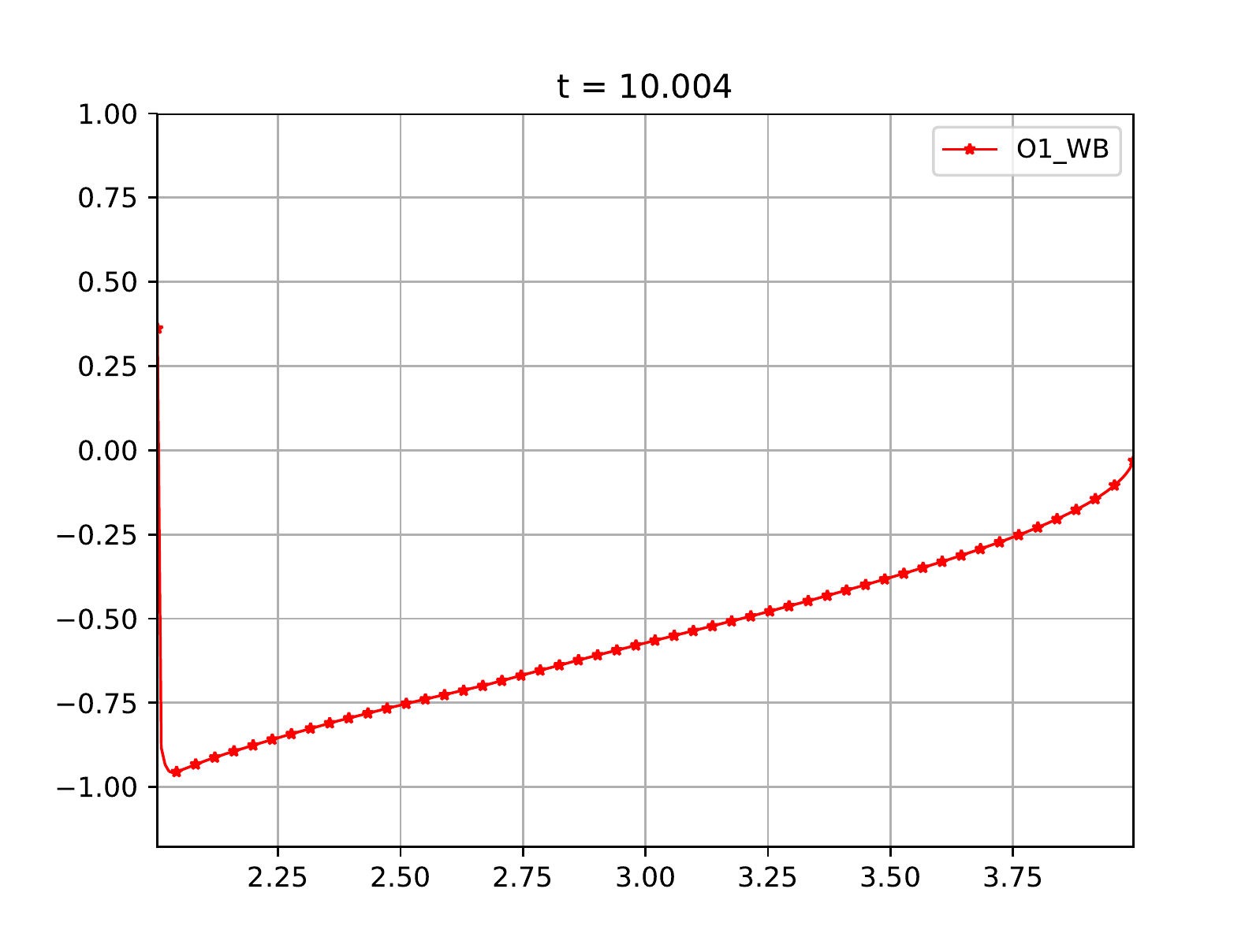}
		\label{fig:ko1_wb_testNonOnePositiveNonSteady_t_10}
	\end{subfigure}
		\quad	\quad	\quad	\quad
	\begin{subfigure}[h]{0.45\textwidth}
		\centering
		\includegraphics[width=1\linewidth]{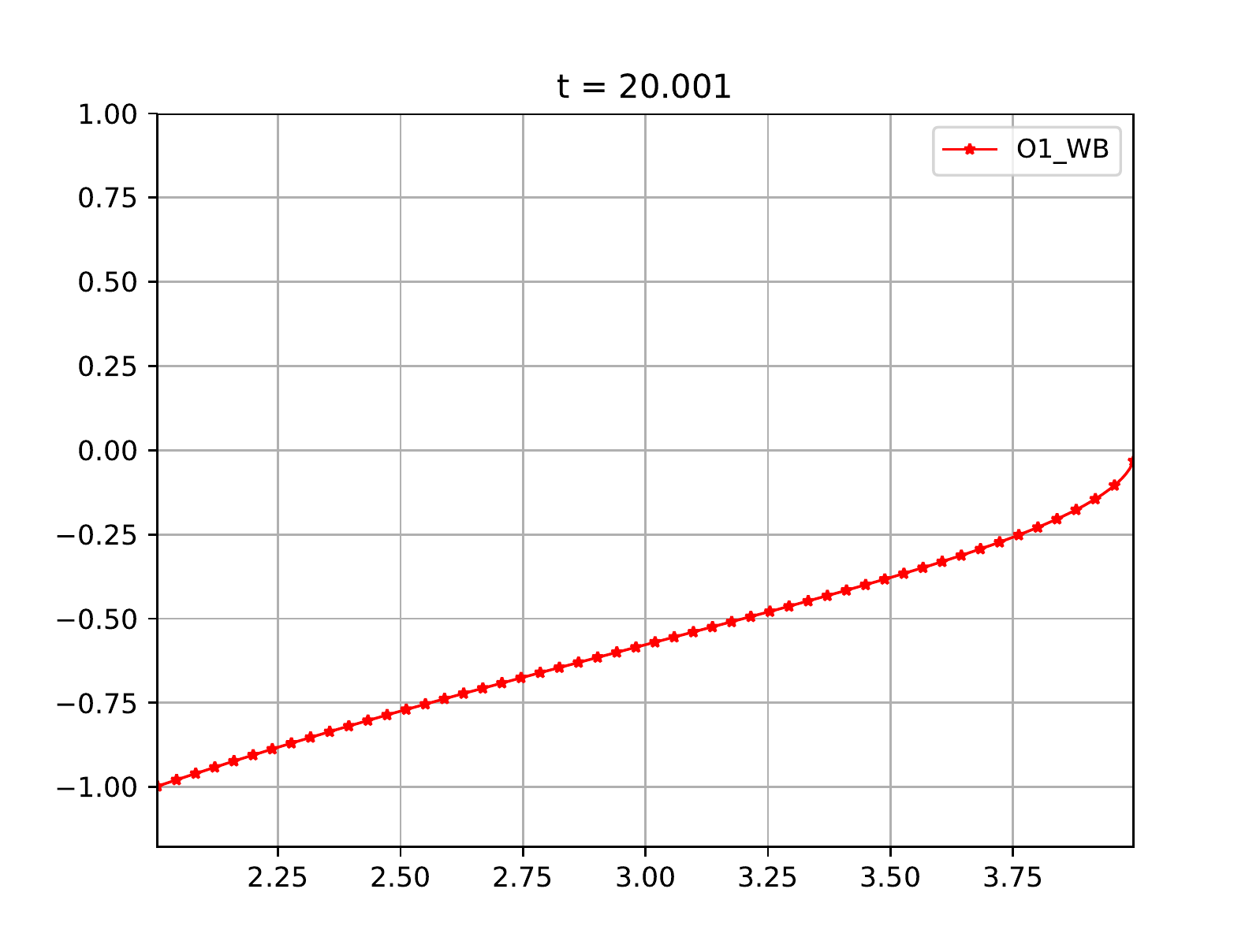}
		\label{fig:ko1_wb_testNonOnePositiveNonSteady_t_20}
	\end{subfigure}
	\caption{Burgers-Schwarzschild model with the initial condition \eqref{testB11}:
	 first-order well-balanced scheme at selected times.}
	\label{fig:ko1_wb_testNonOnePositiveNonSteady}
\end{figure}


\item Initial condition satisfying $v_{0}(2M)<1$ and $v_{0}(L)< 0$: we finally consider the   initial condition
\bel{testB12}
v_{0}(r) = \begin{cases}
0.8,  & \text{  $2<r<2.1$,}\\
\cos(20r)e^{\frac{-1}{(x-2.5)^{2}}}, & \text{ otherwise}.
\end{cases}
\ee
The numerical solution reaches in finite time  the negative stationary solution $v^*$ such that $v^*(2) = -1$ and $v^*(4) = v_0(4)$: see  Figure 
\ref{fig:ko1_wb_testNonOneNegativeNonSteady}.
\begin{figure}[h]
	\begin{subfigure}[h]{0.45\textwidth}
		\centering
		\includegraphics[width=1\linewidth]{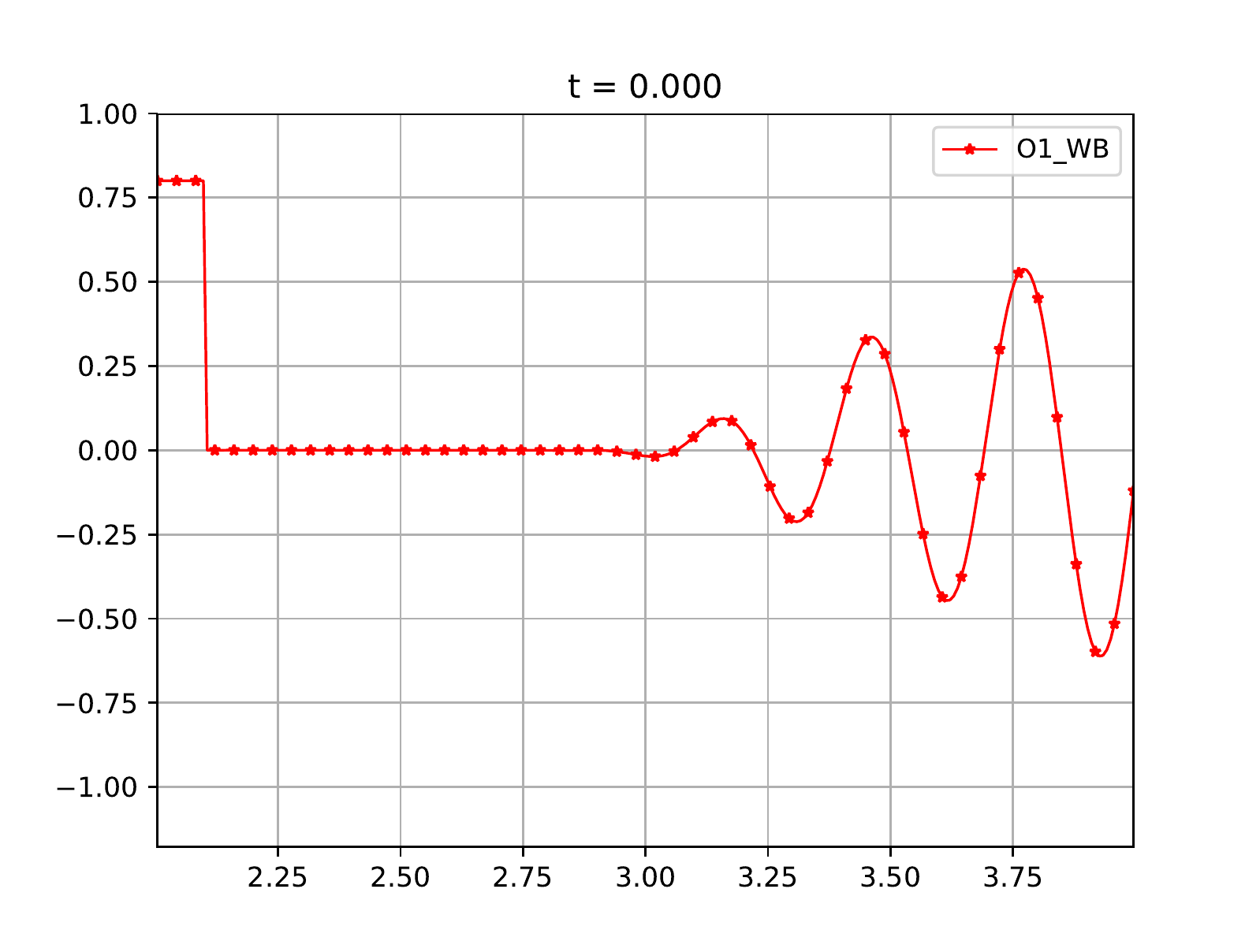}
		\label{fig:ko1_wb_testNonOneNegativeNonSteady_t_0}
	\end{subfigure}
	\begin{subfigure}[h]{0.45\textwidth}
		\centering
		\includegraphics[width=1\linewidth]{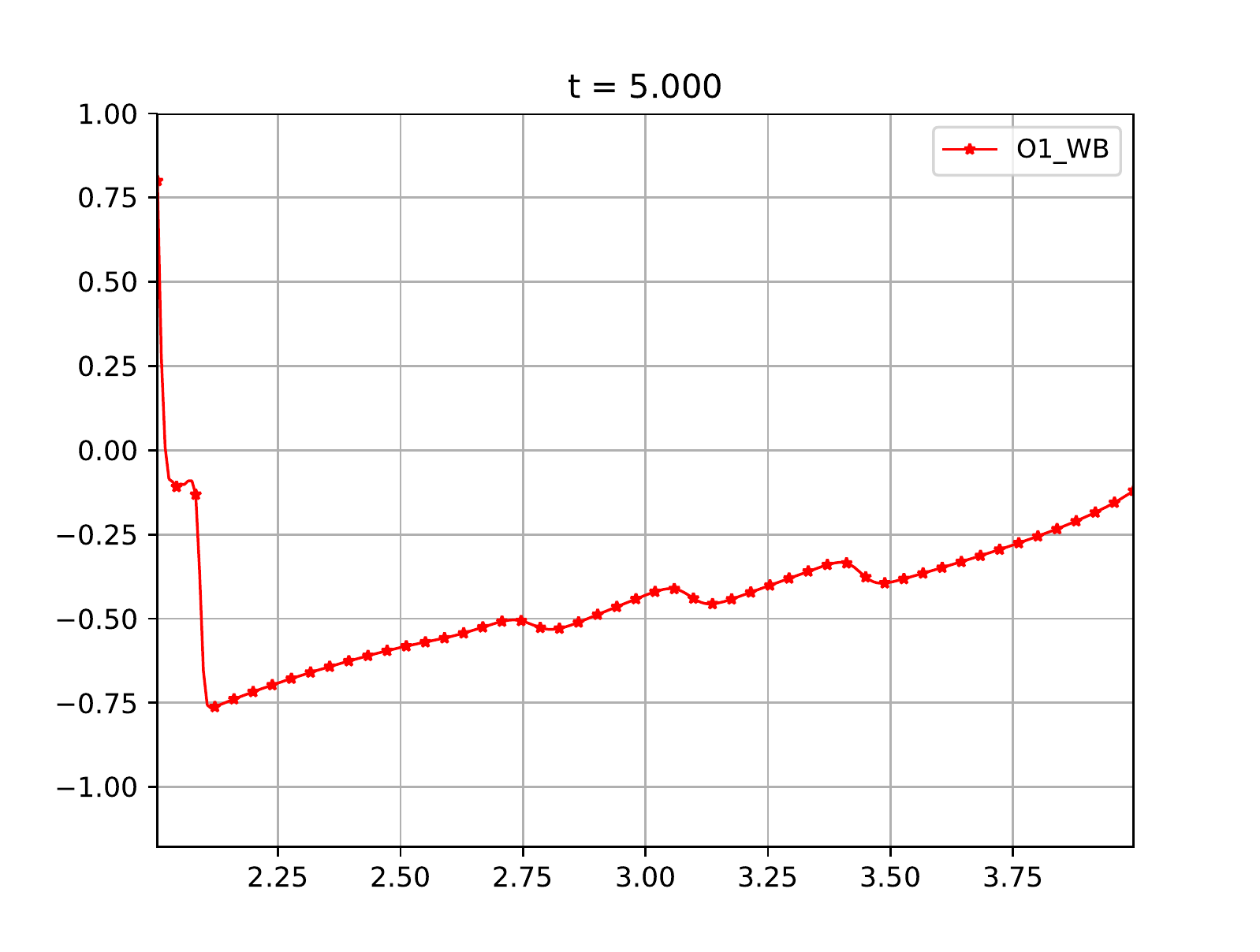}
		\label{fig:ko1_wb_testNonOneNegativeNonSteady_t_5}
	\end{subfigure}
	\begin{subfigure}[h]{0.45\textwidth}
		\centering
		\includegraphics[width=1\linewidth]{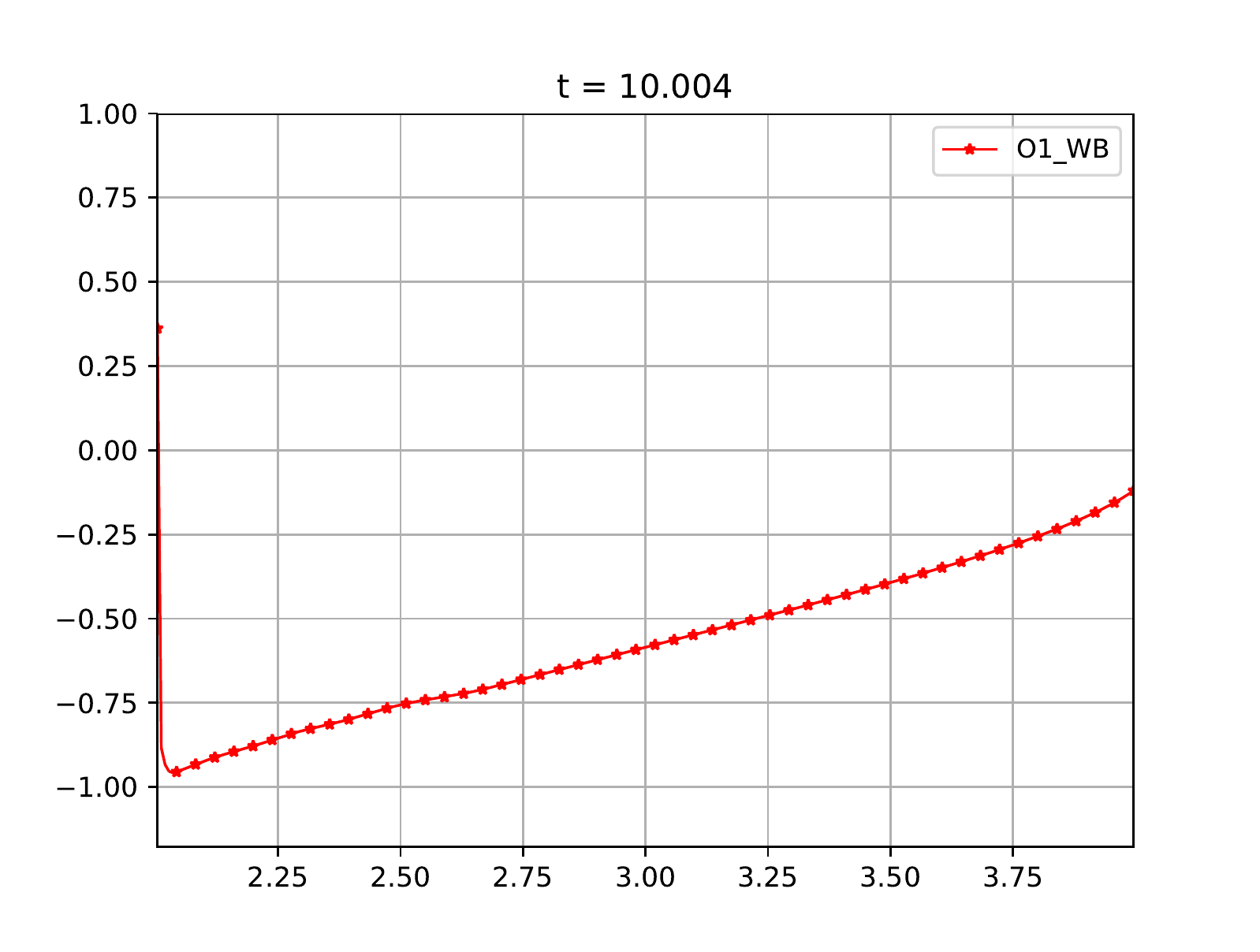}
		\label{fig:ko1_wb_testNonOneNegativeNonSteady_t_10}
	\end{subfigure}
	\quad	\quad	\quad \quad 
	\begin{subfigure}[h]{0.45\textwidth}
		\centering
		\includegraphics[width=1\linewidth]{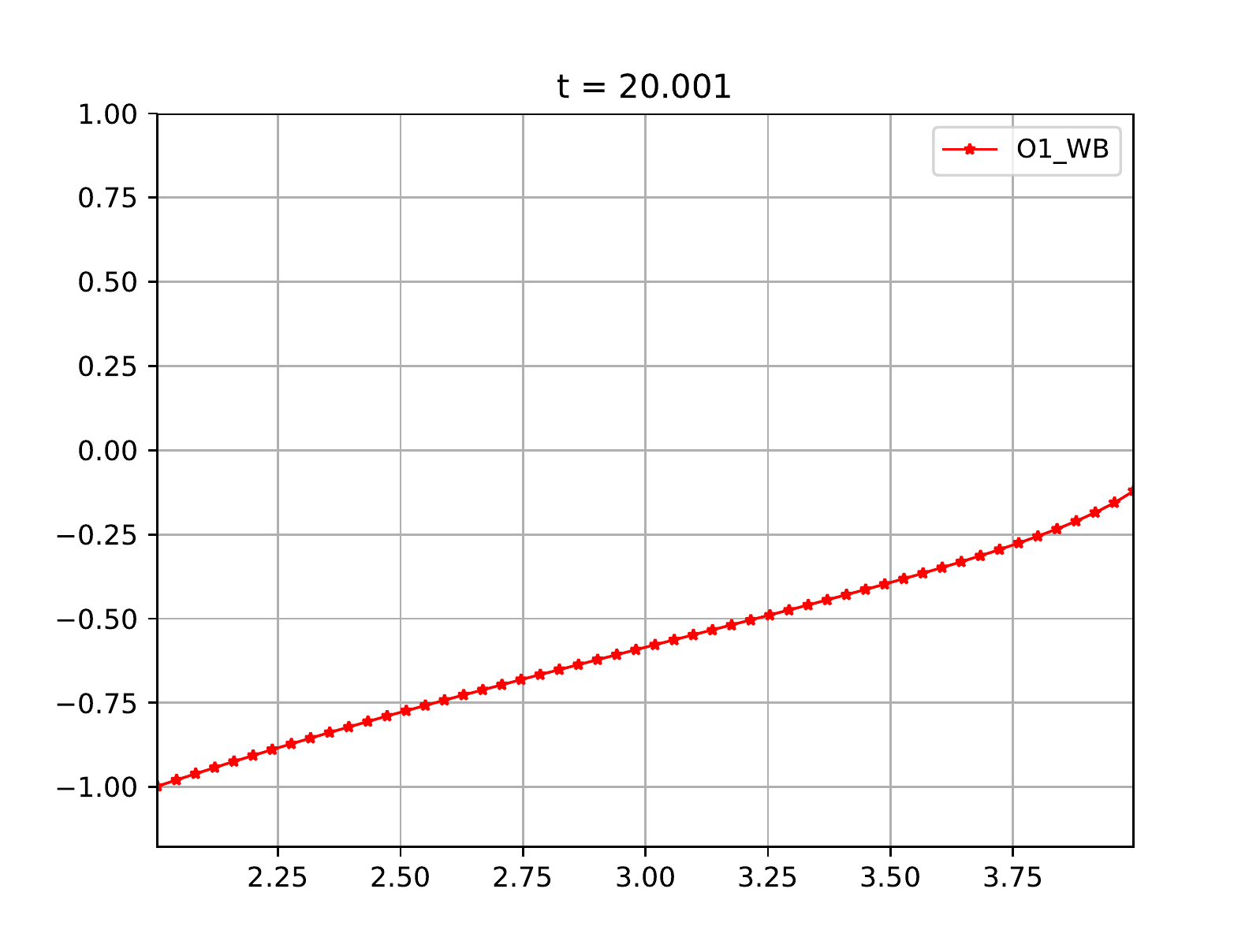}
		\label{fig:ko1_wb_testNonOneNegativeNonSteady_t_20}
	\end{subfigure}
	\caption{Burgers-Schwarzschild model with the initial condition \eqref{testB12}: numerical solution obtained with the first-order well-balanced scheme at selected times.}
	\label{fig:ko1_wb_testNonOneNegativeNonSteady}
\end{figure}
\end{enumerate}

%


\subsection{Main conclusions for the Burgers-Schwarzschild model}\label{subsec_conclusions_Burgers}

From Figures \ref{fig:ko1_ko2_ko3_wb_testSteadyShockLeftPerturbated} to \ref{fig:integralperturbationvsintegralsteadyshockperturbedlabelled} and Table \ref{tab:Areas_for_different_values_of_alpha} we can conclude the following:
\begin{conclusion}
	If  a perturbation $\delta_{v}$ is added to a  steady shock solution of the form
	$$
	v_{0}(r) = \begin{cases}
	\sqrt{1-K_{0}^{2}(1-\frac{2M}{r})}, & \text{  $2M<r<r_{0}$,}\\
	-\sqrt{1-K_{0}^{2}(1-\frac{2M}{r})}, & \text{ otherwise,}
	\end{cases}
	$$
 then the solution reaches at finite time another steady shock solution of the form:  
	$$
	v(r) = \begin{cases}
	\sqrt{1-K_{0}^{2}(1-\frac{2M}{r})}, & \text{ $2M<r<r_{1}$,}\\
	-\sqrt{1-K_{0}^{2}(1-\frac{2M}{r})}, & \text{ othewise.}
	\end{cases}
	$$
	\begin{enumerate}
		\item If $\int_{2M}^{r_{0}} \delta_{v} = 0$ and $\int_{r_{0}}^{\infty} \delta_{v} = 0$, then $r_{1}=r_{0}$, i.e.~the initial stationary solution is recovered.
		\item If $\int_{2M}^{\infty} \delta_{v} = 0$ and $\int_{2M}^{r_{0}} \delta_{v} = -\int_{r_{0}}^{\infty} \delta_{v}$, then $r_{1}\neq r_{0}$ and a different stationary solution is obtained.
		\item If $\int \delta_{v} \neq 0$, then $r_{1}\neq r_{0}$ and a different stationary solution is obtained. In this case the distance between $r_0$ and $r_1$ depends linearly on the amplitude of the perturbation: see  Table \ref{tab:Areas_for_different_values_of_alpha} and Figure \ref{fig:integralperturbationvsintegralsteadyshockperturbedlabelled}.
	\end{enumerate}


\end{conclusion}

In view of Figures \ref{fig:ko1_wb_testOnePositiveNonSteady} to \ref{fig:ko1_wb_testNonOneNegativeNonSteady} we have reached the following.

\begin{conclusion}
	\begin{enumerate}
		\item For a bounded domain $[2M,L]$:
		\begin{enumerate}
			\item If $v_{0}(r)= 1$ for $r\in [2M, 2M+\epsilon)$, with $\epsilon>0$, $v_{0}(L)\geq 0$ and $v_{0}\neq 1$, in finite time  the solution has the form of a right-moving shock that links the stationary solution $v\equiv 1$ and the negative steady solution $v^*$ such that $ v^*(2M) = -1$ and $v^{*}(L) = 0$, that is, 
$
v^{*}_{0}(r)= -\sqrt{1-\frac{1}{1-\frac{2M}{L}}\Big( 1-\frac{2M}{r} \Big)}.
$
			\item If $v_{0}(r)= 1$ for $r\in [2M, 2M+\epsilon)$, with $\epsilon>0$ and $v_{0}(L)=a$, with $a<0$, then in finite time the solution has the form of a right-moving shock that links the stationary solution $v\equiv 1$ and the negative steady solution $v*$  such that $v^*(2M) = -1$ and $v^{*}_{0}(L) = a$, that is,
			$
			 v^{*}_{0}(r)= -\sqrt{1-\frac{1-a^{2}}{1-\frac{2M}{L}}\Big( 1-\frac{2M}{r} \Big)}.
			$
			
			\item If $v_{0}(2M)< 1$ and $v_{0}(L)\geq 0$, then in finite time  the solution coincides  with the negative steady solution such that 
$v^*(2M) = -1$ and $v^{*}_{0}(L) = 0$, that is,
			$
			 v^{*}_{0}(r)= -\sqrt{1-\frac{1}{1-\frac{2M}{L}}\Big( 1-\frac{2M}{r} \Big)}.
			 $ 
			\item If $v_{0}(2M)< 1$ and $v_{0}(L)=a$, with $a<0$, then in finite time the solution coincides with  the negative stationary solution $v^*$  such that $v^*(L) = a$, that is,
			$
			 v^{*}_{0}(r)= -\sqrt{1-\frac{1-a^{2}}{1-\frac{2M}{L}}\Big( 1-\frac{2M}{r} \Big)}.
			 $ 
		\end{enumerate}
		\item For the unbounded domain $[2M,\infty)$ the following conclusions can be drawn  by passing to the limit when $L \to \infty$:
		\begin{enumerate}
			\item If $v_{0}(r)= 1$ for $r\in [2M, 2M+\epsilon)$, with $\epsilon>0$, $\lim_{r\to\infty}v_{0}(r)\geq 0$ and $v_{0}\neq 1$, in finite time  the solution has the form of a right-moving shock that links the stationary solution $v \equiv 1$ and the negative stationary solution
			$
			 v^{*}_{0}(r)= -\sqrt{\frac{2M}{r}},
			 $
corresponding to $K^2 = 1$.
			\item If $v_{0}(r)= 1$ for $r\in [2M, 2M+\epsilon)$, with $\epsilon>0$ and $\lim_{r\to\infty} v_{0}(r)=a$, with $a<0$, then in finite time $t_{0}$ the solution has the form of a  right- moving shock  that links the stationary solution $v \equiv 1$ and the negative stationary solution $v^*$  such that $v^*(2M) = -1$ and  $\lim_{r\to\infty} v^{*}_{0}(r) = a$, that is,
			$
			 v^{*}_{0}(r)= -\sqrt{1-(1-a^{2})\Big( 1-\frac{2M}{r} \Big)}.
			 $
			\item If $v_{0}(2M)< 1$ and $\lim_{r\to\infty} v_{0}(r)\geq 0$, then the solution converges as $t \to \infty$ to the negative stationary solution $v^*$  such that 
$v^*(2M) = -1$ and $\lim_{r\to\infty} v^{*}_{0}(r) = 0$, that is, 
			$ v^{*}_{0}(r)= -\sqrt{\frac{2M}{r}}.$
			\item If $v_{0}(2M)< 1$ and $\lim_{r\to\infty} v_{0}(r)=a$, with $a<0$, then  the solution converges as $t \to \infty$ to the negative stationary solution $v^*$
such that $v^*(2M) = -1$ and  $\lim_{r\to\infty} v^{*}_{0}(r) = a$, that is, 
			$v^{*}_{0}(r)= -\sqrt{1-(1-a^{2})\Big( 1-\frac{2M}{r} \Big)}.$
		\end{enumerate}
	\end{enumerate} 
\end{conclusion}


\section{Euler-Schwarzschild model: designing the numerical algorithm}
\label{section--5} 

\subsection{Preliminaries}

\begin{figure}
		\centering
		\includegraphics[width=.7\linewidth]{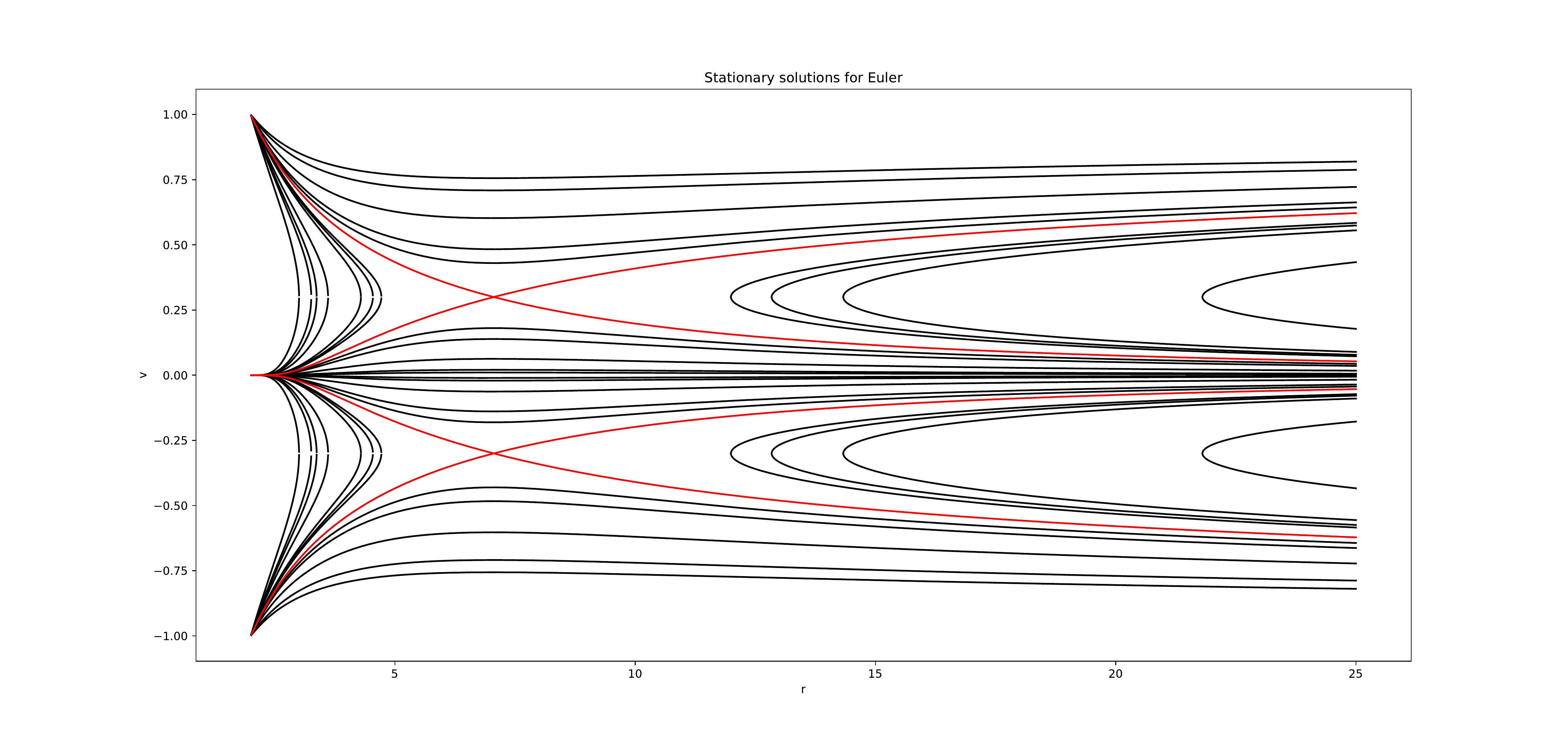}
		\caption{Euler-Schwarzschild model with $k = 0.3$: $v$ variable for some stationary solutions.}
		\label{fig:steadysolutionsforeuler}
\end{figure}

In the case of the Euler-Schwarzschild  equations \eqref{eq:Euler}, the stationary solutions are implicitly given by the equations:
	\bel{eq:steady_state_Euler} 
	\begin{array}{l}
		\displaystyle \frac{\sign(v)(1-v^{2})|v|^{\frac{2k^{2}}{1-k^{2}}}r^{\frac{4k^{2}}{1-k^{2}}}}{\Big( 1-\frac{2M}{r} \Big)} = C_{1},
		\qquad 
		\displaystyle r(r-2M)\rho\frac{v}{1-v^{2}}=C_{2},
	\end{array}
	\ee
	where 
	$C_1, C_2$ are constants. 
	The pair  $(v, \rho)$ of a stationary solution satisfies the following 
	ODE system analyzed first in \cite{PLF-SX1,PLF-SX2}: 
\bel{eq:v_derivative_steady_state_Euler}
\frac{dv}{dr} = v\frac{(1-v^{2})(1-k^{2})}{r(r-2M)}\left(\frac{2k^{2}}{1-k^{2}}(r-2M)-M\right)\Big/ (v^{2}-k^{2}),
\ee
\bel{eq:rho_derivative_steady_state_Euler}
\frac{d\rho}{dr} = -\frac{2(r-M)}{r(r-2M)}\rho - \rho\frac{(1+v^{2})(1-k^{2})}{r(r-2M)}\left(\frac{2k^{2}}{1-k^{2}}(r-2M)-M\right)\Big/ (v^{2}-k^{2}). 
\ee
 Figure \ref{fig:steadysolutionsforeuler} shows the graph of $v$ for some of them. When these functions are defined in $(2M,\infty)$,
 they have a  maximum or a minimum  in
 \bel{eq:r_critical_steady_state_Euler}
 r_{c} = \frac{M(1-k^{2})}{2k^{2}} +2M,
 \ee
 that comes from solving $\frac{dv}{dr} =0$. In Figure \ref{fig:steadysolutionsforeuler} the red stationary solutions are those such that at $r=r_{c}$ the take the value $v=\pm k$.

Given two constants $C_1$ and $C_2$, in order  to compute the stationary solution given by \eqref{eq:steady_state_Euler} in a point $r=a$,  the following nonlinear system has to be solved:
\bel{eq:system} 
	\sign(v)(1-v^{2})|v|^{\frac{2k^{2}}{1-k^{2}}} 
	= \displaystyle \frac{\left(1-\frac{2M}{a}\right)}{a^{\frac{4k^{2}}{1-k^{2}}}}C_1,\\ \\
	\displaystyle \rho=\frac{1-v^{2}}{va(a-2M)}C_2.
\ee
It is enough thus to solve, if it is possible, the nonlinear equation
\bel{eq:eq}
g(v) = K_a,
\ee
with
\be g(v) =\sign(v)(1-v^{2})|v|^{\frac{2k^{2}}{1-k^{2}}}, \quad v \in [-1,1],
\qquad
K_a = \frac{\left(1-\frac{2M}{a}\right)}{a^{\frac{4k^{2}}{1-k^{2}}}}C_1,
\ee
to compute $v$. Once this equation is solved, $\rho$ is computed using the second equation of \eqref{eq:system}. 

It can be easily checked that $g$ satisfies:
$$
 -(1-k^{2})k^{\frac{2k^{2}}{1-k^{2}}} = g(-k) \leq g(v) \leq g(k) = (1-k^{2})k^{\frac{2k^{2}}{1-k^{2}}}, 
 \qquad v \in [-1,1].
$$
Moreover, $g$ is strictly monotone in $[-1, -k)$, $(-k, k)$, and $(k,1]$. As a consequence we have the following conclusion. ((For further properties, see the study in \cite{PLF-SX1,PLF-SX2}.) 
\begin{itemize}

\item If  $|K_a| > g(k)$, the equation \eqref{eq:eq} has no solution, i.e.~a stationary solution given by $C_1$ and $C_2$ cannot be defined at $r = a$.

\item If $|K_a| = g(k)$ then the equation \eqref{eq:eq} has only one solution ($k$ if $K_a > 0$, $-k$ if $K_a < 0$). Therefore, \eqref{eq:system} has only one solution
that is a sonic state.

\item Otherwise, \eqref{eq:eq} has two possible solutions. Therefore there are two states $(\rho, v)$ that solve  \eqref{eq:system}, one supersonic and one subsonic.

\end{itemize}

\noindent For the sake of clarity, together with the representation 
$$V = [\rho (1+k^{2}v^{2})/(1-v^{2}), \rho v (1+k^{2})/(1-v^{2})]^T, $$
 for the states, we will use  
$\widetilde V= [\rho, v]^T.$
Here, $V$ can be easily computed from $\widetilde V$ and, given $V$,  $\widetilde V$ is also easily computed by (\ref{eq:primitive_variables}) that comes from solving a second-degree equation.


\subsection{First-order method}

If the midpoint rule is used to compute the initial averages, given a family of cell values $\widetilde V_i$, in the first step of the well-balanced  reconstruction procedure
one has to find, if it is possible,  a stationary solution $\widetilde V^*_i$ defined in $[r_{i-\frac{1}{2}}$, $r_{i+\frac{1}{2}}]$ such that
$
\widetilde  V_i^*(r_i) = \widetilde  V_i = [\rho_i, v_i]^T.
$ 
Obviously such a stationary solution would correspond to the choice of constants:
  \begin{eqnarray}
    C_{i,1} & = & \displaystyle \frac{\sign(v_{i})(1-v_{i}^{2})|v_{i}|^{\frac{2k^{2}}{1-k^{2}}}r_{i}^{\frac{4k^{2}}{1-k^{2}}}}{\left(1-\frac{2M}{r_{i}}\right)},\label{Ci1212}
    \qquad
    C_{i,2}  =   \displaystyle r_{i}(r_{i}-2M)\rho_{i}\frac{v_{i}}{1-v_{i}^{2}}. 
  \end{eqnarray}
According to the discussion above, the corresponding stationary solution is defined in $r_{i \pm \frac{1}{2}}$ provided that:
 \bel{eq:Econdexfo}
|K_{i \pm \frac{1}{2}}| \leq g(k),
\qquad
K_{i \pm\frac{1}{2}} = \left(1-\frac{2M}{r_{i \pm \frac{1}{2}}}\right)r_{i \pm \frac{1}{2}}^{-\frac{4k^{2}}{1-k^{2}}}C_{i,1}.
\ee
When $|K_{i \pm \frac{1}{2}}| < g(k)$ there are two possible values for $\widetilde V^*_i(r_{i \pm \frac{1}{2}})$, one subsonic and one supersonic. Therefore, a criterion is needed to select one or the other. The following criterion will be used here:
\begin{itemize}
\item if $\widetilde V_i$ is not sonic, then the state whose regime (sub or supersonic) is the same  as $\widetilde V_i$ is selected for  $\widetilde  V^*_i(r_{i \pm \frac{1}{2}}).$
\item if $\widetilde  V_i$ is sonic, then the state whose regime is the same as $\widetilde  V_{i+1}$  is selected for  $\widetilde  V^*_i(r_{i + \frac{1}{2}})$ and the state whose regime is the same as $\widetilde  V_{i-1}$ is selected for  $\widetilde V^*_i(r_{i - \frac{1}{2}})$.
\end{itemize}
Observe that this criterion aims to preserve the regime of the given cell values. 

If  condition \eqref{eq:Econdexfo}  is satisfied, then the numerical method \eqref{eq:nummetfo} is used. Otherwise the standard  trivial reconstruction is considered.

The expression of the semi-discrete first-order method is then as follows:
\bel{eq:semi_disc-merge_E}
\frac{dV_{i}}{dt}= -\frac{1}{\Delta r}\left(F_{i+\frac{1}{2}}-F_{i-\frac{1}{2}} - S_i \right),
\ee
where 
\bel{eq:source-merge_E}
S_ i 
= \begin{cases}
F(V_{i}^{*}(r_{i+\frac{1}{2}}), r_{i+\frac{1}{2}}) - F(V_{i}^{*}(r_{i-\frac{1}{2}}), r_{i-\frac{1}{2}})), & \text{ if \eqref{eq:Econdexfo} holds,} 
\\
S(V_i, r_i), & \text{ otherwise.} 
\end{cases}
\ee
The forward Euler method is used again for the time discretization.


\subsection{Second-order method}

Let us use again the midpoint rule to compute cell averages and  the minmod reconstruction operator. 
The stationary solution sought at the first stage of the well-balanced reconstruction  procedure is again characterized by the constants \eqref{Ci1212}.
This time, this stationary solution has to be computed at the points  $r_{i-1}$, $r_{i -\frac{1}{2}}$, $r_{i+\frac{1}{2}}$, $r_{i+1}$ so that the following condition has to be satisfied:
\bel{eq:Econdexso}
|K_{i+j}| \leq g(k),\quad j = -1, -\frac{1}{2}, \frac{1}{2}, 1,
\ee
where $K_{i\pm \frac{1}{2}}$ are given by \eqref{eq:Econdexfo} and
\bel{Ki+1}
K_{i \pm 1} = \left(1-\frac{2M}{r_{i \pm 1}}\right)r_{i \pm 1}^{-\frac{4k^{2}}{1-k^{2}}}C_{i,1}.
\ee

If this condition is satisfied,  the reconstruction is defined as follows:
$$
\mathbb{P}_{i}(r) = V_{i}^{*}(r) + \text{minmod}\left(\frac{W_{i+1}-W_{i}}{\Delta r}, \frac{W_{i+1}-W_{i-1}}{2\Delta r},\frac{W_{i}-W_{i-1}}{\Delta r}\right)(r-r_{i}),
$$
where the minmod function is applied component by component and
$
W_j = V_j - V_i^*(r_j)$ for $j = i-1, i, i+1$.
Observe that the \textit{conserved} variables $V$ are used in the reconstruction procedure.
On the other hand, if \eqref{eq:Econdexso} is not satisfied, then the standard MUSCL reconstruction is applied:
$$
\mathbb{Q}_{i}(r) = V_i + \text{minmod}\left(\frac{V_{i+1}-V_{i}}{\Delta r}, \frac{V_{i+1}-V_{i-1}}{2\Delta r},\frac{V_{i}-V_{i-1}}{\Delta r}\right)(r-r_{i}). 
$$
The expression of the numerical method is given again by \eqref{eq:semi_disc-merge_E}-\eqref{eq:source-merge_E} with the difference that the second-order reconstructions are used now to compute
the numerical fluxes. The TVDRK2 method is used now to discretize the equations in time. 


\subsection{Third-order method}

Although it will not be implemented in the present paper, let us briefly describe the first step of a third-order well-balanced reconstruction procedure based on the  two-point Gauss quadrature i order to compute averages:.
It consists on 
finding $C_1$ and $C_2$ such that
\be
\frac{1}{2}V^*(r_{i,0}; C_1, C_2) + \frac{1}{2}V^*(r_{i,1}; C_1, C_2) = V_i,
\ee
where $V^*(r; C_1, C_2)$ represents the value at $r$ of a stationary solution characterized by the constants $C_1$ and $C_2$. 


\section{Euler-Schwarzschild model: a numerical study}
\label{section--6} 

\subsection{Preliminaries}

In this section several tests are considered to check the performance of the first- and second-order well-balanced numerical methods  for Euler-Schwarzschild model introduced in the previous section. We consider the spatial  interval $[2M, L]$ with $M=1$ and $L=10$, a 500-point uniform mesh, $k=0.3$ and the CFL number is set to 0.5 again. At $r=2M$ we impose $F_{-\frac{1}{2}}=0$ as boundary condition since $\Big( 1-\frac{2M}{r} \Big)= 0$. At $r=L$ we will use a transmissive boundary condition in the case we are not in a stationary solution or we will expand the stationary solution if we are in one. 

In order to test the dependency of the results on the numerical method, two different first-order numerical fluxes are considered:
the Lax-Friedrichs numerical flux
\begin{equation}\label{LFflux}
F_{i+\frac{1}{2}} = \frac{1}{2}(F(V_{i})+F(V_{i+1}))-\frac{1}{2}\frac{\Delta t}{\Delta x}(V_{i+1}-V_{i}),
\end{equation}
and  a HLL-like numerical flux in PVM form (see \cite{Castro-FN}):
\begin{equation}\label{HLLflux}
F_{i+\frac{1}{2}} = \frac{1}{2}(F(V_{i})+F(V_{i+1}))-\frac{1}{2}\left(\alpha_{0}(V_{i+1}-V_{i})+\alpha_{1}(F(V_{i+1})-F(V_{i}))\right),
\end{equation}
with
\begin{equation}\label{coef}
\alpha_{0}=\frac{\overline{\lambda_{2}}|\overline{\lambda_{1}}| - \overline{\lambda_{1}}|\overline{\lambda_{2}}|}{\overline{\lambda_{2}}- \overline{\lambda_{1}}}, \qquad
 \alpha_{1}=\frac{|\overline{\lambda_{2}}| - |\overline{\lambda_{1}}|}{\overline{\lambda_{2}}- \overline{\lambda_{1}}},
\end{equation}
where $\overline{\lambda_{1}}$ and $\overline{\lambda_{2}}$ are the eigenvalues of some intermediate  matrix $J_{i+\frac{1}{2}}$ of the form
\begin{equation}\label{Jac}
J_{i+\frac{1}{2}} =  \left(1- \frac{2M}{r_{i+\frac{1}{2}}}\right)  \left[ \begin{array}{cc} 0 &  1\\ \displaystyle \frac{k^2 - v_m^2}{1 - k^2 v_m^2} & \displaystyle \frac{2(1 - k^2)v_m}{1 - k^2 v_m^2} \end{array}\right]
\end{equation}
where $v_m$ is some intermediate value between $v_i^n$ and $v_{i+1}^n$. 

Given two states $V_L$ and $V_R$, in order to choose an adequate intermediate value $v_m$, we look for $v$ such that the following Roe-type property is satisfied:
\bel{roe}
\aligned
&
 \left[ \begin{array}{cc} 0 &  1\\ \displaystyle \frac{k^2 - v^2}{1 - k^2 v^2} & \displaystyle \frac{2(1 - k^2)v}{1 - k^2 v^2} \end{array}\right]\cdot  (V_R - V_L ) = \widehat F_R - \widehat F_L,
\qquad
\quad
 \widehat F_\alpha  =  \left(\begin{array}{c}
\displaystyle \frac{1+k^{2}}{1-v_\alpha^{2}}\rho_\alpha v_\alpha \\
\displaystyle \frac{v_\alpha^{2}+k^{2}}{1-v_\alpha^{2}}\rho_\alpha
\end{array}\right),\quad \alpha = L, R,
\endaligned
\end{equation}
i.e.~the factor $(1- 2M/r)$ is neglected for simplicity. Due to the form of the matrix, it is enough to find $v$ such that
$$
\frac{k^2 - v^2}{1 - k^2 v^2}( V_{1, R} - V_{1,L}) +  \frac{2(1 - k^2)v}{1 - k^2 v^2}  ( V_{2,R} - V_{2,L}) = F_{2,R} - F_{2,L}.
$$
This equality is equivalent to a second-order degree equation for $v$, namely 
$
\alpha v^2  + \beta v + \gamma = 0,
$
where
\begin{eqnarray*}
	& & \alpha = \rho_R(1- v_L^2) - \rho_L(1- v_R^2), \\
	& & \beta =- 2 \left(\rho_Rv_R(1- v_L^2) - \rho_Lv_L(1- v_R^2)\right),
	\\
	& & \gamma = \rho_Rv_R^2(1- v_L^2) - \rho_Lv_L^2(1- v_R^2).
\end{eqnarray*}
Since the discriminant 
$
D = \rho_L \rho_R (1- v_L^2)(1- v_R^2) (v_R - v_L)^2.
$
 is always positive, there are always two real solutions:
\begin{eqnarray*}
	v_{\pm} = \frac{\rho_R v_R(1-v_L^2) - \rho_L v_L (1-v_R^2) \pm |v_R-v_L| \sqrt{\rho_L \rho_R (1- v_L^2)(1- v_R^2)}}{\rho_R(1-v_L^2)-\rho_L(1-v_R^2)}
\end{eqnarray*}
and it can be proven that:
\begin{itemize}

\item  if $v_L<v_R$, then $v_{-} \in (v_L,v_R)$ and $v_{+} \notin (v_L,v_R)$, so that we will take $v_m=v_{-}$;

\item   if $v_L>v_R$ then $v_{-} \notin (v_R,v_L)$ and $v_{+} \notin (v_R,v_L)$,  so that we will take $v_m=v_{+}$.

\end{itemize}

Finally, in the case $\alpha=0$ and $V_R\neq V_L$, we take $v_m = -\frac{\gamma}{\beta}$ and in the case $||V_R- V_L||_\infty < \epsilon$ we take $v_m=\frac{v_L+v_R}{2}$. 

Once $v_m$ has been chosen, the expression of $\overline{\lambda_j}$, $j=1,2$ in \eqref{coef} is as follows:
\begin{eqnarray*}
& \overline{\lambda_{1}} = \lambda_{1}\left(v_m\right) = \left(1-\frac{2M}{r_{i+\frac{1}{2}}}\right)\frac{v_m-k}{1-k^{2}v_m}, \qquad
\quad
 \overline{\lambda_{2}} = \lambda_{2}\left(v_m\right) = \left(1-\frac{2M}{r_{i+\frac{1}{2}}}\right)\frac{v_m+k}{1+k^{2}v_m}.
\end{eqnarray*}
Since for a $2 \time 2$-systems HLL and Roe methods are equivalent and the intermediate value chosen  to compute the wave speeds satisfies a Roe-type property, this numerical flux will be called Roe-type numerical flux in what follows. 
The proposed numerical  method will be compared with those based on the same numerical flux and the standard first- and second-order reconstructions.


\subsection{Stationary solutions} 

\paragraph{{{\red Positive stationary solution}}}

\begin{figure}[h]
	\begin{subfigure}[h]{0.5\textwidth}
		\centering
		\includegraphics[width=1\linewidth]{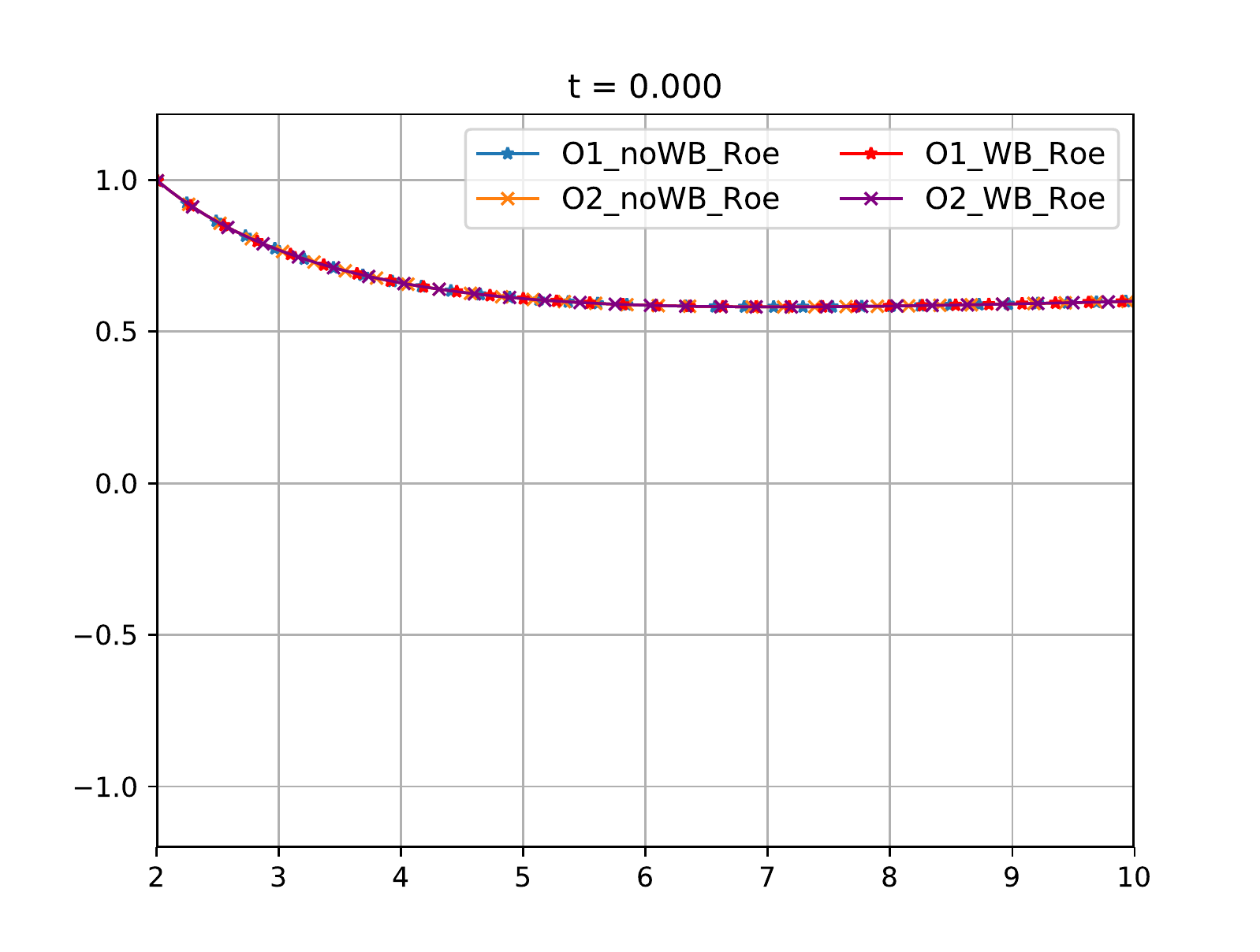}
		\label{fig:Euler_ko1_ko2_WB_vs_noWB_testWB1_t_0_hepse14}
	\end{subfigure}
	\begin{subfigure}[h]{0.5\textwidth}
		\centering
		\includegraphics[width=1\linewidth]{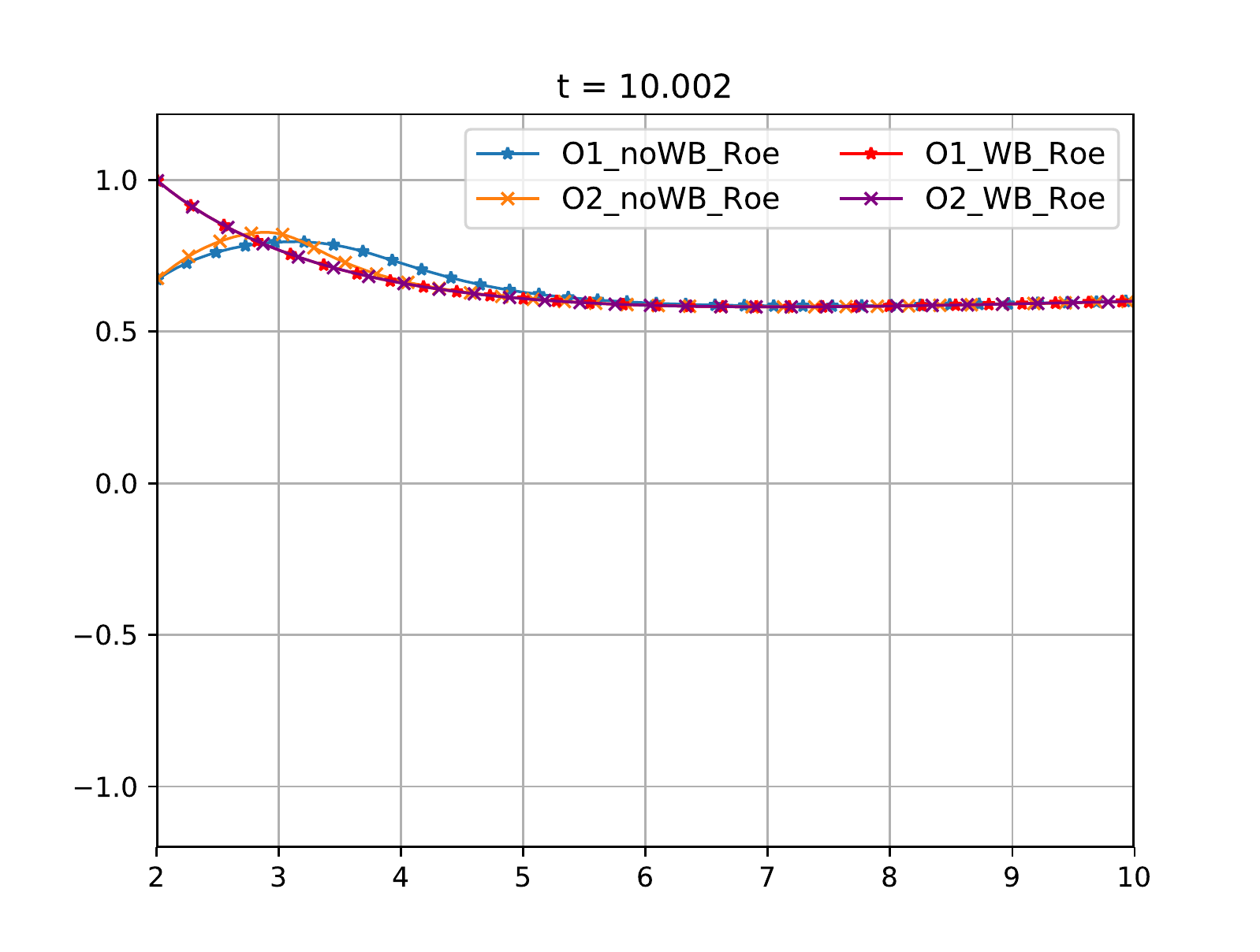}
		\label{fig:Euler_ko1_ko2_WB_vs_noWB_testWB1_t_10_hepse14}
	\end{subfigure}
	\begin{subfigure}[h]{0.5\textwidth}
		\centering
		\includegraphics[width=1\linewidth]{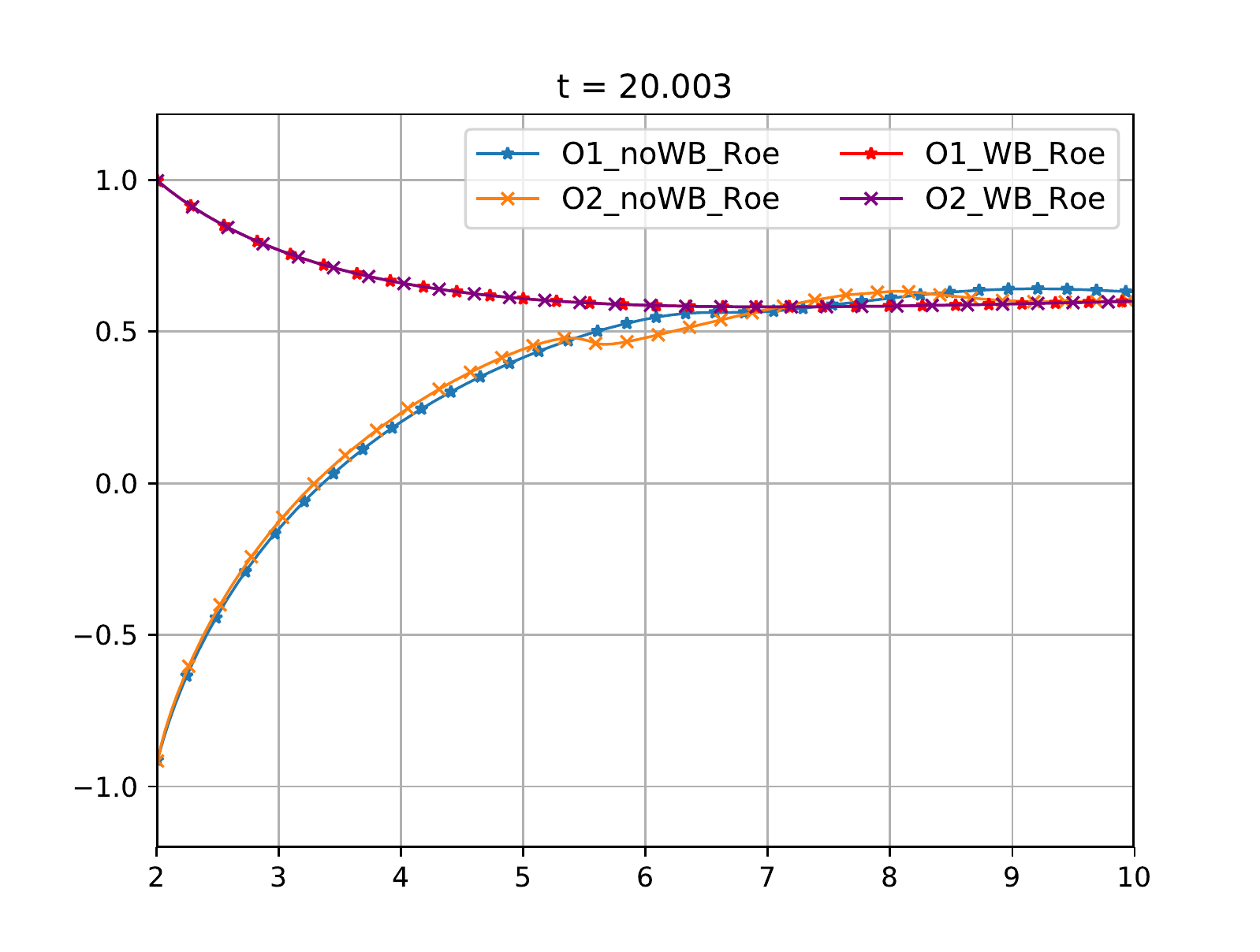}
		\label{fig:Euler_ko1_ko2_WB_vs_noWB_testWB1_t_20_hepse14}
	\end{subfigure}
	\begin{subfigure}[h]{0.5\textwidth}
		\centering
		\includegraphics[width=1\linewidth]{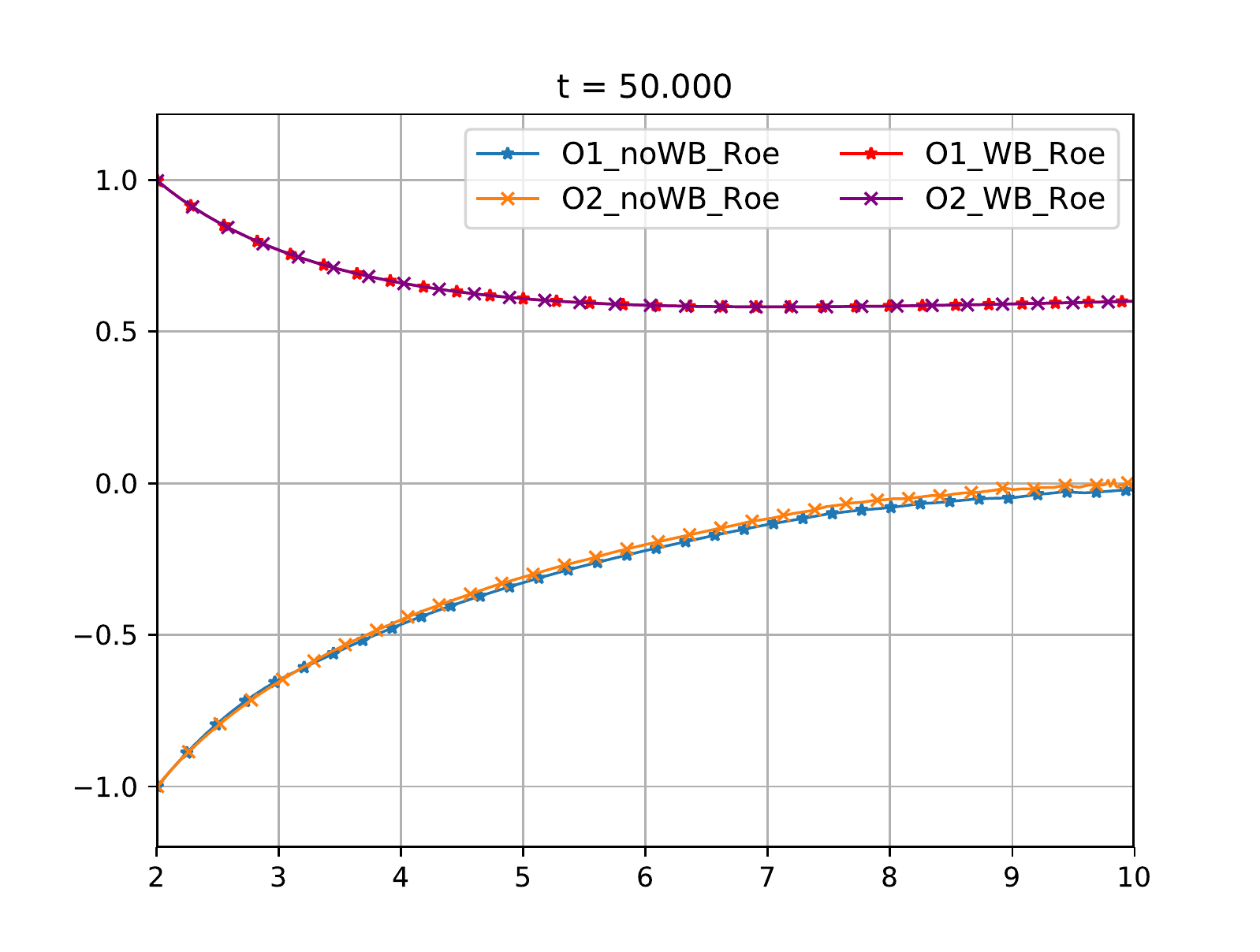}
		\label{fig:Euler_ko1_ko2_WB_vs_noWB_testWB1_t_50_hepse14}
	\end{subfigure}
	\caption{Euler-Schwarzschild model with the initial condition \eqref{testE1}:  first- and second-order well-balanced and non-well-balanced methods at selected times for the variable $v$.}
	\label{fig:Euler_ko1_ko2_WB_vs_noWB_testWB1_hepse14}
\end{figure}

\begin{figure}[h]
	\begin{subfigure}[h]{0.5\textwidth}
		\centering
		\includegraphics[width=1\linewidth]{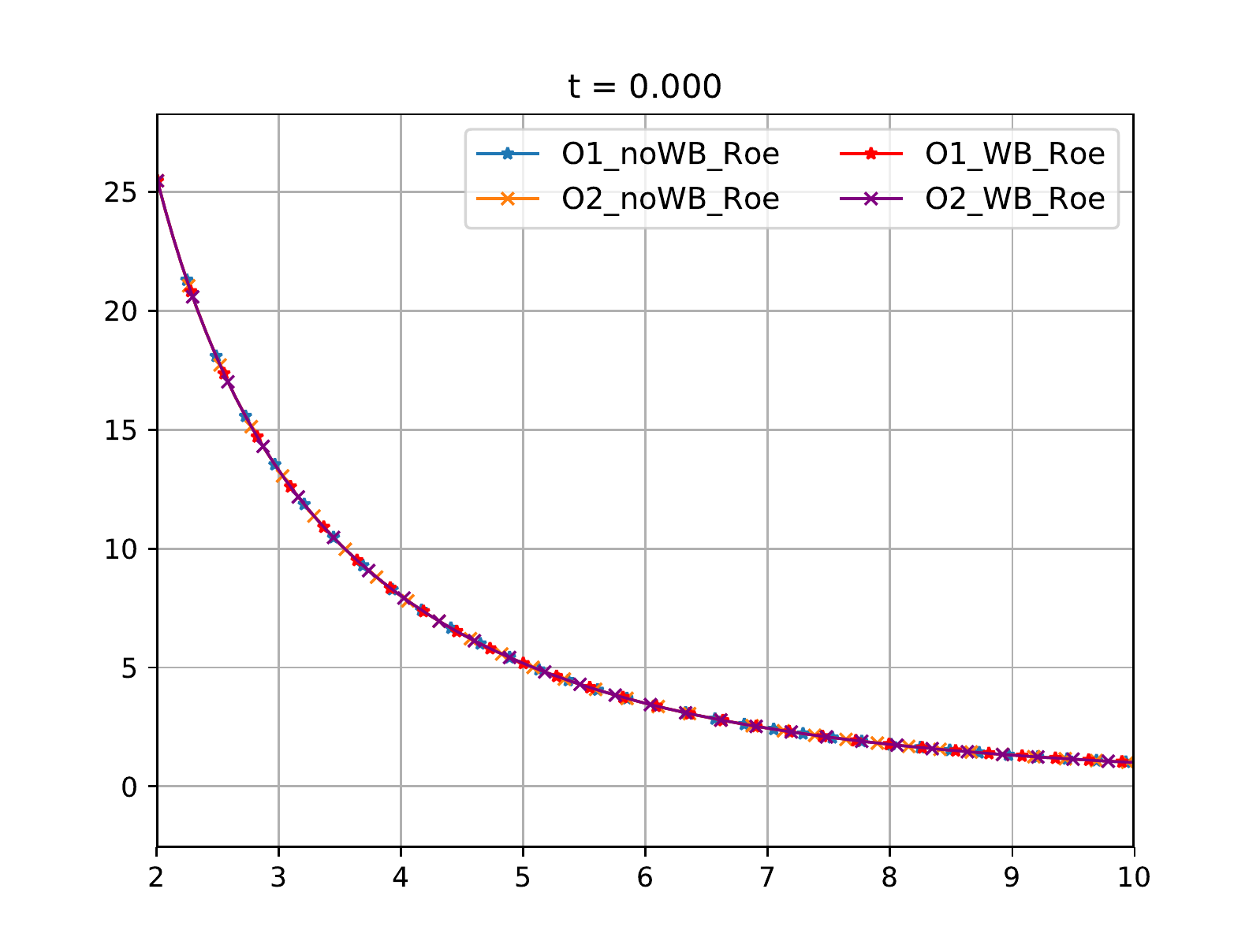}
		\label{fig:Euler_ko1_ko2_WB_vs_noWB_testWB1_t_0_hepse14_rho}
	\end{subfigure}
	\begin{subfigure}[h]{0.5\textwidth}
		\centering
		\includegraphics[width=1\linewidth]{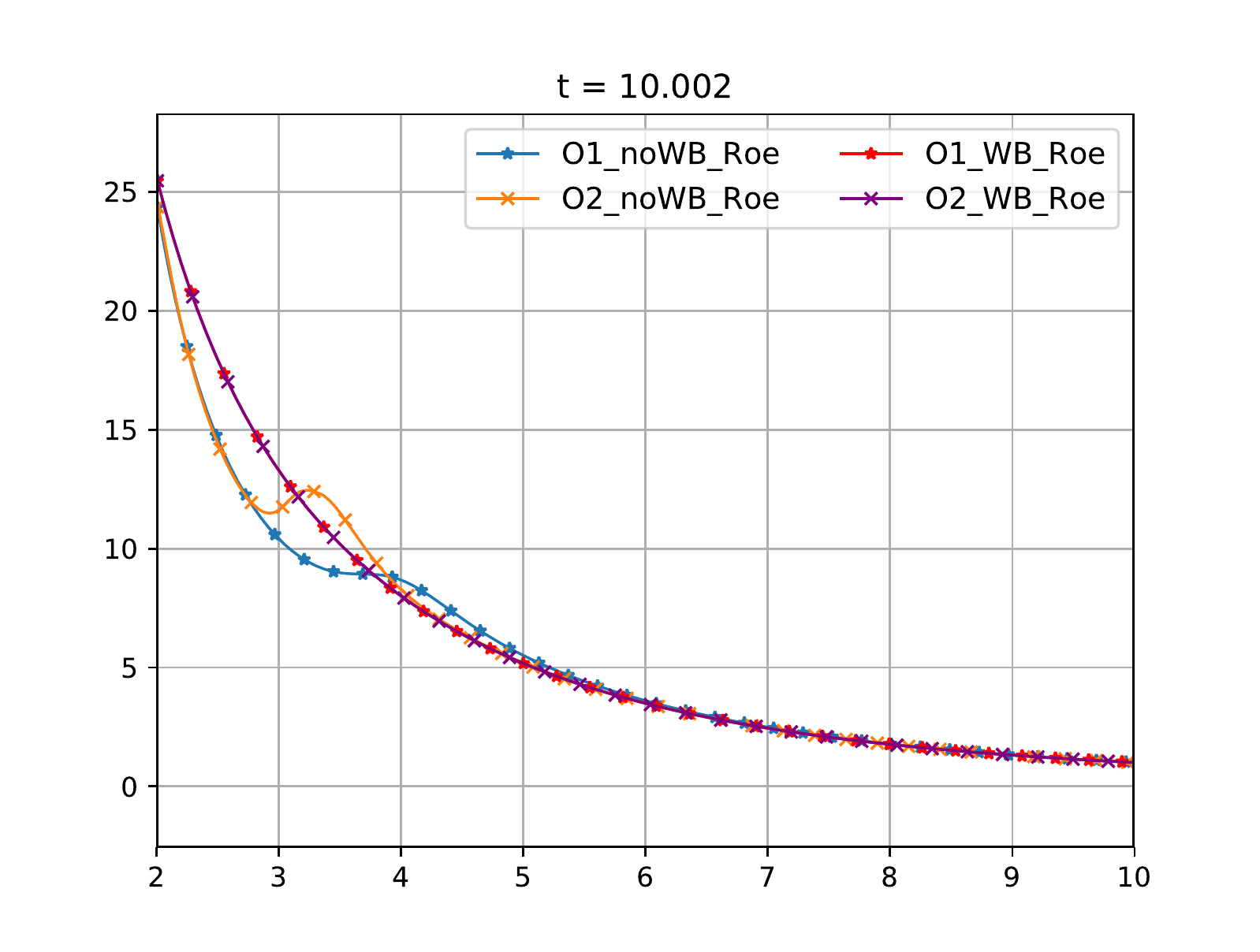}
		\label{fig:Euler_ko1_ko2_WB_vs_noWB_testWB1_t_10_hepse14_rho}
	\end{subfigure}
	\begin{subfigure}[h]{0.5\textwidth}
		\centering
		\includegraphics[width=1\linewidth]{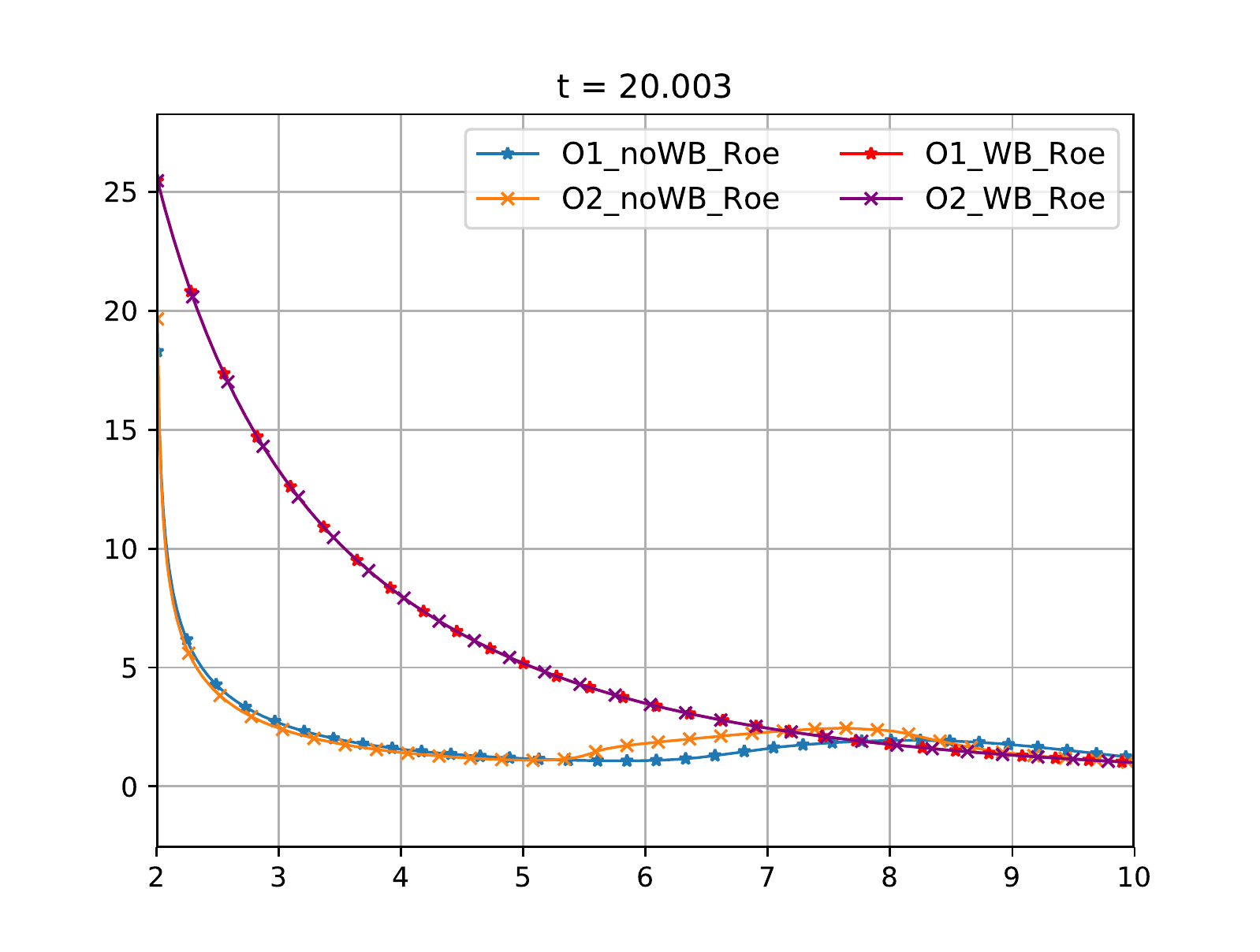}
		\label{fig:Euler_ko1_ko2_WB_vs_noWB_testWB1_t_20_hepse14_rho}
	\end{subfigure}
	\begin{subfigure}[h]{0.5\textwidth}
		\centering
		\includegraphics[width=1\linewidth]{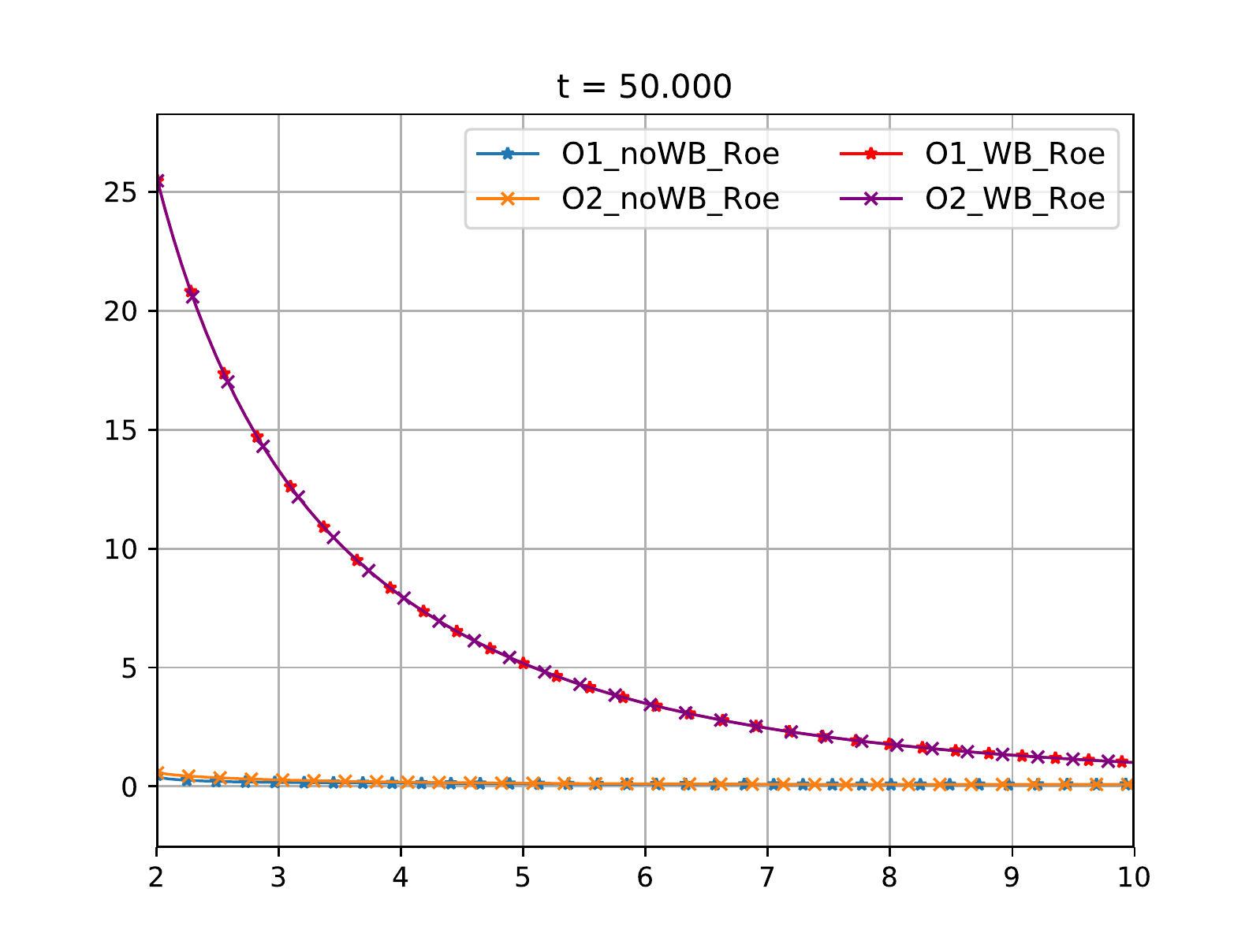}
		\label{fig:Euler_ko1_ko2_WB_vs_noWB_testWB1_t_50_hepse14_rho}
	\end{subfigure}
	\caption{Euler-Schwarzschild model with the initial condition \eqref{testE1}:  first- and second-order well-balanced and non-well-balanced methods at selected times for the variable $\rho$.}
	\label{fig:Euler_ko1_ko2_WB_vs_noWB_testWB1_hepse14_rho}
\end{figure}

We take as initial condition the positive supersonic stationary solution satisfying
\bel{testE1}
 \rho^{*}(10)=1, \quad v^{*}(10)=0.6.
\ee
 Table \ref{tab:Error_TestE1} shows the error in $L^1$ norm between the numerical solution
at time $ t = 50$ for the well-balanced and non-well-balanced methods using the Roe-type numerical flux.
  Figures \ref{fig:Euler_ko1_ko2_WB_vs_noWB_testWB1_hepse14} and \ref{fig:Euler_ko1_ko2_WB_vs_noWB_testWB1_hepse14_rho} compare the numerical solutions obtained with the well-balanced and the non-well-balanced methods: as it happened for the Burgers-Schwarzschild model, the numerical solutions obtained with non-well-balanced methods depart from the initial steady state.
  
  \begin{table}[ht]
  	\centering
  	\begin{tabular}{|c|c|c|c|c|}
  		\hline 
  		Scheme (500 cells) & Error $v$ (1st) & Error $\rho$ (1st) & Error $v$ (2nd) & Error $\rho$ (2nd) \\ 
  		\hline 
  		Well-balanced & 3.34E-13 & 5.61E-12 & 3.43E-13 & 7.12E-12  \\ 
  		\hline 
  		Non well-balance & 0.94 & 5.79 & 0.93 & 5.75 \\ 
  		
  		\hline 
  	\end{tabular} 
	  	\caption{Well-balanced versus non-well-balanced schemes: $L^{1}$ errors at time $t=50$ for the Euler-Schwarzschild model with the initial condition (\ref{testE1}).}

  	\label{tab:Error_TestE1}
  \end{table}

\paragraph{{{\red Negative stationary solution}}}

\begin{figure}[h]
	\begin{subfigure}[h]{0.32\textwidth}
		\centering
		\includegraphics[width=1\linewidth]{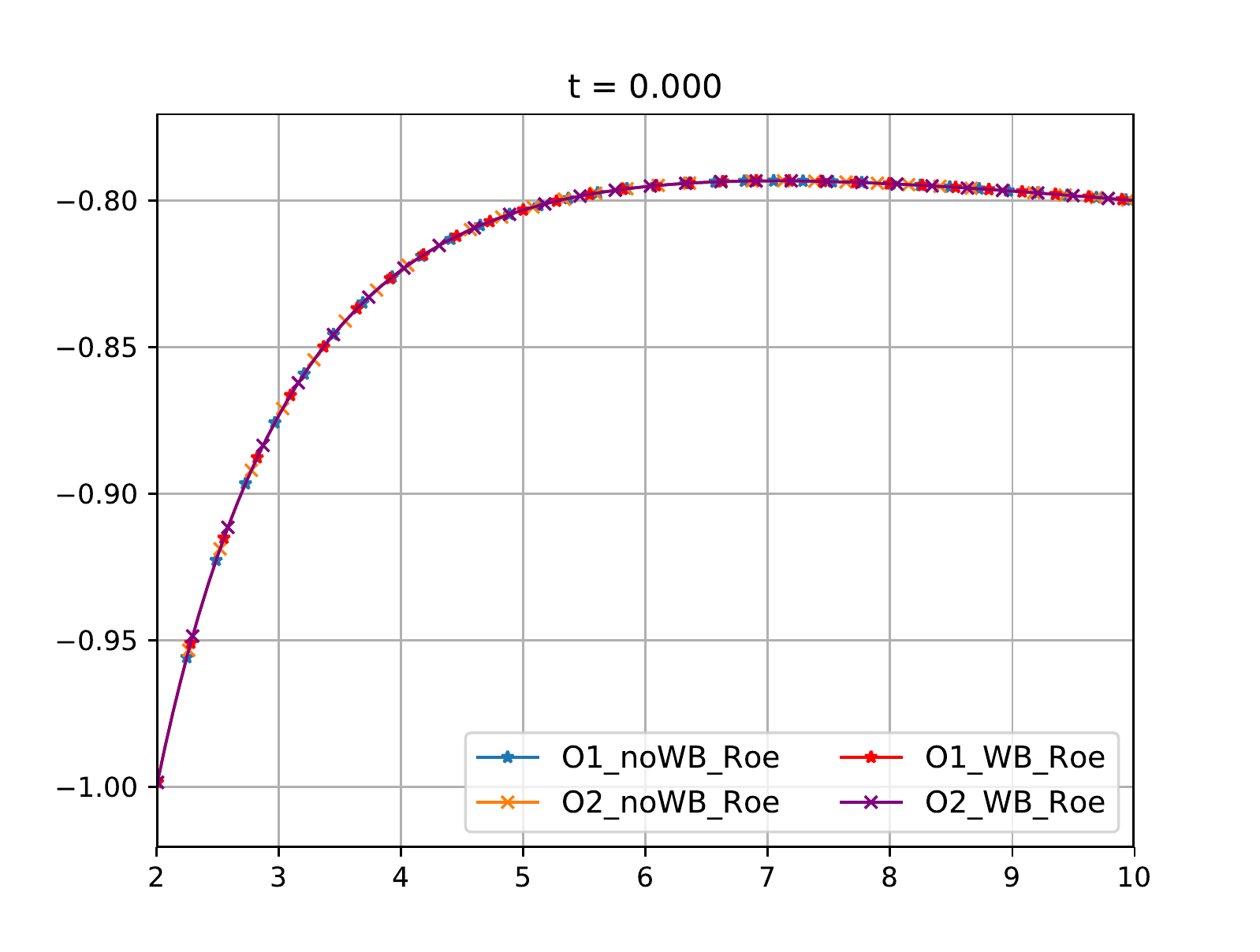}
		\label{fig:Euler_ko1_ko2_WB_vs_noWB_testWB2_t_0_hepse14}
	\end{subfigure}
	\begin{subfigure}[h]{0.32\textwidth}
		\centering
		\includegraphics[width=1\linewidth]{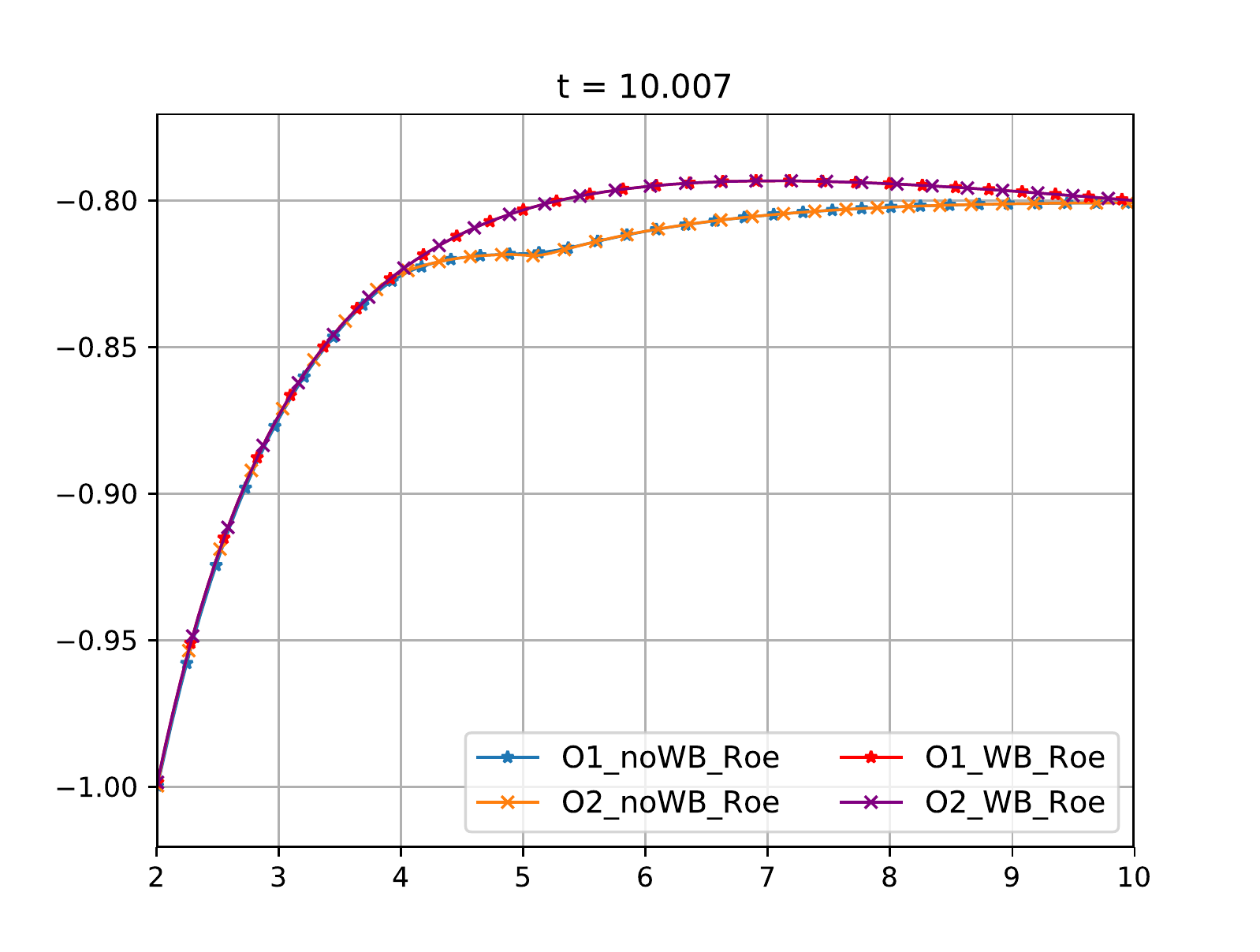}
		\label{fig:Euler_ko1_ko2_WB_vs_noWB_testWB2_t_10_hepse14}
	\end{subfigure}
	\begin{subfigure}[h]{0.32\textwidth}
		\centering
		\includegraphics[width=1\linewidth]{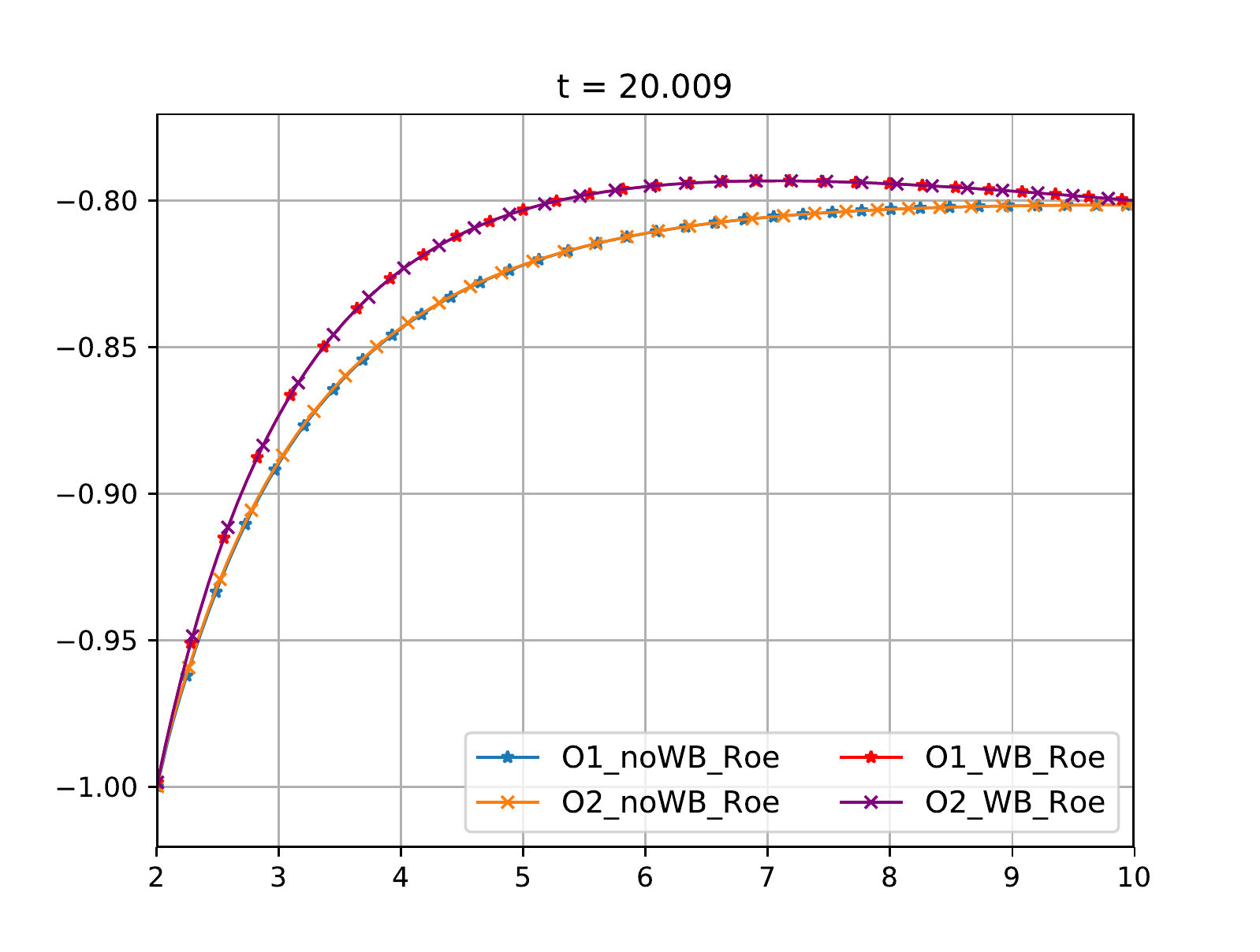}
		\label{fig:Euler_ko1_ko2_WB_vs_noWB_testWB2_t_20_hepse14}
	\end{subfigure}
	\begin{subfigure}[h]{0.32\textwidth}
		\centering
		\includegraphics[width=1\linewidth]{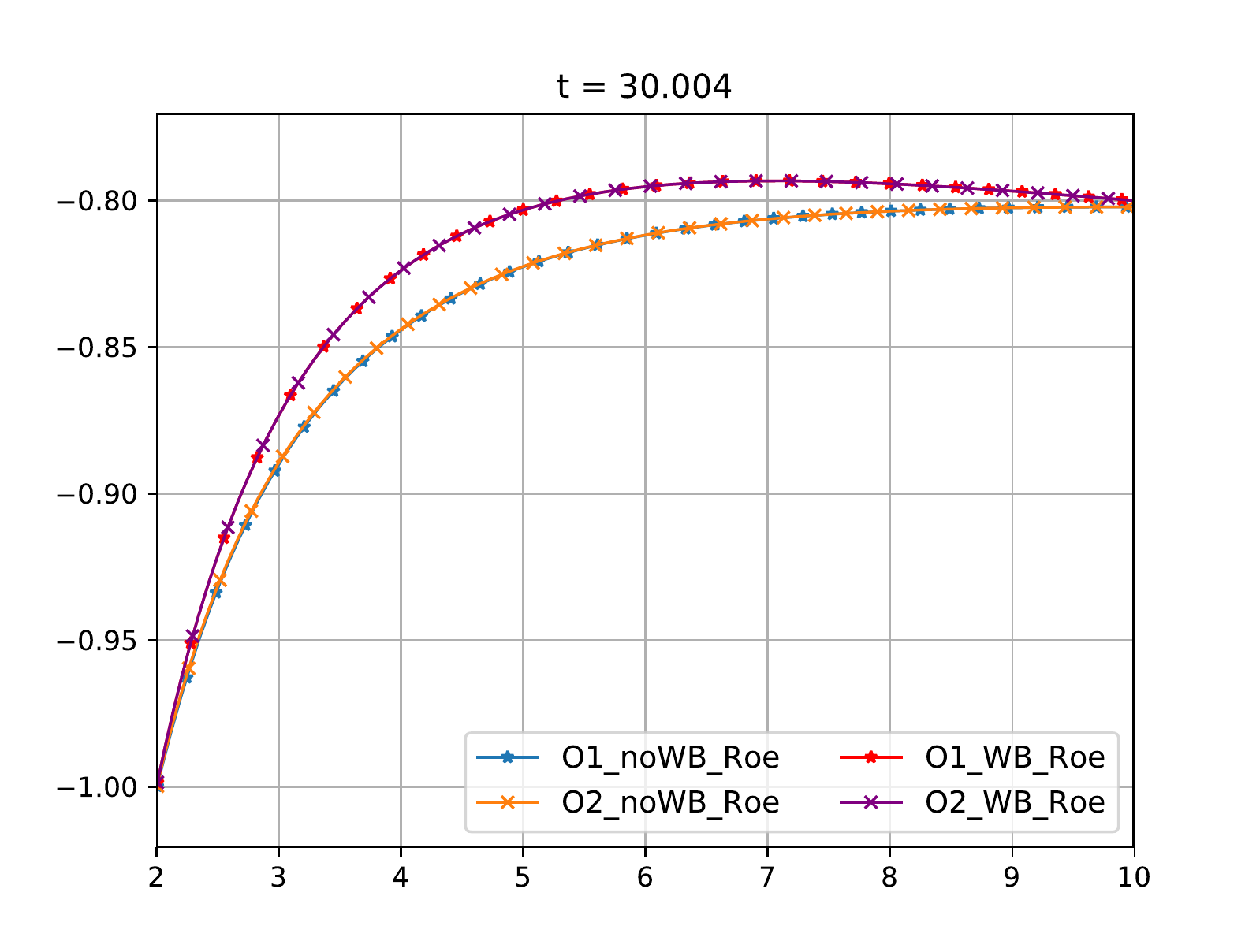}
		\label{fig:Euler_ko1_ko2_WB_vs_noWB_testWB2_t_30_hepse14}
	\end{subfigure}
\begin{subfigure}[h]{0.32\textwidth}
	\centering
	\includegraphics[width=1\linewidth]{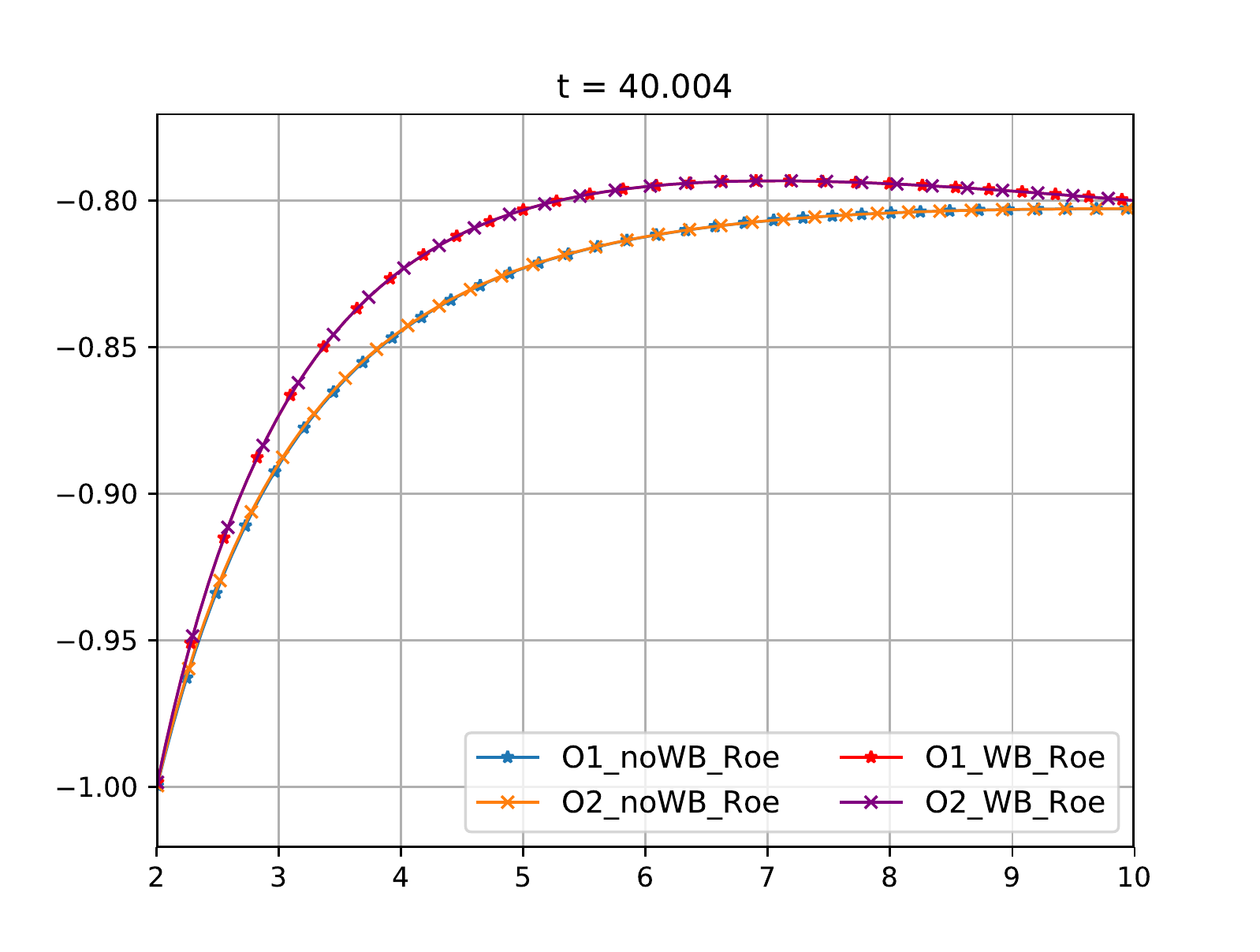}
	\label{fig:Euler_ko1_ko2_WB_vs_noWB_testWB2_t_40_hepse14}
\end{subfigure}
\begin{subfigure}[h]{0.32\textwidth}
	\centering
	\includegraphics[width=1\linewidth]{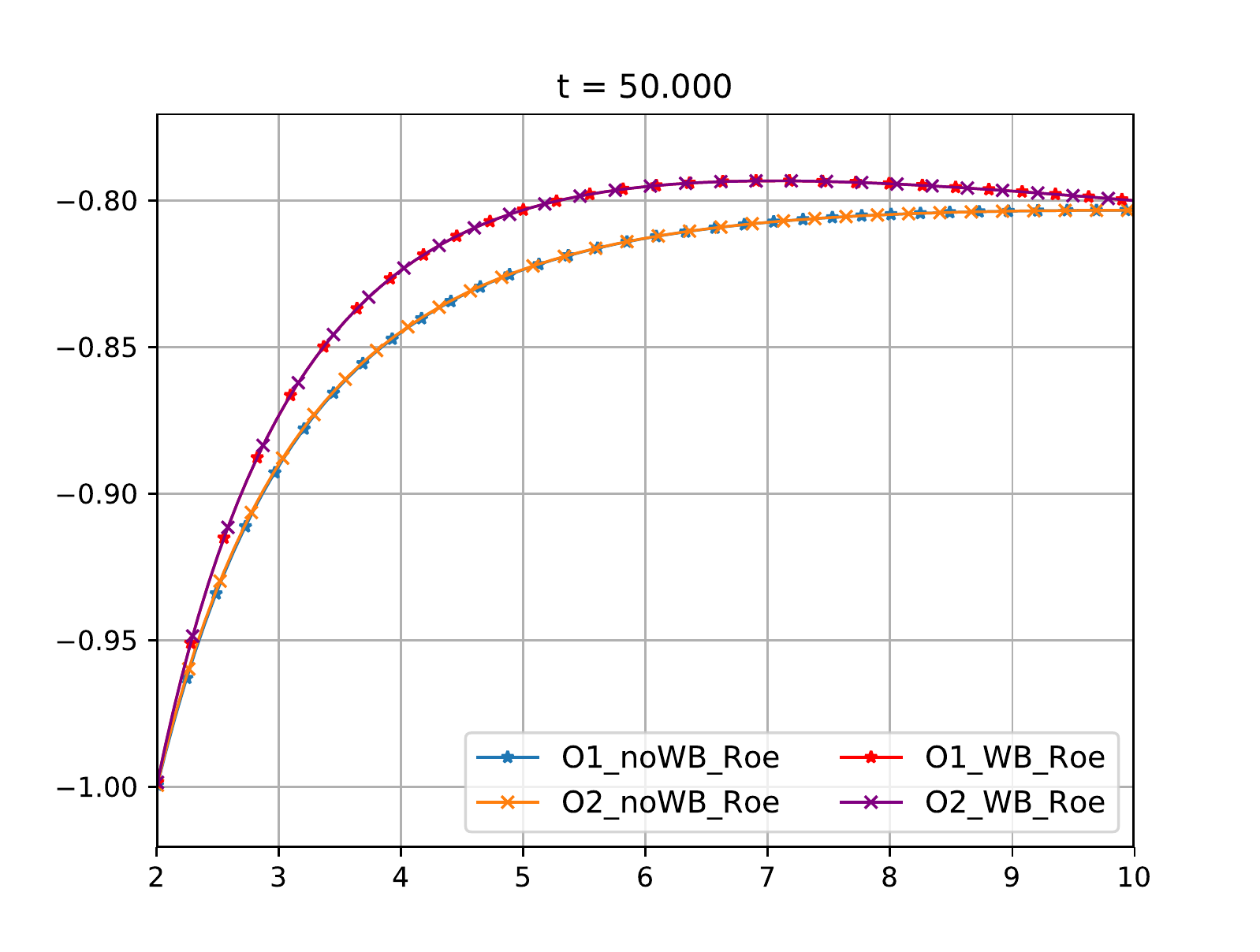}
	\label{fig:Euler_ko1_ko2_WB_vs_noWB_testWB2_t_50_hepse14}
\end{subfigure}
	\caption{Euler-Schwarzschild model with the initial condition the stationary solution satisfying \eqref{testE2}:  first- and second-order well-balanced and non-well-balanced methods at selected times for the variable $v$.}
	\label{fig:Euler_ko1_ko2_WB_vs_noWB_testWB2_hepse14}
\end{figure}

\begin{figure}[h]
	\begin{subfigure}[h]{0.32\textwidth}
		\centering
		\includegraphics[width=1\linewidth]{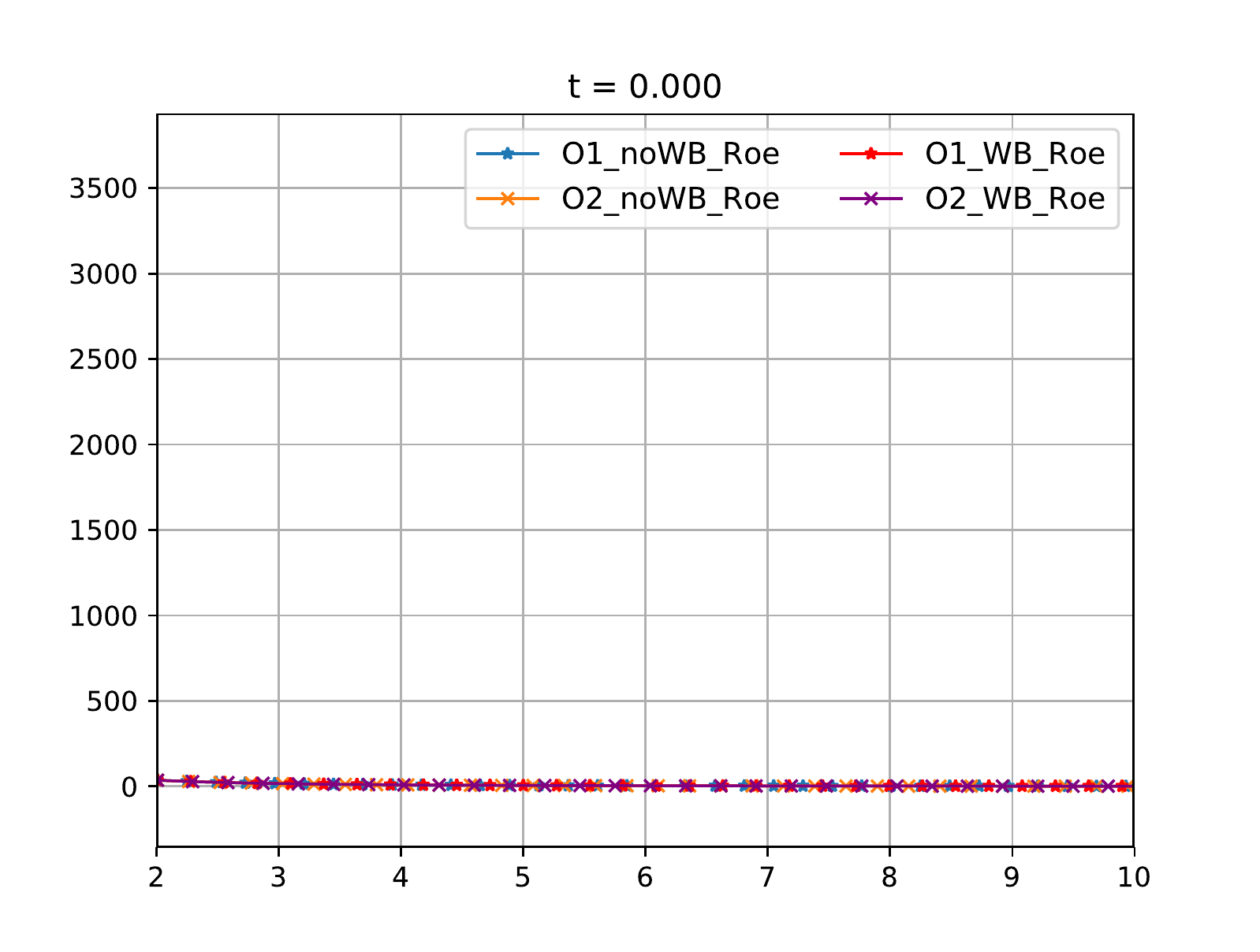}
		\label{fig:Euler_ko1_ko2_WB_vs_noWB_testWB2_t_0_hepse14_rho}
	\end{subfigure}
	\begin{subfigure}[h]{0.32\textwidth}
		\centering
		\includegraphics[width=1\linewidth]{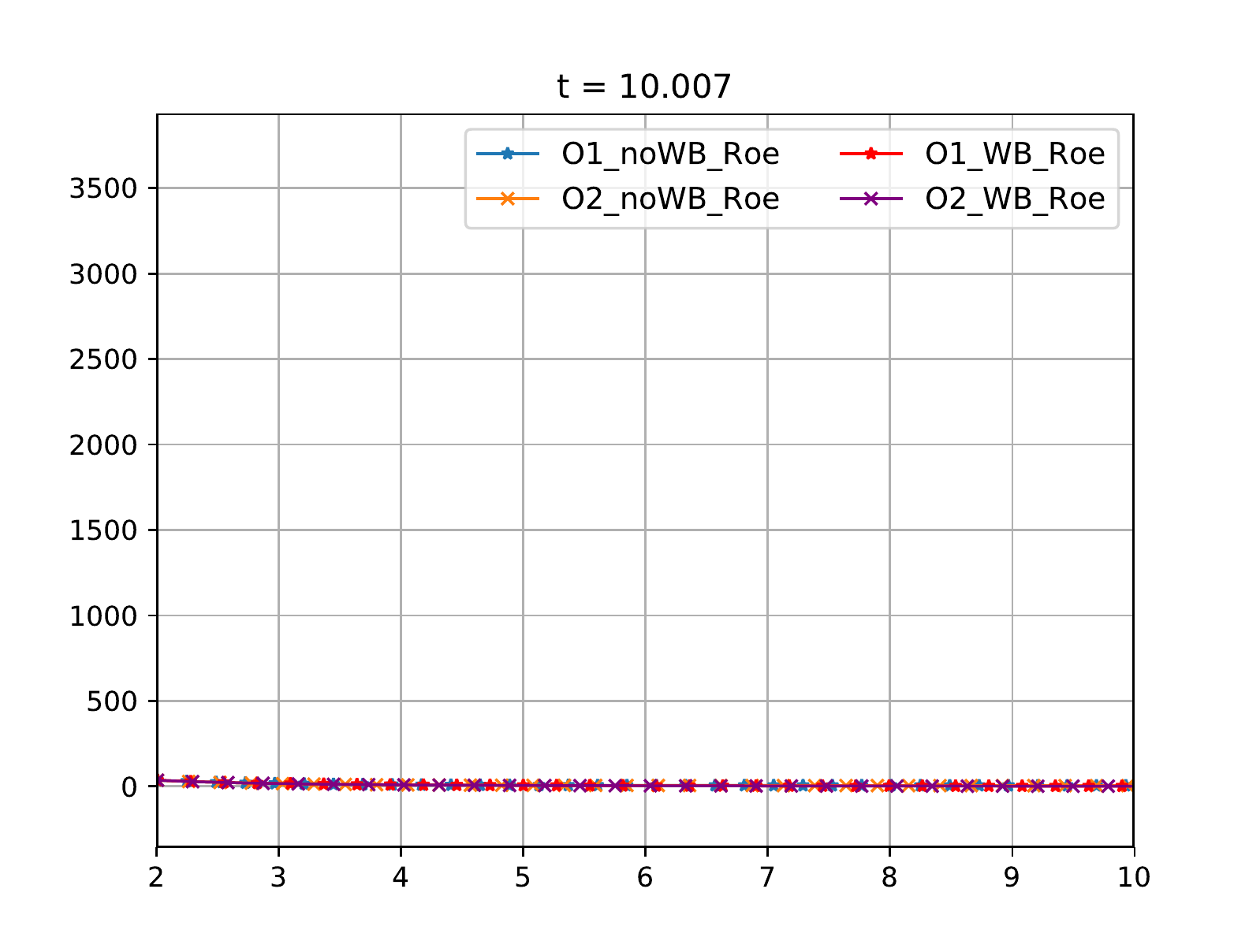}
		\label{fig:Euler_ko1_ko2_WB_vs_noWB_testWB2_t_10_hepse14_rho}
	\end{subfigure}
	\begin{subfigure}[h]{0.32\textwidth}
		\centering
		\includegraphics[width=1\linewidth]{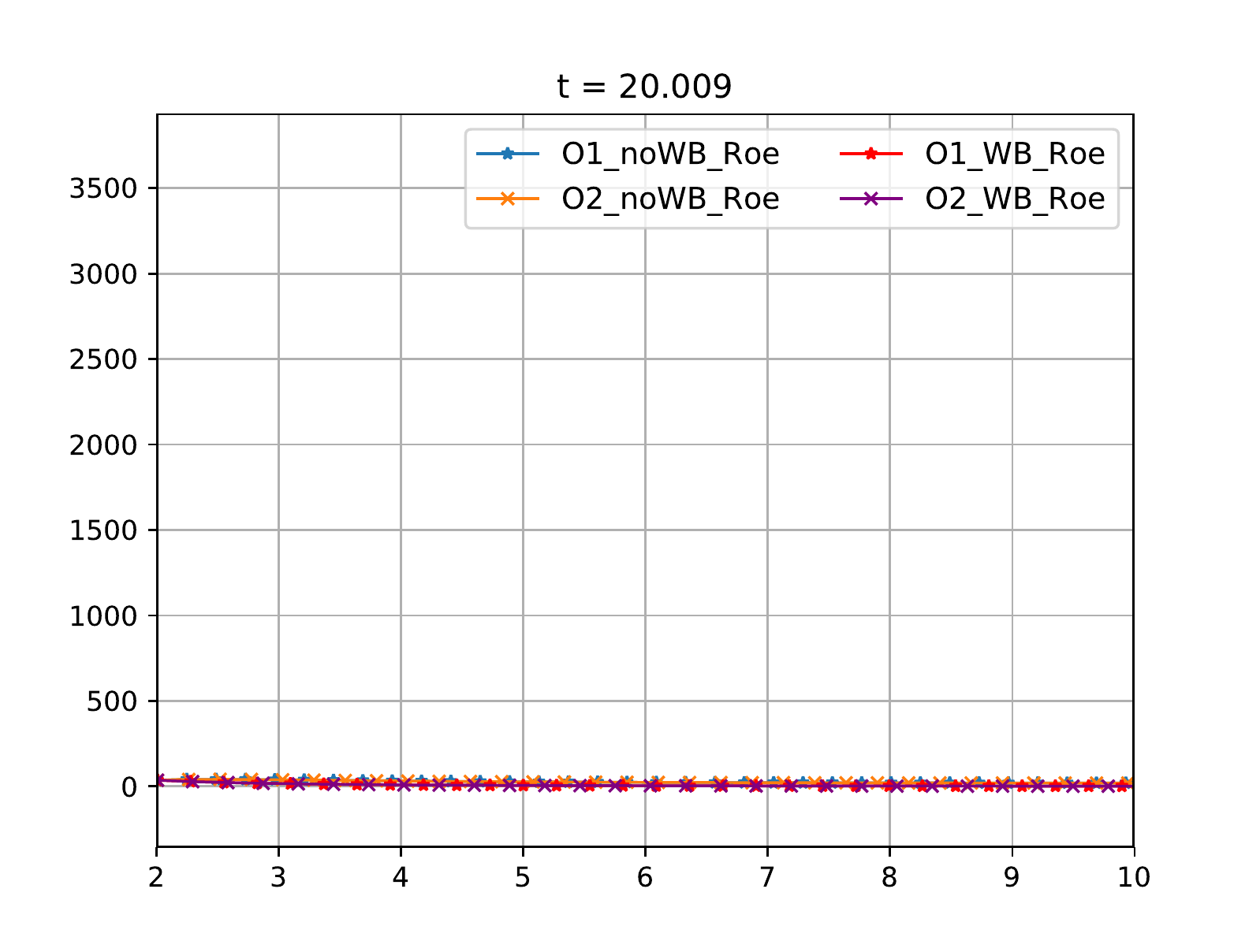}
		\label{fig:Euler_ko1_ko2_WB_vs_noWB_testWB2_t_20_hepse14_rho}
	\end{subfigure}
	\begin{subfigure}[h]{0.32\textwidth}
		\centering
		\includegraphics[width=1\linewidth]{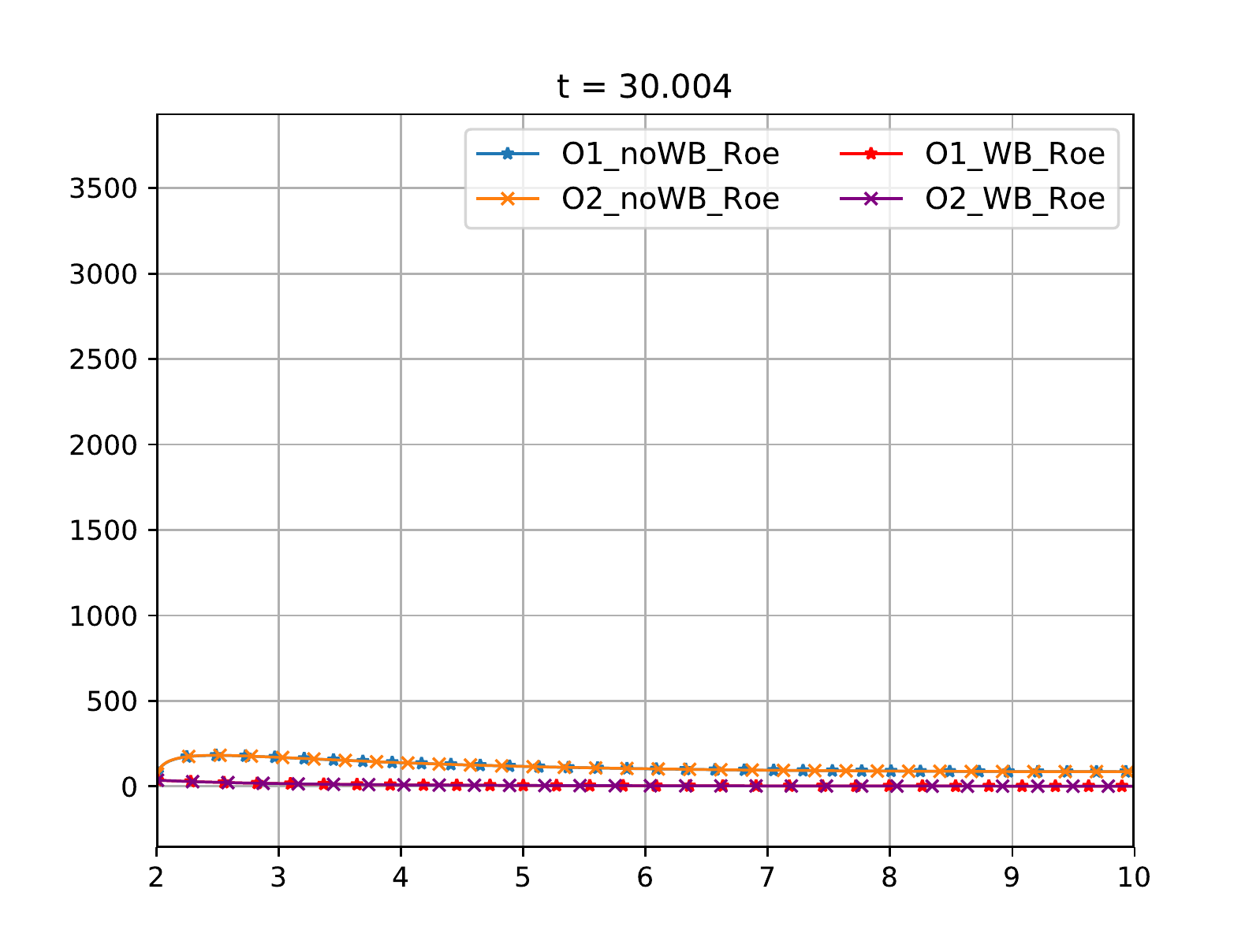}
		\label{fig:Euler_ko1_ko2_WB_vs_noWB_testWB2_t_30_hepse14_rho}
	\end{subfigure}
\begin{subfigure}[h]{0.32\textwidth}
	\centering
	\includegraphics[width=1\linewidth]{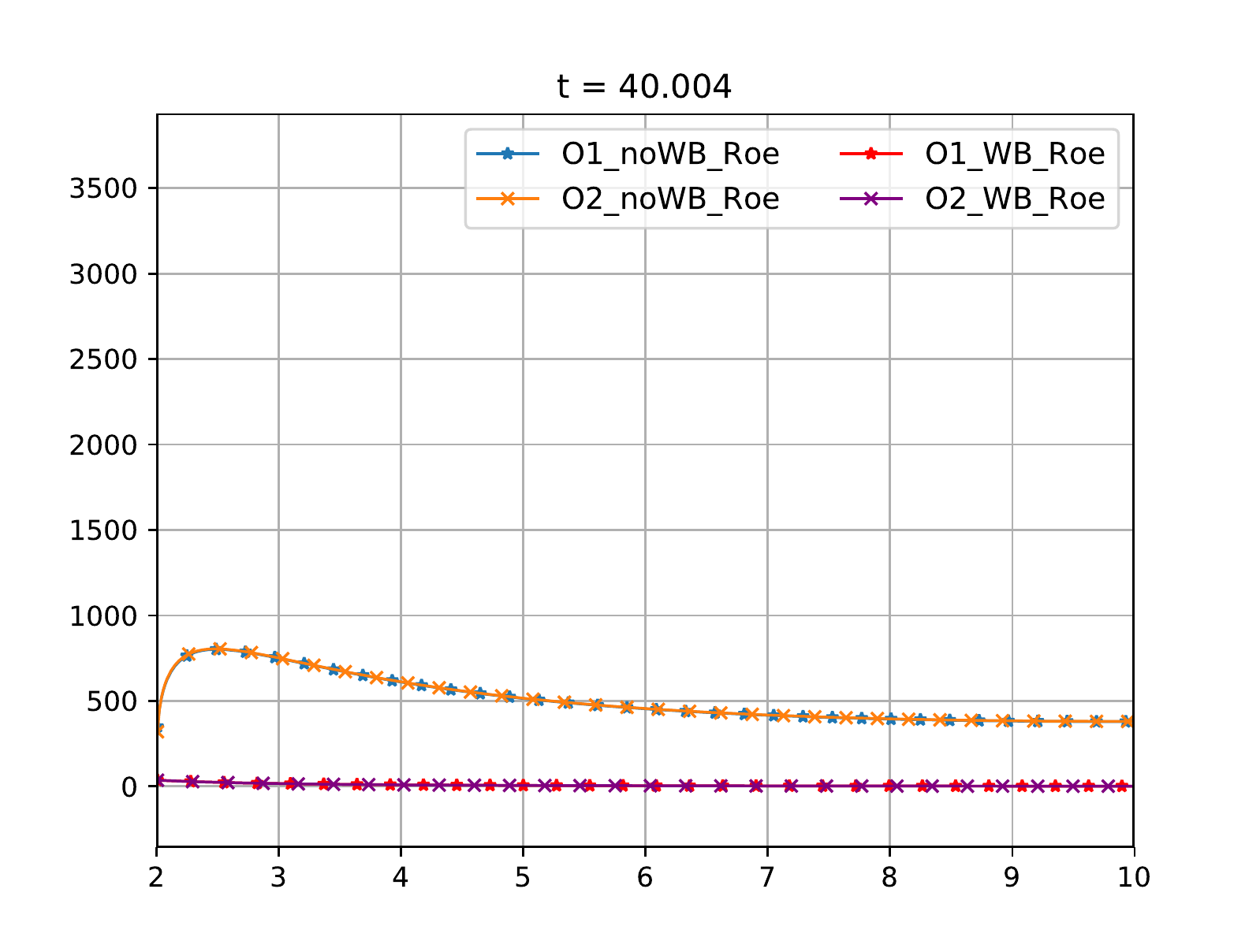}
	\label{fig:Euler_ko1_ko2_WB_vs_noWB_testWB2_t_40_hepse14_rho}
\end{subfigure}
\begin{subfigure}[h]{0.32\textwidth}
	\centering
	\includegraphics[width=1\linewidth]{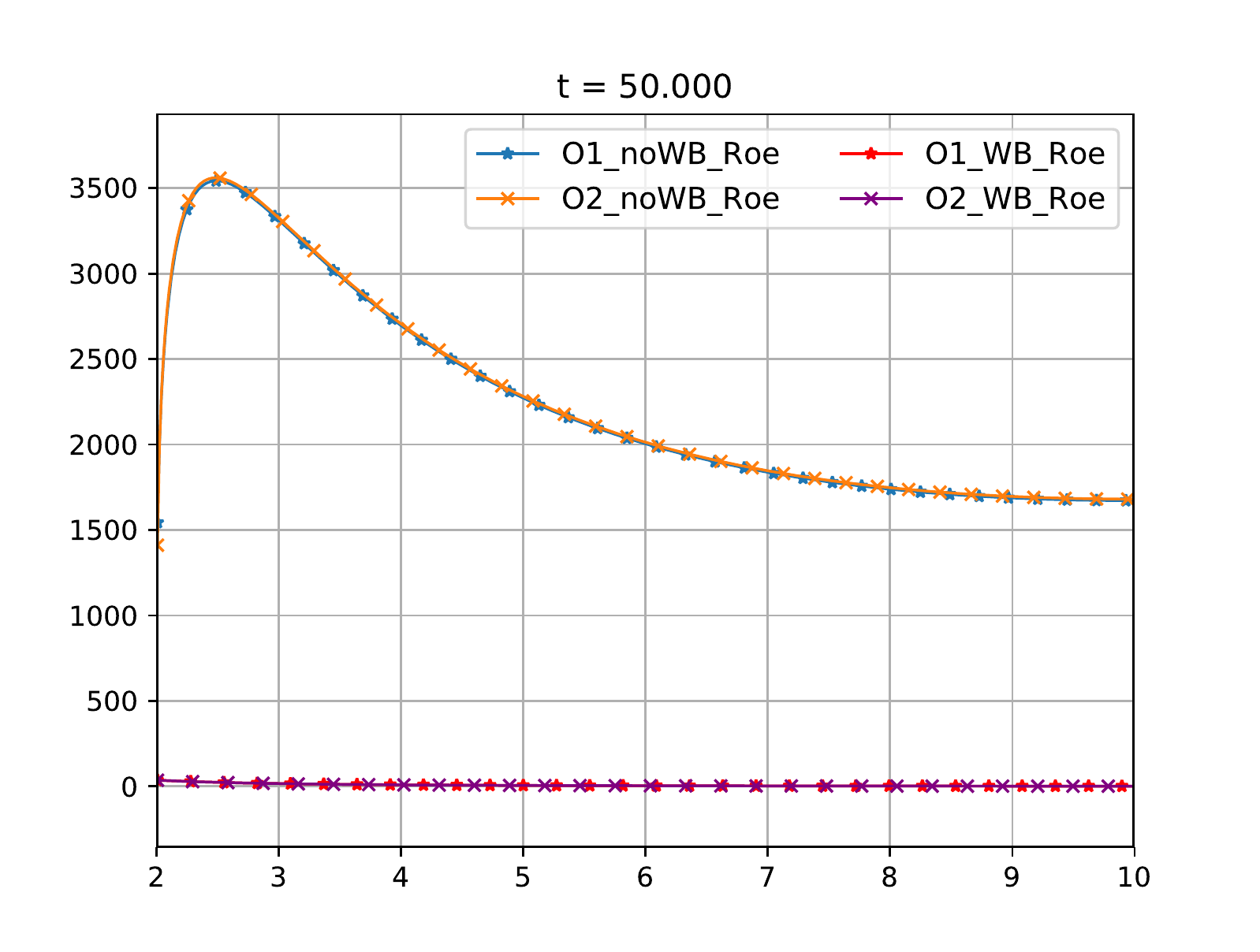}
	\label{fig:Euler_ko1_ko2_WB_vs_noWB_testWB2_t_50_hepse14_rho}
\end{subfigure}
	\caption{Euler-Schwarzschild model with the initial condition  \eqref{testE2}:  first- and second-order well-balanced and non-well-balanced methods at selected times for the variable $\rho$.}
	\label{fig:Euler_ko1_ko2_WB_vs_noWB_testWB2_hepse14_rho}
\end{figure}

Let us consider now as initial condition the negative supersonic stationary solution $V^*$ that satisfies
\bel{testE2}
\rho^{*}(10)=1, \quad v^{*}(10)=-0.8.
\ee
Table \ref{tab:Error_TestE2} shows the error in $L^1$ norm between the numerical solution
at time $ t = 50$.
 Figures \ref{fig:Euler_ko1_ko2_WB_vs_noWB_testWB2_hepse14}, \ref{fig:Euler_ko1_ko2_WB_vs_noWB_testWB2_hepse14_rho} show the difference between the numerical results given by well-balanced and non-well-balanced methods using the Roe-type numerical flux. Again the numerical solutions obtained with non-well-balanced methods depart from the initial steady state.
 
 \begin{table}[ht]
 	\centering
 	\begin{tabular}{|c|c|c|c|c|}
 		\hline 
 		Scheme (500 cells) & Error $v$ (1st) & Error $\rho$ (1st) & Error $v$ (2nd) & Error $\rho$ (2nd) \\ 
 		\hline 
 		Well-balanced & 1.54E-15 & 7.02E-13 & 1.35E-15 & 5.01E-13 \\ 
 		\hline 
 		Non well-balance & 0.01 & 2240.72 & 0.01 & 2250.77 \\ 
 		
 		\hline 
 	\end{tabular} 
	 	\caption{Well-balanced versus non-well-balanced schemes: $L^{1}$ errors at time $t=50$ for the Burgers-Schwarzschild model with the initial condition  (\ref{testE2})}

 	\label{tab:Error_TestE2}
 \end{table}

\paragraph{{{\red Discontinuous stationary entropy weak solution}}}

\begin{figure}[h]
	\begin{subfigure}[h]{0.32\textwidth}
		\centering
		\includegraphics[width=1\linewidth]{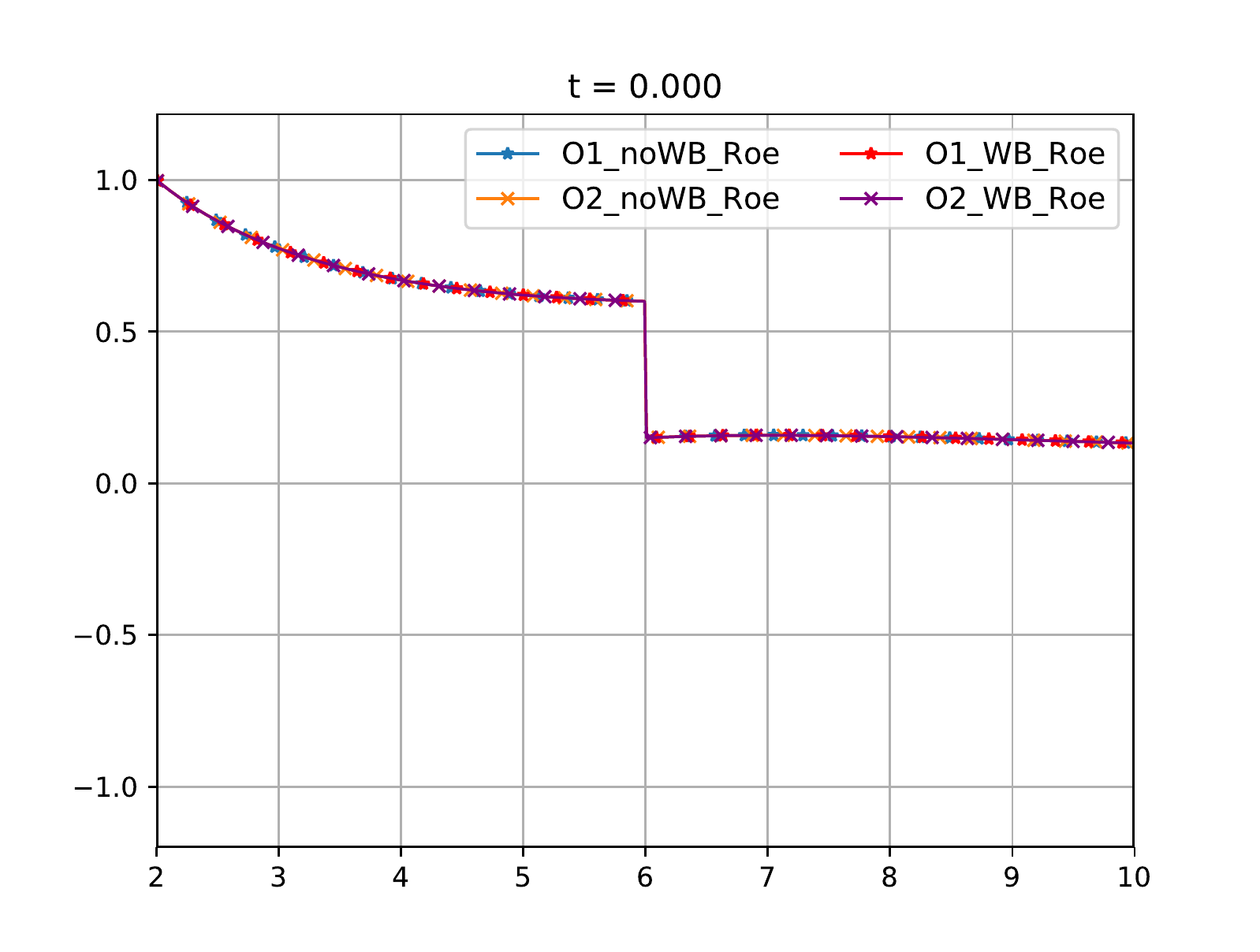}
		\label{fig:Euler_ko1_ko2_WB_vs_noWB_testWB3_t_0_hepse14}
	\end{subfigure}
	\begin{subfigure}[h]{0.32\textwidth}
		\centering
		\includegraphics[width=1\linewidth]{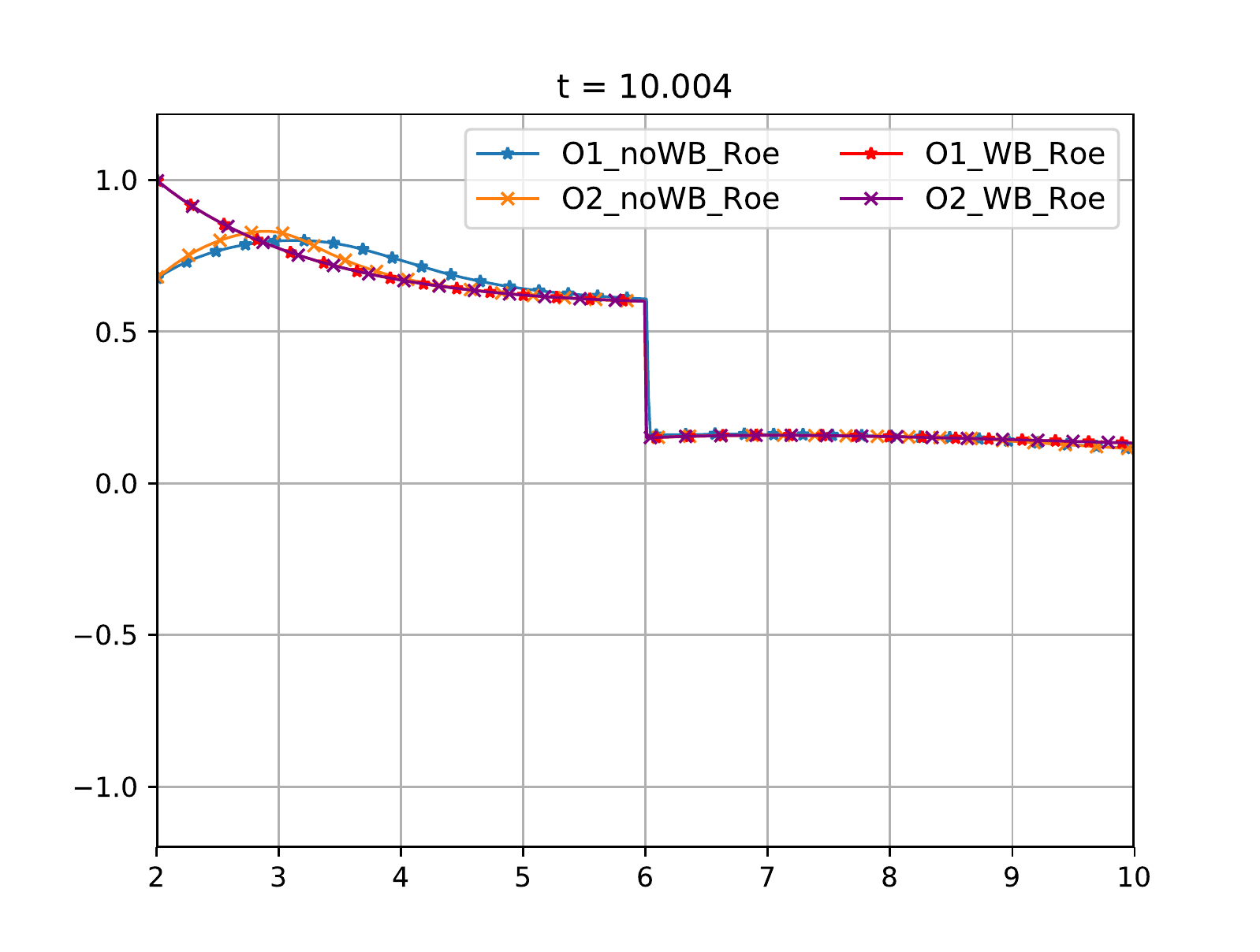}
		\label{fig:Euler_ko1_ko2_WB_vs_noWB_testWB3_t_10_hepse14}
	\end{subfigure}
	\begin{subfigure}[h]{0.32\textwidth}
		\centering
		\includegraphics[width=1\linewidth]{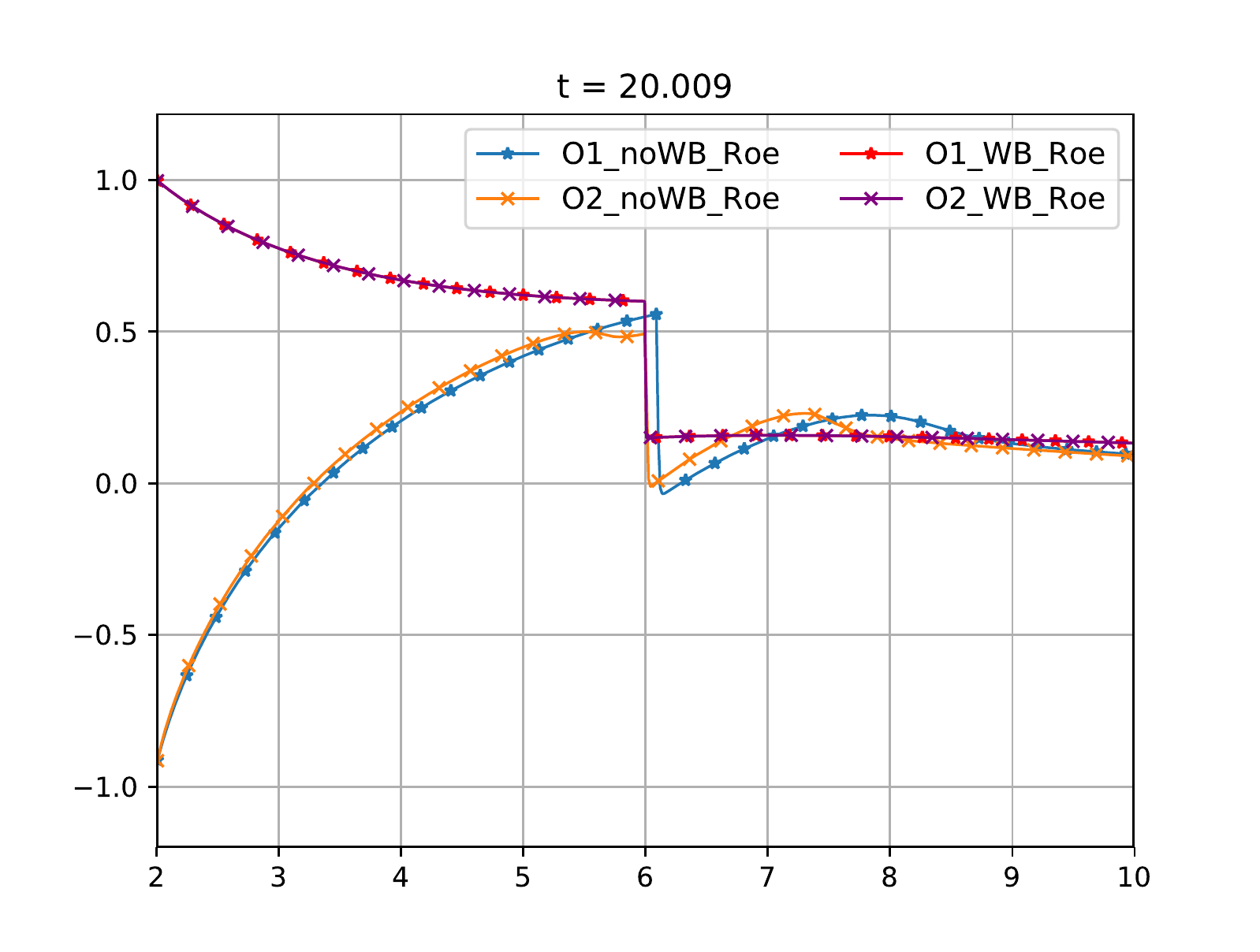}
		\label{fig:Euler_ko1_ko2_WB_vs_noWB_testWB3_t_20_hepse14}
	\end{subfigure}
	\begin{subfigure}[h]{0.32\textwidth}
		\centering
		\includegraphics[width=1\linewidth]{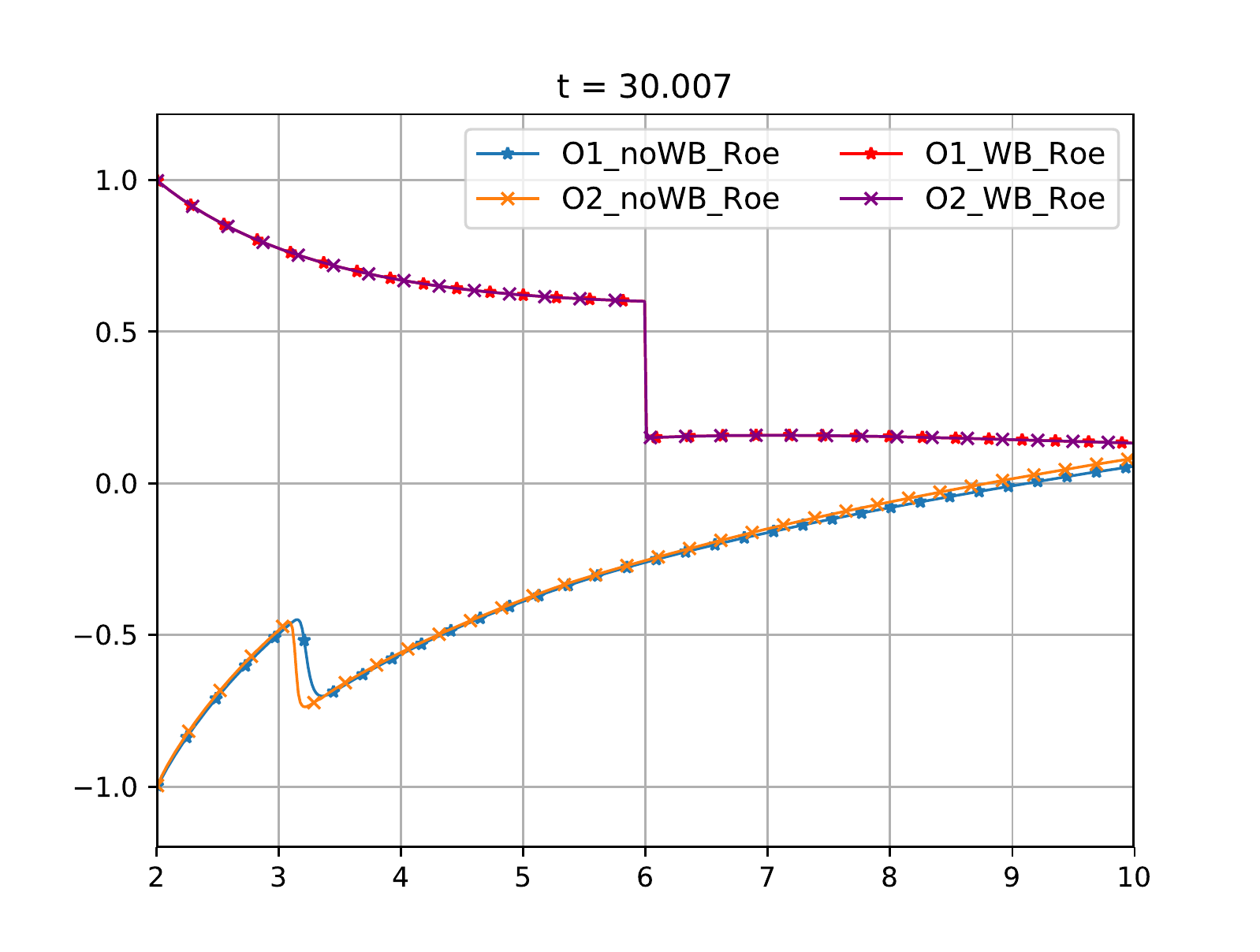}
		\label{fig:Euler_ko1_ko2_WB_vs_noWB_testWB3_t_30_hepse14}
	\end{subfigure}
	\begin{subfigure}[h]{0.32\textwidth}
		\centering
		\includegraphics[width=1\linewidth]{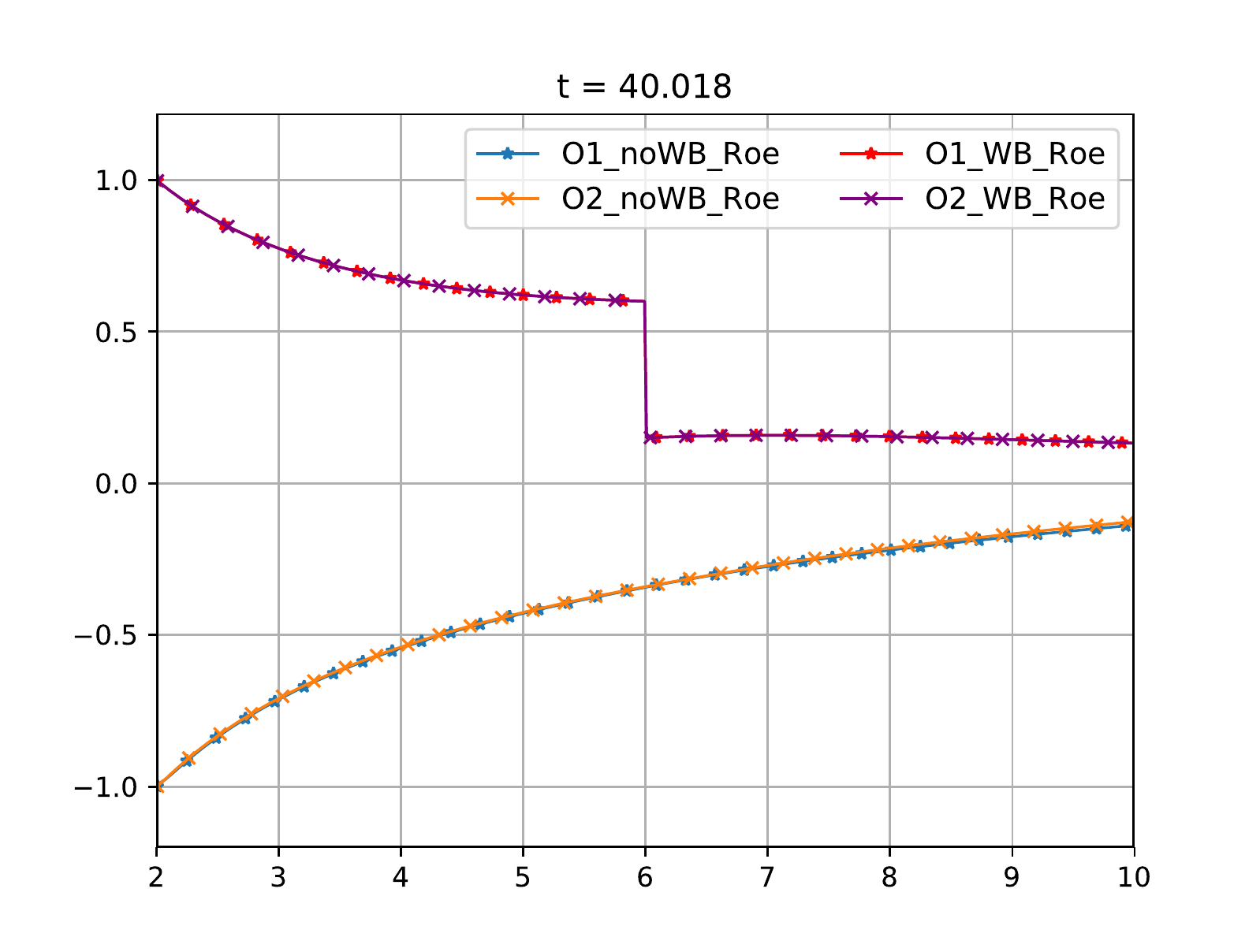}
		\label{fig:Euler_ko1_ko2_WB_vs_noWB_testWB3_t_40_hepse14}
	\end{subfigure}
	\begin{subfigure}[h]{0.32\textwidth}
		\centering
		\includegraphics[width=1\linewidth]{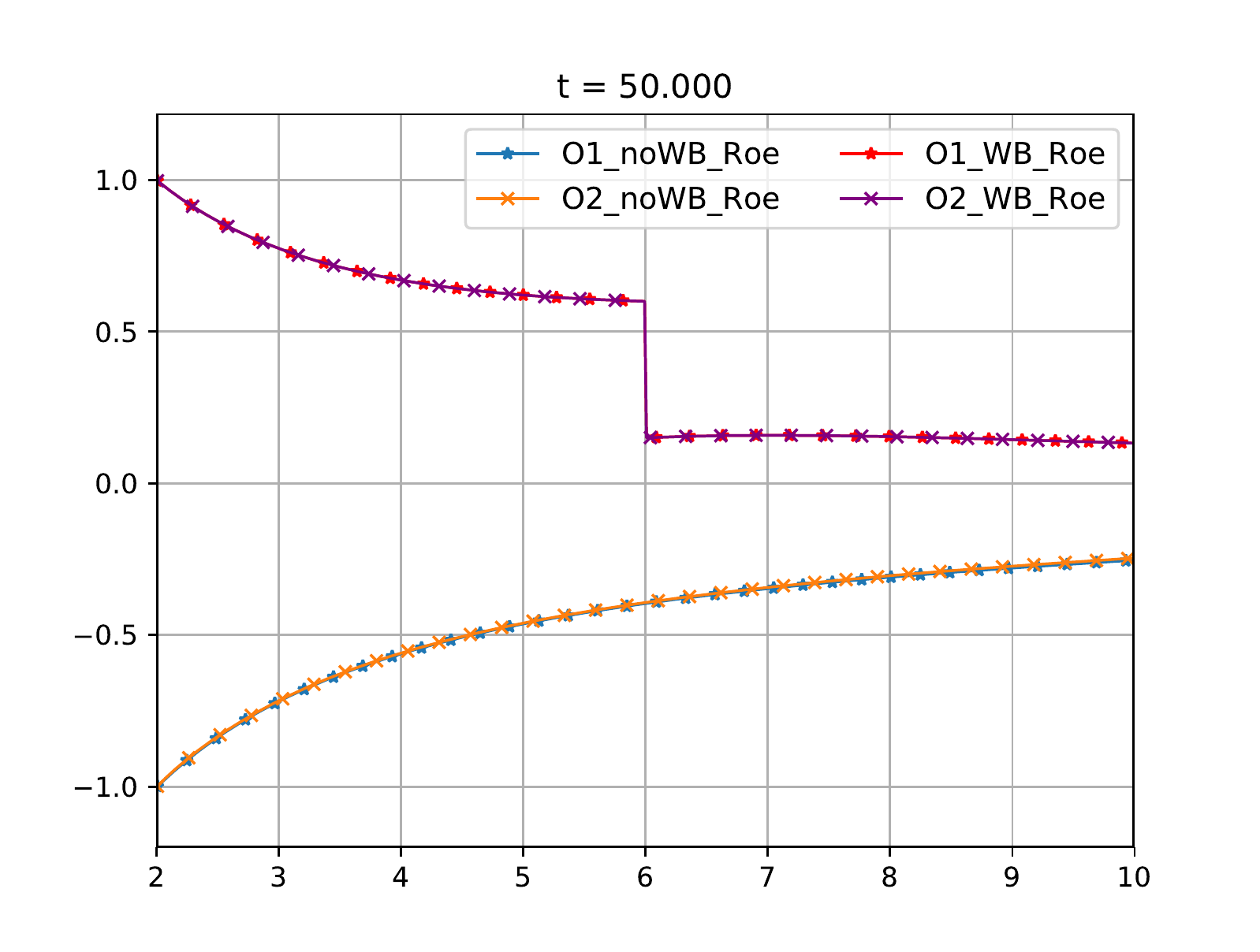}
		\label{fig:Euler_ko1_ko2_WB_vs_noWB_testWB3_t_50_hepse14}
	\end{subfigure}
	\caption{Euler-Schwarzschild model with the initial condition \eqref{eq:testE3a}:  first- and second-order well-balanced and non-well-balanced methods at selected times for the variable $v$.}
	\label{fig:Euler_ko1_ko2_WB_vs_noWB_testWB3_hepse14}
\end{figure}

\begin{figure}[h]
	\begin{subfigure}[h]{0.32\textwidth}
		\centering
		\includegraphics[width=1\linewidth]{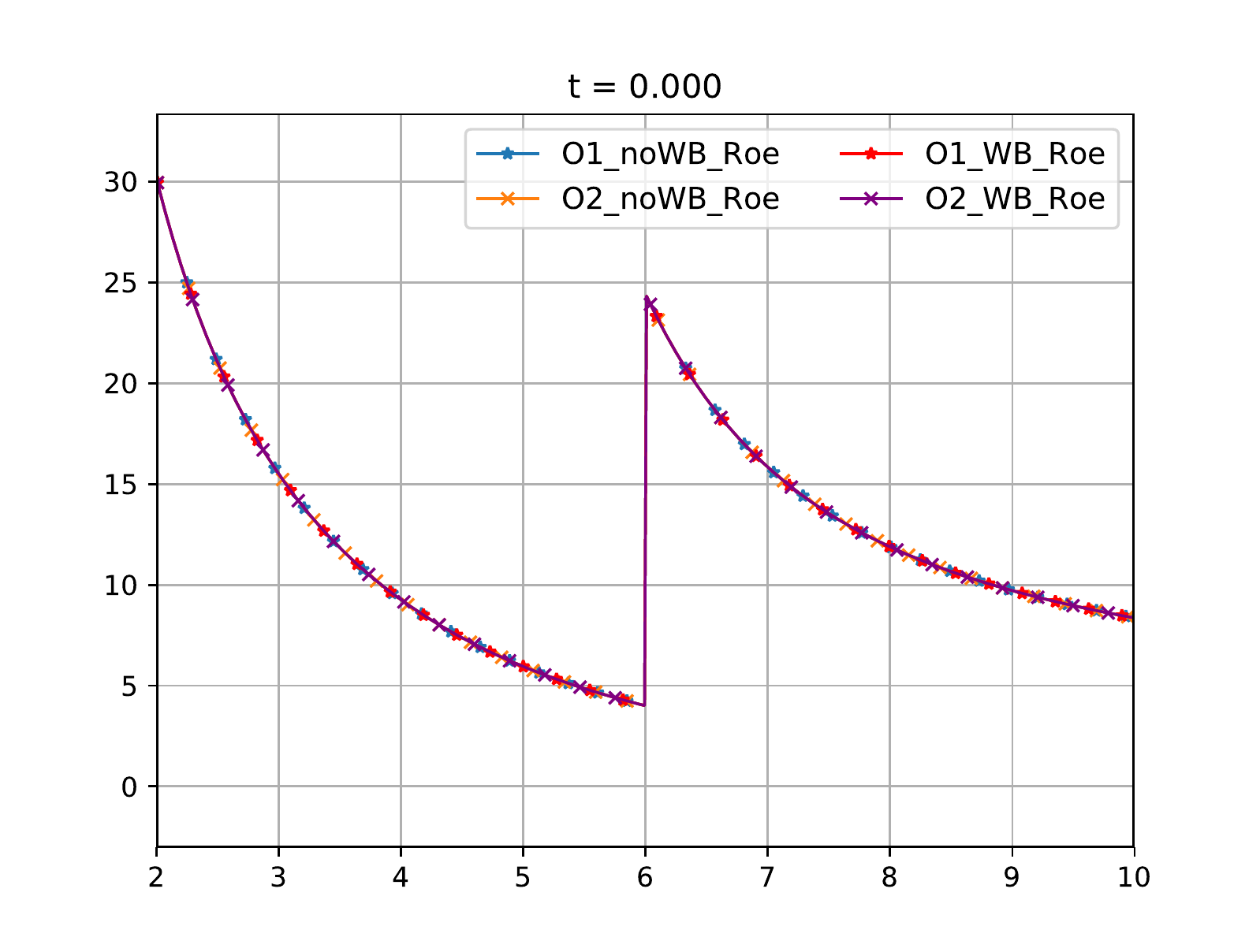}
		\label{fig:Euler_ko1_ko2_WB_vs_noWB_testWB3_t_0_hepse14_rho}
	\end{subfigure}
	\begin{subfigure}[h]{0.32\textwidth}
		\centering
		\includegraphics[width=1\linewidth]{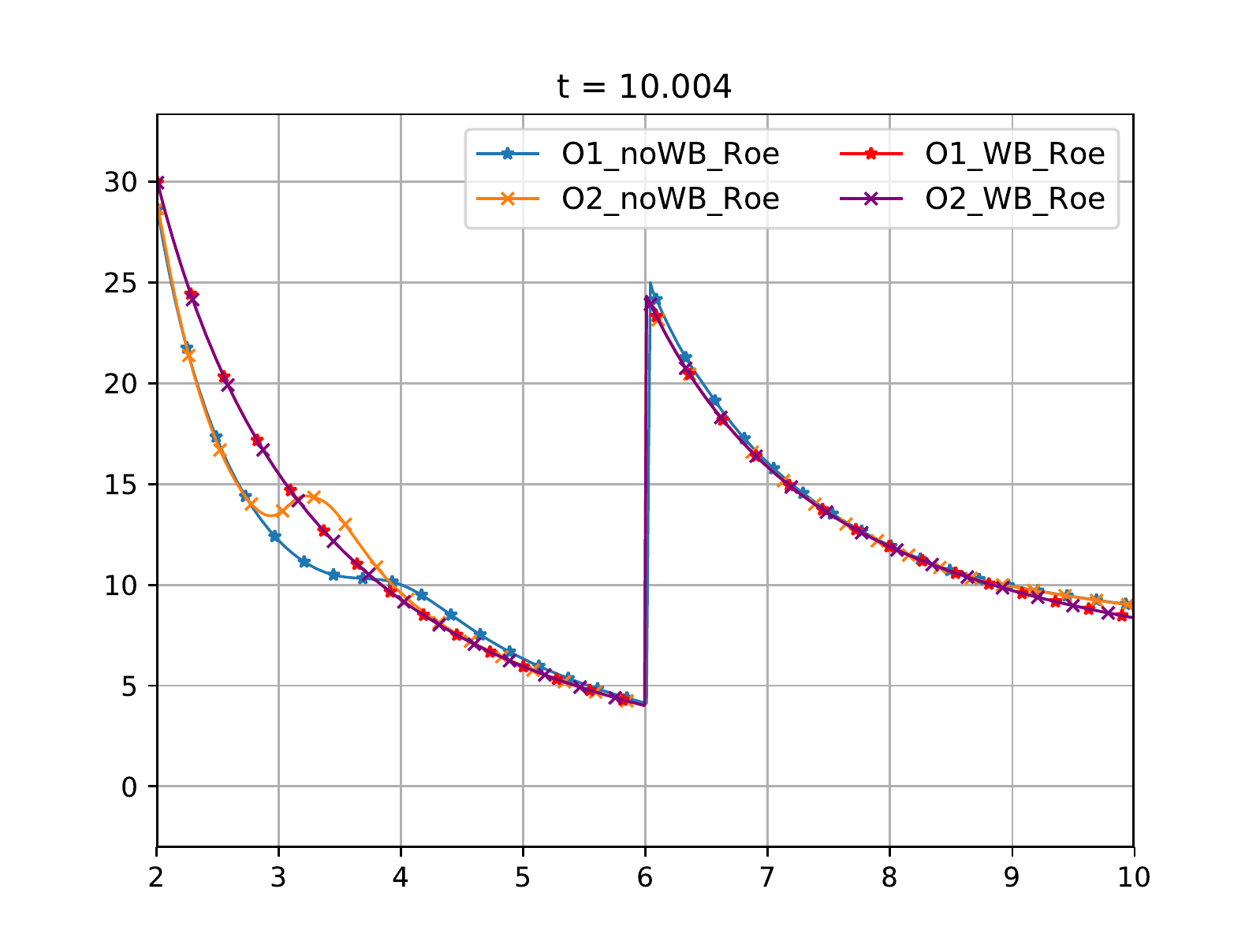}
		\label{fig:Euler_ko1_ko2_WB_vs_noWB_testWB3_t_10_hepse14_rho}
	\end{subfigure}
	\begin{subfigure}[h]{0.32\textwidth}
		\centering
		\includegraphics[width=1\linewidth]{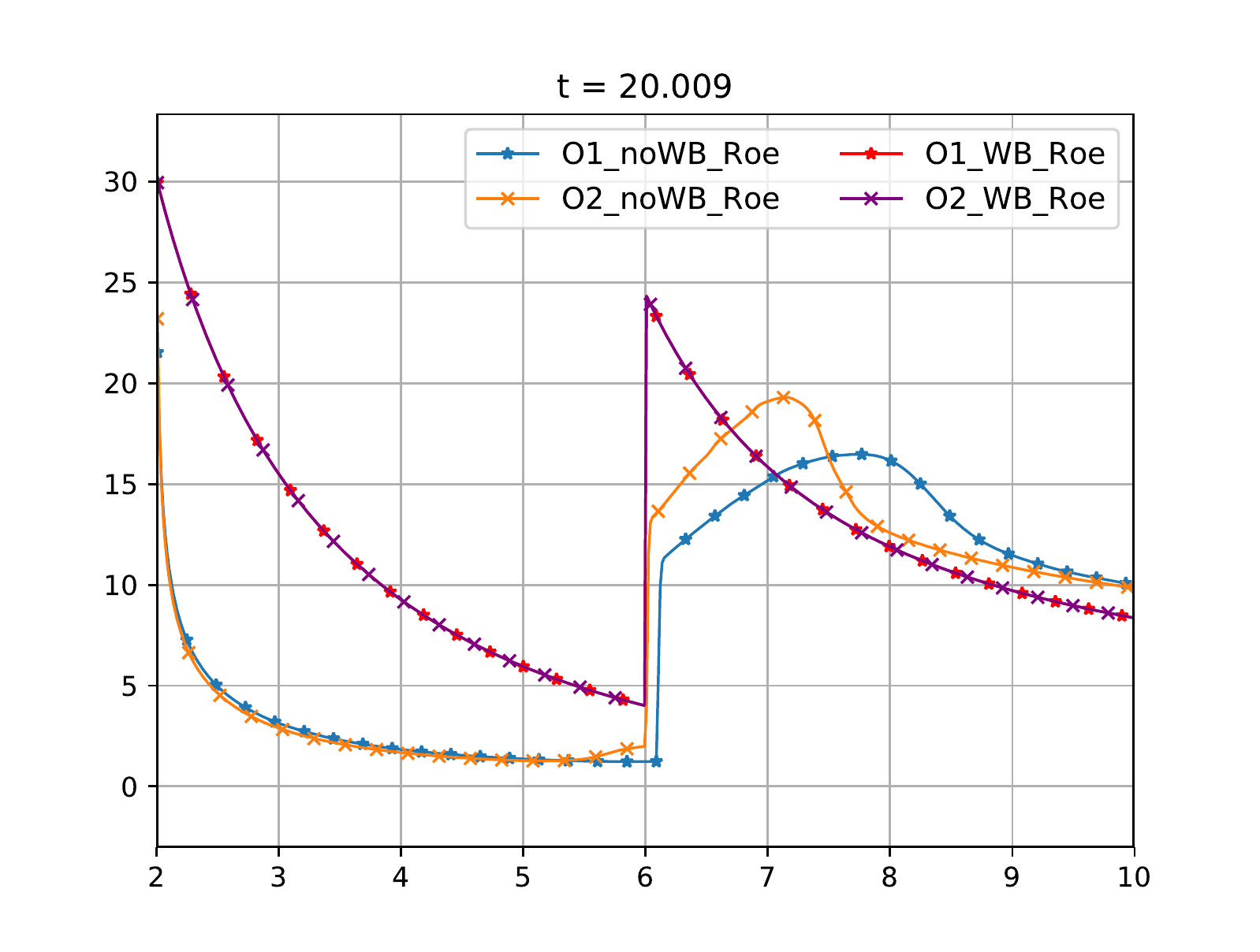}
		\label{fig:Euler_ko1_ko2_WB_vs_noWB_testWB3_t_20_hepse14_rho}
	\end{subfigure}
	\begin{subfigure}[h]{0.32\textwidth}
		\centering
		\includegraphics[width=1\linewidth]{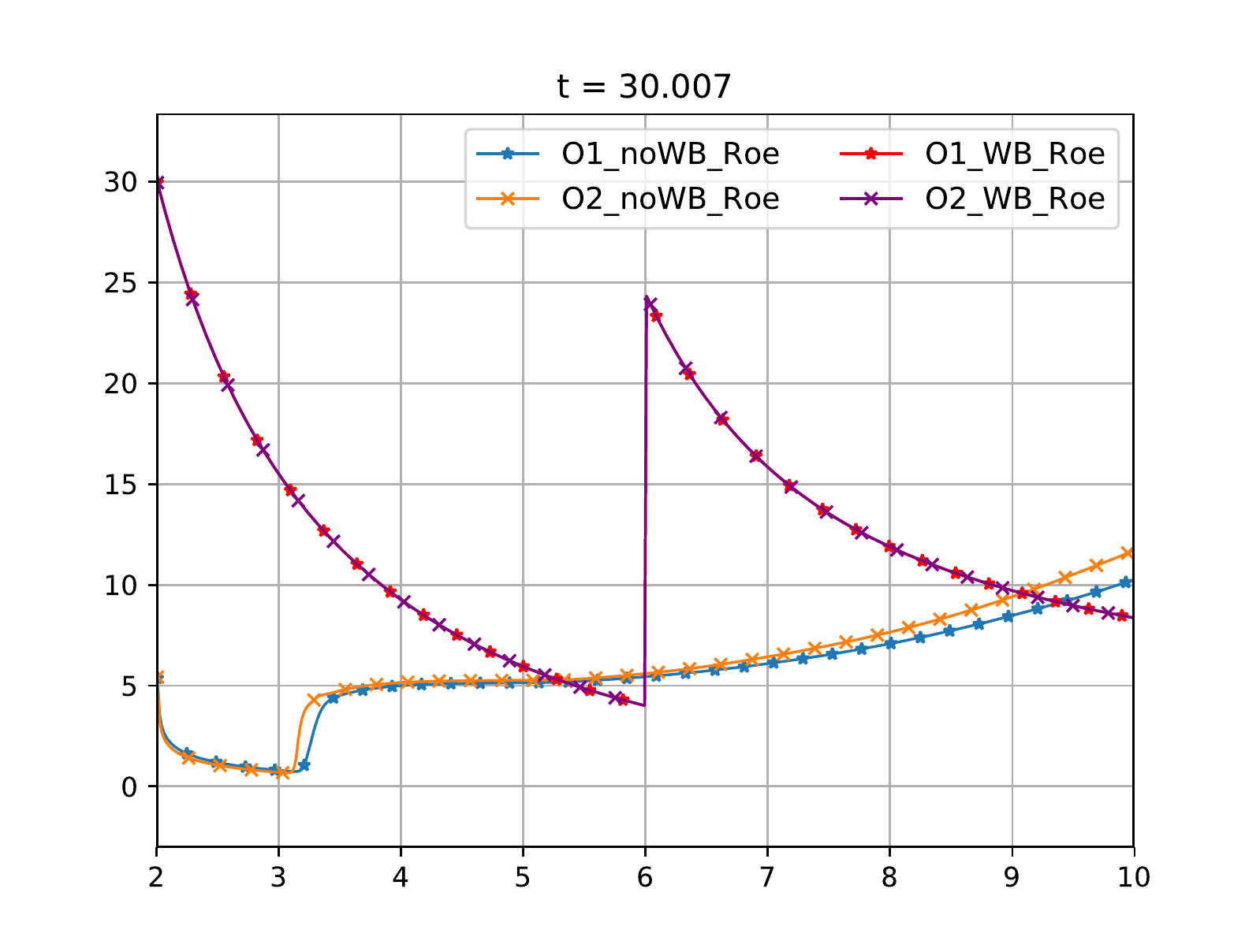}
		\label{fig:Euler_ko1_ko2_WB_vs_noWB_testWB3_t_30_hepse14_rho}
	\end{subfigure}
	\begin{subfigure}[h]{0.32\textwidth}
		\centering
		\includegraphics[width=1\linewidth]{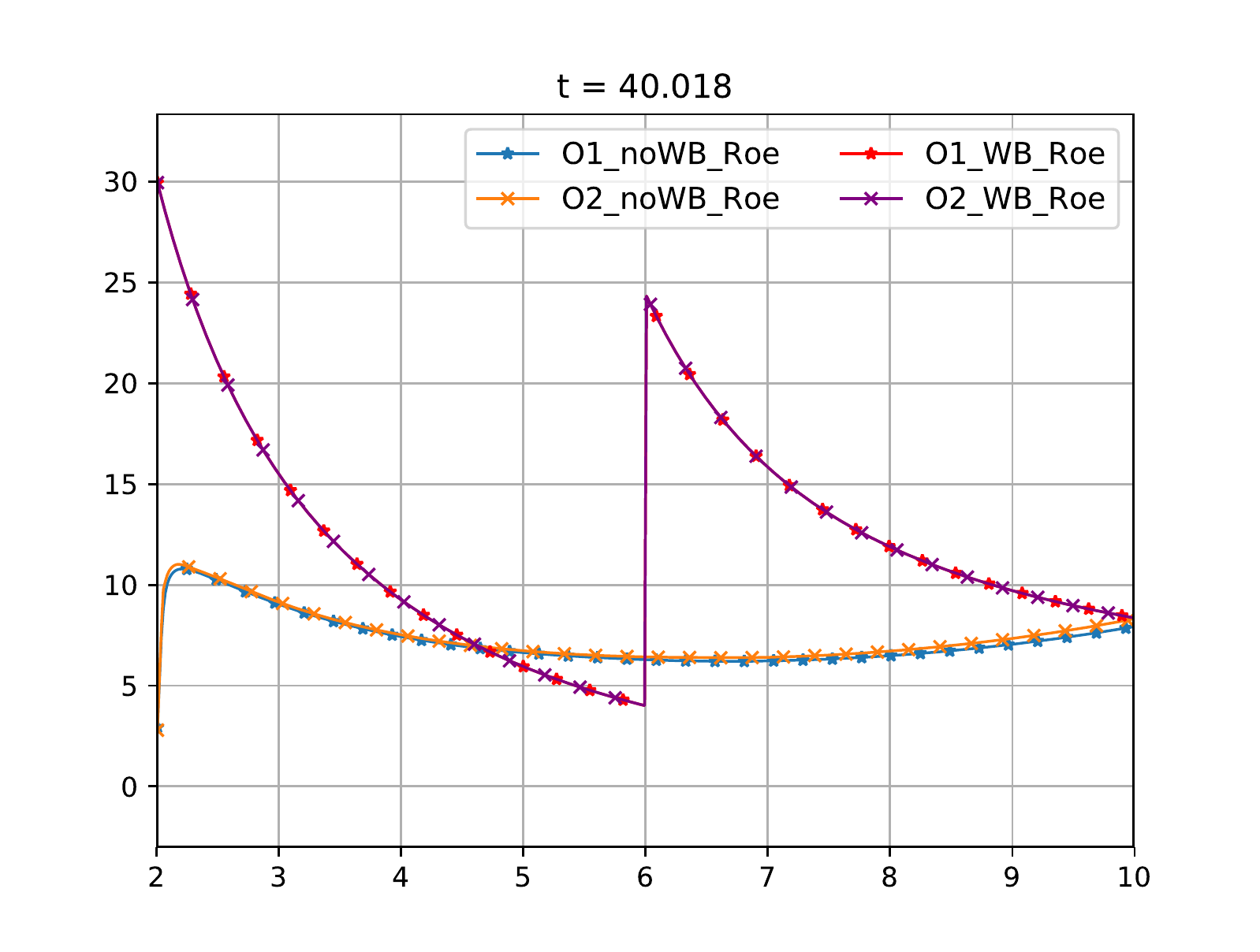}
		\label{fig:Euler_ko1_ko2_WB_vs_noWB_testWB3_t_40_hepse14_rho}
	\end{subfigure}
	\begin{subfigure}[h]{0.32\textwidth}
		\centering
		\includegraphics[width=1\linewidth]{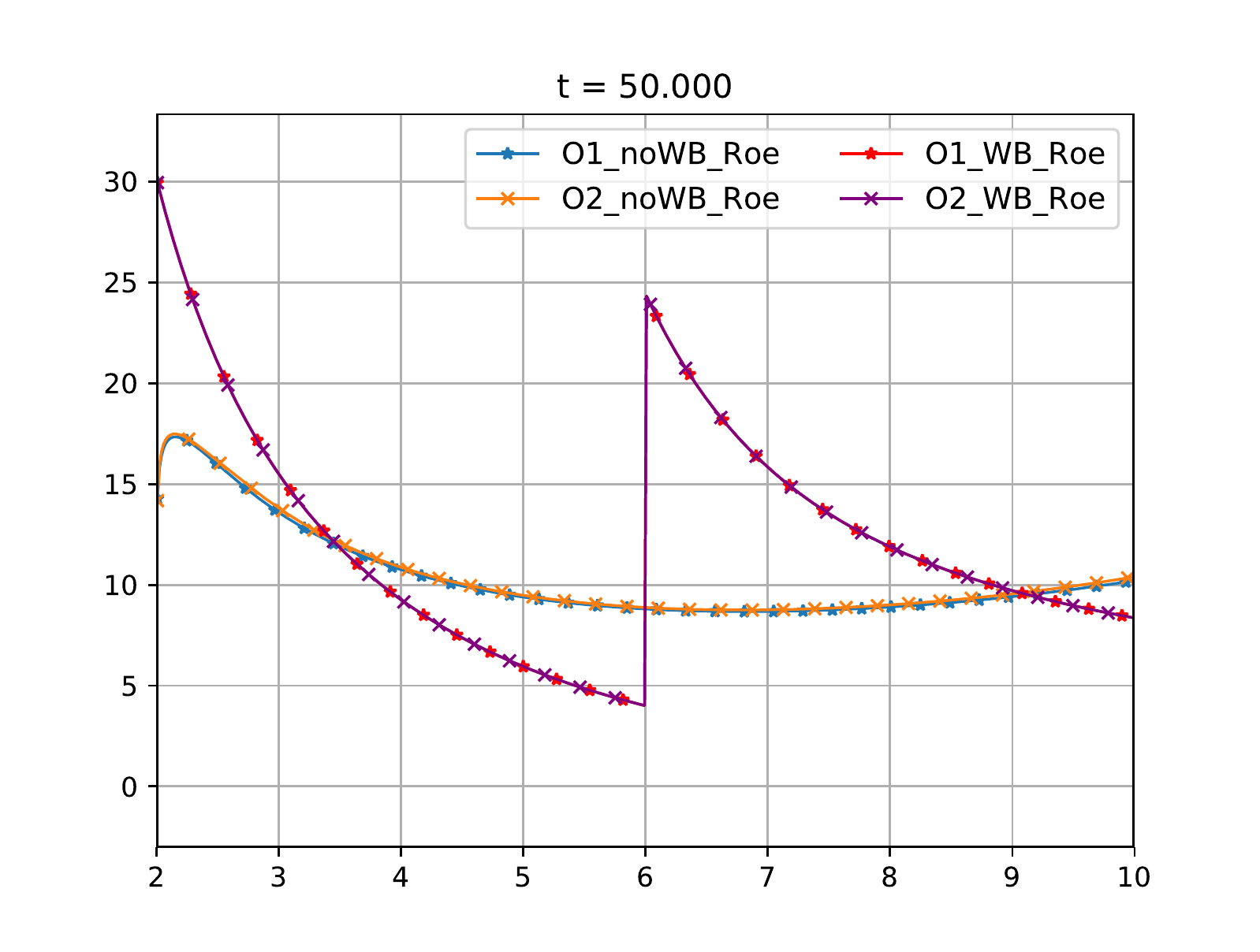}
		\label{fig:Euler_ko1_ko2_WB_vs_noWB_testWB3_t_50_hepse14_rho}
	\end{subfigure}
	\caption{Euler-Schwarzschild model with the initial condition \eqref{eq:testE3a}:  first- and second-order well-balanced and non-well-balanced methods at selected times for the variable $\rho$.}
	\label{fig:Euler_ko1_ko2_WB_vs_noWB_testWB3_hepse14_rho}
\end{figure}

We consider finally the initial condition
\bel{eq:testE3a}
V_{0}(r) = \begin{cases}
V_{-}^{*}(r), & \text{ $r\leq 6$},\\
V_{+}^{*}(r), & \text{ otherwise,}
\end{cases}
\ee 
where $V_{-}^{*}(r)$ is the supersonic stationary solution such that 
\bel{eq:testE3b}
\rho_{-}^{*}(6)=4, \quad v_{-}^{*}(6)=0.6
\ee
and $V_{+}^{*}(r)$ is the subsonic one such that 
\bel{eq:testE3c}
\rho_{+}^{*}(6)=\frac{\rho_{-}^{*}(6)(v_{-}^{*}(6)^{2}-k^{4})}{k^{2}(1-v_{-}^{*}(6)^{2})}, \quad v_{+}^{*}(6)=\frac{k^{2}}{v_{-}^{*}(6)}.
\ee
$V_0$ is an entropy weak stationary solution of the system: see \cite{PLF-SX1,PLF-SX2}.
Table \ref{tab:Error_TestE3} shows the error in $L^1$ norm between the numerical solution
at time $ t = 50$ and Figures \ref{fig:Euler_ko1_ko2_WB_vs_noWB_testWB3_hepse14}, \ref{fig:Euler_ko1_ko2_WB_vs_noWB_testWB3_hepse14_rho} show the comparison of the numerical results obtained with well-balanced and non-well-balanced methods at selected times.
On the other hand, the numerical results of this section put on evidence, as for the Burgers-Schwarzschild  system, the relevance of using well-balanced methods for the Euler-Schwarzschild model.

\begin{table}[ht]
	\centering
	\begin{tabular}{|c|c|c|c|c|}
		\hline 
		Scheme (500 cells) & Error $v$ (1st) & Error $\rho$ (1st) & Error $v$ (2nd) & Error $\rho$ (2nd) \\ 
		\hline 
		Well-balanced & 2.20E-13 & 1.25E-11 & 1.92E-13 & 1.03E-11 \\ 
		\hline 
		Non well-balance & 0.89 & 3.94 & 0.89 & 3.92 \\ 
		
		\hline 
	\end{tabular} 
		\caption{Well-balanced versus non-well-balanced schemes: $L^{1}$ errors at time $t=50$ for the Burgers-Schwarzschild model with the initial condition (\ref{eq:testE3a})}

	\label{tab:Error_TestE3}
\end{table}

\subsection{Perturbation of a regular stationary solution}

\begin{figure}[h]
	\begin{subfigure}[h]{0.32\textwidth}
		\centering
		\includegraphics[width=1\linewidth]{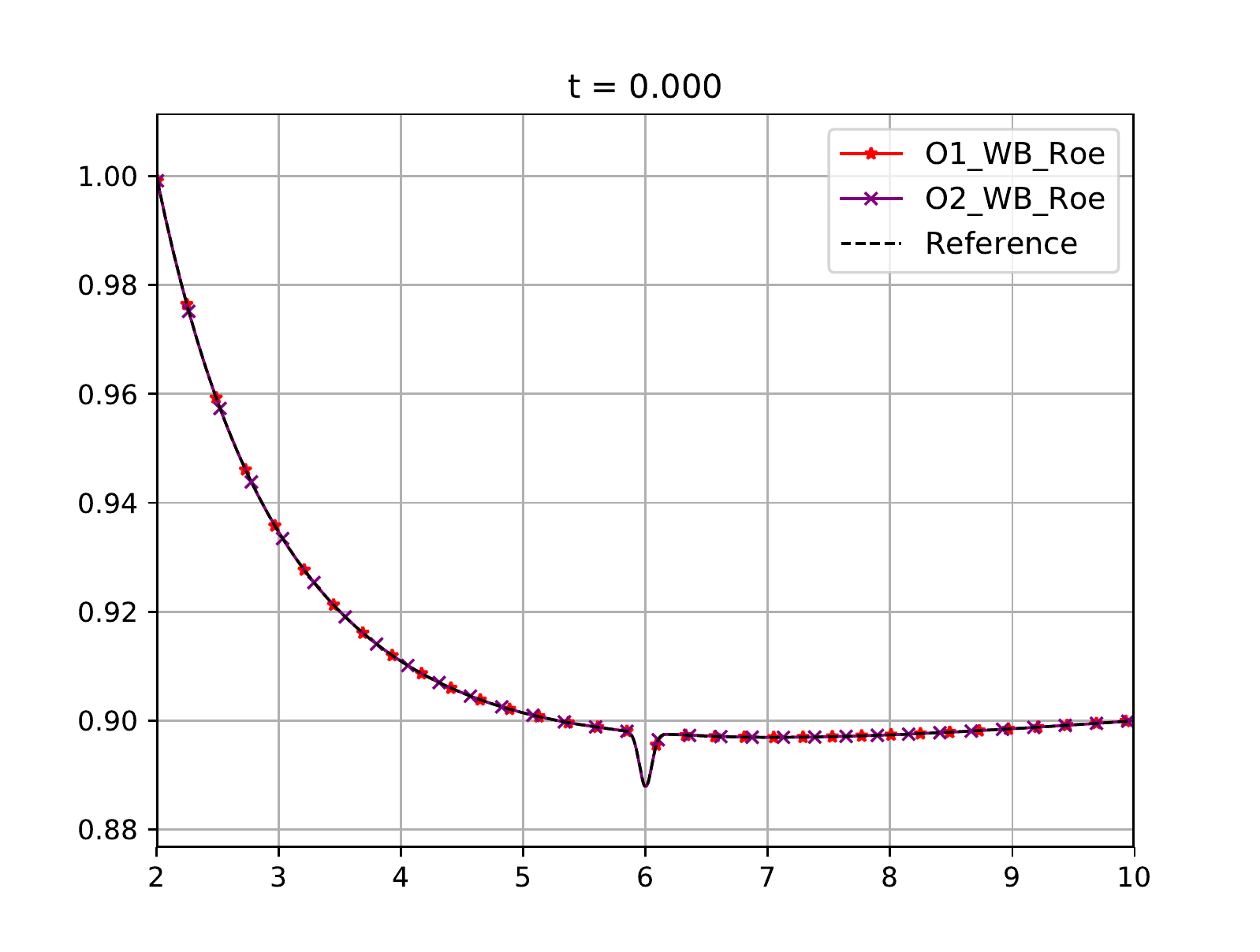}
		\label{fig:Euler_ko1_ko2_WB_testPerturbedWB1_t_0_hepse14}
	\end{subfigure}
	\begin{subfigure}[h]{0.32\textwidth}
		\centering
		\includegraphics[width=1\linewidth]{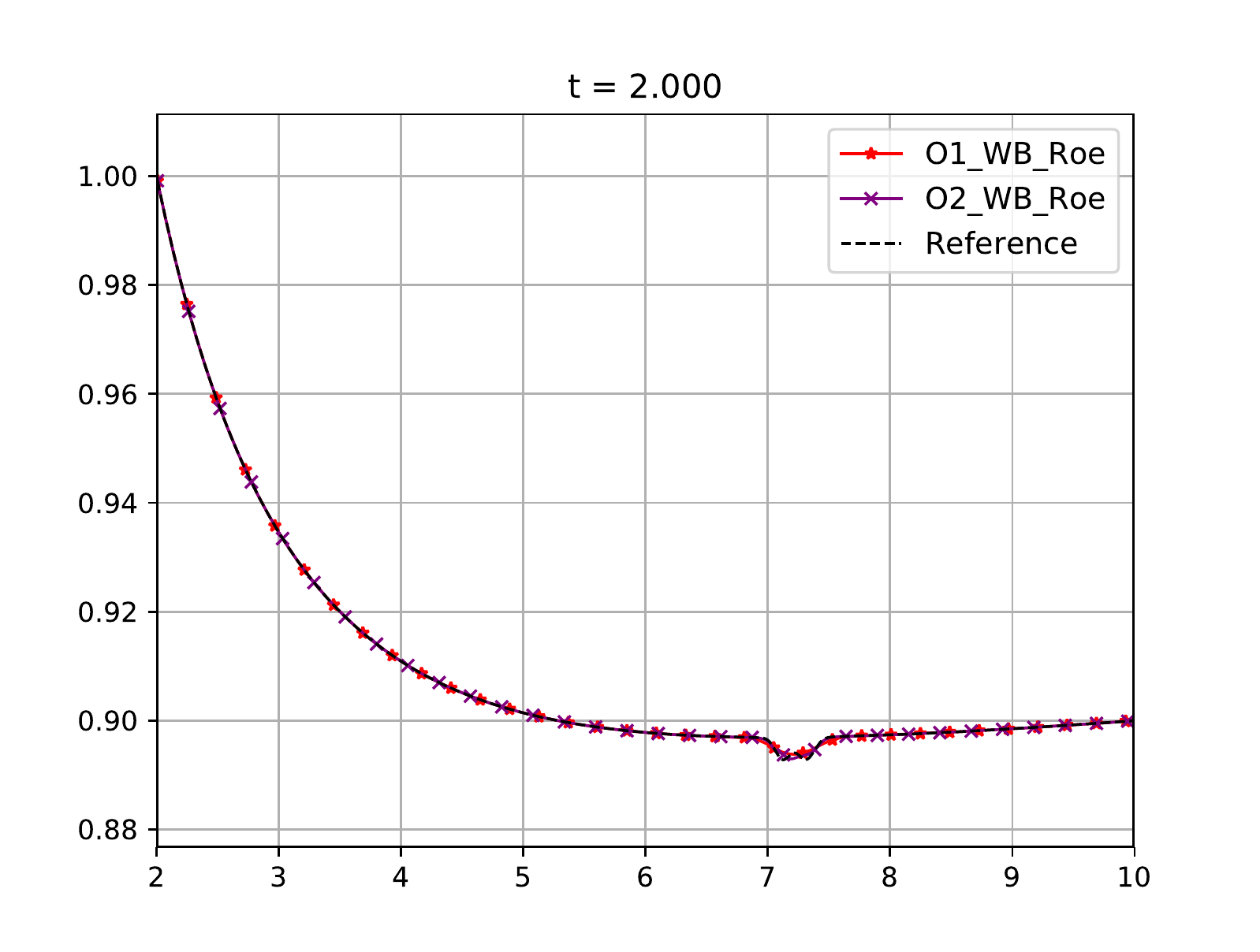}
		\label{fig:Euler_ko1_ko2_WB_testPerturbedWB1_t_2_hepse14}
	\end{subfigure}
	\begin{subfigure}[h]{0.32\textwidth}
		\centering
		\includegraphics[width=1\linewidth]{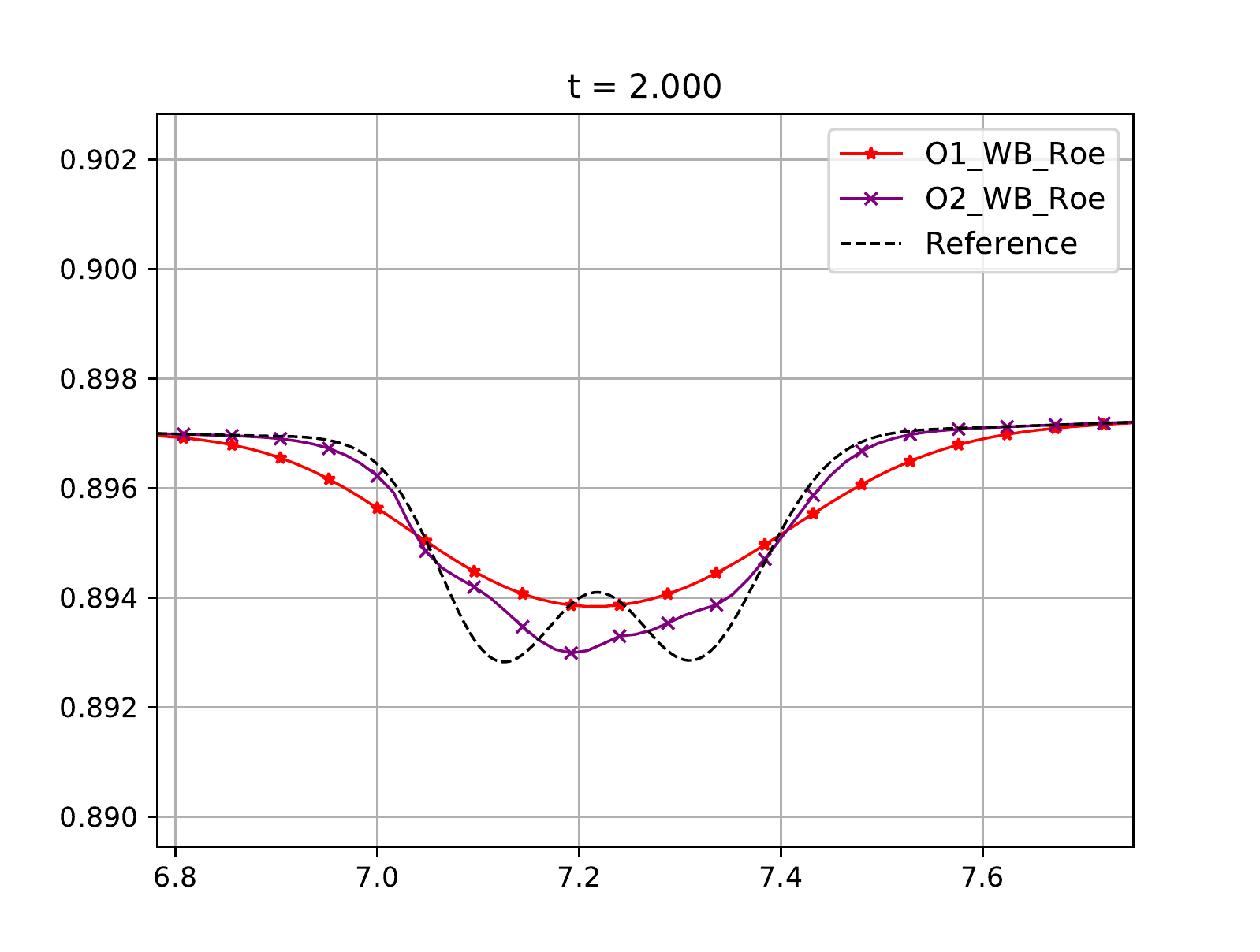}
		\caption{Zoom.}
		\label{fig:Euler_ko1_ko2_WB_testPerturbedWB1_t_2_hepse14_zoom}
	\end{subfigure}
	\begin{subfigure}[h]{0.32\textwidth}
		\centering
		\includegraphics[width=1\linewidth]{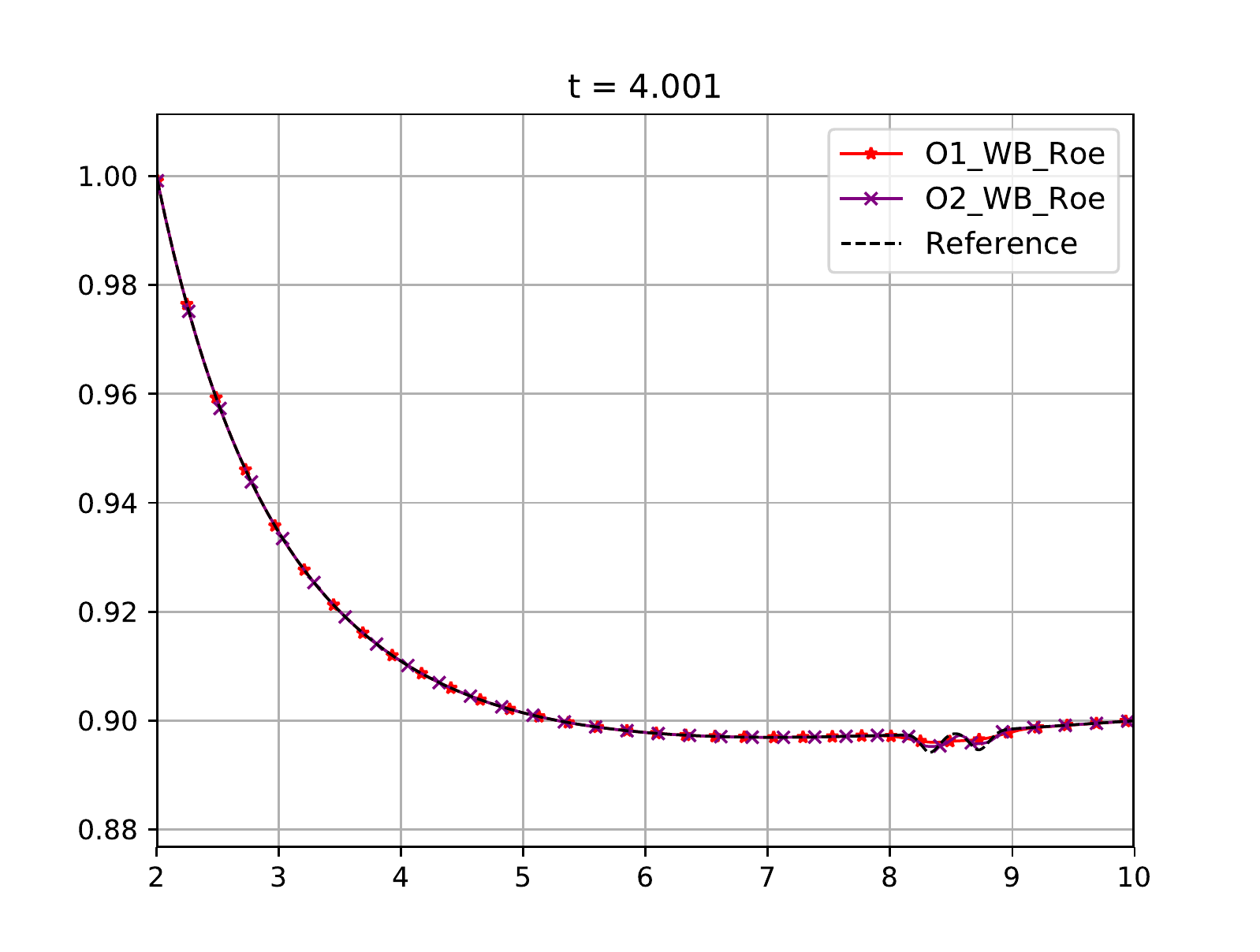}
		\label{fig:Euler_ko1_ko2_WB_testPerturbedWB1_t_4_hepse14}
	\end{subfigure}
	\begin{subfigure}[h]{0.32\textwidth}
		\centering
		\includegraphics[width=1\linewidth]{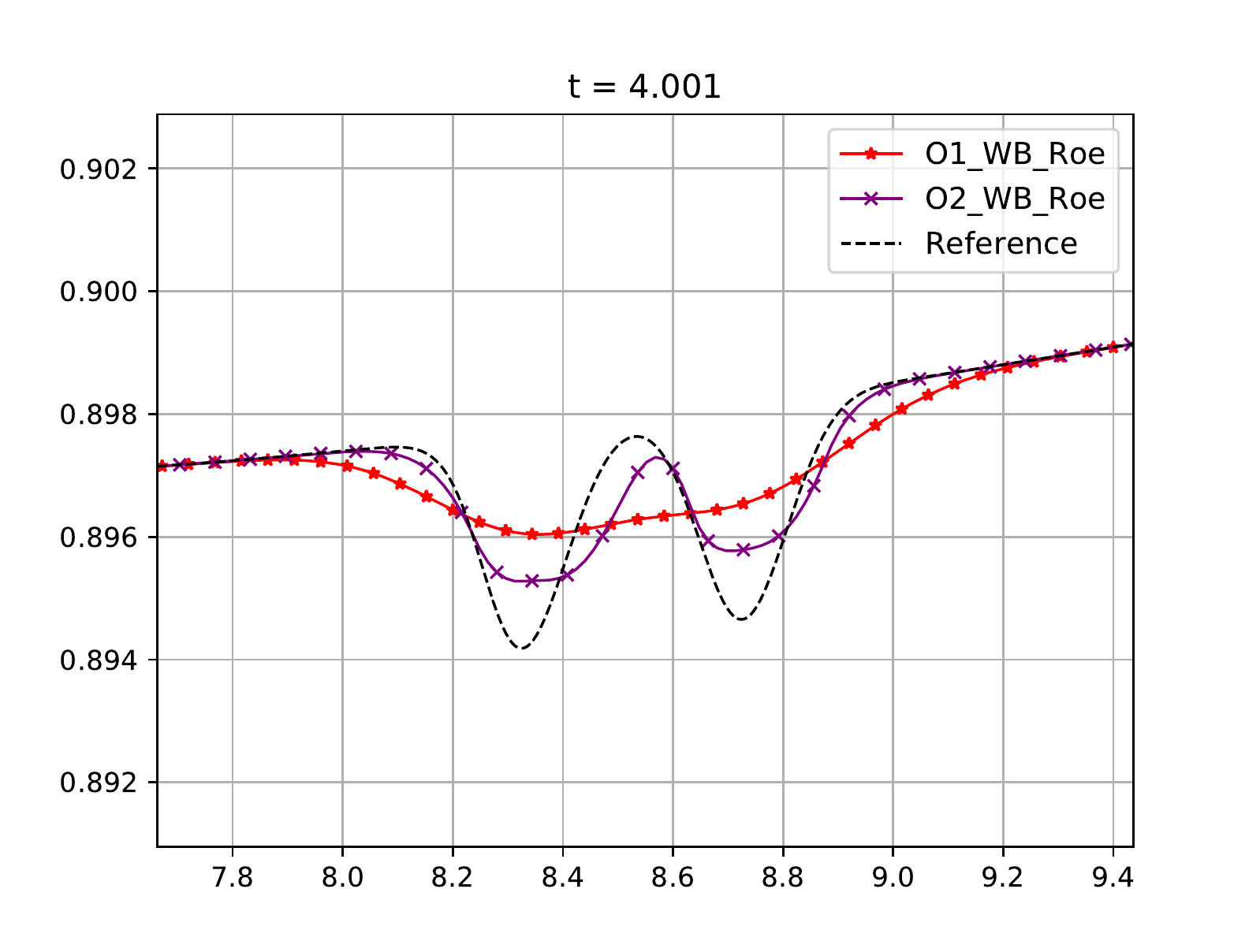}
		\caption{Zoom.}
		\label{fig:Euler_ko1_ko2_WB_testPerturbedWB1_t_4_hepse14_zoom}
	\end{subfigure}
	\begin{subfigure}[h]{0.32\textwidth}
		\centering
		\includegraphics[width=1\linewidth]{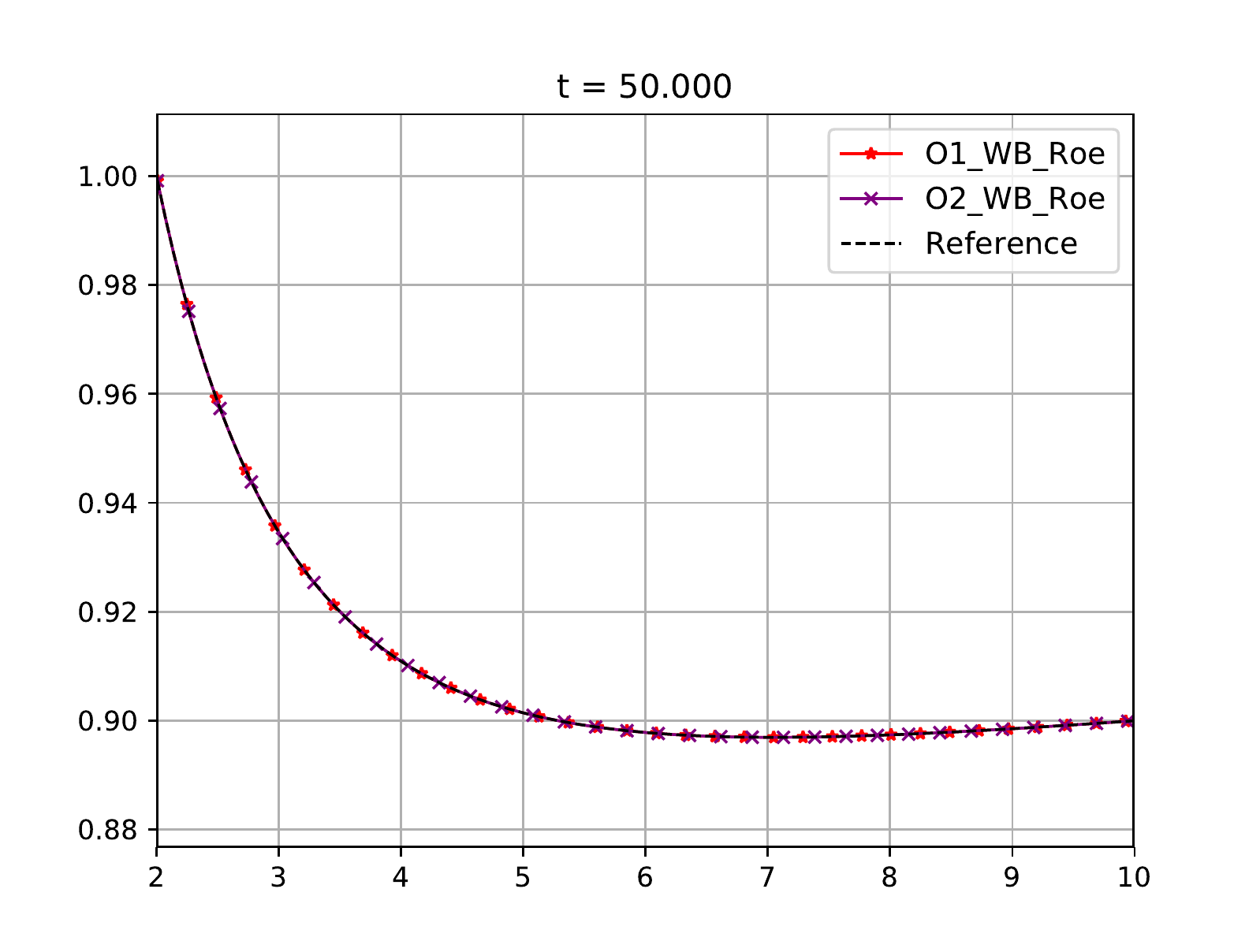}
		\label{fig:Euler_ko1_ko2_WB_testPerturbedWB1_t_50_hepse14}
	\end{subfigure}
	\caption{Euler-Schwarzschild model with the initial condition \eqref{eq:testE4a}:  first- and second-order well-balanced at selected times for the variable $v$.}
	\label{fig:Euler_ko1_ko2_WB_testPerturbedWB1_hepse14}
\end{figure}

In this test we consider the initial condition
\bel{eq:testE4a}
\widetilde V_0(r) =\widetilde V^*(r) + \delta(r),
\ee
where $\widetilde V^*$ is the supersonic stationary solution 
\bel{eq:testE4b}
\rho^{*}(10)=1, \qquad v^{*}(10)=0.9
\ee
and
\bel{eq:testE4c}
\delta(r) = [\delta_{v}(r),
\delta_{\rho}(r)]^T = \begin{cases}
[-0.01e^{-200(r-6)^{2}}, 0]^T, & \text{ $ 5<r<7$,}\\
[0, 0]^T, & \text{ otherwise.}
\end{cases}
\ee
It can be observed in Figure \ref{fig:Euler_ko1_ko2_WB_testPerturbedWB1_hepse14} that the stationary solution $V^*$ is recovered once the perturbation has left the domain. In this Figure, the numerical results obtained with the first- and second-order well-balanced methods are compared with a reference solution computed using the first-order well-balanced method with  a 5000-point mesh (the Roe-type numerical flux is used again). As expected, the second-order method is less diffusive.


\subsection{Perturbation of a steady shock solution}

\paragraph{{{\red Left-hand perturbation.}}}

\begin{figure}[h]
	\begin{subfigure}[h]{0.32\textwidth}
		\centering
		\includegraphics[width=1\linewidth]{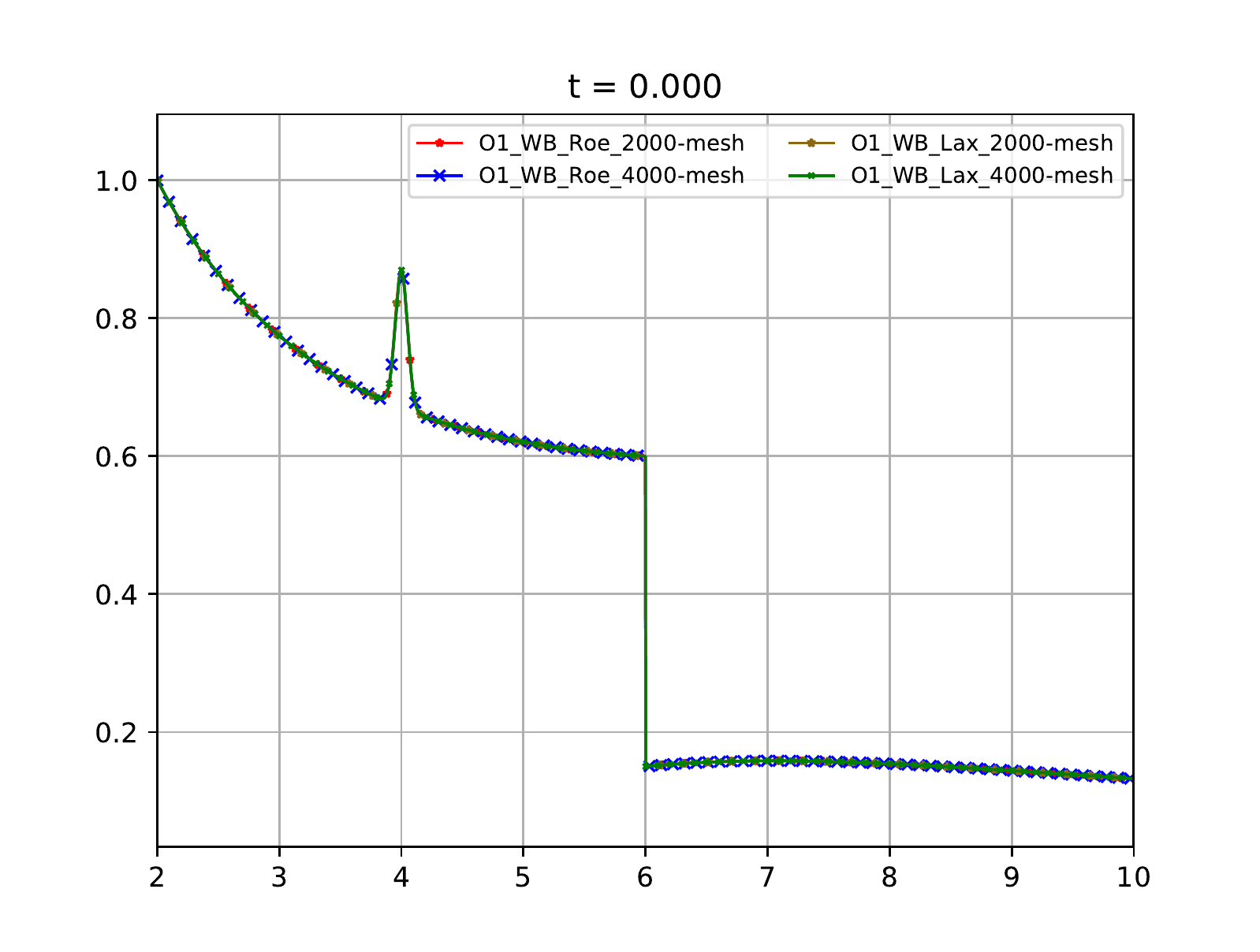}
		\label{fig:Euler_ko1_WB_testPerturbedWB3_HLL_vs_Lax_t_0}
	\end{subfigure}
	\begin{subfigure}[h]{0.32\textwidth}
		\centering
		\includegraphics[width=1\linewidth]{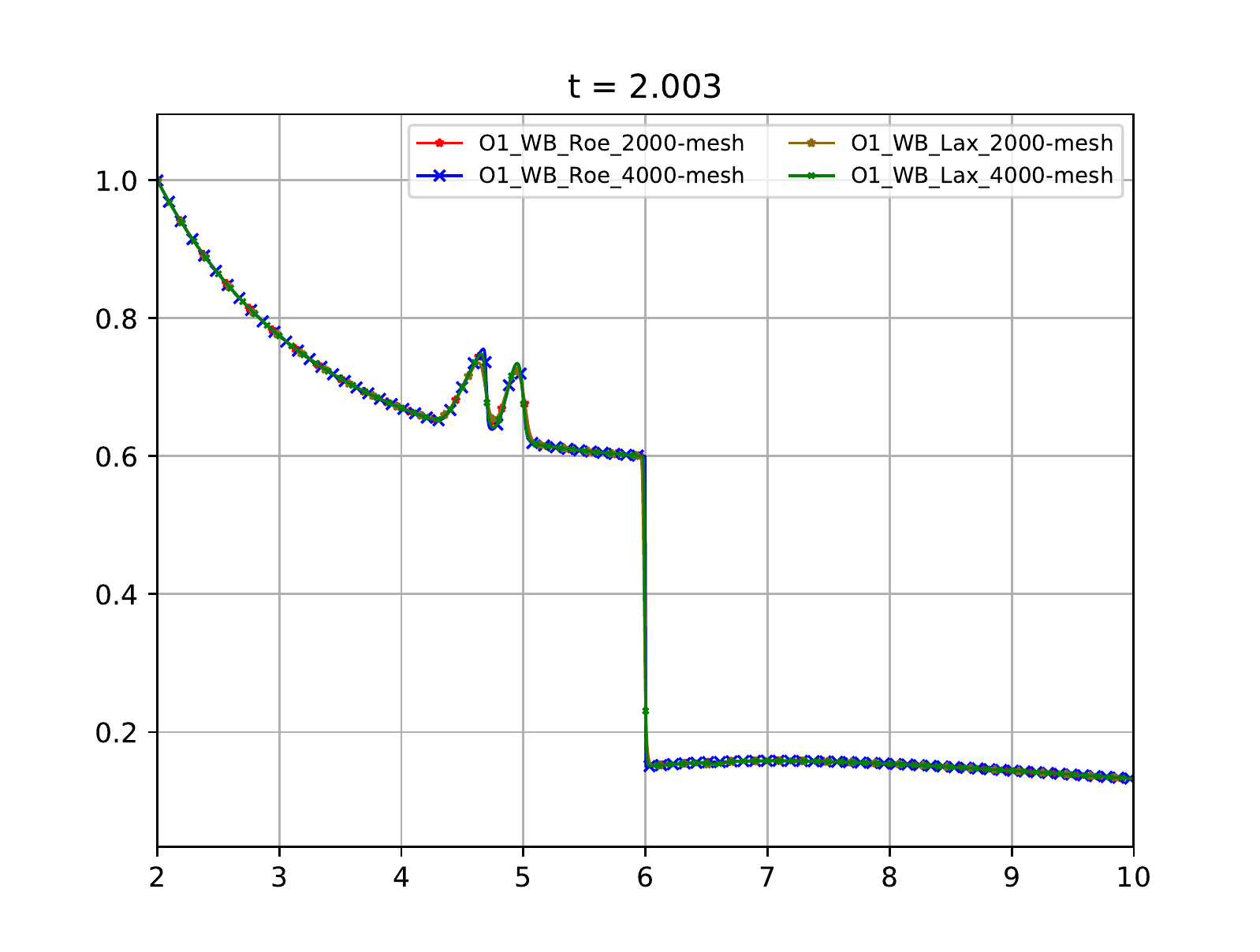}
		\label{fig:Euler_ko1_WB_testPerturbedWB3_HLL_vs_Lax_t_2}
	\end{subfigure}
	\begin{subfigure}[h]{0.32\textwidth}
		\centering
		\includegraphics[width=1\linewidth]{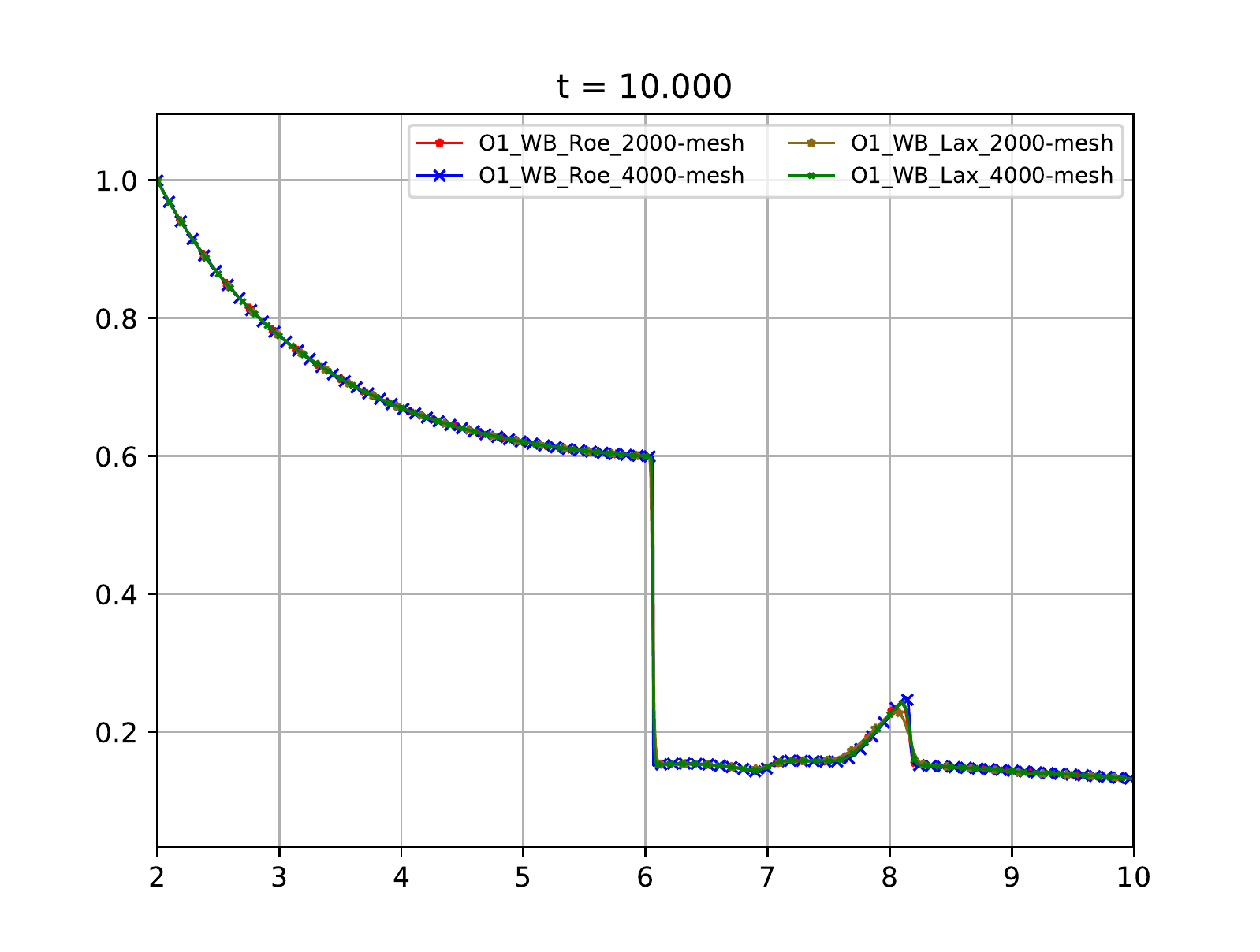}
		\label{fig:Euler_ko1_WB_testPerturbedWB3_HLL_vs_Lax_t_10}
	\end{subfigure}
	\begin{subfigure}[h]{0.32\textwidth}
		\centering
		\includegraphics[width=1\linewidth]{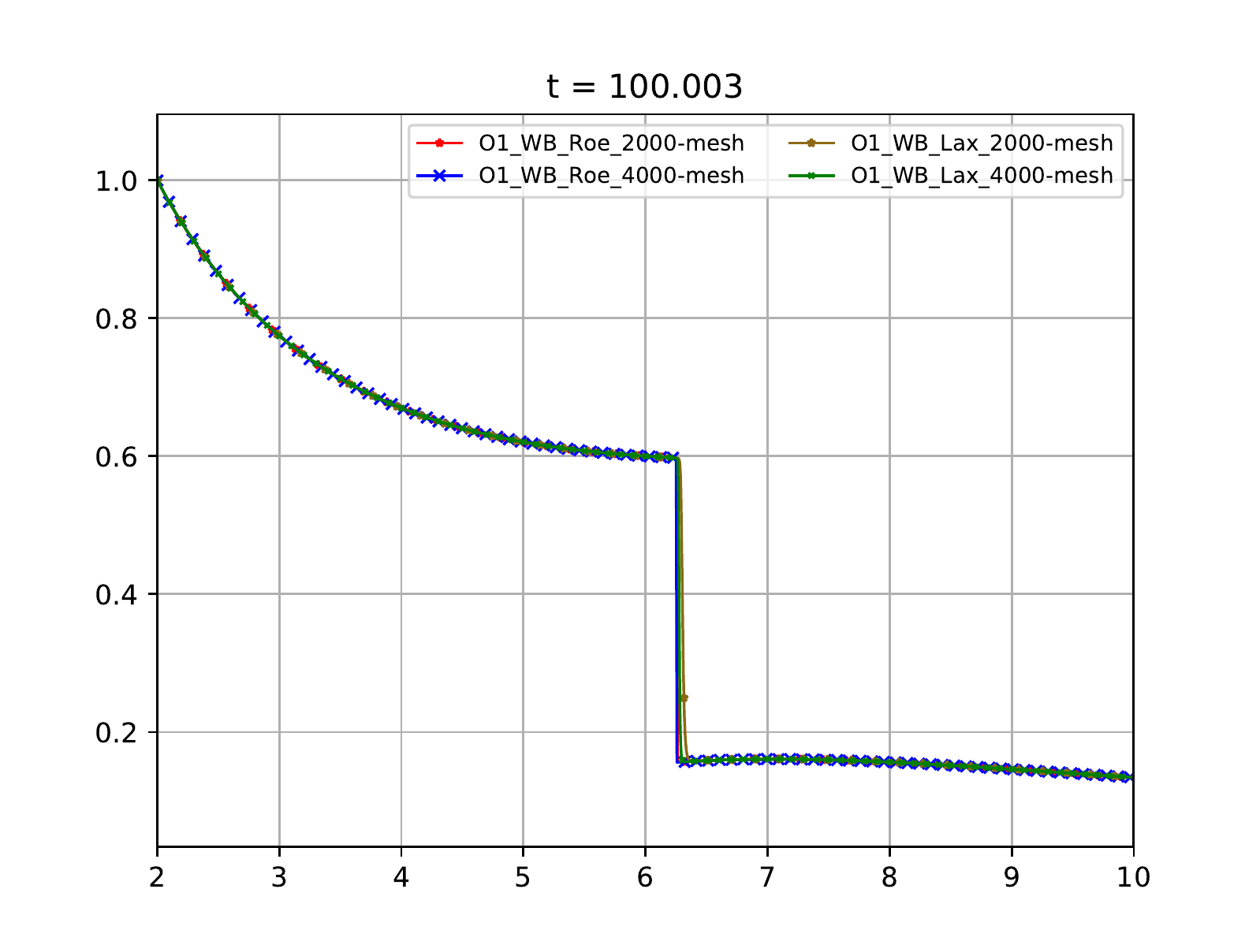}
		\label{fig:Euler_ko1_WB_testPerturbedWB3_HLL_vs_Lax_t_100}
	\end{subfigure}
	\begin{subfigure}[h]{0.32\textwidth}
		\centering
		\includegraphics[width=1\linewidth]{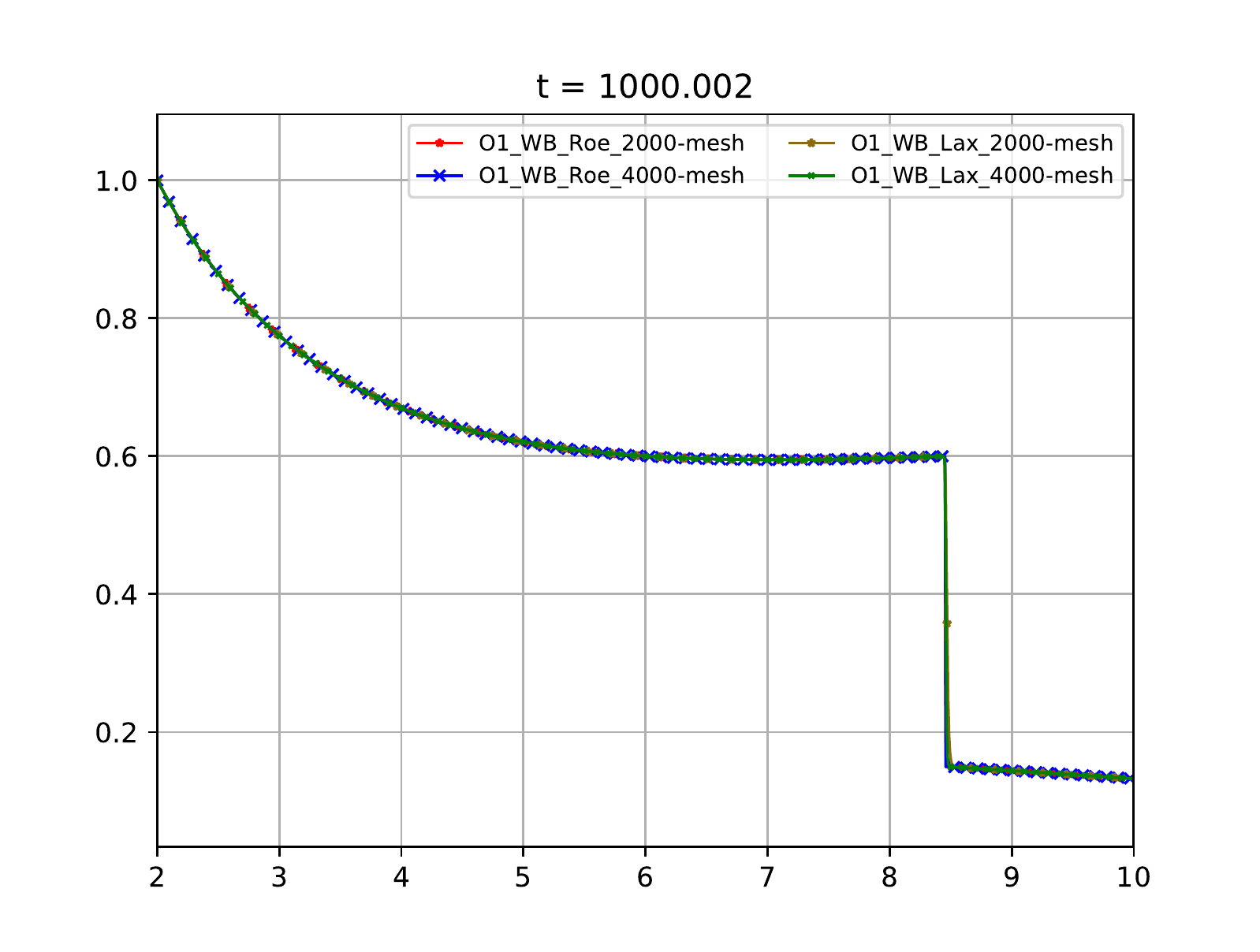}
		\label{fig:Euler_ko1_WB_testPerturbedWB3_HLL_vs_Lax_t_1000}
	\end{subfigure}
	\begin{subfigure}[h]{0.32\textwidth}
		\centering
		\includegraphics[width=1\linewidth]{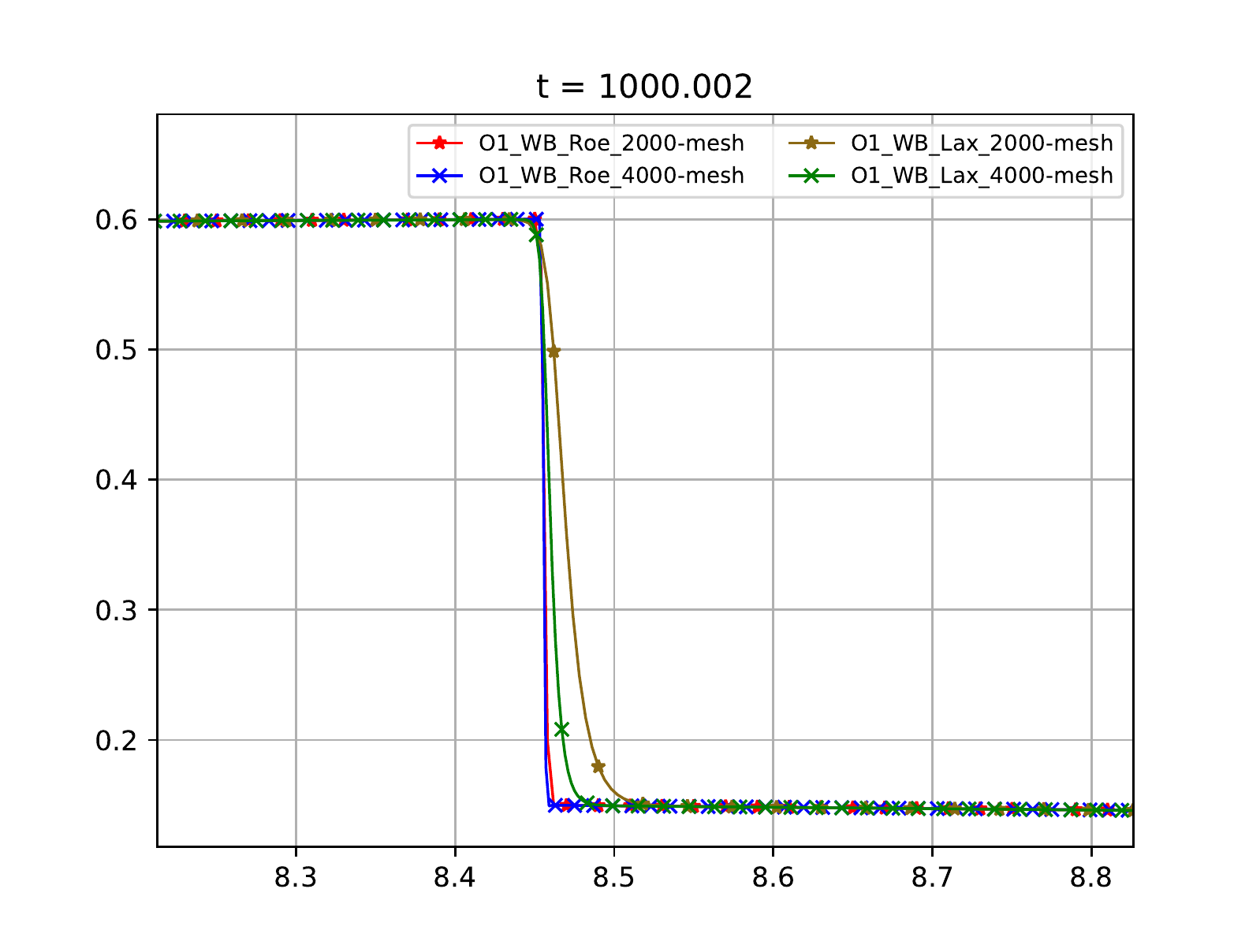}
		\caption{Zoom.}
		\label{fig:Euler_ko1_WB_testPerturbedWB3_HLL_vs_Lax_t_1000_zoom}
	\end{subfigure}
	\caption{Euler-Schwarzschild model with the initial condition \eqref{eq:testE5a}: comparison between the first-order well-balanced method with different meshes using the Roe-type and the Lax numerical fluxes at selected times for the variable $v$.}
	\label{fig:Euler_ko1_WB_testPerturbedWB3_HLL_vs_Lax}
\end{figure}

\begin{figure}[h]
	\begin{subfigure}[h]{0.32\textwidth}
		\centering
		\includegraphics[width=1\linewidth]{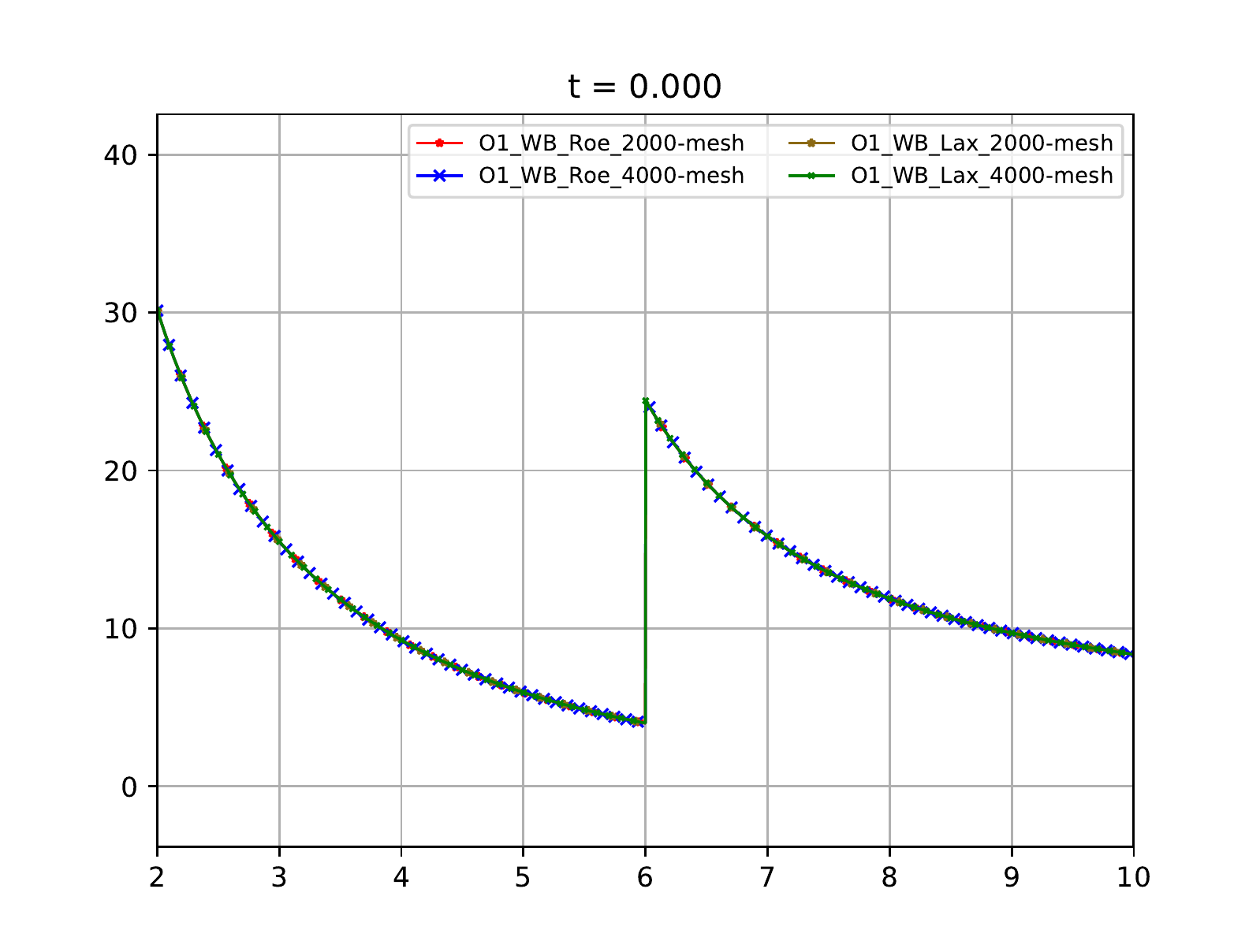}
		\label{fig:Euler_ko1_WB_testPerturbedWB3_HLL_vs_Lax_t_0_rho}
	\end{subfigure}
	\begin{subfigure}[h]{0.32\textwidth}
		\centering
		\includegraphics[width=1\linewidth]{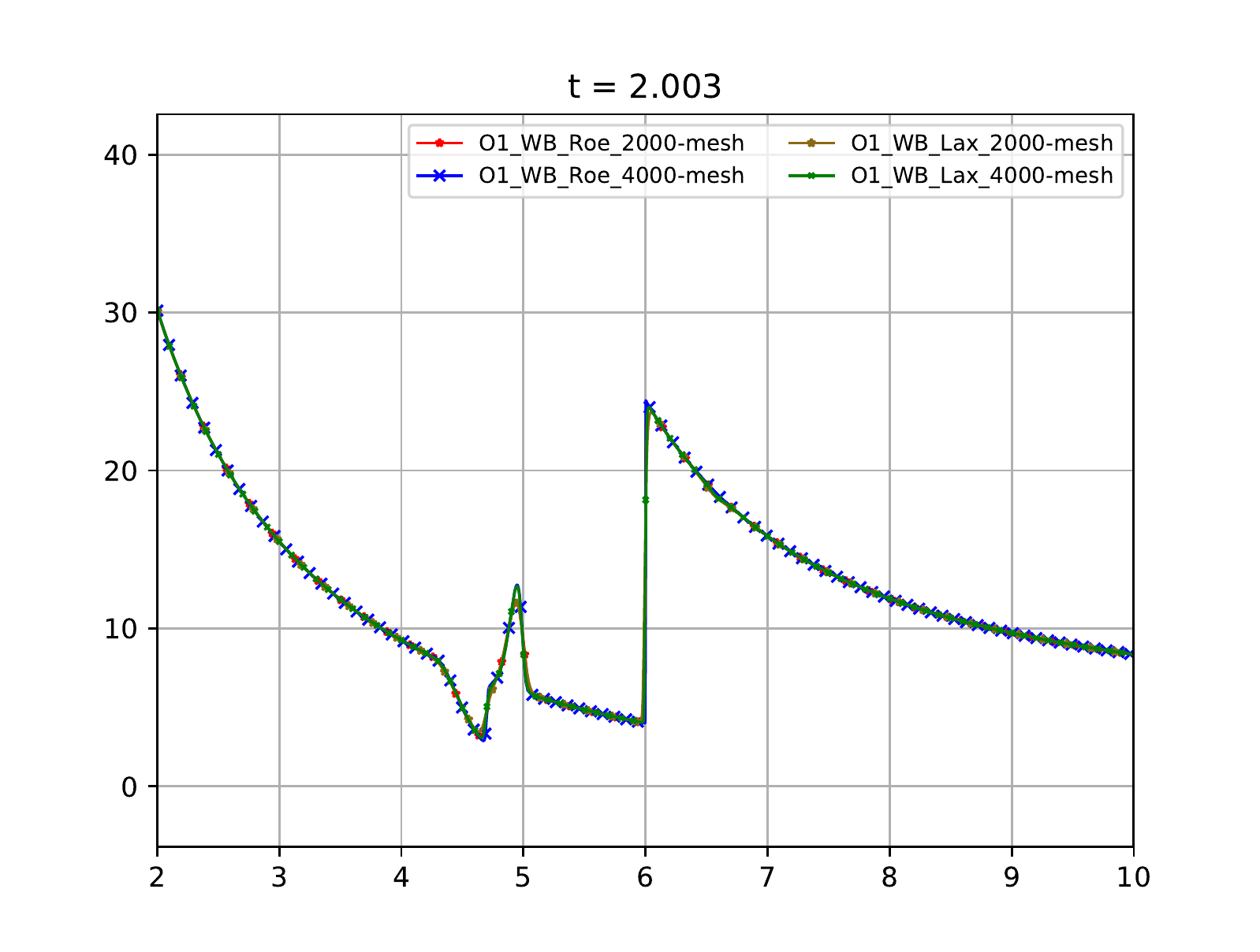}
		\label{fig:Euler_ko1_WB_testPerturbedWB3_HLL_vs_Lax_t_2_rho}
	\end{subfigure}
	\begin{subfigure}[h]{0.32\textwidth}
		\centering
		\includegraphics[width=1\linewidth]{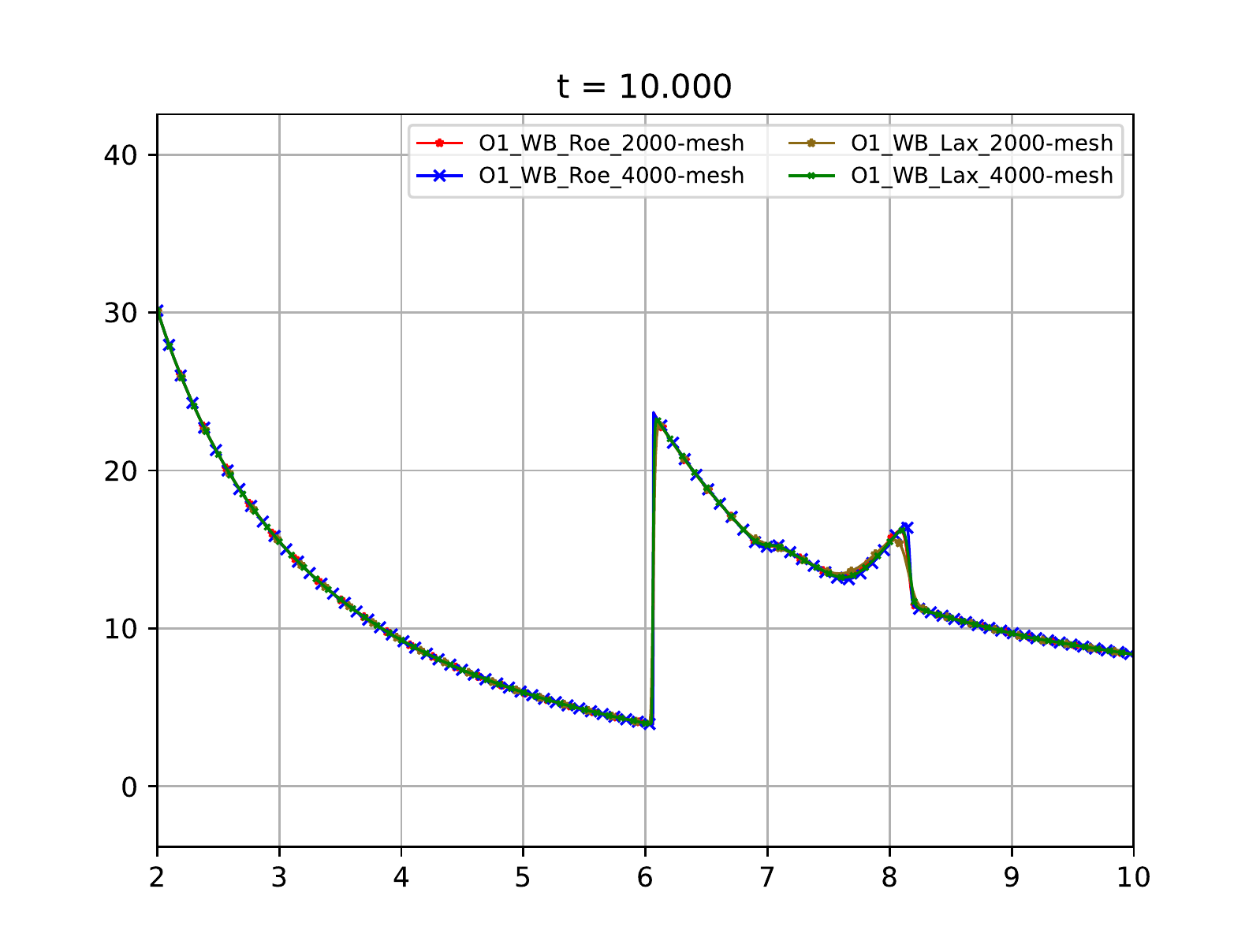}
		\label{fig:Euler_ko1_WB_testPerturbedWB3_HLL_vs_Lax_t_10_rho}
	\end{subfigure}
	\begin{subfigure}[h]{0.32\textwidth}
		\centering
		\includegraphics[width=1\linewidth]{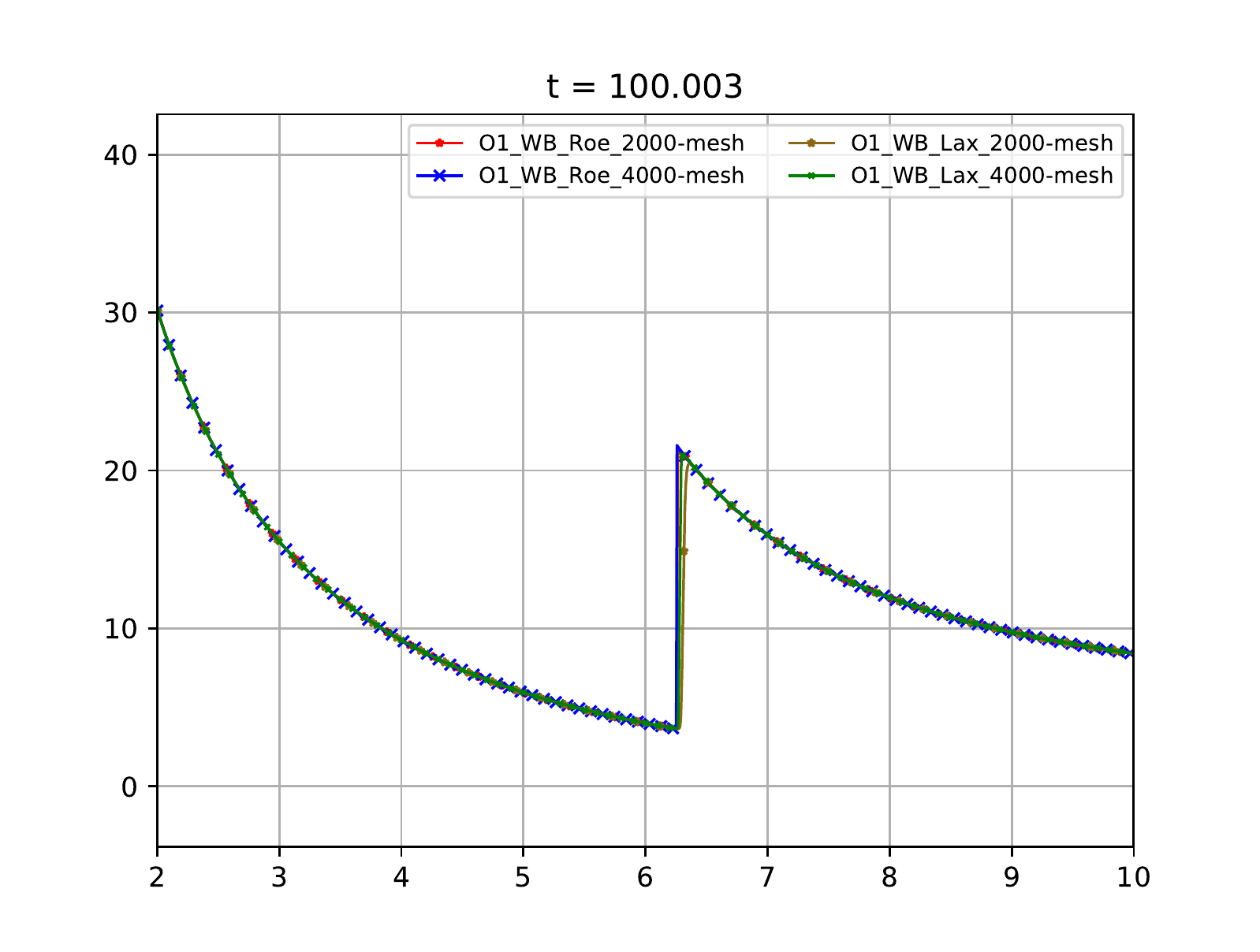}
		\label{fig:Euler_ko1_WB_testPerturbedWB3_HLL_vs_Lax_t_100_rho}
	\end{subfigure}
	\begin{subfigure}[h]{0.32\textwidth}
		\centering
		\includegraphics[width=1\linewidth]{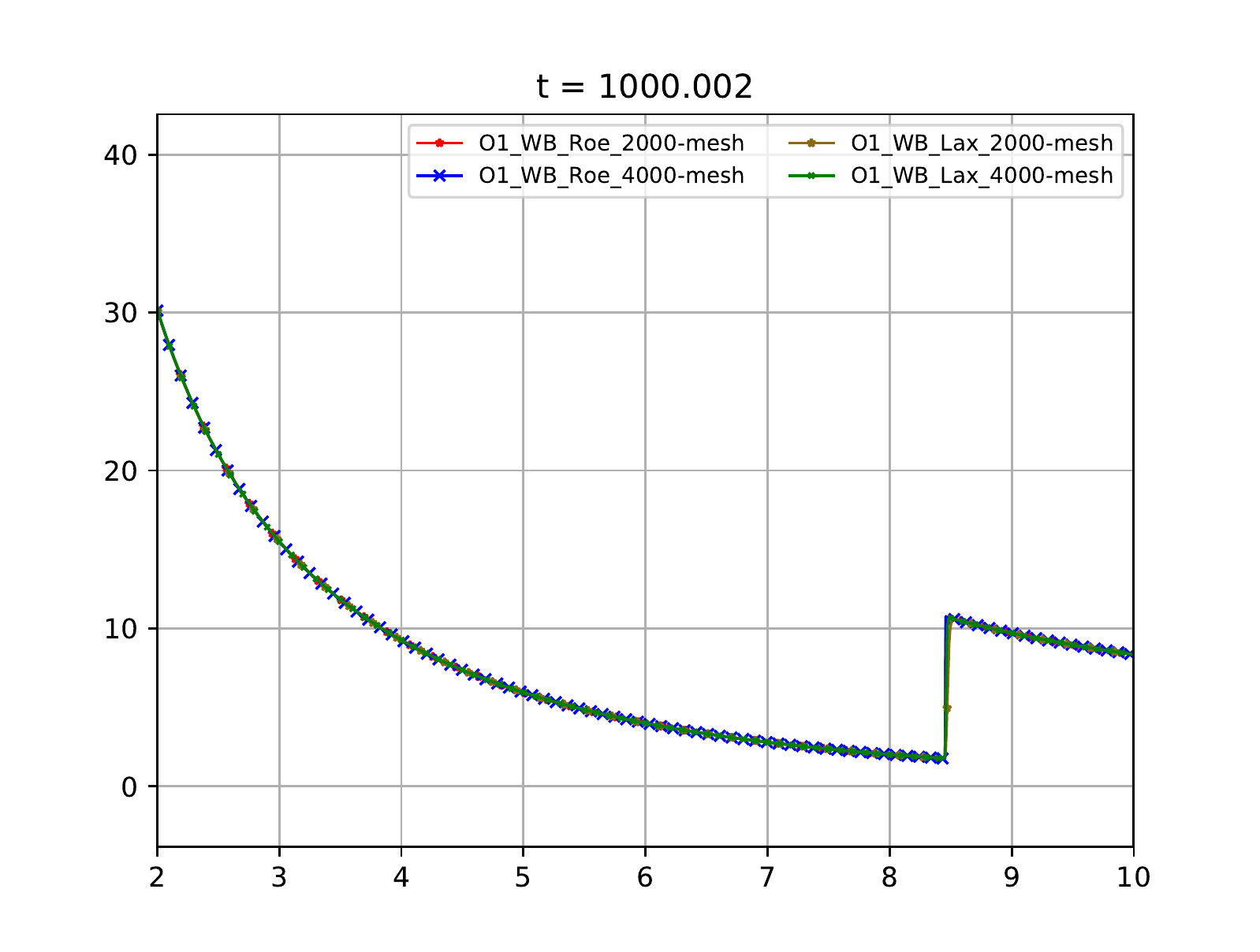}
		\label{fig:Euler_ko1_WB_testPerturbedWB3_HLL_vs_Lax_t_1000_rho}
	\end{subfigure}
	\begin{subfigure}[h]{0.32\textwidth}
		\centering
		\includegraphics[width=1\linewidth]{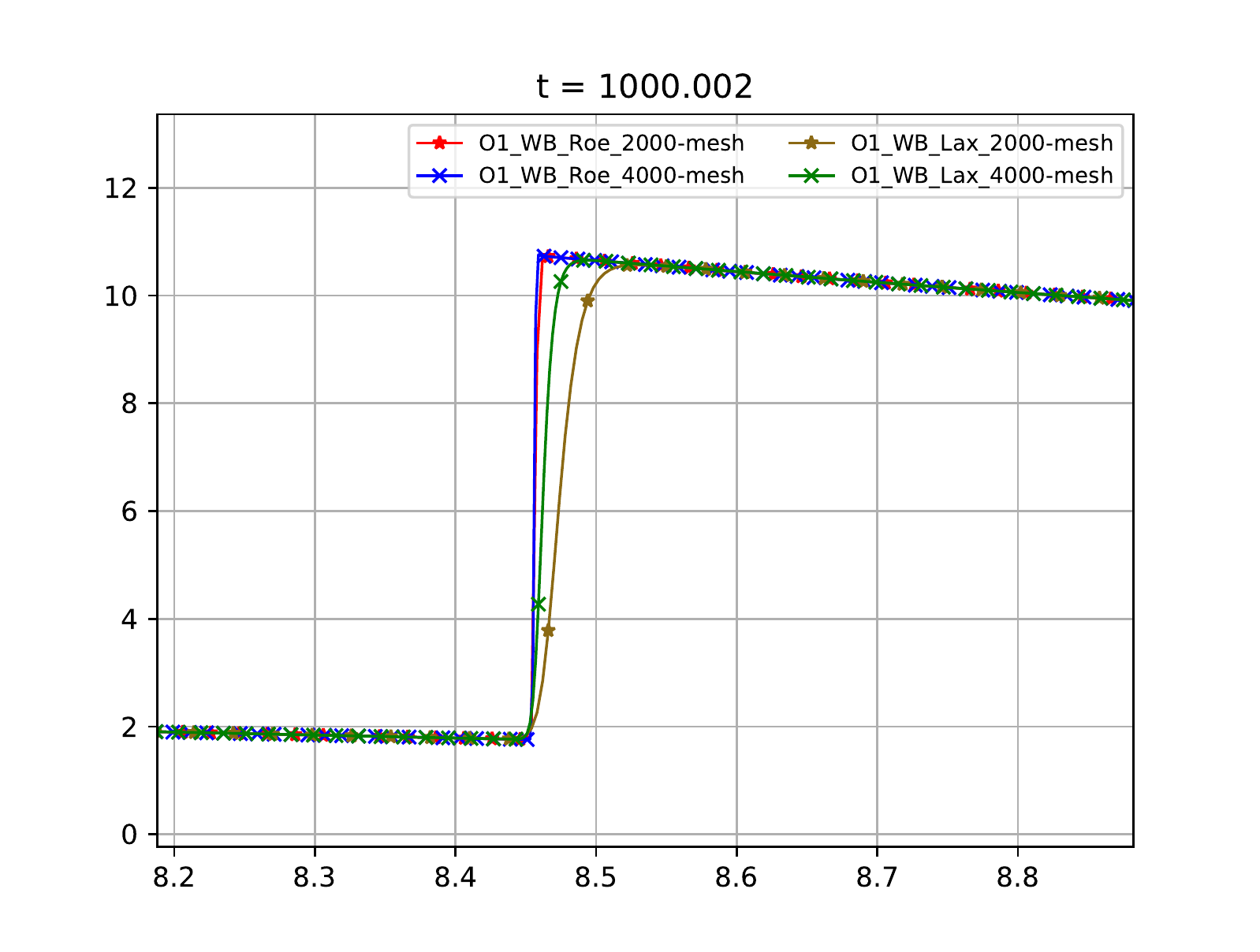}
		\caption{Zoom.}
		\label{fig:Euler_ko1_WB_testPerturbedWB3_HLL_vs_Lax_t_1000_zoom_rho}
	\end{subfigure}
	\caption{Euler-Schwarzschild model with the initial condition \eqref{eq:testE5a}: first-order well-balanced method with different meshes using the Roe-type and the Lax numerical fluxes at selected times for the variable $\rho$.}
	\label{fig:Euler_ko1_WB_testPerturbedWB3_HLL_vs_Lax_rho}
\end{figure}

\begin{figure}[h]
	\centering
	\includegraphics[width=.6\linewidth]{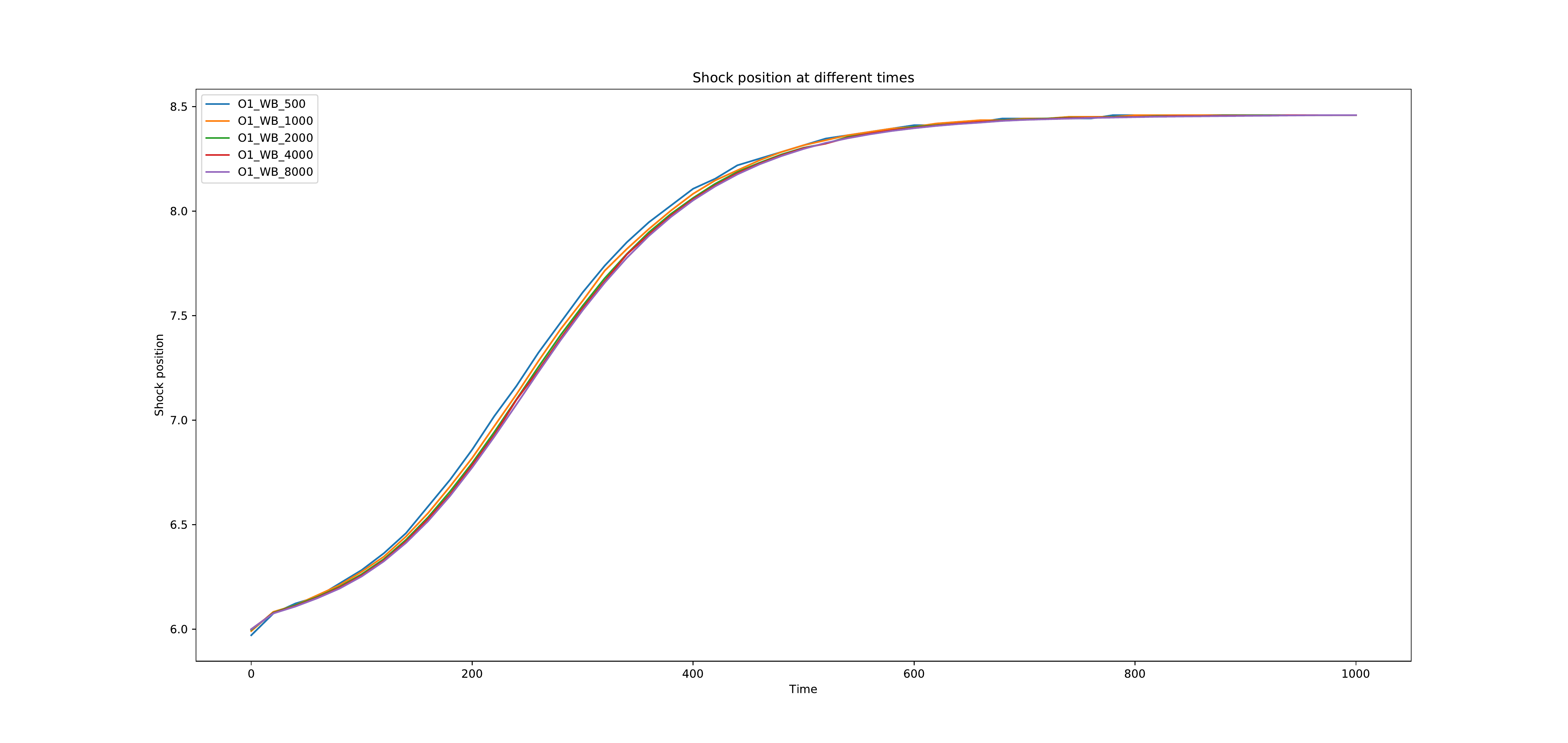}
	\caption{Euler-Schwarzschild model taking as initial condition \eqref{eq:testE5a}:  evolution of the shock position with time obtained with the first-order well-balanced method using the Roe-type numerical flux with different meshes.}
	\label{fig:eulerko1wbshockpositiondifferentmeshes}
\end{figure}

We consider the initial condition 
\bel{eq:testE5a}
\widetilde V_{0}(r) = 
\widetilde{V}^*(r) + \delta_L(r),
\ee
where
\bel{eq:testE5b}
\delta_L(r) = [\delta_{v,L}(r),
\delta_{\rho,L}(r)]^T = \begin{cases}
	[0.2e^{-200(r-4)^{2}}, 0]^T, & \text{ $ 3<r<5$,}\\
	[0, 0]^T, & \text{ otherwise,}
\end{cases}
\ee
and $V^*(r)$ is the stationary solution given by \eqref{eq:testE3a}-\eqref{eq:testE3c}.
In Figures \ref{fig:Euler_ko1_WB_testPerturbedWB3_HLL_vs_Lax} and \ref{fig:Euler_ko1_WB_testPerturbedWB3_HLL_vs_Lax_rho} the numerical results obtained with the first- and second-order well-balanced methods using the Lax-Friedrichs and the Roe-type numerical methods with different meshes are compared. As it happened for the Burgers-Schwarzschild model, the location of the stationary shock changes after the passage of the wave generated by the perturbation. Nevertheless in this case the displacement of the shock is slower. Different numerical methods have been applied to check the dependency of the motion on the scheme: although the evolution of the shock slightly depends on the number of points of the mesh,  all the numerical solutions capture the same final location of the shock. In Figure \ref{fig:eulerko1wbshockpositiondifferentmeshes} the evolution of the shock given by the first-order WB method with different number of cells are compared. The location of the shock at every time step has been detected by using the condition $\displaystyle \frac{v_i-v_{i-1}}{v_{i+1}-v_{i}} \geq 0.8$.


\paragraph{\red Right-hand perturbation.}

Let us consider now   two different  initial conditions: on the one hand
\bel{eq:testE6a}
\widetilde V_{0,1}(r) = 
\widetilde{V}^*(r) + \delta_R(r),
\ee
where $\widetilde{V}^*(r) $ is again the discontinuous stationary solution given by \eqref{eq:testE3a}-\eqref{eq:testE3c} 
and
\bel{eq:testE6b}
\delta_R(r) = [\delta_{v,R}(r),
\delta_{\rho,R}(r)]^T = \begin{cases}
	[-0.05e^{-200(r-8)^{2}}, 0]^T, & \text{ $ 7<r<9$,}\\
	[0, 0]^T, & \text{ otherwise.}
\end{cases}
\ee 
On the other hand, 
\bel{eq:testE6aa}
\widetilde V_{0,2}(r) = 
\widetilde{V}_2^*(r) + \delta_R(r),
\ee
where $\delta_R$ is given again by \eqref{eq:testE6a} and $\widetilde{V}_2^*(r)$ is the steady shock of the form (\ref{eq:testE3a}) 
satisfying 
\bel{eq:testE3bbis}
\rho_{-}^{*}(6)=5, \quad v_{-}^{*}(6)=0.6.
\ee
Observe that the definition of $v$ is identical for both stationary solutions but $\rho$ is different.

\begin{figure}[h]
	\begin{subfigure}[h]{0.32\textwidth}
		\centering
		\includegraphics[width=1\linewidth]{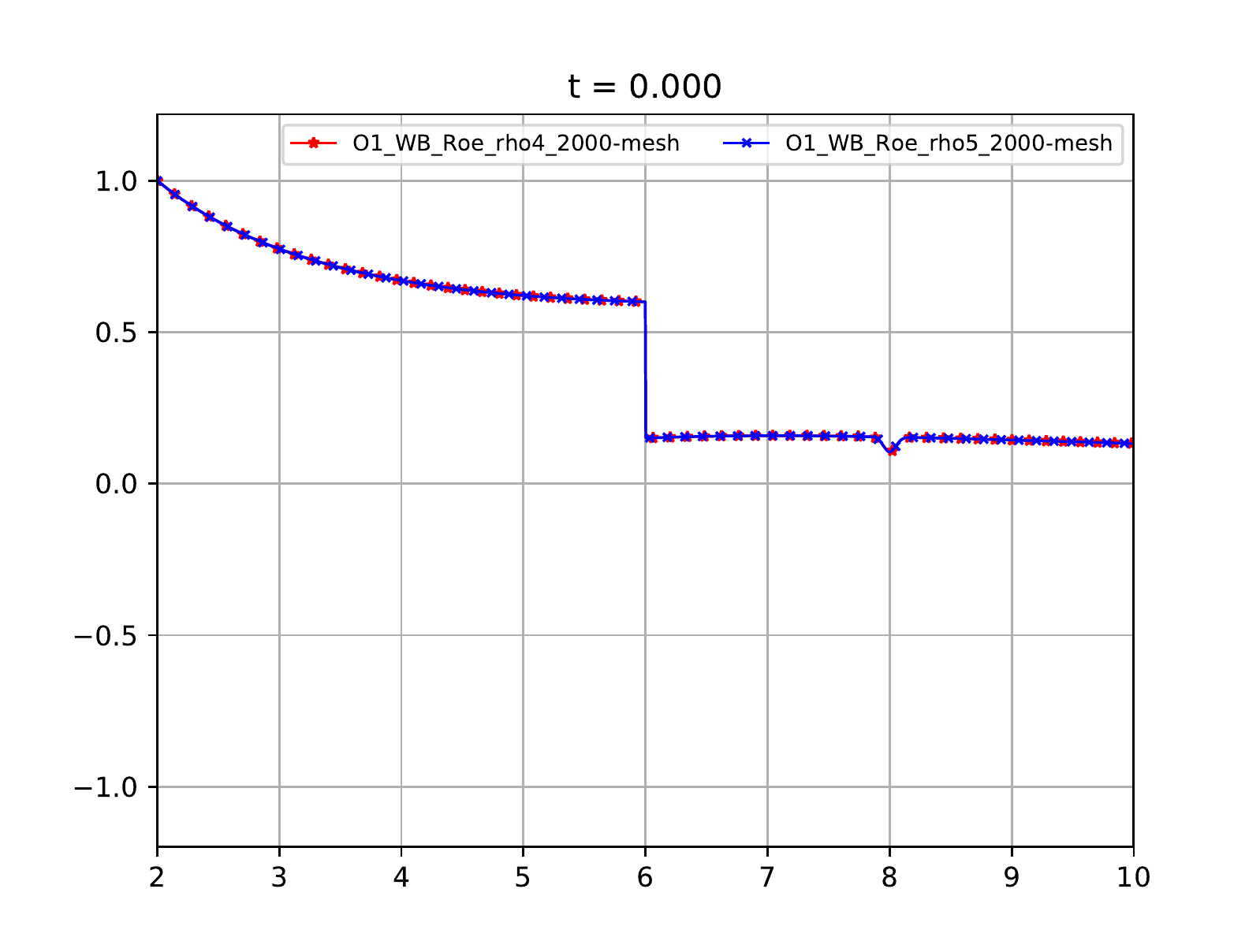}
		\label{fig:Euler_ko1_WB_testPerturbedWB3_negative_perturbation_t_0}
	\end{subfigure}
	\begin{subfigure}[h]{0.32\textwidth}
		\centering
		\includegraphics[width=1\linewidth]{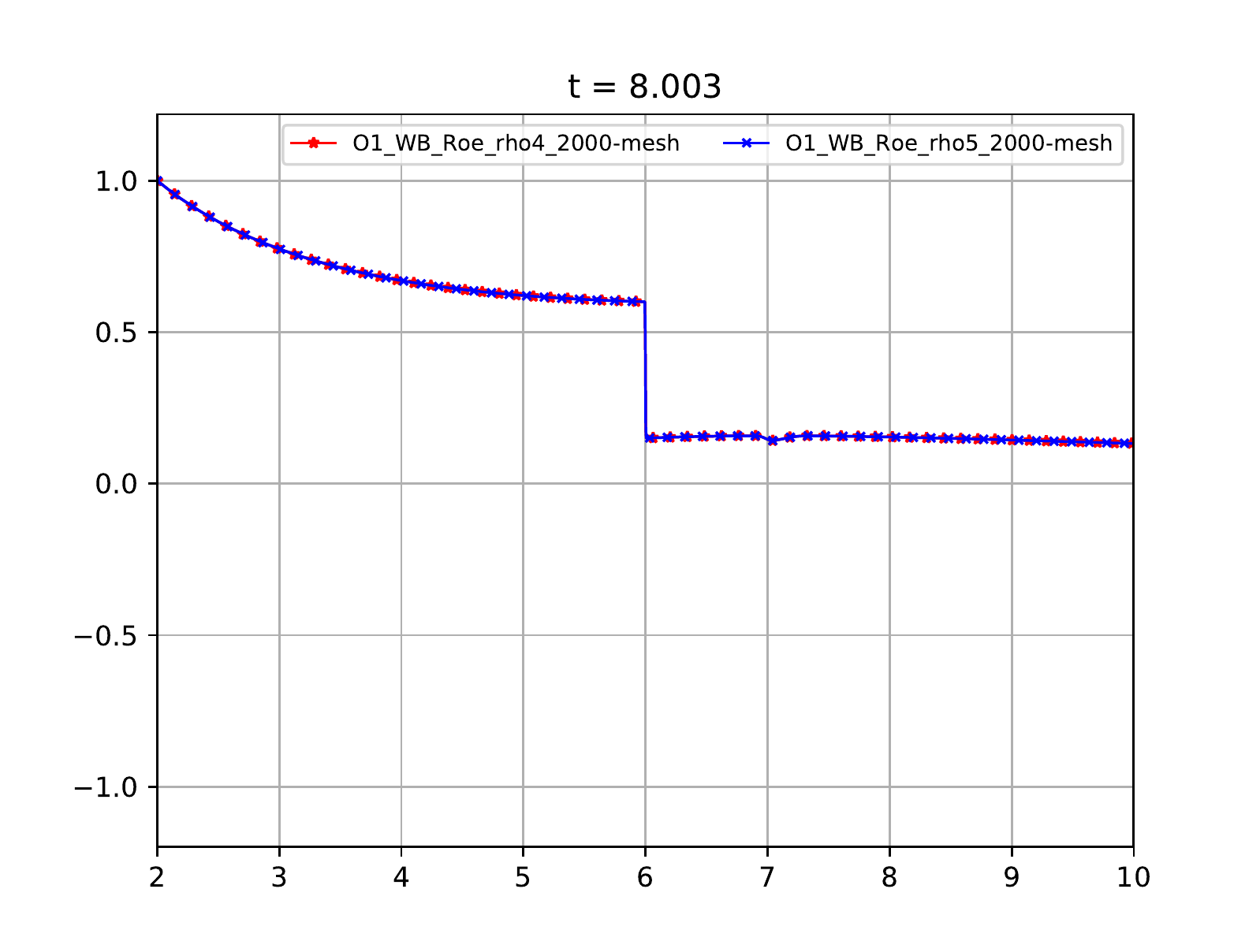}
		\label{fig:Euler_ko1_WB_testPerturbedWB3_negative_perturbation_t_8}
	\end{subfigure}
	\begin{subfigure}[h]{0.32\textwidth}
		\centering
		\includegraphics[width=1\linewidth]{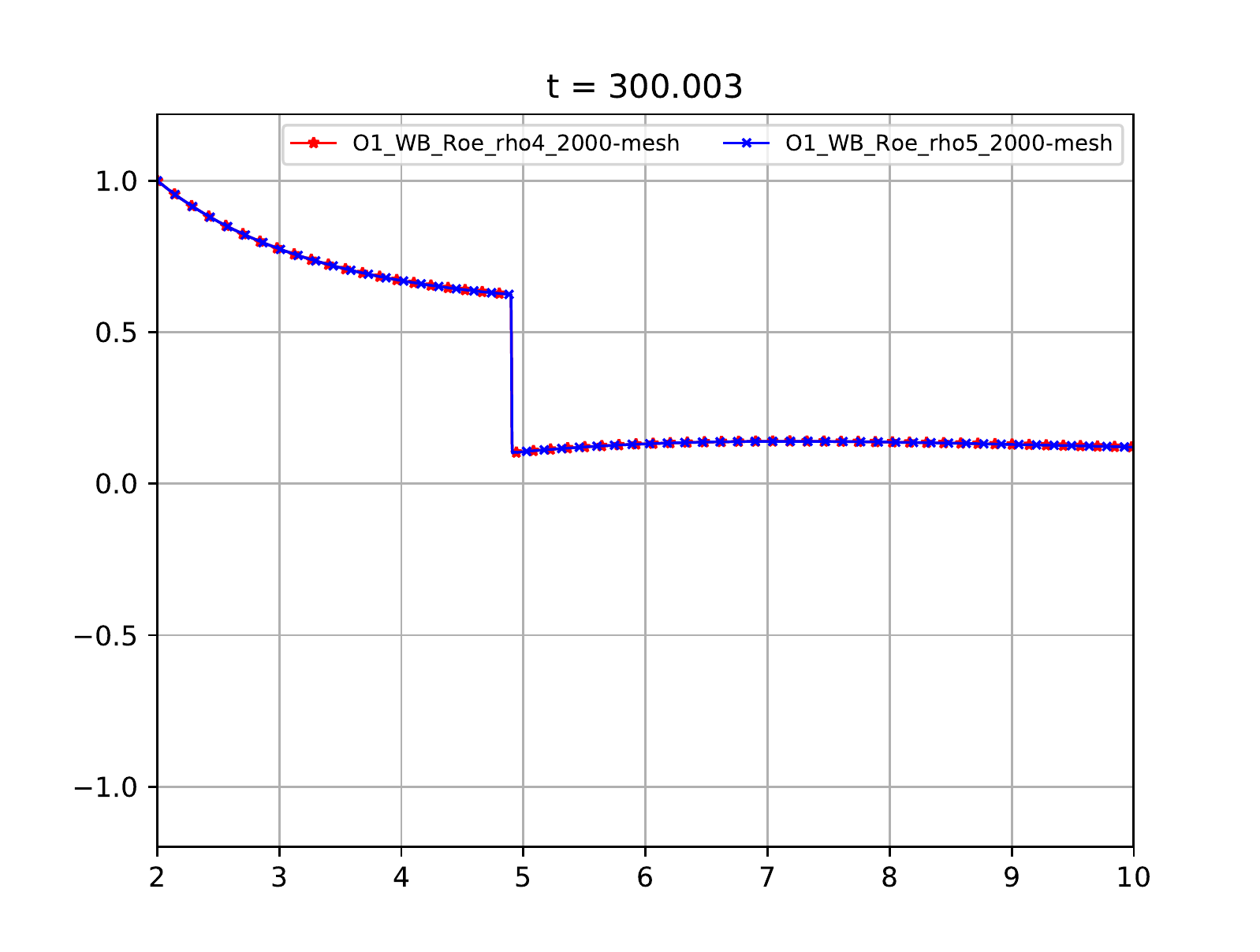}
		\label{fig:Euler_ko1_WB_testPerturbedWB3_negative_perturbation_t_300}
	\end{subfigure}
	\begin{subfigure}[h]{0.32\textwidth}
		\centering
		\includegraphics[width=1\linewidth]{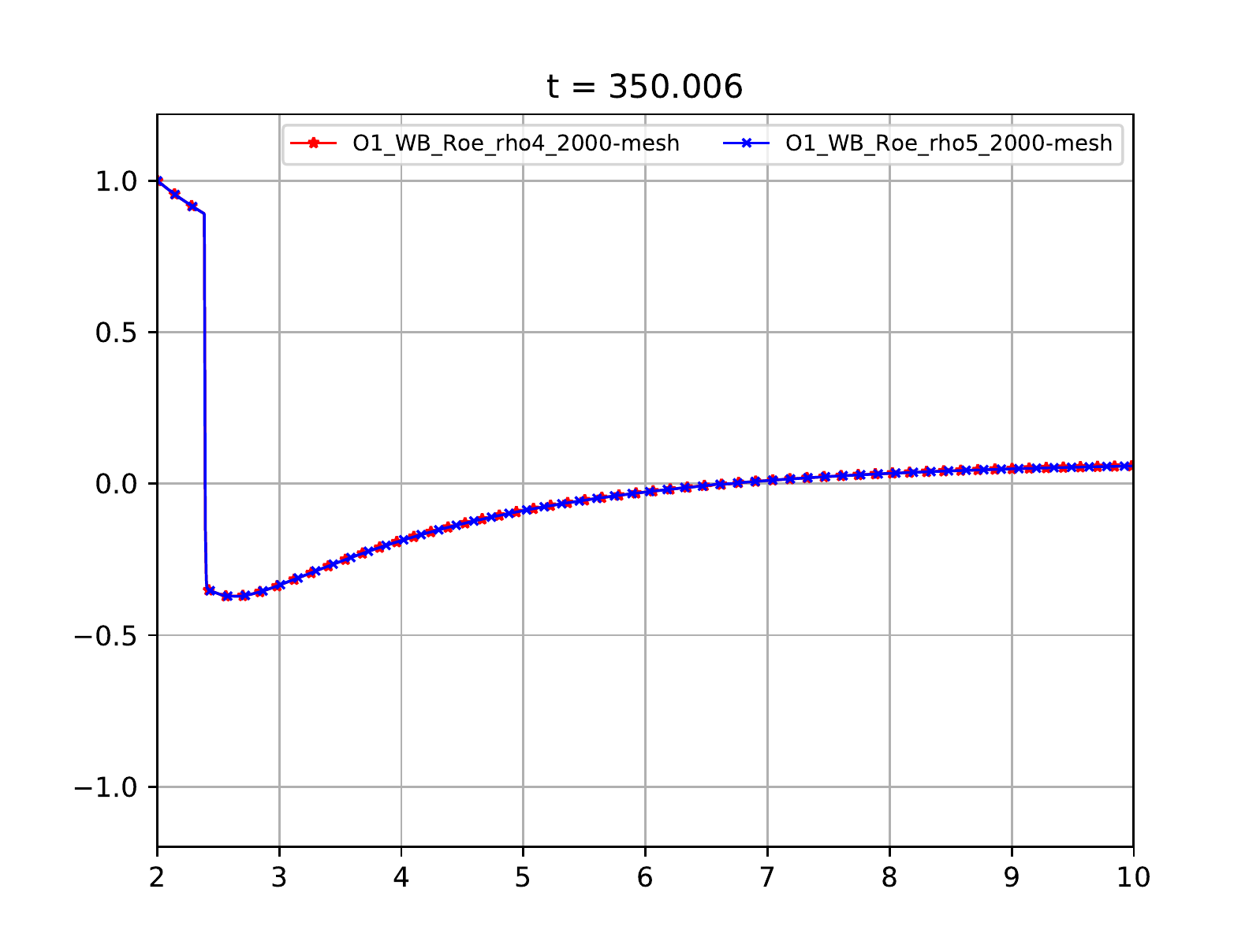}
		\label{fig:Euler_ko1_WB_testPerturbedWB3_negative_perturbation_t_350}
	\end{subfigure}
	\begin{subfigure}[h]{0.32\textwidth}
		\centering
		\includegraphics[width=1\linewidth]{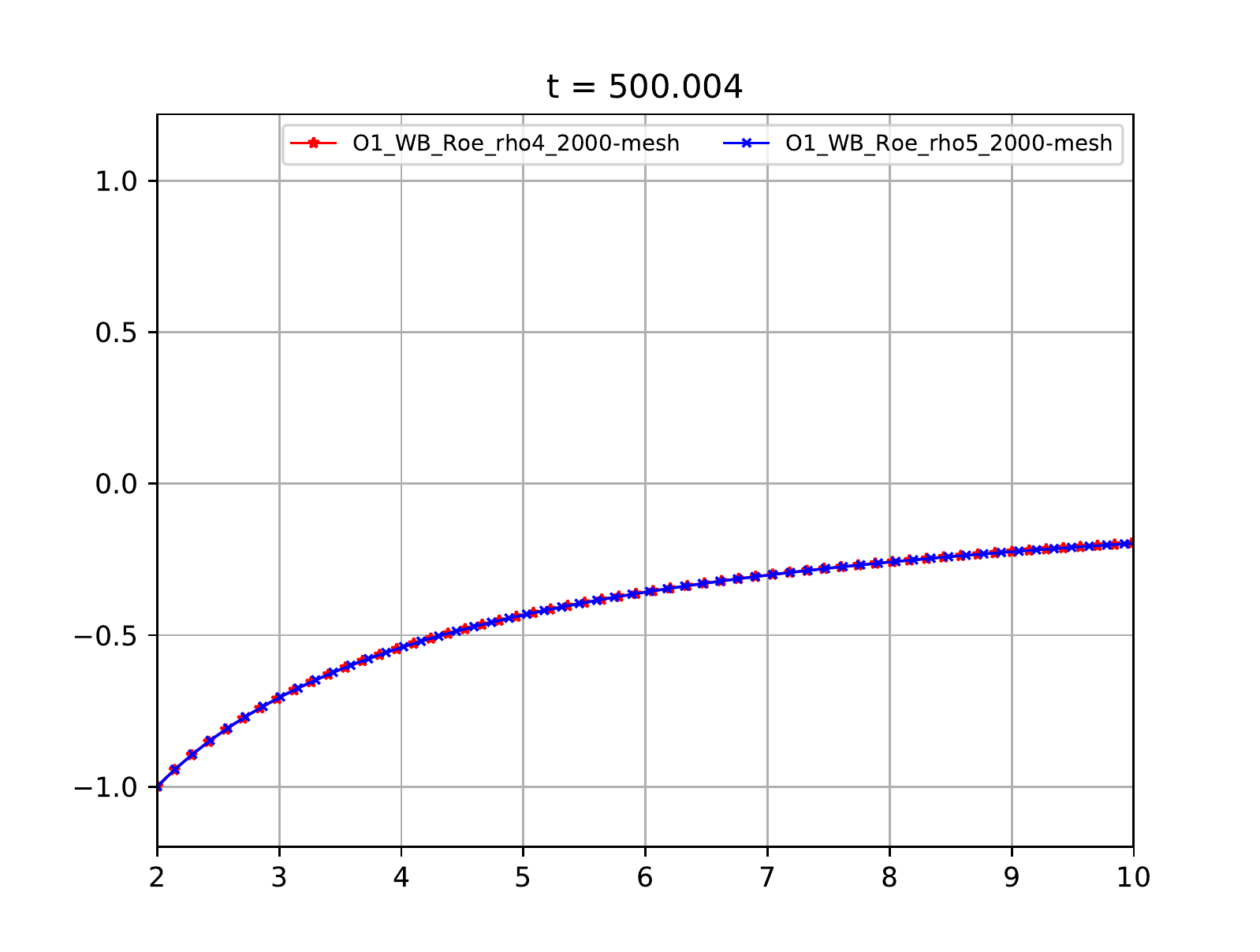}
		\label{fig:Euler_ko1_WB_testPerturbedWB3_negative_perturbation_t_500}
	\end{subfigure}
	\begin{subfigure}[h]{0.32\textwidth}
		\centering
		\includegraphics[width=1\linewidth]{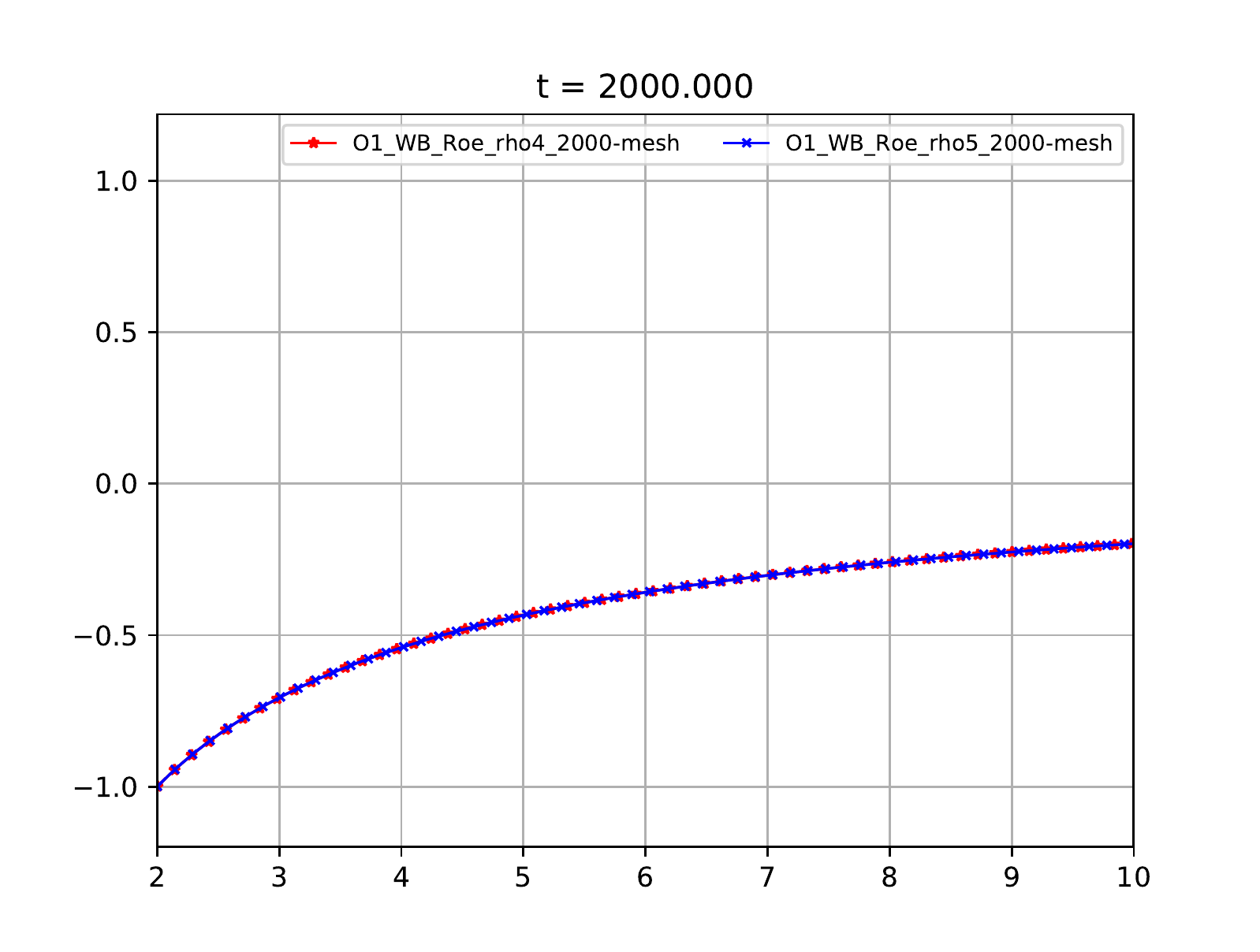}
		\label{fig:Euler_ko1_WB_testPerturbedWB3_negative_perturbation_t_2000}
	\end{subfigure}
	\caption{Euler-Schwarzschild model with the initial conditions \eqref{eq:testE6a} and \eqref{eq:testE6aa}: first-order well-balanced method with a 2000-point mesh using the Roe-type numerical flux at selected times for the variable $v$: the numerical solutions coincide.}
	\label{fig:Euler_ko1_WB_testPerturbedWB3_negative_perturbation}
\end{figure}

\begin{figure}[h]
	\begin{subfigure}[h]{0.32\textwidth}
		\centering
		\includegraphics[width=1\linewidth]{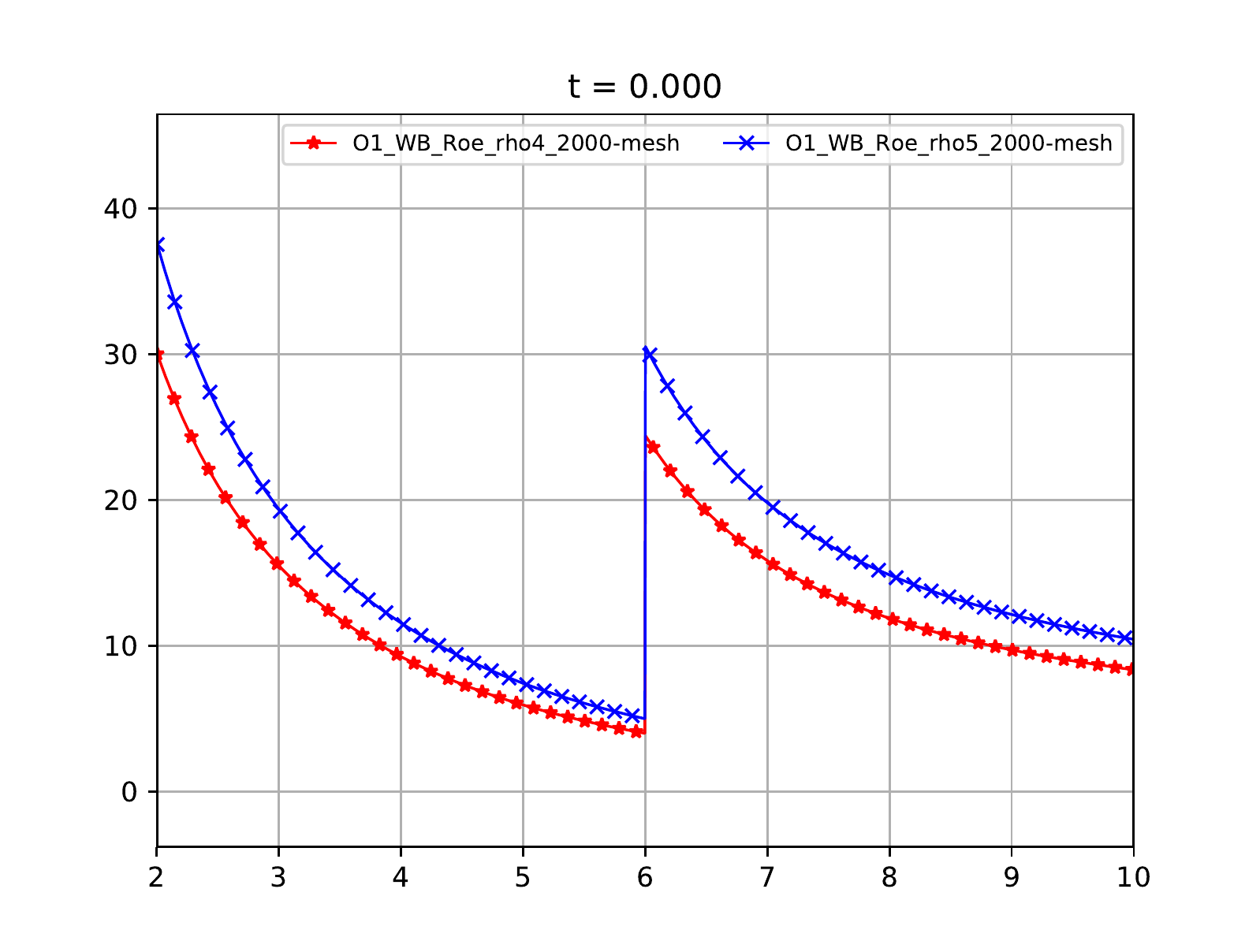}
		\label{fig:Euler_ko1_WB_testPerturbedWB3_negative_perturbation_t_0_rho}
	\end{subfigure}
	\begin{subfigure}[h]{0.32\textwidth}
		\centering
		\includegraphics[width=1\linewidth]{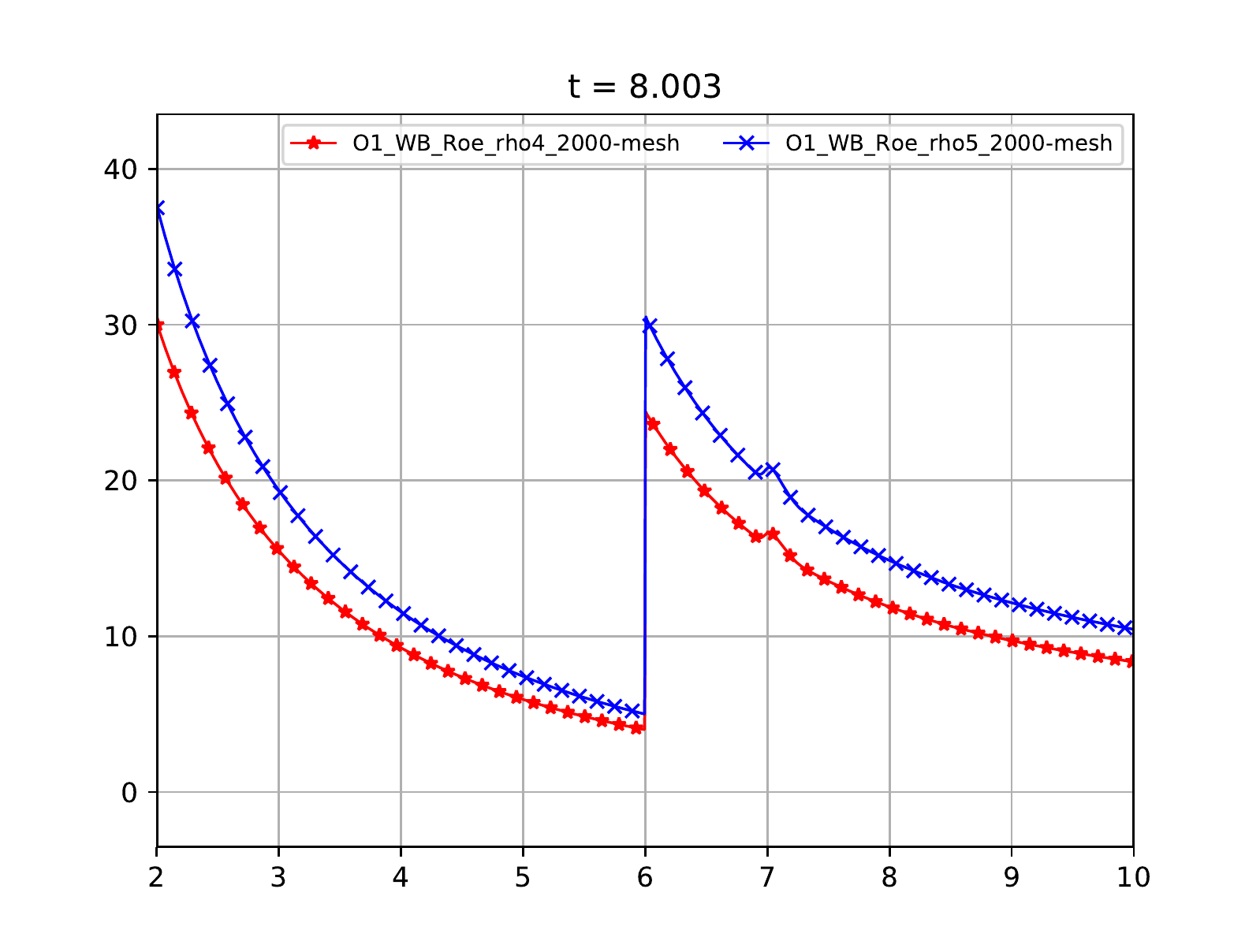}
		\label{fig:Euler_ko1_WB_testPerturbedWB3_negative_perturbation_t_8_rho}
	\end{subfigure}
	\begin{subfigure}[h]{0.32\textwidth}
		\centering
		\includegraphics[width=1\linewidth]{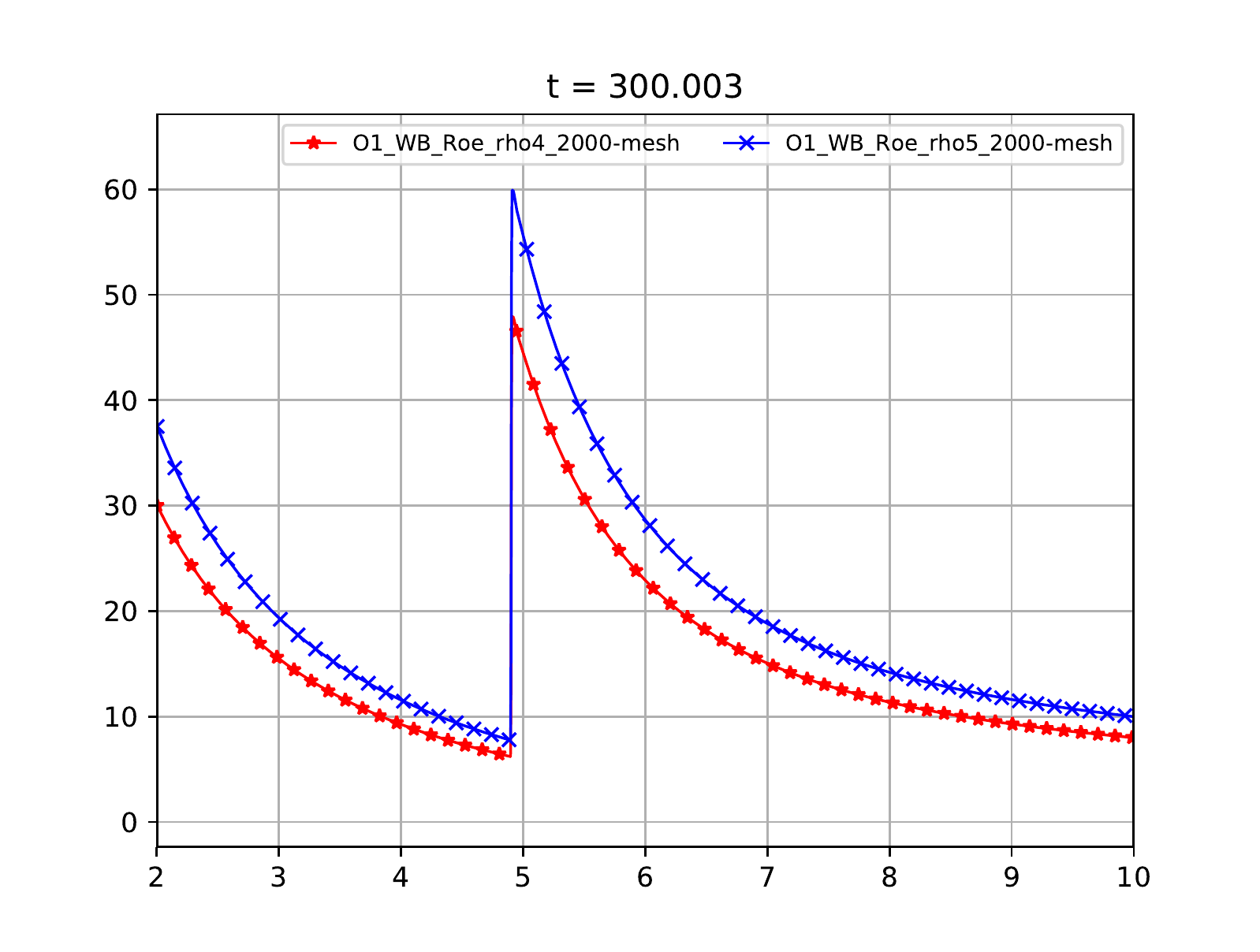}
		\label{fig:Euler_ko1_WB_testPerturbedWB3_negative_perturbation_t_300_rho}
	\end{subfigure}
	\begin{subfigure}[h]{0.32\textwidth}
		\centering
		\includegraphics[width=1\linewidth]{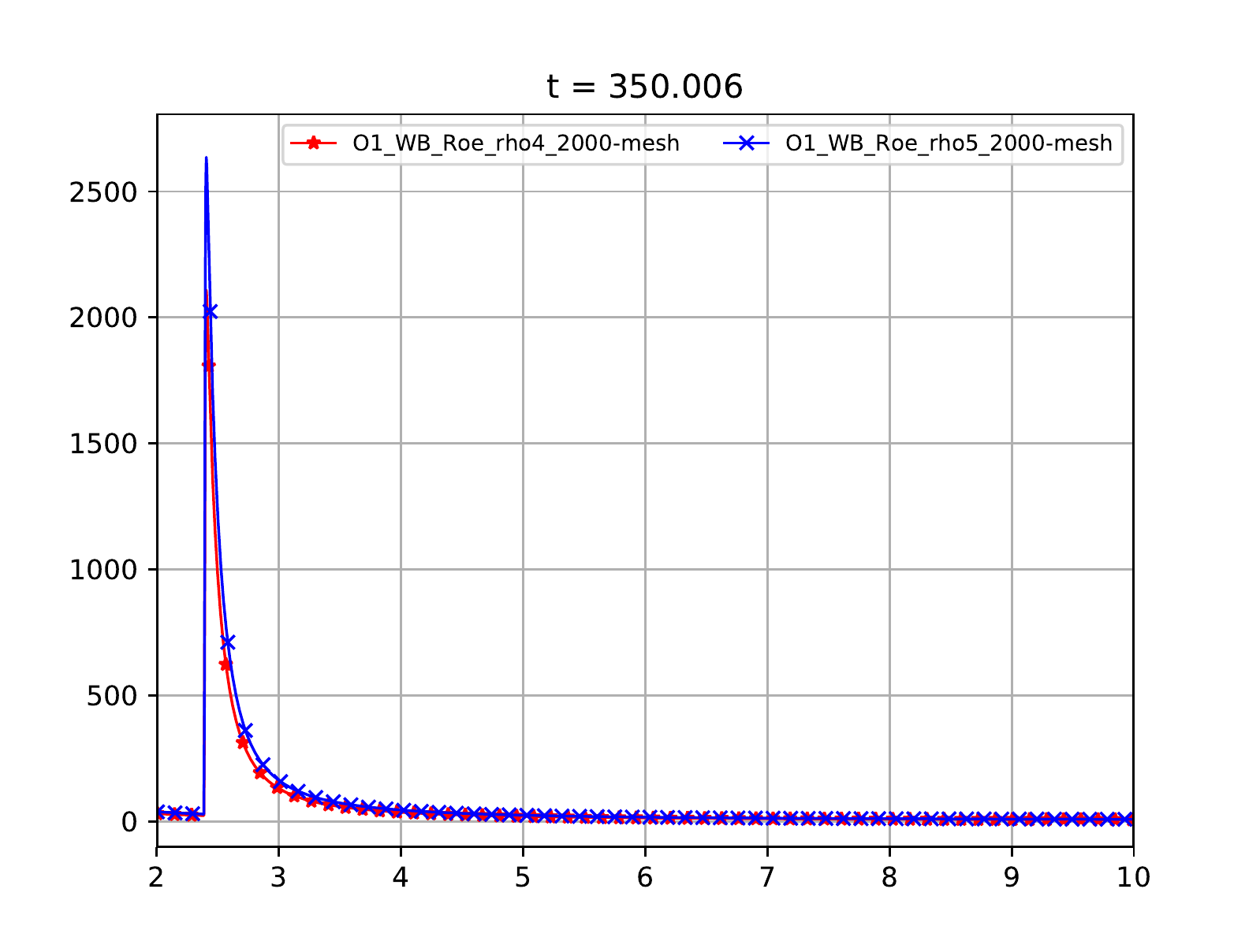}
		\label{fig:Euler_ko1_WB_testPerturbedWB3_negative_perturbation_t_350_rho}
	\end{subfigure}
	\begin{subfigure}[h]{0.32\textwidth}
		\centering
		\includegraphics[width=1\linewidth]{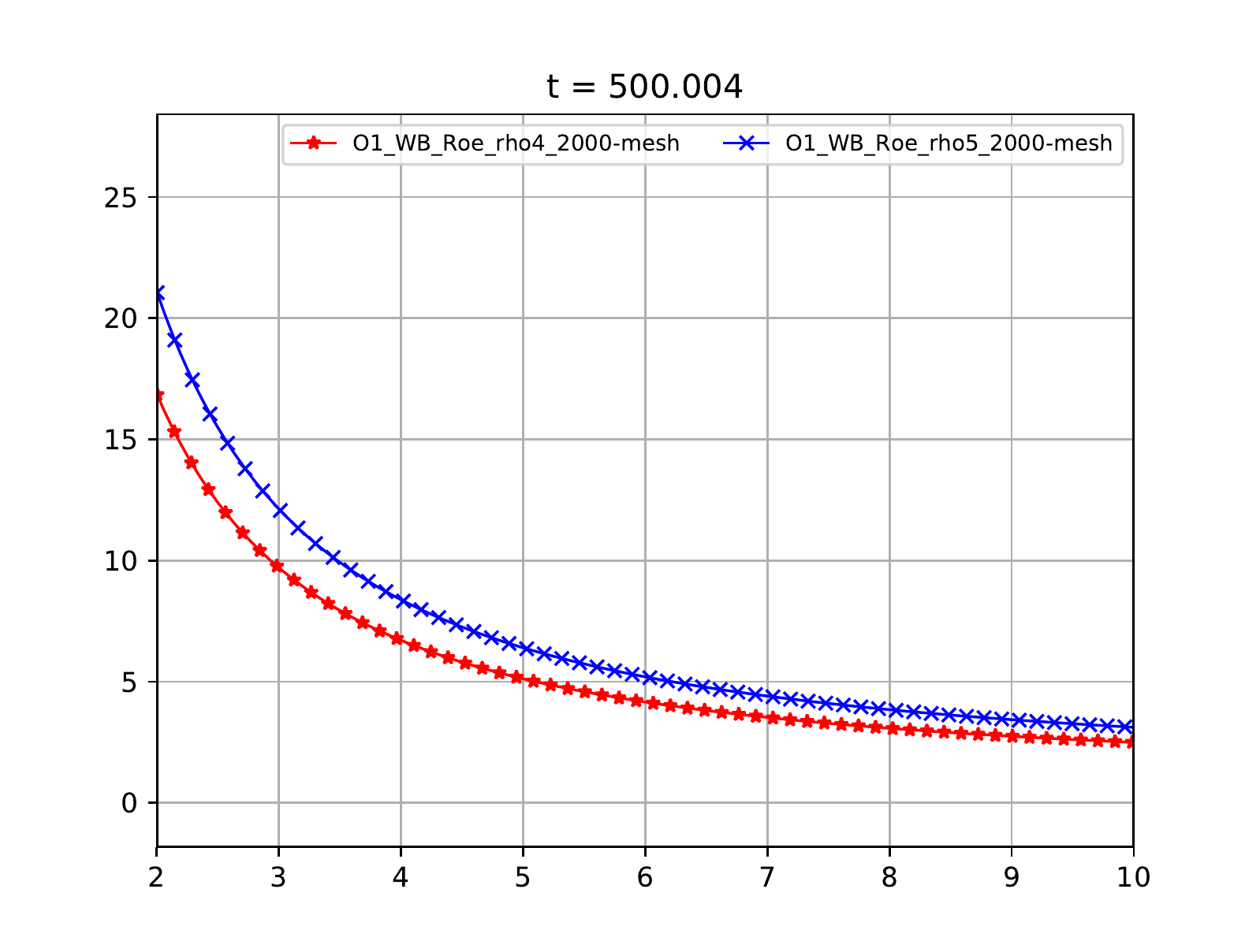}
		\label{fig:Euler_ko1_WB_testPerturbedWB3_negative_perturbation_t_500_rho}
	\end{subfigure}
	\begin{subfigure}[h]{0.32\textwidth}
		\centering
		\includegraphics[width=1\linewidth]{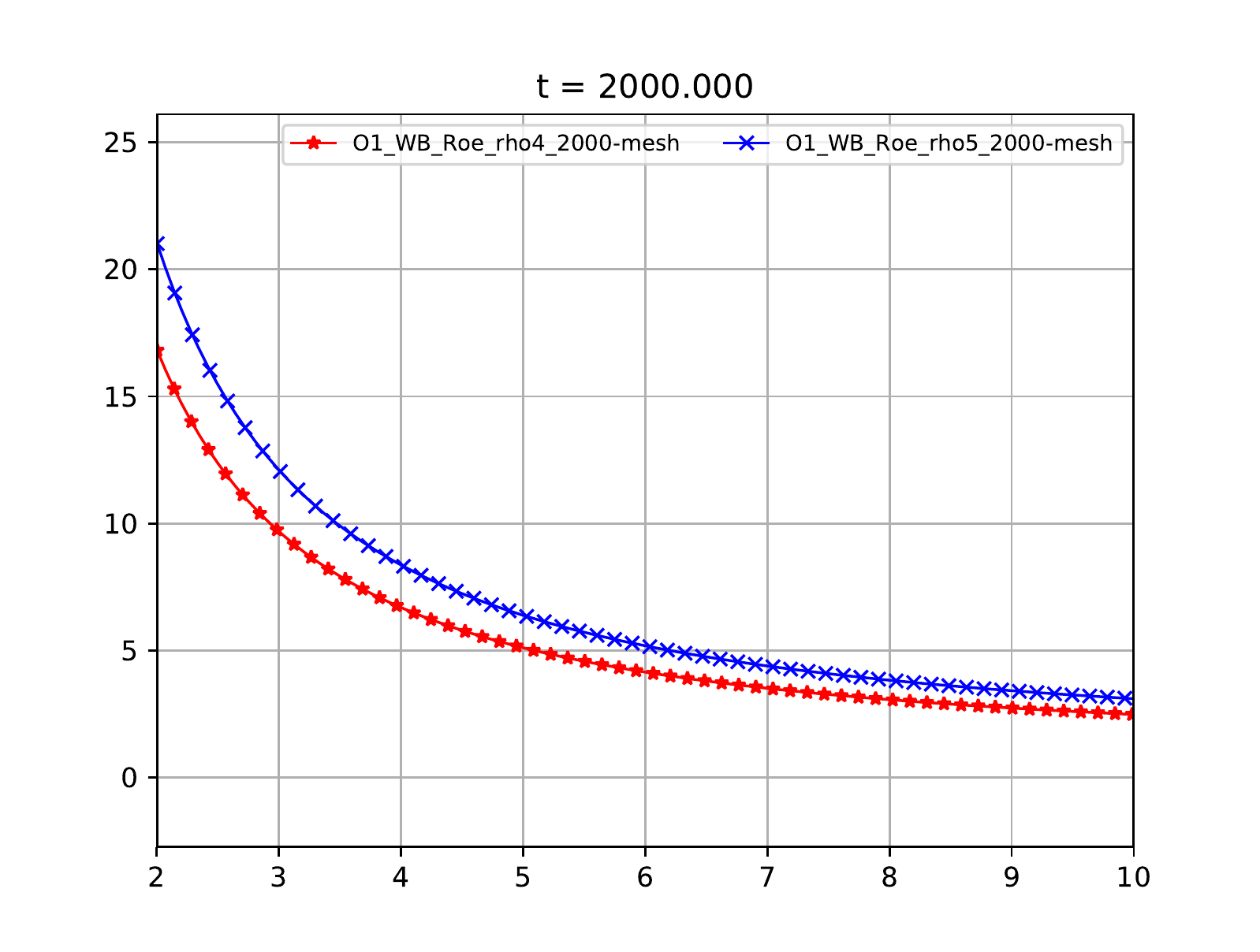}
		\label{fig:Euler_ko1_WB_testPerturbedWB3_negative_perturbation_t_2000_rho}
	\end{subfigure}
	\caption{Euler-Schwarzschild model with the initial conditions \eqref{eq:testE6a} and \eqref{eq:testE6aa}: first-order well-balanced method with a 2000-point mesh using the Roe-type numerical flux at selected times for the variable $\rho$.}
	\label{fig:Euler_ko1_WB_testPerturbedWB3_negative_perturbation_rho}
\end{figure}

After the passage of the perturbation, the  shock starts moving leftward  and, in both cases, the numerical solution converges to a smooth transonic  stationary solution of the form:
\bel{eq:testE6_solution}
V^*(r) = \begin{cases}
	V_{-}^{*}(r), & \text{ $r\leq r_c$},\\
	V_{+}^{*}(r), & \text{ otherwise,}
\end{cases}
\ee 
where $r_c$ is given by \eqref{eq:r_critical_steady_state_Euler};  $V_{-}^{*}(r)$  and $V_{+}^{*}(r)$  are respectively a subsonic and a supersonic stationary solution  satisfying
$v_{\pm}^{*}(r_c)=-k$: see Figure \ref{fig:Euler_ko1_WB_testPerturbedWB3_negative_perturbation}. Nevertheless, the limits in time  of the approximations of $\rho$ are different:
see Figure \ref{fig:Euler_ko1_WB_testPerturbedWB3_negative_perturbation_rho}. Observe that, in the Euler-Schwarzschild model (\ref{eq:steady_state_Euler}), 
there are infinitely many stationary solutions with the same function  $v$ and  different $\rho$.\\

\paragraph{{{\red Relation between the perturbation and the displacement of the shock}}}

\begin{figure}[h]
	\begin{subfigure}[h]{0.32\textwidth}
		\centering
		\includegraphics[width=1\linewidth]{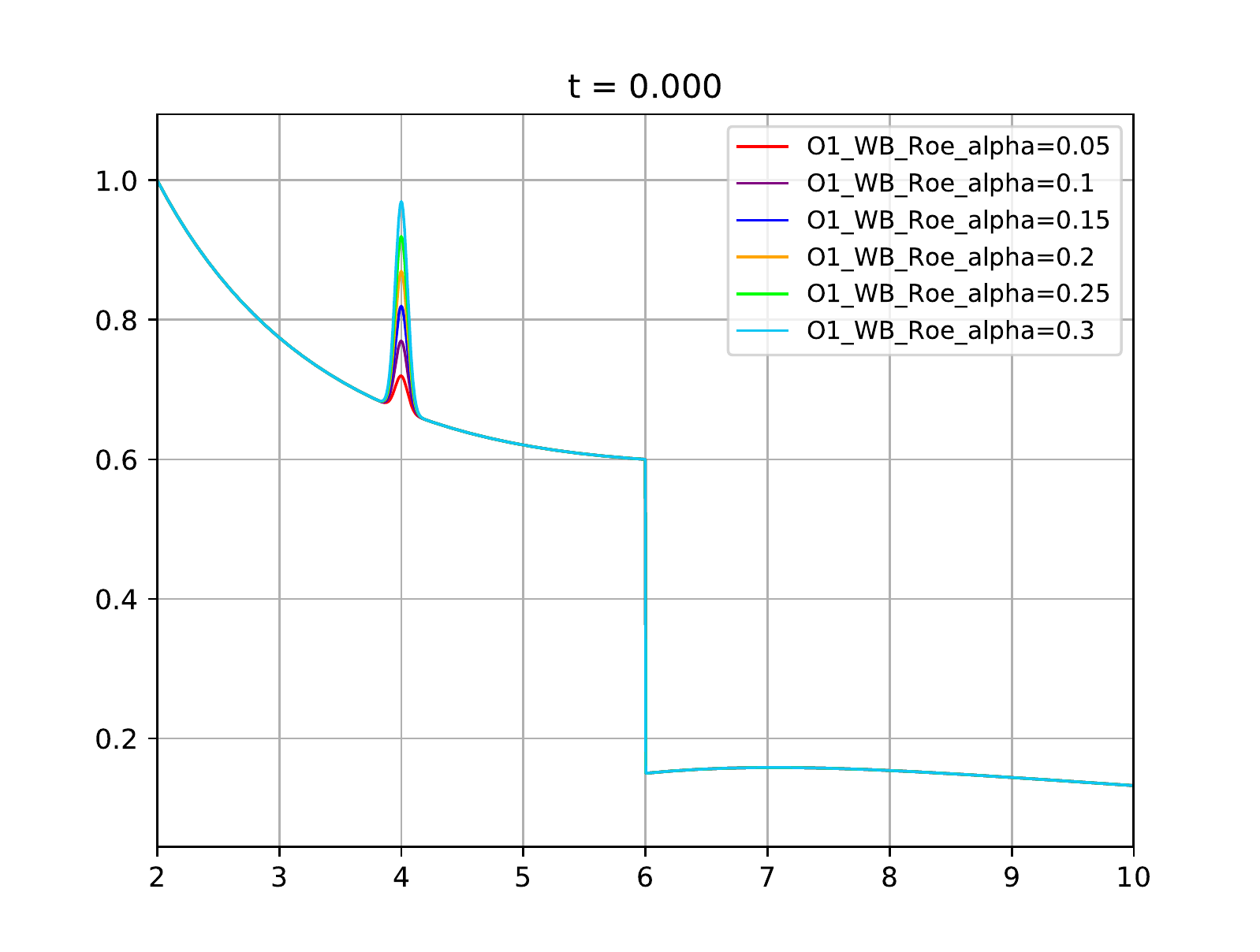}
		\label{fig:Euler_ko1_WB_testPerturbedWB3_alpha_comparison_t_0}
	\end{subfigure}
\begin{subfigure}[h]{0.32\textwidth}
	\centering
	\includegraphics[width=1\linewidth]{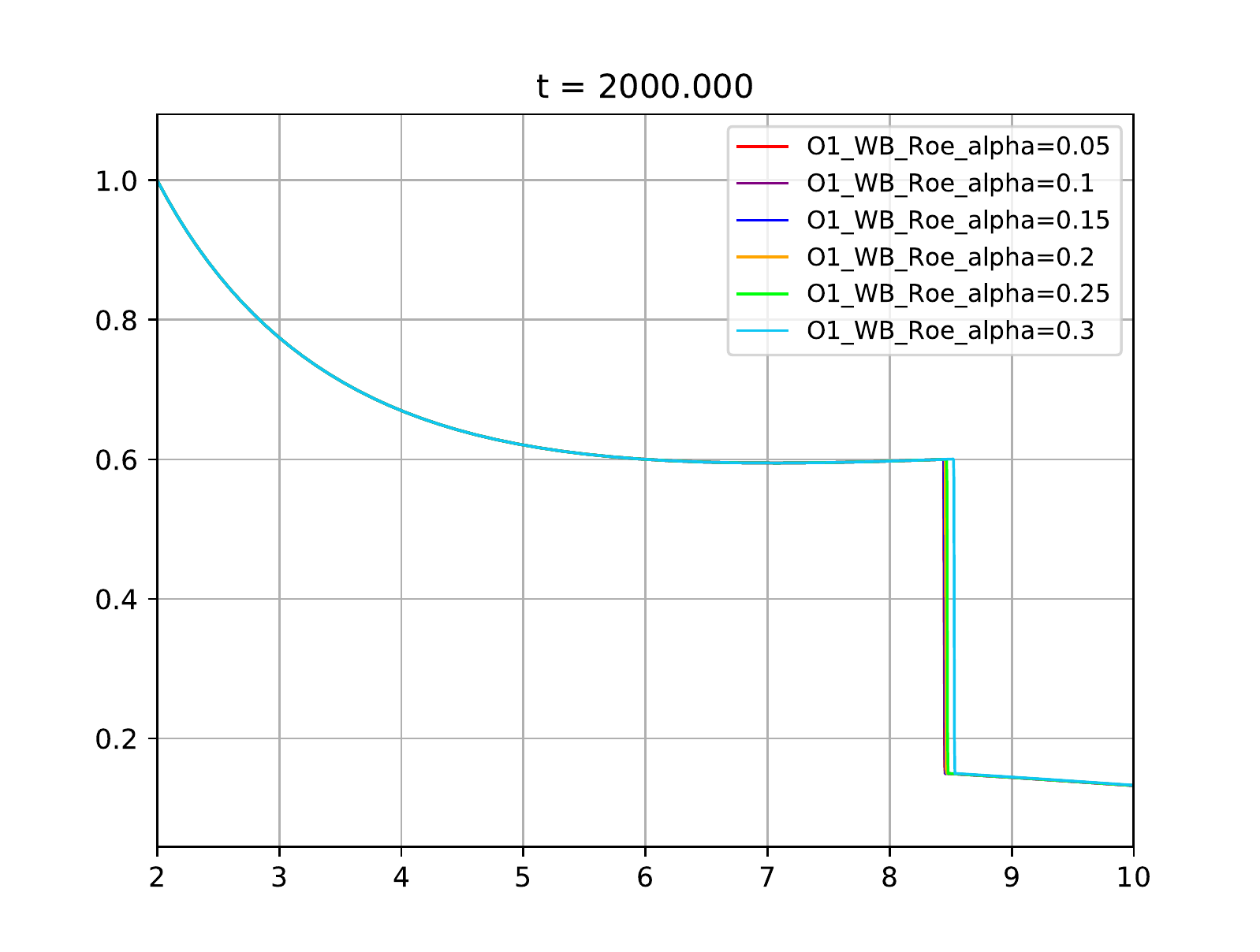}
	\label{fig:Euler_ko1_WB_testPerturbedWB3_alpha_comparison_t_2000}
\end{subfigure}
\begin{subfigure}[h]{0.32\textwidth}
	\centering
	\includegraphics[width=1\linewidth]{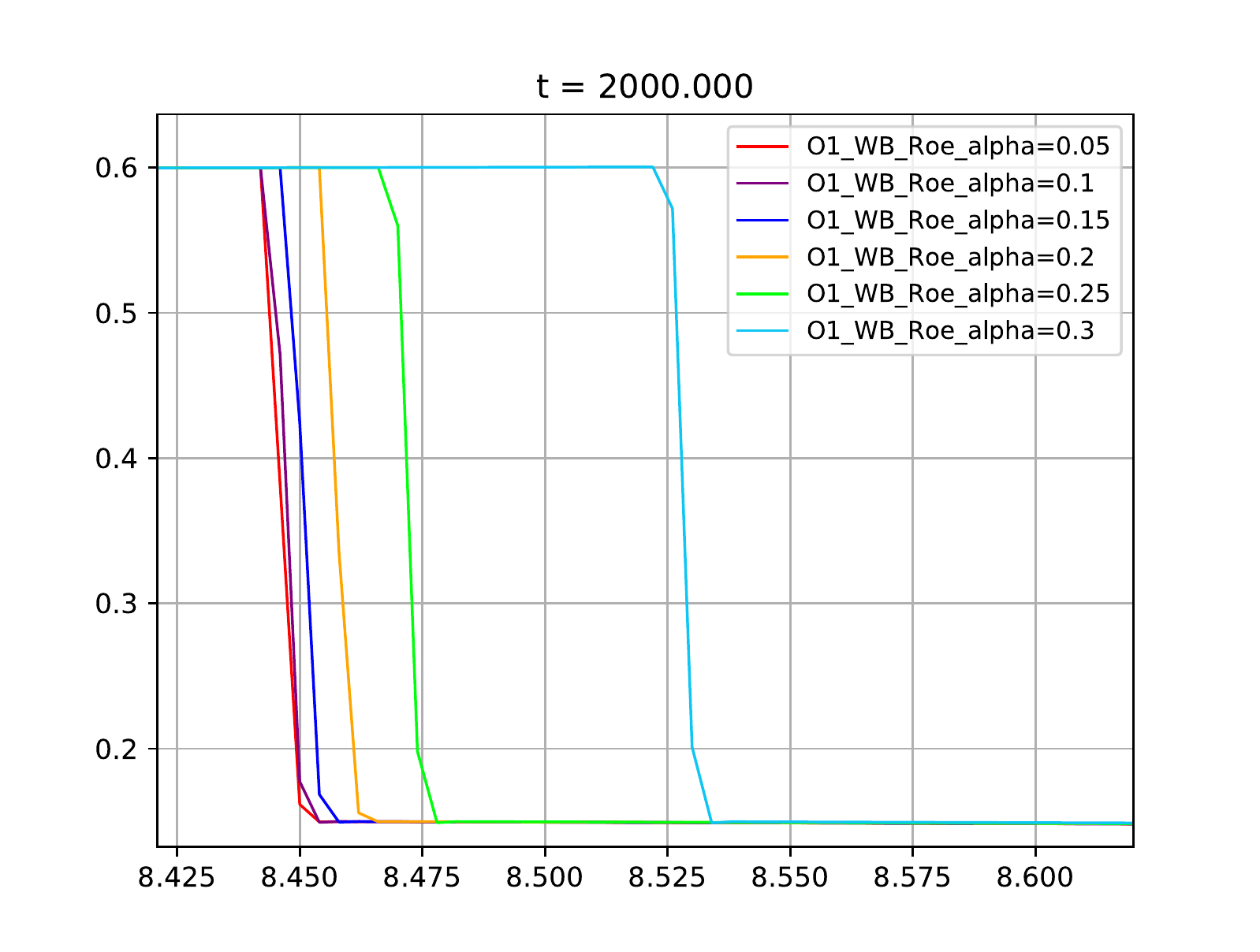}
	\caption{Zoom.}
	\label{fig:Euler_ko1_WB_testPerturbedWB3_alpha_comparison_t_2000_zoom}
\end{subfigure}
	\caption{Euler-Schwarzschild model with the initial condition \eqref{eq:testE7a}: first-order well-balanced method taking different values of $\alpha$ for variable $v$.}
	\label{fig:Euler_ko1_WB_testPerturbedWB3_alpha_comparison}
\end{figure}

\begin{figure}[h]
	\begin{subfigure}[h]{0.32\textwidth}
		\centering
		\includegraphics[width=1\linewidth]{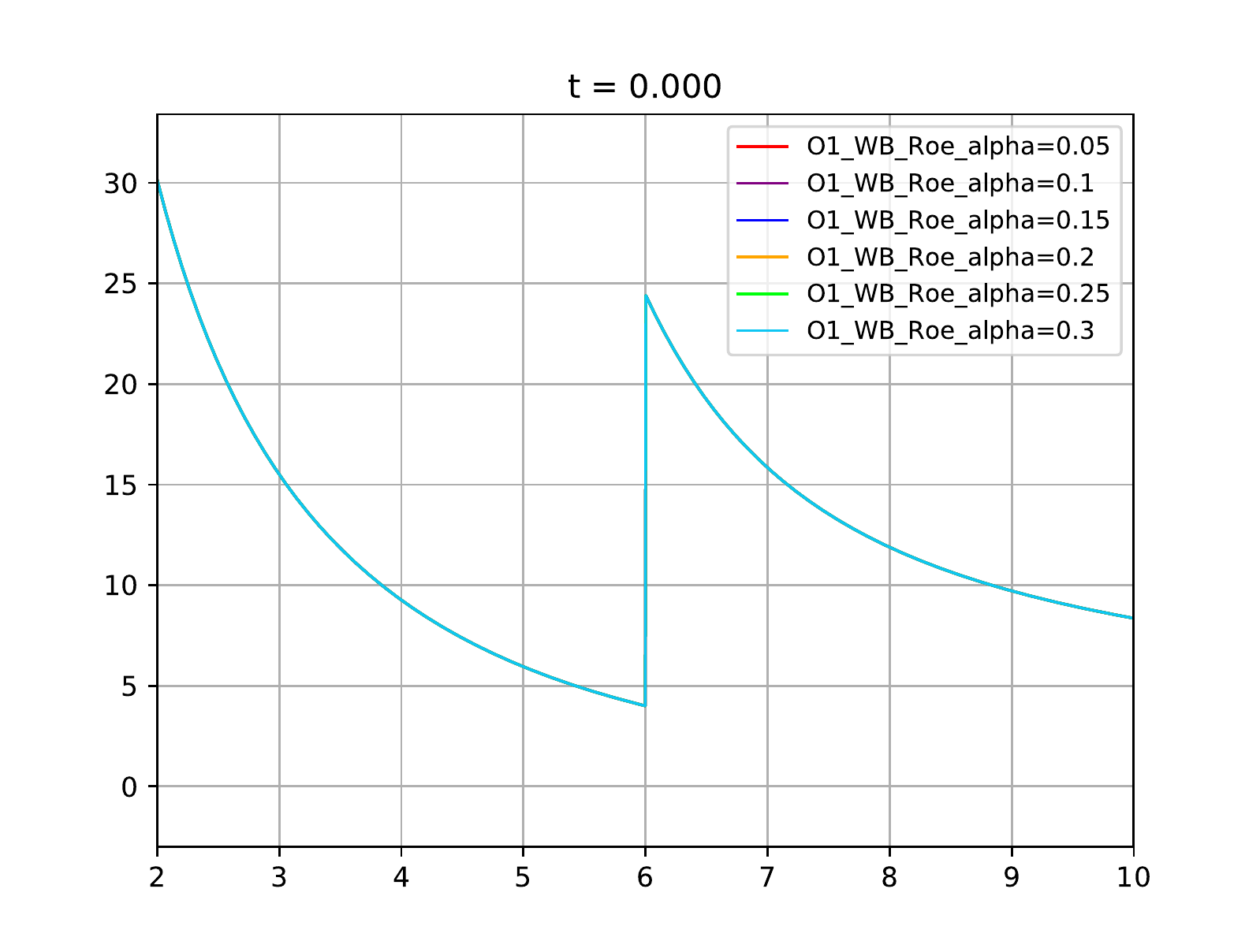}
		\label{fig:Euler_ko1_WB_testPerturbedWB3_alpha_comparison_t_0_rho}
	\end{subfigure}
	\begin{subfigure}[h]{0.32\textwidth}
		\centering
		\includegraphics[width=1\linewidth]{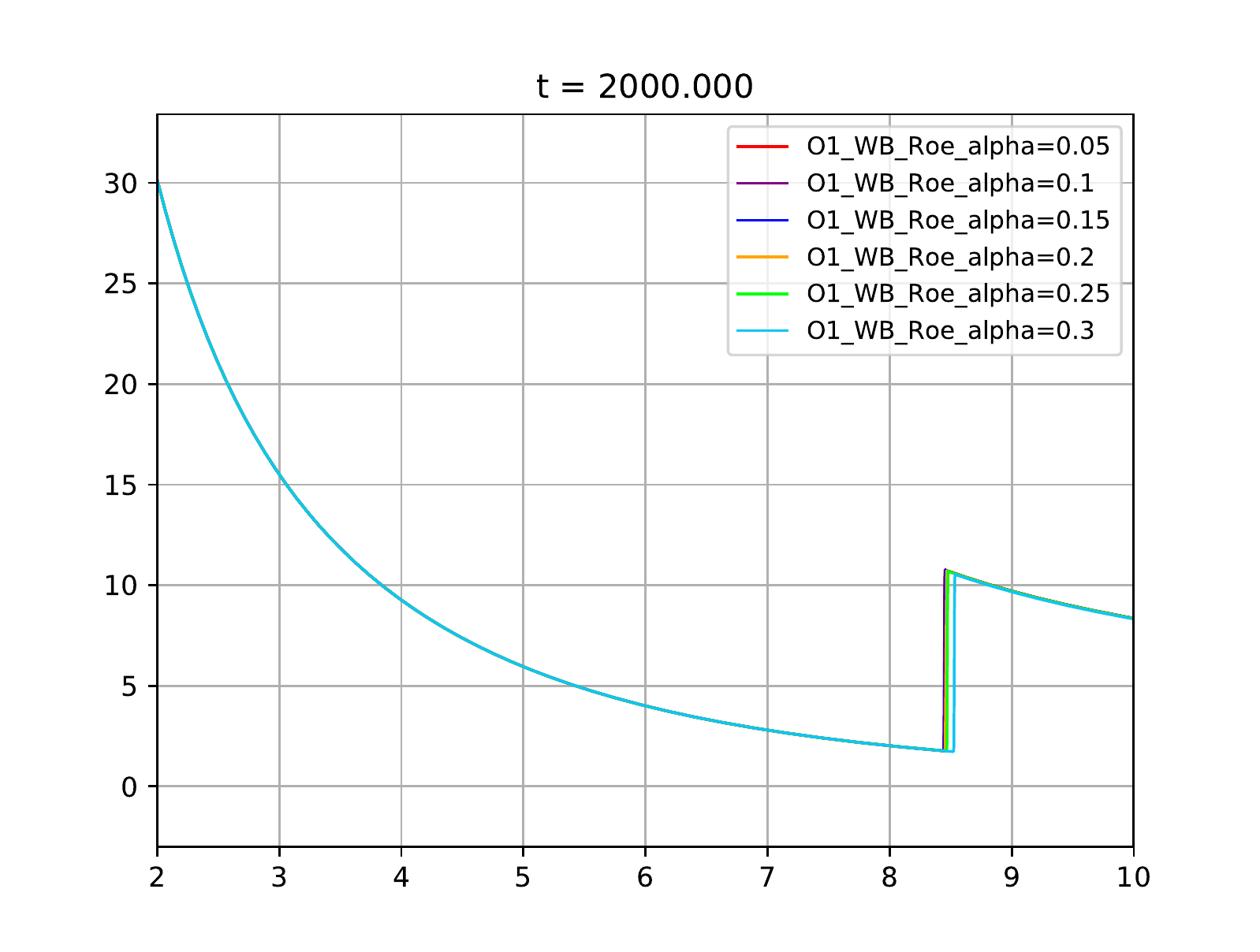}
		\label{fig:Euler_ko1_WB_testPerturbedWB3_alpha_comparison_t_2000_rho}
	\end{subfigure}
	\begin{subfigure}[h]{0.32\textwidth}
		\centering
		\includegraphics[width=1\linewidth]{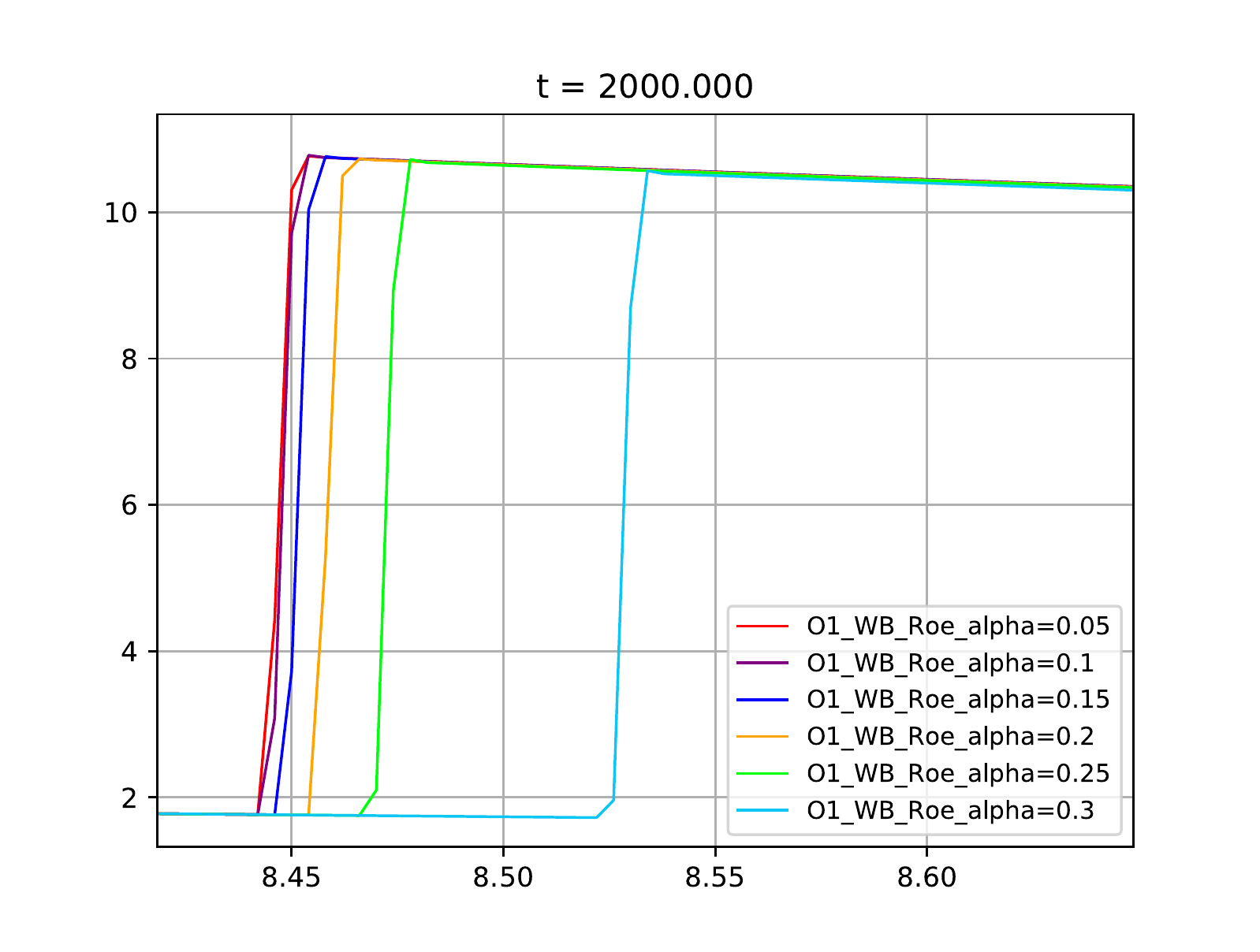}
		\caption{Zoom.}
		\label{fig:Euler_ko1_WB_testPerturbedWB3_alpha_comparison_t_2000_zoom_rho}
	\end{subfigure}
	\caption{Euler-Schwarzschild model with the initial condition \eqref{eq:testE7a}: first-order well-balanced method taking different values of $\alpha$ for variable $\rho$.}
	\label{fig:Euler_ko1_WB_testPerturbedWB3_alpha_comparison_rho}
\end{figure}

In order to study the relationship between the amplitude of the perturbation and the distance between the initial and the final shock locations, we consider the family of initial conditions:
\bel{eq:testE7a}
\widetilde V_{0}(r) = 
\widetilde{V}^*(r) + \delta(\alpha, r),
\ee
where $\widetilde V^*$ is the steady shock solution given by \eqref{eq:testE3a}-\eqref{eq:testE3c}  and 
\bel{eq:testE7b}
\delta(\alpha, r) = [\delta_{v}(\alpha, r),
\delta_{\rho}(\alpha, r)]^T = \begin{cases}
	[\alpha e^{-200(r-4)^{2}}, 0]^T, & \text{ $ 3<r<5$,}\\
	[0, 0]^T, & \text{ otherwise,}
\end{cases}
\ee
with $\alpha>0$. In this case we will also use the Roe-type numerical flux and a 2000-point uniform mesh.
Figures \ref{fig:Euler_ko1_WB_testPerturbedWB3_alpha_comparison} and \ref{fig:Euler_ko1_WB_testPerturbedWB3_alpha_comparison_rho} show the numerical solution for different values of $\alpha$ and we observe that depending on the amplitude of the perturbation the numerical solutions converge in time to  different steady shock solutions.

\begin{figure}
	\centering
	\includegraphics[width=.65\linewidth]{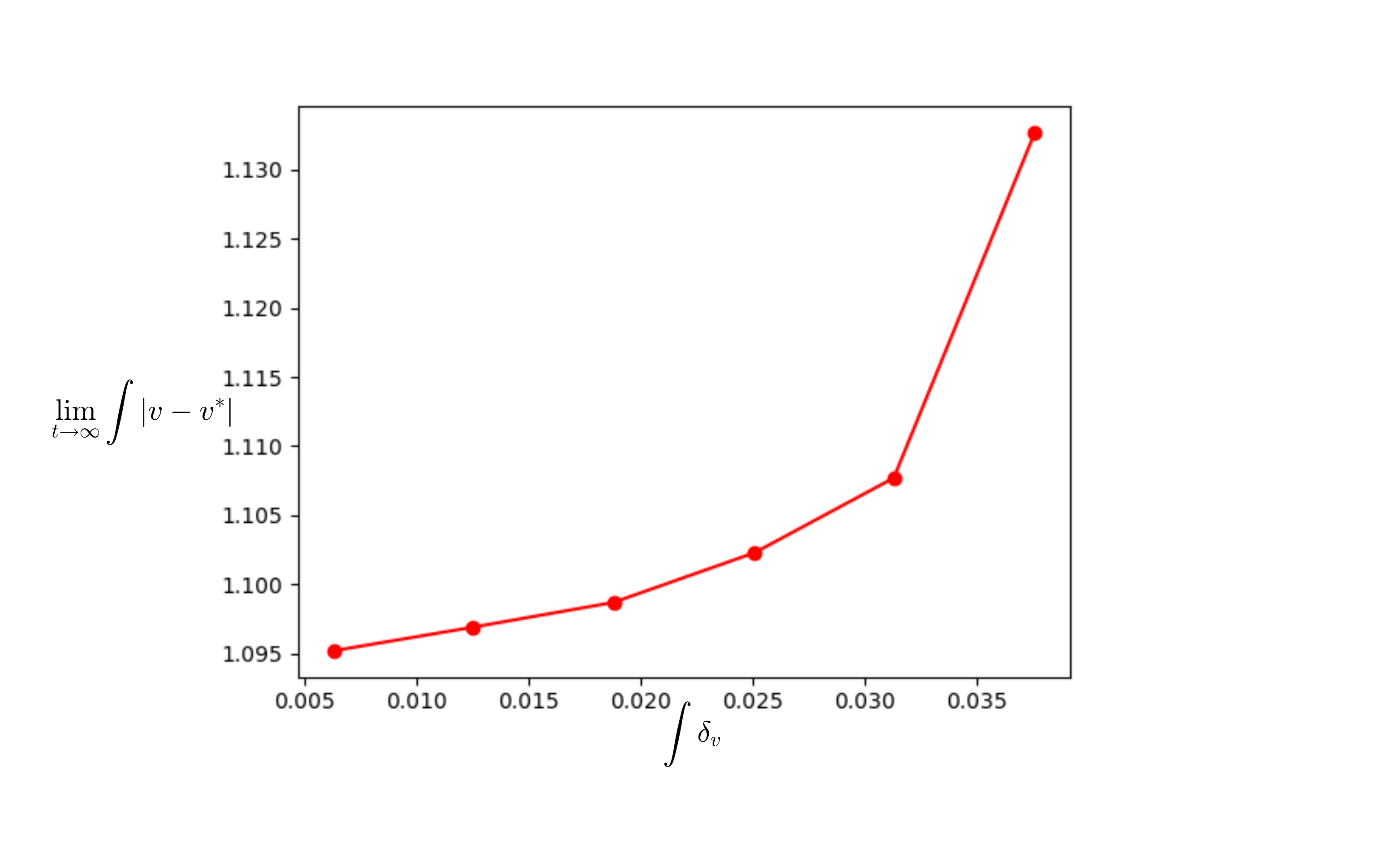}
	\caption{Euler-Schwarzschild model with the initial condition \eqref{eq:testE7a}: values of $\lim_{t\to\infty} \int |v -v^{*}|$ as a function of
		$\int \delta_{v}$.}
	\label{fig:eulerintegralperturbationvsintegralsteadyshockperturbedlabelled}
\end{figure}

The amplitude of the perturbation is measured with 
$
\int \delta_v(\alpha, r) \,dr
$
and the distance between the shocks are measured by
$
\lim_{t \to \infty} \int |v(r,t) - v^*(r)| \, dr,
$
as we did for the Burgers-Schwarzschild model.  Table \ref{tab:Euler_Areas_for_different_values_of_alpha} and Figure \ref{fig:eulerintegralperturbationvsintegralsteadyshockperturbedlabelled} show the relationship between those magnitudes: the displacement of the shock seems to grow exponentially with the amplitude.

\begin{table}[h]
	\centering
	\begin{tabular}{||c|c|c||} 
		\hline
		$\alpha$ & $\int \delta_{v}$ & $\lim_{t\to\infty} \int |v -v^{*}|$\\ 
		\hline
		0.05 & 0.0063 & 1.0952 \\ 
		\hline
		0.1 & 0.0125 & 1.0969 \\ 
		\hline
		0.15 & 0.0188 & 1.0987 \\
		\hline
		0.2 & 0.0251 & 1.1023 \\
		\hline
		0.25 & 0.0313 & 1.1077  \\
		\hline
		0.3 & 0.0376 & 1.1327 \\ 
		\hline
	\end{tabular}
	\caption{Euler-Schwarzschild model with the initial condition \eqref{eq:testE7a}: measures of the perturbation and the shock displacement for 
		different values of $\alpha$}
	\label{tab:Euler_Areas_for_different_values_of_alpha}
\end{table}

Let us finally consider a family of initial conditions that generate  leftward displacement of the initial steady shock:
\bel{eq:testE8a}
\widetilde V_{0}(r) = 
\widetilde{V}^*(r) + \delta(\beta, r),
\ee
where $\widetilde V^*$ is again the steady shock solution given by \eqref{eq:testE3a}- \eqref{eq:testE3c} and 

\begin{figure}[h]
	\begin{subfigure}[h]{0.49\textwidth}
		\centering
		\includegraphics[width=0.9\linewidth]{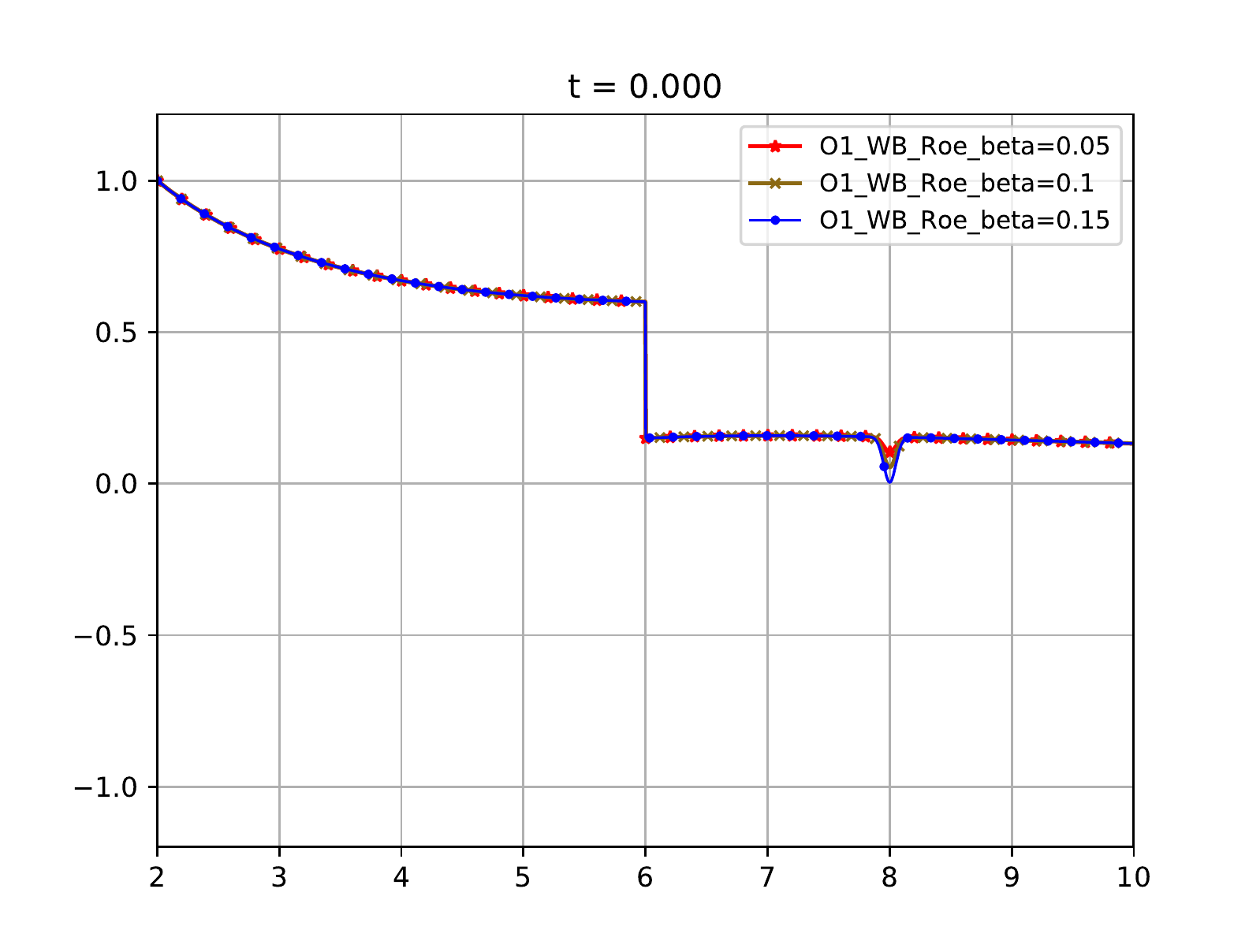}
		\label{fig:Euler_ko1_WB_testPerturbedWB3_beta_comparison_t_0}
	\end{subfigure}
	\begin{subfigure}[h]{0.49\textwidth}
		\centering
		\includegraphics[width=0.9\linewidth]{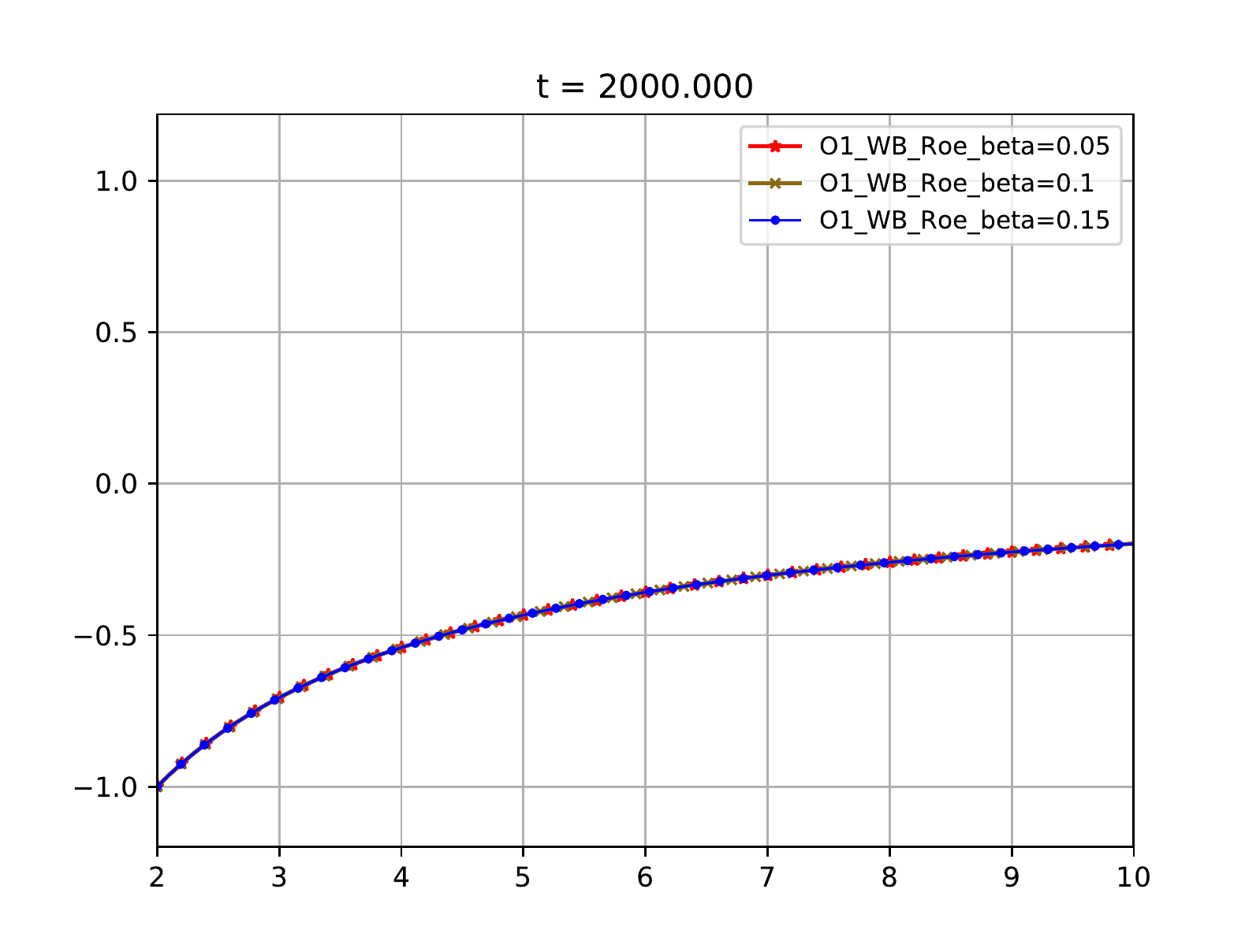}
		\label{fig:Euler_ko1_WB_testPerturbedWB3_beta_comparison_t_2000}
	\end{subfigure}
	\caption{Euler-Schwarzschild model with the initial condition \eqref{eq:testE8a}: first-order well-balanced method taking different values of $\beta$ for variable $v$.}
	\label{fig:Euler_ko1_WB_testPerturbedWB3_beta_comparison}
\end{figure}

\begin{figure}[h]
	\begin{subfigure}[h]{0.49\textwidth}
		\centering
		\includegraphics[width=0.9\linewidth]{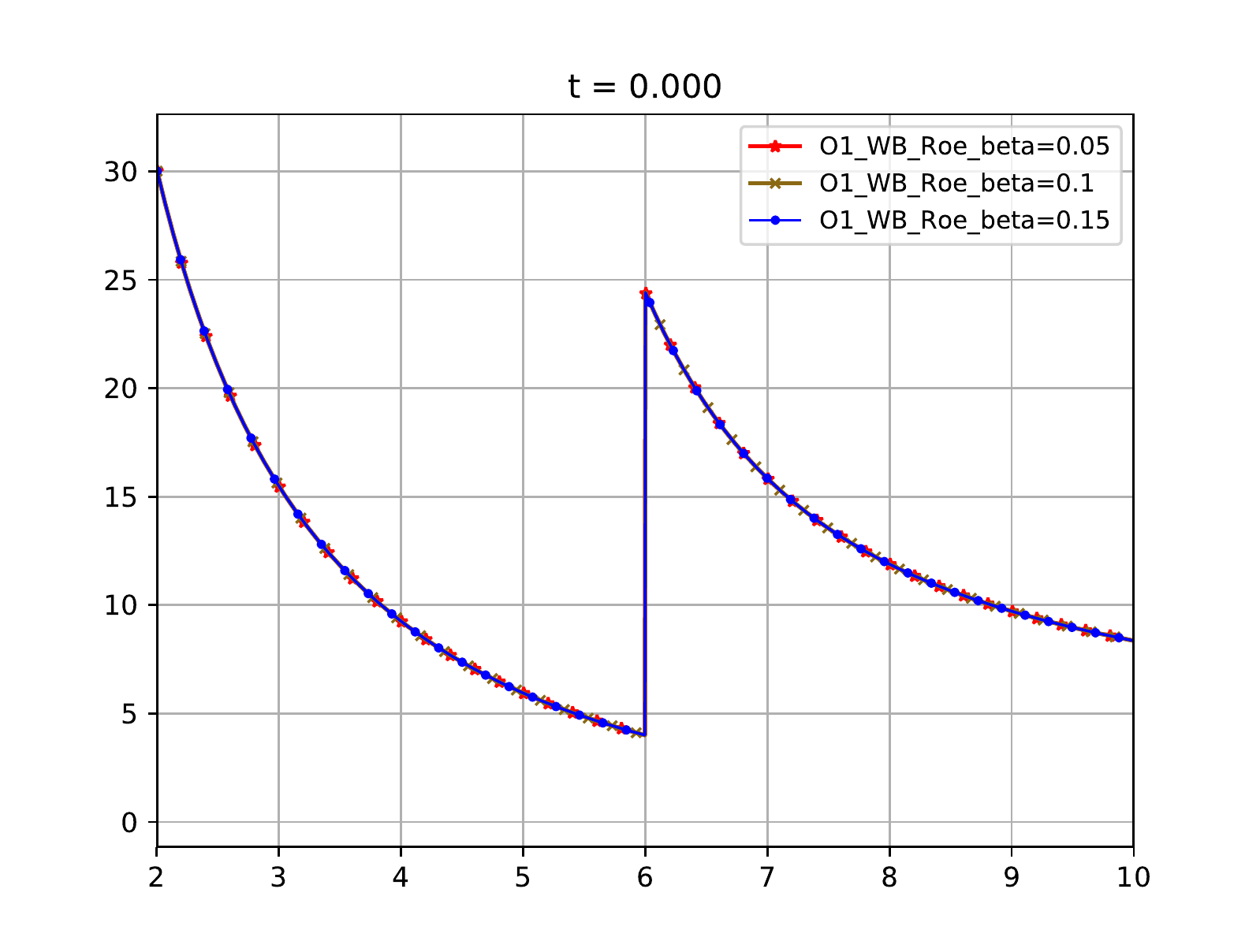}
		\label{fig:Euler_ko1_WB_testPerturbedWB3_beta_comparison_t_0_rho}
	\end{subfigure}
	\begin{subfigure}[h]{0.49\textwidth}
		\centering
		\includegraphics[width=0.9\linewidth]{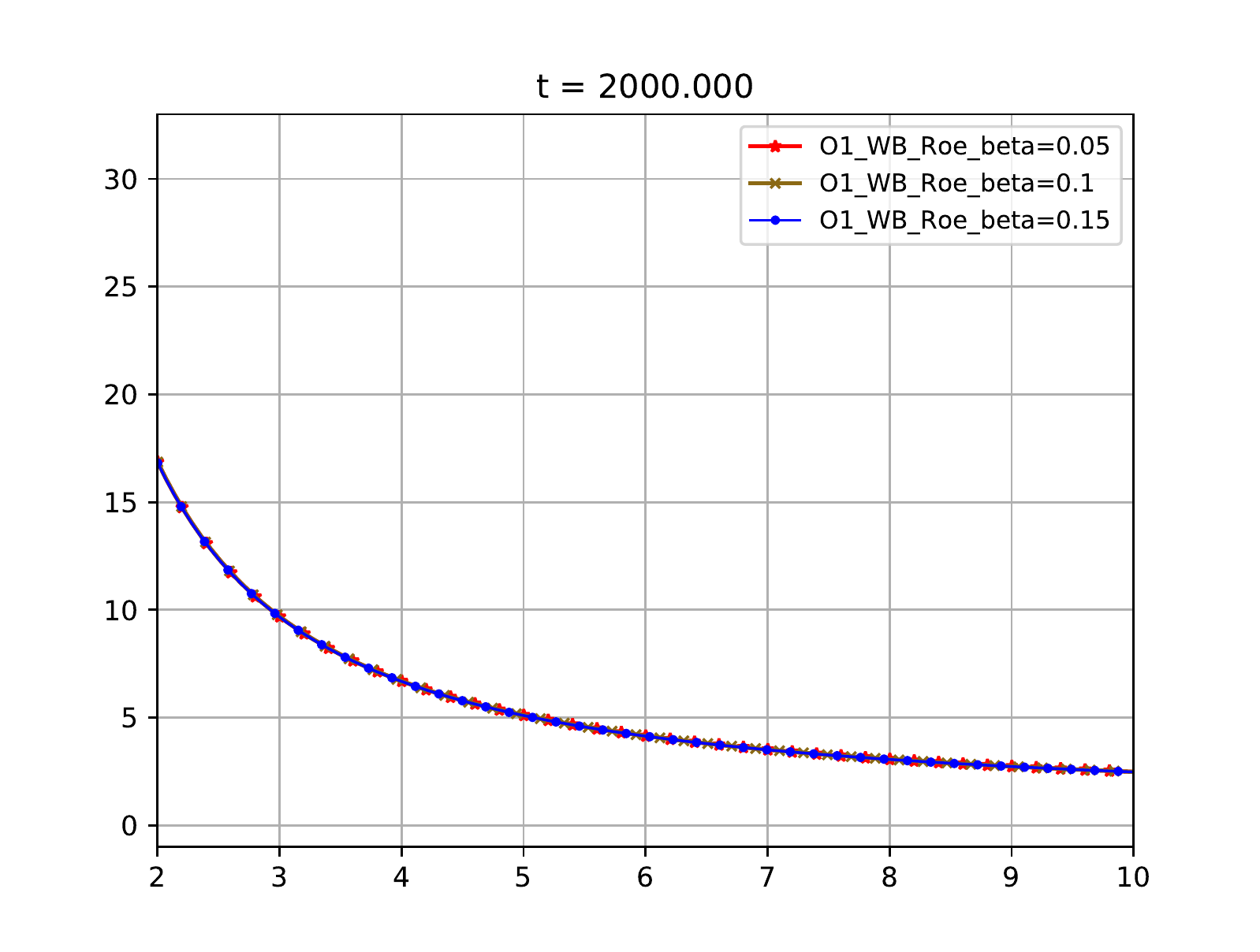}
		\label{fig:Euler_ko1_WB_testPerturbedWB3_beta_comparison_t_2000_rho}
	\end{subfigure}
	\caption{Euler-Schwarzschild model with the initial condition \eqref{eq:testE8a}: first-order well-balanced method taking different values of $\beta$ for variable $\rho$.}
	\label{fig:Euler_ko1_WB_testPerturbedWB3_beta_comparison_rho}
\end{figure}

\bel{eq:testE8b}
\delta(\beta, r) = [\beta{v}(\beta, r),
\delta_{\rho}(\beta, r)]^T = \begin{cases}
	[\beta e^{-200(r-8)^{2}}, 0]^T, & \text{ $7<r<8$,}\\
	[0, 0]^T, & \text{ otherwise,}
\end{cases}
\ee
with $\beta<0$. In this case we will use the Roe-type numerical flux and a 2000-point uniform mesh.
Figures \ref{fig:Euler_ko1_WB_testPerturbedWB3_beta_comparison} and \ref{fig:Euler_ko1_WB_testPerturbedWB3_beta_comparison_rho} show the numerical solution for different values of $\beta$ and we observe that it converges to same stationary solution without depending on the perturbation.\\

\subsection{Main conclusions for the Euler-Schwarzschild model}\label{subsec_conclusions_Euler}

In view of Figure \ref{fig:Euler_ko1_ko2_WB_testPerturbedWB1_hepse14} we arrive at the following observation.

\begin{conclusion}
	If a smooth stationary solution of the Euler system (\ref{eq:Euler}) is perturbed, the solution is restored once the wave generated by the perturbation  goes away.
\end{conclusion}

In view of Figures \ref{fig:Euler_ko1_WB_testPerturbedWB3_HLL_vs_Lax}-\ref{fig:Euler_ko1_WB_testPerturbedWB3_beta_comparison_rho} and Table \ref{tab:Euler_Areas_for_different_values_of_alpha} we arrive at the following observation.

\begin{conclusion}
Consider a perturbation $\delta= (\delta_v, \delta_\rho)$ added to a steady shock solution of the form
	$$
	V_{0}(r) =  \begin{cases}
	V_{-}^{*}(r), & \text{ $r\leq r_0$},\\
	V_{+}^{*}(r), & \text{ otherwise.}
	\end{cases}
	$$
	\begin{enumerate}
		\item If the perturbation moves the steady shock to the right, then a different stationary solution of the form
		$$
		V(r) =  \begin{cases}
		V_{-}^{*}(r), & \text{ $r\leq r_1$},\\
		V_{+}^{*}(r), & \text{ otherwise,}
		\end{cases}
		$$
		with $r_0 \neq r_1$, is obtained, and the distance between $r_0$ and $r_1$ seems to depend exponentially on the amplitude of the perturbation: see Table \ref{tab:Euler_Areas_for_different_values_of_alpha} and Figure \ref{fig:eulerintegralperturbationvsintegralsteadyshockperturbedlabelled}.
		\item If the perturbation moves the steady shock to the left, then a steady shock solution of the form \eqref{eq:testE6_solution} is obtained.
	\end{enumerate}
\end{conclusion}


\paragraph{Acknowledgments} 

Part of this paper was written during the Academic year 2018-2019 when PLF was a visiting fellow at the Courant Institute for Mathematical Sciences, New York University. The authors were supported by an Innovative Training Network (ITN) under the grant 642768 ``ModCompShock''. The research of CP and EPG was also partially supported by the Spanish Government(SG), the European Regional Development Fund (ERDF), the Regional Government of Andalusia (RGA), and the University of M\'alaga (UMA) through the projects RTI2018-096064-B-C21 (SG-ERDF), UMA18-Federja-161 (RGA-ERDF-UMA), and P18-RT-3163 (RGA-ERDF).
 


\end{document}